\documentclass[12pt]{amsart}

\usepackage[top=30mm, bottom=30mm, left=30mm, right=30mm]{geometry}
\usepackage{subfig}
\usepackage{bm}
\usepackage{caption}
\usepackage{commath}

\usepackage{wrapfig}
\usepackage{pinlabel, color}
\usepackage{multirow}

\usepackage[hyperindex, pdftex,colorlinks=true,
                      pdfstartview=FitV,
                      linkcolor=blue,
                      citecolor=blue,
                      urlcolor=blue,
          ]{hyperref}
\usepackage{ragged2e} 

\usepackage{array}

\usepackage[section]{placeins}

\newcommand{\HH}{\mathbb{H}}
\newcommand{\EE}{\mathbb{E}}
\newcommand{\RR}{\mathbb{R}}

\newif\ifimages
\imagestrue

\newif\iftables
\tablestrue

\thanks{The second author was supported in part by National Science Foundation grant DMS-1308767.}
\begin{document}

\title[Visualizing Hyperbolic Honeycombs]{Visualizing Hyperbolic Honeycombs}

\author[Nelson]{Roice Nelson}
\address{\hskip-\parindent
        Austin, TX USA}
\email{roice3@gmail.com}

\author[Segerman]{Henry Segerman}
\address{\hskip-\parindent
        Department of Mathematics \\
        Oklahoma State University \\
        Stillwater, OK USA}
\email{segerman@math.okstate.edu}


\begin{abstract}
We explore visual representations of tilings corresponding to Schl\"afli symbols. In three dimensions, we call these tilings ``honeycombs''. Schl\"afli symbols encode, in a very efficient way, regular tilings of spherical, euclidean and hyperbolic spaces in all dimensions. In three dimensions, there are only a finite number of spherical and euclidean honeycombs, but infinitely many hyperbolic honeycombs. Moreover, there are only four hyperbolic honeycombs with material vertices and material cells (the cells are entirely inside of hyperbolic space), eleven with ideal vertices or cells (the cells touch the boundary of hyperbolic space in some way), and all others have either hyperideal vertices or hyperideal cells (the cells go outside of the boundary of hyperbolic space in some way). We develop strategies for visualizing honeycombs in all of these categories, either via rendered images or 3D prints. High resolution images are available at \htmladdnormallink{hyperbolichoneycombs.org}{http://hyperbolichoneycombs.org}.
\end{abstract}

\maketitle

\begin{quote}
``Art attempts to portray in limited space, and song in time, what is otherwise unlimited.'' - Jacob Dix~\cite{Dix2013}
\end{quote}

\section{Introduction}
\label{Sec:Intro}

Schl\"afli symbols are an efficient system of notation for regular spherical, euclidean and hyperbolic tilings. These symbols can also represent regular polytopes, which convert to spherical tilings under radial projection. Schl\"afli symbols are defined recursively, as follows. A symbol of the form $\{p\}$ denotes the regular $p$-gon. A symbol of the form $\{p,q\}$ denotes a two-dimensional tiling of $p$-gons in which $q$ such polygons meet at each vertex. For example, $\{4,3\}$ denotes the tiling of squares, with three arranged around each vertex, also known as the cube. A Schl\"afli symbol of the form $\{p,q,r\}$ denotes a three dimensional tiling of \emph{cells} with Schl\"afli symbol $\{p,q\}$, with $r$ such cells meeting at each edge. For example, $\{4,3,4\}$ denotes the tiling of three-dimensional euclidean space by cubes.  In three dimensions, we also call these tilings \emph{honeycombs}. 

In this paper, we describe ways to visualize all $\{p,q,r\}$ honeycombs, where each of $p,q$ and $r$ is an integer greater than two, or $\infty$. 

The honeycombs of three-dimensional spherical space, corresponding under radial projection to the six regular 4-dimensional polytopes, have been visualized in many ways. The tiling of euclidean space by cubes, $\{4,3,4\}$, is the only regular euclidean honeycomb. Some hyperbolic honeycombs have previously been visualized in the software \emph{Curved Spaces}~\cite{curved_spaces}, by Weeks, and in the video \emph{Not Knot}~\cite{not_knot}, by Epstein et al., with visualization techniques of the latter described in depth in subsequent papers~\cite{gunn1993discrete, phillips1992visualizing}. In these, hyperbolic space is shown in the ``in-space'' view, and the emphasis is on visualizing hyperbolic manifolds rather than honeycombs. The in-space view works well for honeycombs with smaller cells, corresponding to smaller values of $p,q$ and $r$. As we will see in Section \ref{Sec:hyperidealization}, this breaks down once the cells become so big that they cross the boundary of hyperbolic space. At this point, we switch to viewing the pattern of intersections of cells with the boundary, allowing us to visualize all remaining honeycombs. The technique of drawing an intersection with the boundary of hyperbolic space has been used effectively by Vladimir Bulatov~\cite{bulatovbending}, producing intriguing visualizations of hyperbolic symmetries. We should also mention that many beautiful images drawn on the boundary of hyperbolic space have been generated -- for example the limit sets of Kleinian groups~\cite{jos_leys,Mumford2002}. 

Two-dimensional hyperbolic geometry has been featured in works of art since Escher's \emph{Circle Limit} series of prints. Three-dimensional hyperbolic space is less well explored, but we feel that there is much potential. Our images include, as a side effect, apparent two-dimensional hyperbolic tilings, as well as many other new features.  We will see that the sheer variety of possible images makes them accessible to creative exploration.  Interesting artistic challenges also arise in creating sculptures based on hyperbolic honeycombs.

We warm up by looking at two-dimensional tilings corresponding to Schl\"afli symbols of the form $\{p,q\}$, and view these in spherical, euclidean and hyperbolic spaces (Section \ref{Sec:2d}). We also note the relation between Schl\"afli symbols and duality, and how to determine which geometry a tiling lives in from its symbol. We move on to three-dimensional honeycombs, again noting the relation with duality and how to determine the geometry of a honeycomb (Section \ref{Sec:3d}). We then introduce tilings and honeycombs with ideal vertices and ideal cells (Section \ref{Sec:ideal}). At this point, we can provide a visual map of the different kinds of geometries and vertex/cell types of honeycombs (Section \ref{Sec:3d_map}). We then introduce tilings and honeycombs with hyperideal vertices and hyperideal cells (Section \ref{Sec:hyperidealization}). With all classes of honeycombs described, we are ready to draw pictures (in various forms) of the honeycombs on the boundary of the infinite ``cube'' of three-dimensional Schl\"afli symbols, and on surfaces within the cube corresponding to honeycombs with ideal vertices or ideal cells (Section \ref{Sec:schlafli_cube}). We display artworks based on our investigations, both 3D printed sculpture and 2D renders (Section \ref{Sec:artwork}). Finally, we give some directions for future work (Section \ref{Sec:future}).

This paper is intended to give a tour of the qualitative concepts and images of these hyperbolic honeycombs, without giving too many of the precise geometric details of the construction of the images and models. We give some of these details in Appendix \ref{Sec:implementation}. We also discuss the question of the uniqueness of a honeycomb or tiling corresponding to a Schl\"afli symbol in Appendix \ref{Sec:uniqueness}.

For background information on hyperbolic space, see, e.g. \emph{Euclidean and Non-Euclidean Geometries: Development and History} by Greenberg~\cite{greenberg1993euclidean}. For more technical material, see \emph{Visual Complex Analysis} by Needham~\cite{Needham199902}, and the code we used to generate the images in this paper~\cite{hyp_honeycombs_code}.

\section{Two dimensions}
\label{Sec:2d}
\subsection{Geometries}

We start by looking at some examples of two-dimensional tilings given by Schl\"afli symbols, see Figure \ref{Fig:43_44_45}. For most cases, there is a uniquely defined regular tiling for every Schl\"afli symbol (see Appendices \ref{Sec:honeycomb_construction} and \ref{Sec:uniqueness} for details beyond the discussion here).  Immediately there is a question of how to draw pictures of these tilings. The $\{4,4\}$ tiling is very familiar, and can be drawn on the page without any distortion, but the other two tilings cannot. The $\{4,3\}$ tiling is otherwise known as the radial projection of the cube: three squares meet around each vertex, tiling the sphere with six squares. 

Figure \ref{Fig:43} shows the stereographic projection~\cite[pp. 140-148]{Needham199902} of the $\{4,3\}$ tiling of the sphere to the plane. Note that stereographic projection is a one-to-one map from the sphere to the plane with one point added at infinity. Here, five of the six faces of the cube are mapped inside of the ``clover leaf'' shape, and the sixth is mapped outside, covering the entire rest of the plane and the point at infinity. In Figure \ref{Fig:45} we draw the $\{4,5\}$ tiling on the Poincar\'e disk model of $\HH^2$~\cite[pp. 315-319]{Needham199902}.

All three of these tilings live in \emph{isotropic} spaces (and we draw projections of them to the page). An isotropic space looks the same at every point, looking in any direction. More precisely, there is a symmetry of the space that takes any pair of a point and a direction based at that point to any other such pair. In general, we will interpret Schl\"afli symbols as describing tilings of (simply connected) isotropic spaces: in two dimensions these are the sphere, $S^2$, the euclidean plane, $\EE^2$, and the hyperbolic plane, $\HH^2$.

\begin{figure}[htbp]
\centering 
\hspace{-.7cm}
\subfloat[$\{4,3\}$]
{
\ifimages
\includegraphics[width=0.3\textwidth]{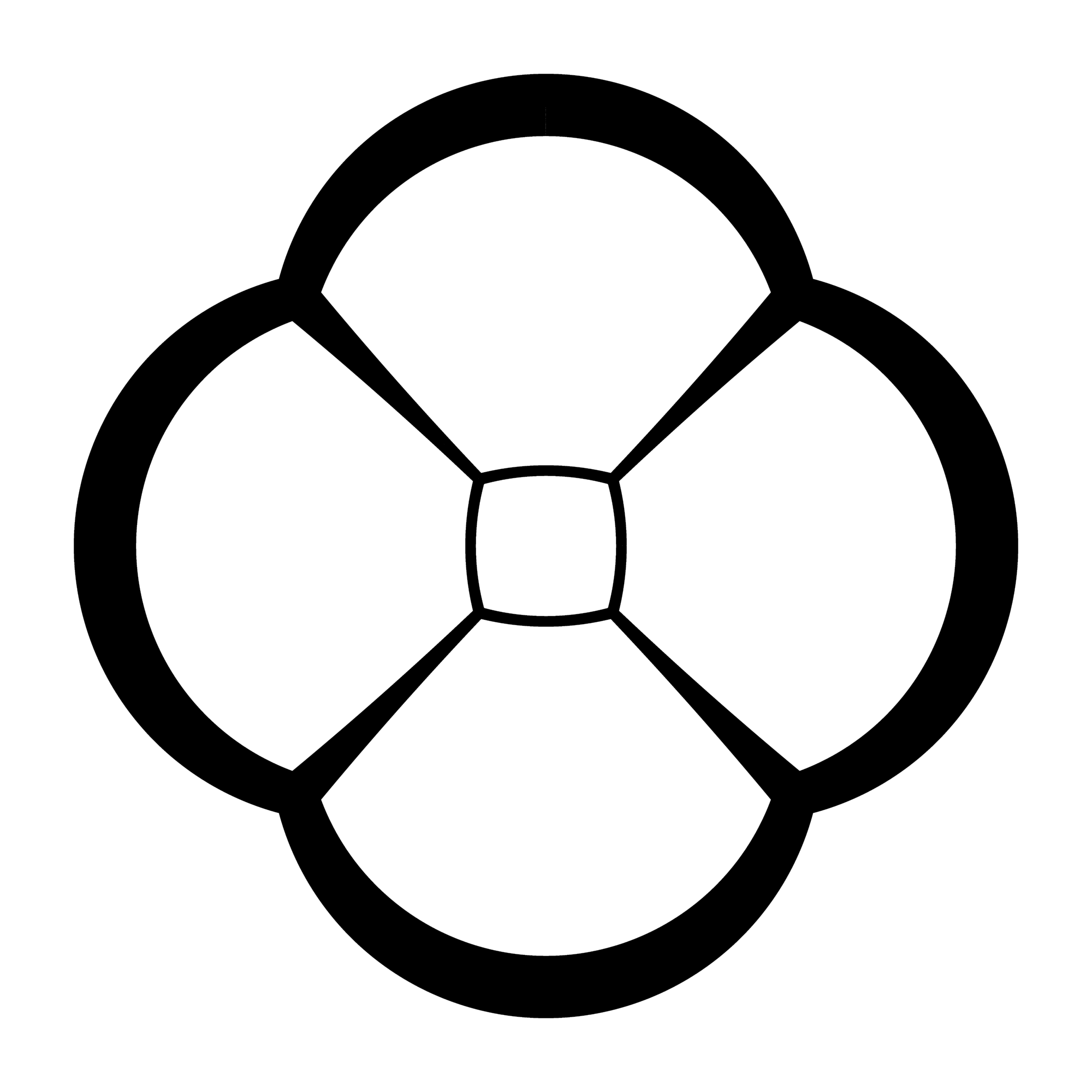}
\fi
\label{Fig:43}
}
\hspace{.1cm}
\subfloat[$\{4,4\}$]
{
\ifimages
\includegraphics[width=0.3\textwidth]{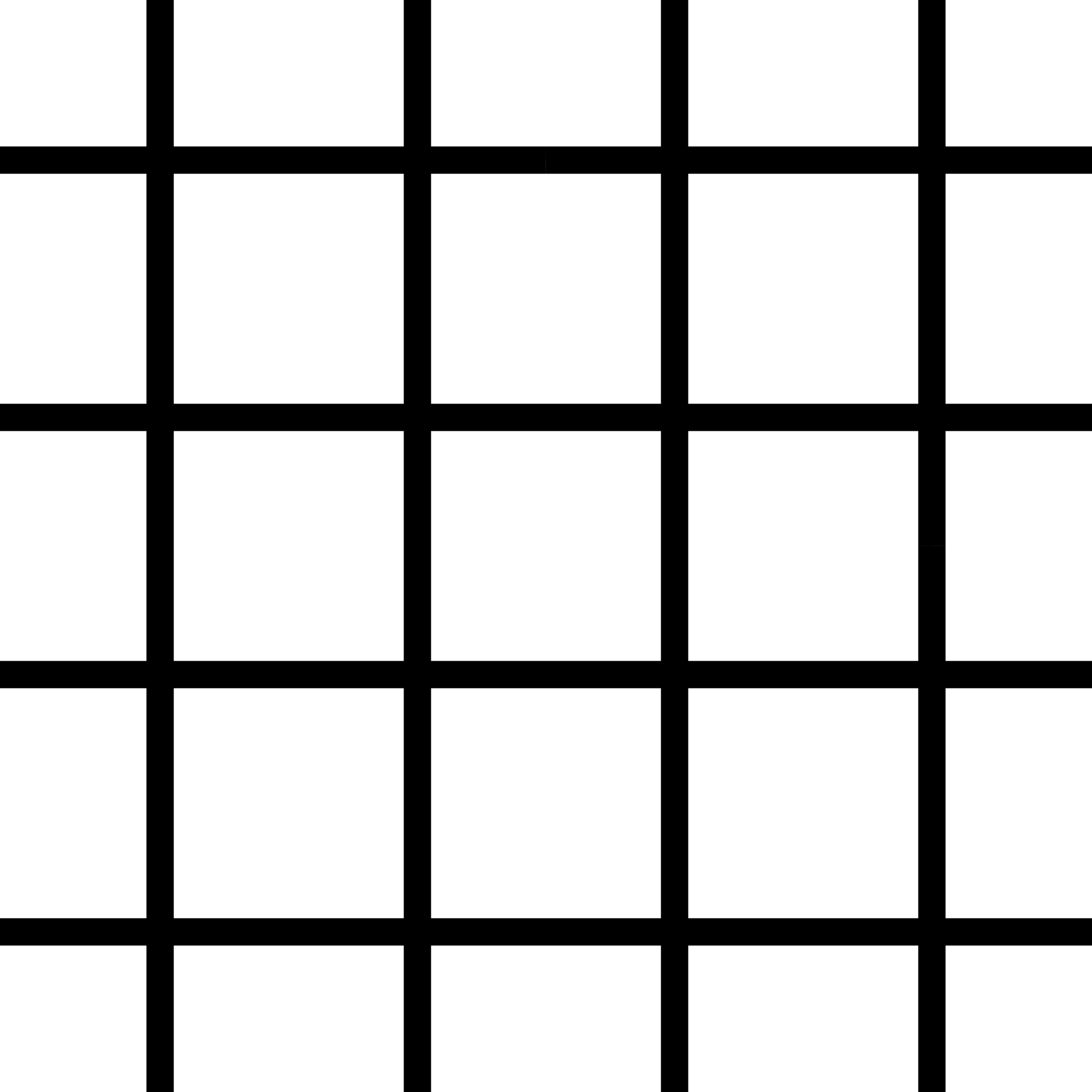}
\fi
\label{Fig:44}
}
\hspace{.5cm}
\subfloat[$\{4,5\}$]
{
\ifimages
\includegraphics[width=0.3\textwidth]{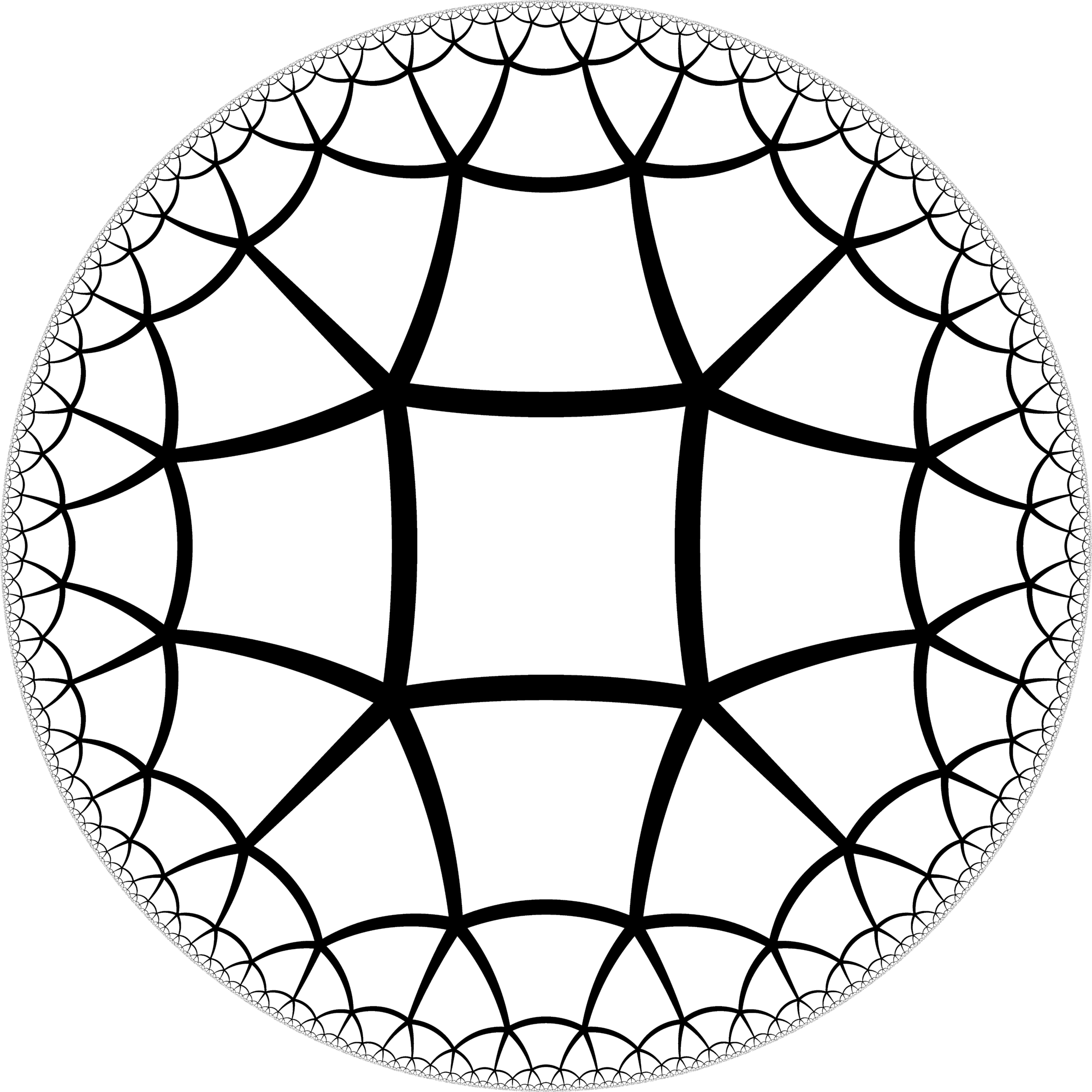}
\fi
\label{Fig:45}
}
\caption{Two dimensional tilings by squares.}
\label{Fig:43_44_45}
\end{figure}

\subsection{Duality}

\begin{figure}[htbp]
\centering 
\hspace{-.7cm}
\subfloat[$\{3,4\}$]
{
\ifimages
\includegraphics[width=0.3\textwidth]{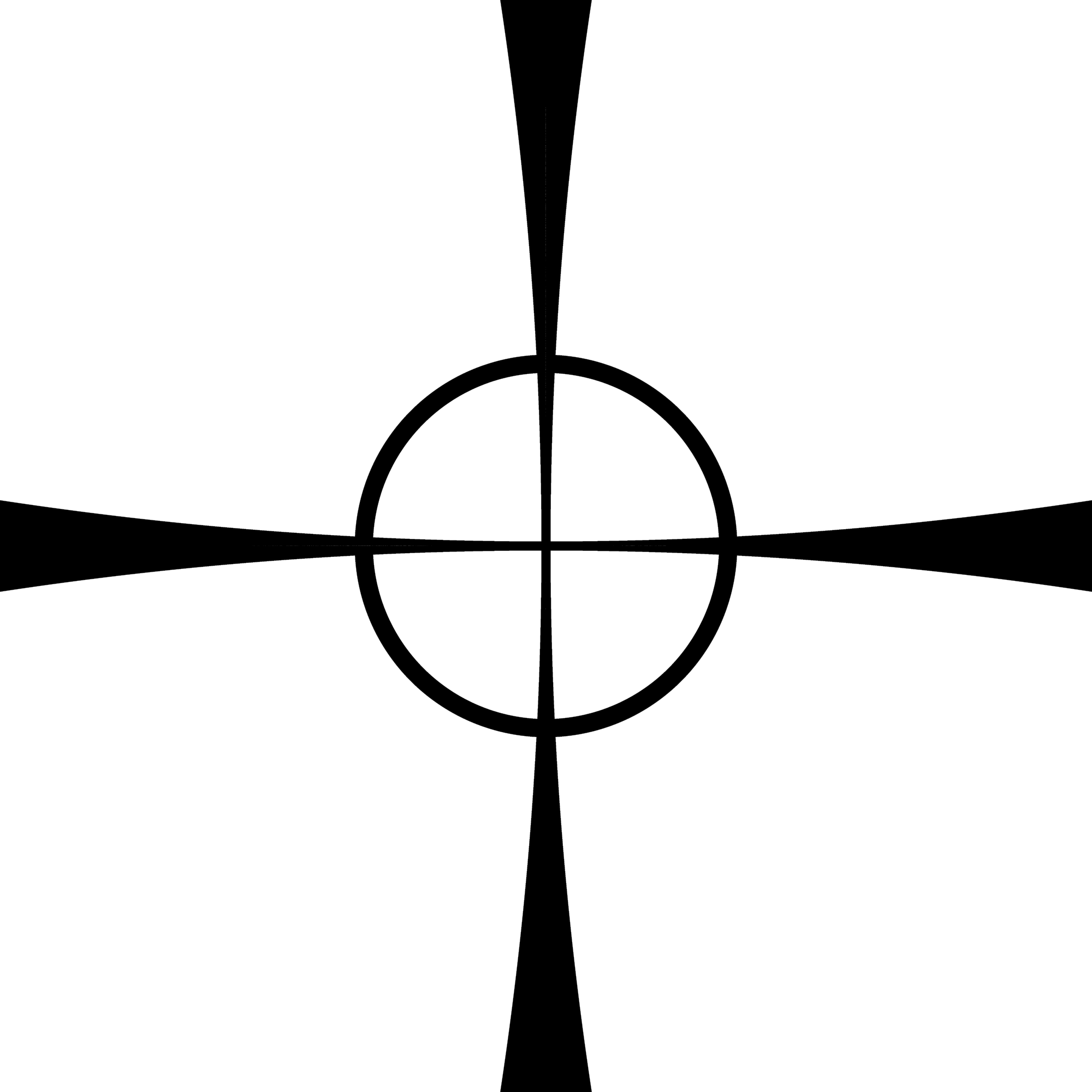}
\fi
\label{Fig:34}
}
\hspace{.1cm}
\subfloat[$\{4,4\}$]
{
\ifimages
\includegraphics[width=0.3\textwidth]{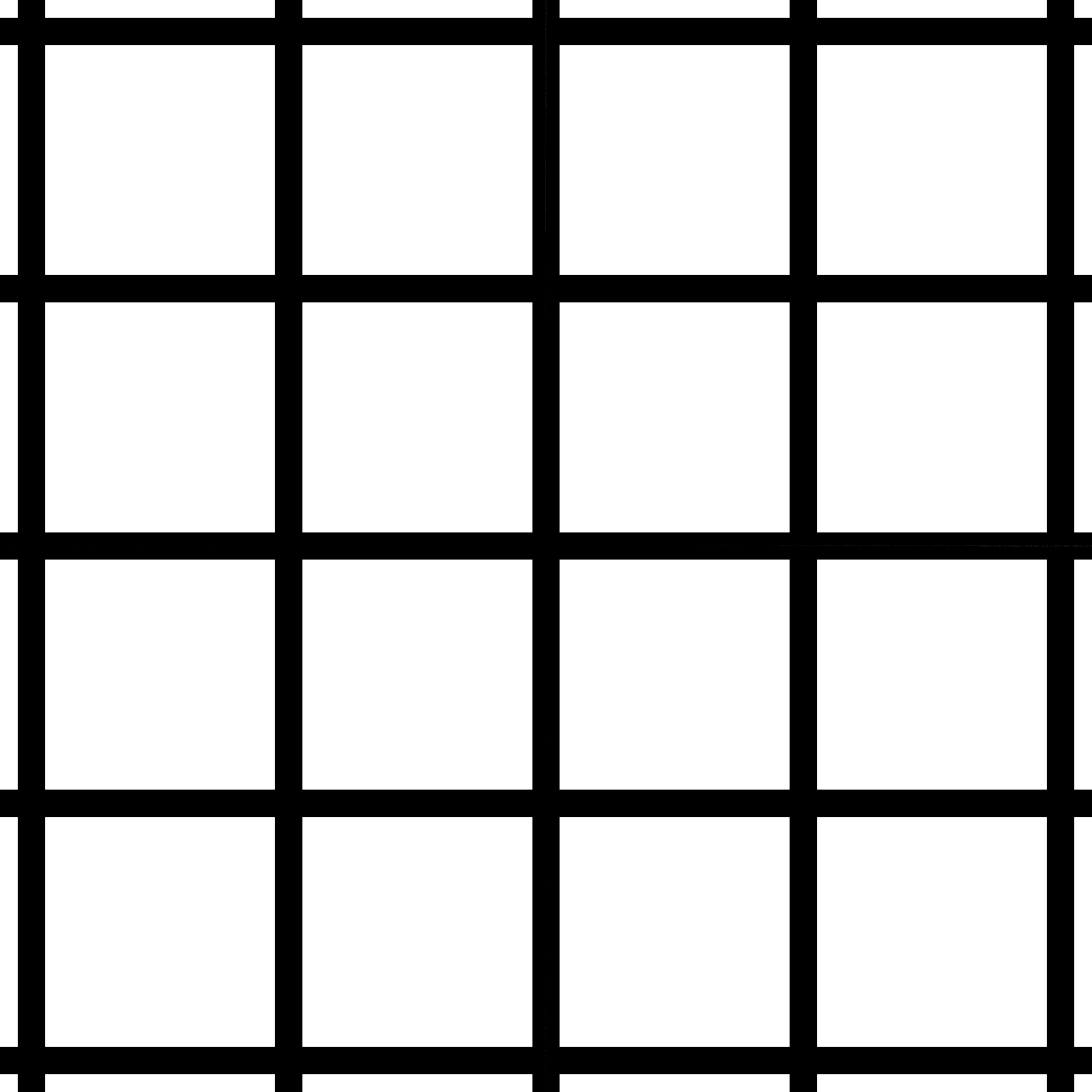}
\fi
\label{Fig:44_dual}
}
\hspace{.5cm}
\subfloat[$\{5,4\}$]
{
\ifimages
\includegraphics[width=0.3\textwidth]{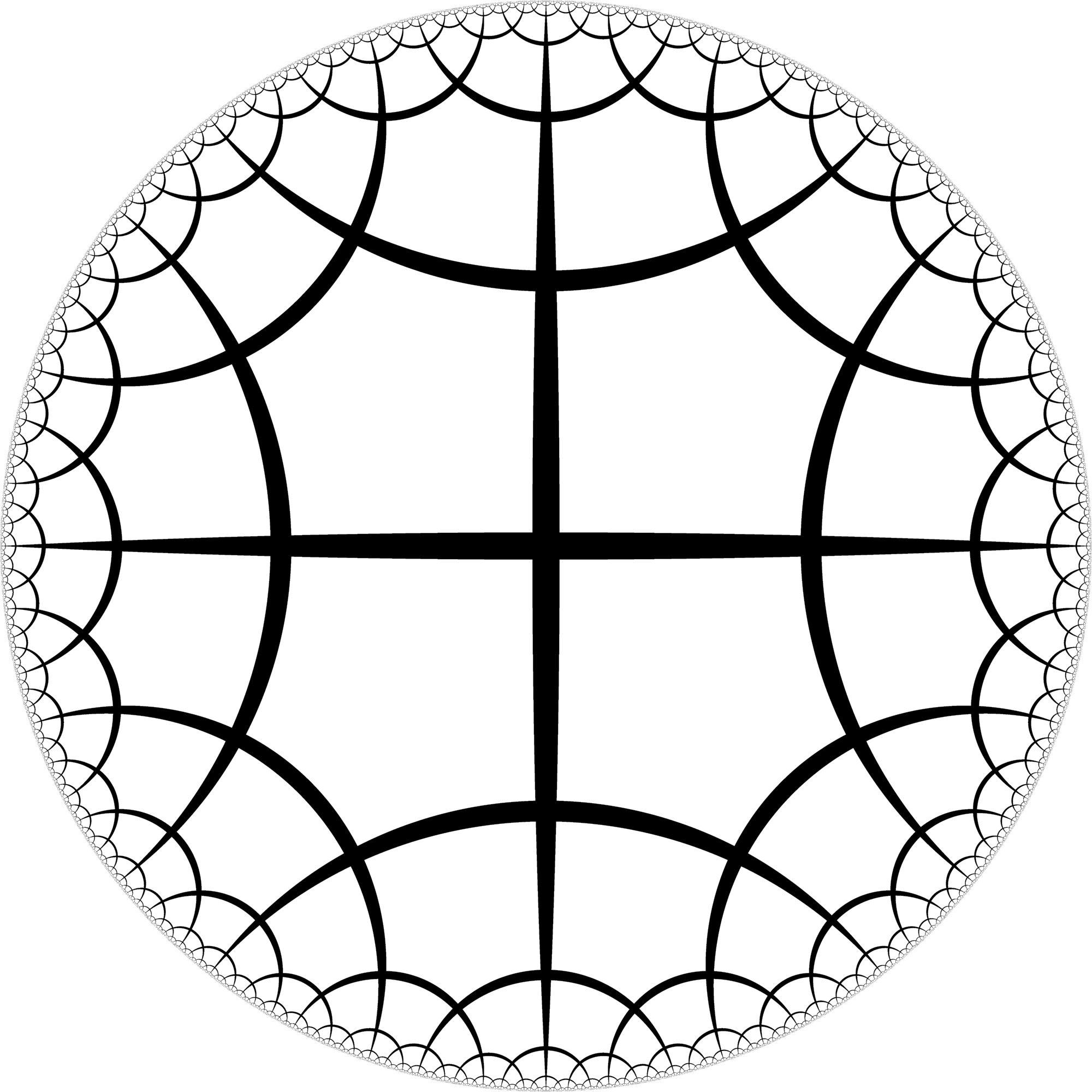}
\fi
\label{Fig:54}
}
\caption{Dual two-dimensional tilings.}
\label{Fig:34_44_54}
\end{figure}

As is well-known, the regular polyhedra relate to each other by duality.  For example, given a cube (which we think of as the spherical $\{4,3\}$ tiling in Figure \ref{Fig:43}), we can construct the dual as follows. We plot a new vertex at the center of each face of the cube. Whenever two of these vertices are separated by an edge of the cube, we connect them by a new edge. These new vertices and edges give the dual octahedron (thought of as the spherical $\{3,4\}$ tiling in Figure \ref{Fig:34} -- note that one of the vertices is mapped to infinity under stereographic projection).  If we start with an octahedron and perform the same operations, we get back the cube. Similarly the dodecahedron and icosahedron are duals of each other, and the tetrahedron is self-dual. In terms of the Schl\"afli symbols, we get the dual of a polyhedron by reversing its Schl\"afli symbol. This relation extends from polyhedra (tilings of the sphere) to tilings of the euclidean plane ($\{3,6\}$ is dual to $\{6,3\}$ and $\{4,4\}$ is self-dual) and tilings of the hyperbolic plane.

For both spherical and hyperbolic spaces, in projecting to the euclidean plane of the page, we make a choice. In particular Figures \ref{Fig:43} and \ref{Fig:45} show \emph{tile-centered} pictures of the $\{4,3\}$ and $\{4,5\}$ tilings, while Figures \ref{Fig:34} and \ref{Fig:54} show \emph{vertex-centered} pictures of the $\{3,4\}$ and $\{5,4\}$ tilings, so called because a tile (respectively, vertex) of the tiling is at the center of the stereographic projection or Poincar\'e disk model on the page.

\subsection{Which two-dimensional Schl\"afli symbols produce tilings of which geometries?}

Suppose we try to tile the euclidean plane with polygons as given to us by a Schl\"afli symbol $\{p,q\}$. Then $q$ $p$-gons meet at a vertex, each of which has internal angle $(1-2/p)\pi$, for a total angle of $q(1-2/p)\pi$. We will be able to tile the euclidean plane if this value is equal to $2\pi$, or equivalently if 

\begin{equation}
(p-2)(q-2)=4
\label{2d_equation}
\end{equation}

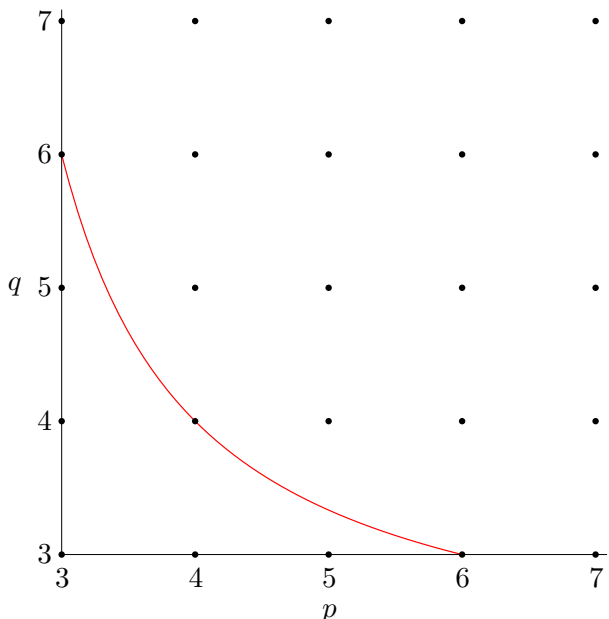
\begin{figure}[htbp]
\centering 
\labellist
\small\hair 2pt
\pinlabel 3 at 35 15
\pinlabel 4 at 147 15
\pinlabel 5 at 259 15
\pinlabel 6 at 371 15
\pinlabel 7 at 483 15
\pinlabel $p$ at 259 -15

\pinlabel 3 at 20 35 
\pinlabel 4 at 20 147
\pinlabel 5 at 20 259
\pinlabel 6 at 20 371
\pinlabel 7 at 20 483

\pinlabel $q$ at -5 259 

\endlabellist

\includegraphics[width=0.5\textwidth]{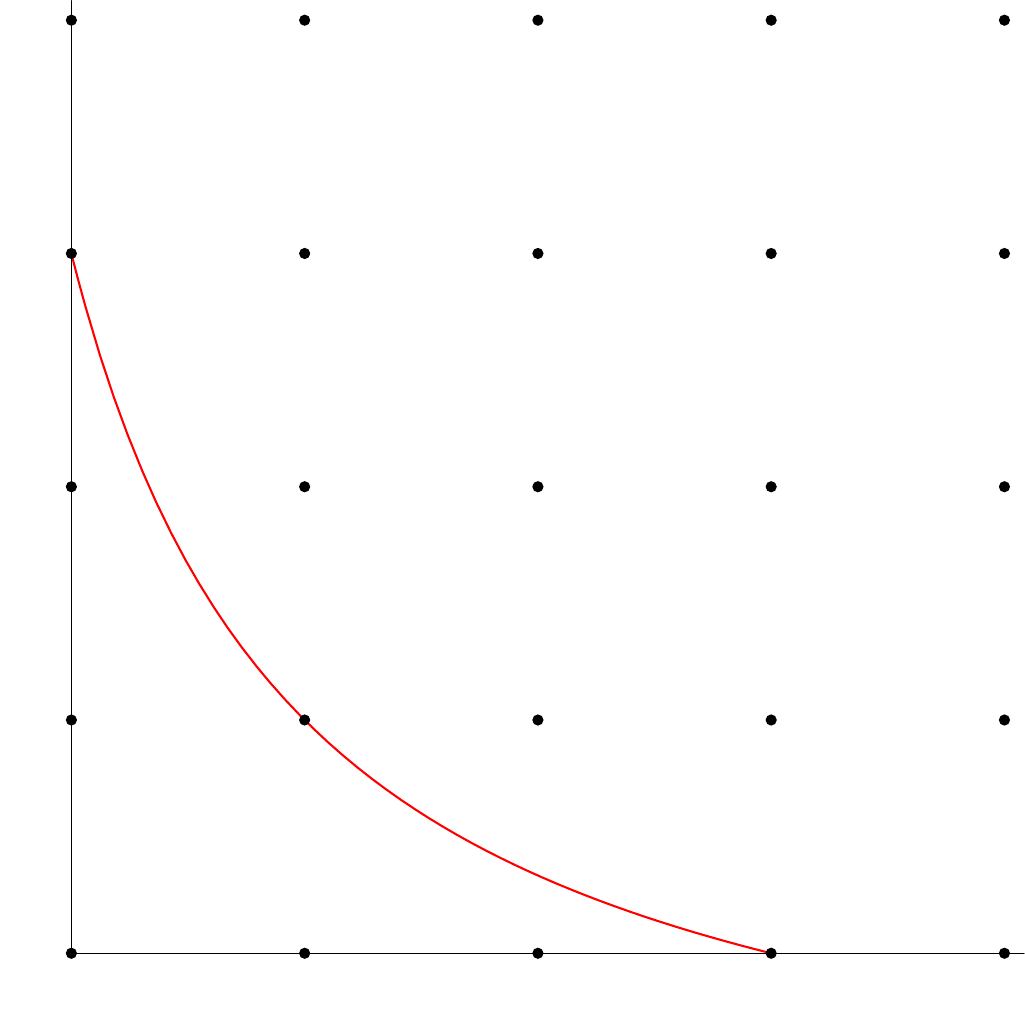}
\caption{Geometries associated to the two-dimensional Schl\"afli symbols: below the curve $(p-2)(q-2)=4,$ the integer points correspond to tilings of the sphere, on the curve to tilings of the euclidean plane, and above the curve to tilings of the hyperbolic plane.}
\label{Fig:2D_Schlafli_map}
\end{figure}

If the expression on the left hand side is less than $4$, then we instead get a polyhedron~\cite[p. 5]{coxeter1973}, which as before we view as a tiling of the sphere. If this value is greater than $4$, then we can tile the hyperbolic plane~\cite[pp. 155-156]{Coxeter1954}. See Figure \ref{Fig:2D_Schlafli_map}.

\section{Three dimensions}
\label{Sec:3d}
\subsection{Geometries}

\begin{figure}[htbp]
\centering 
\hspace{-.7cm}
\subfloat[$\{4,3,3\}$]
{
\ifimages
\includegraphics[width=0.3\textwidth]{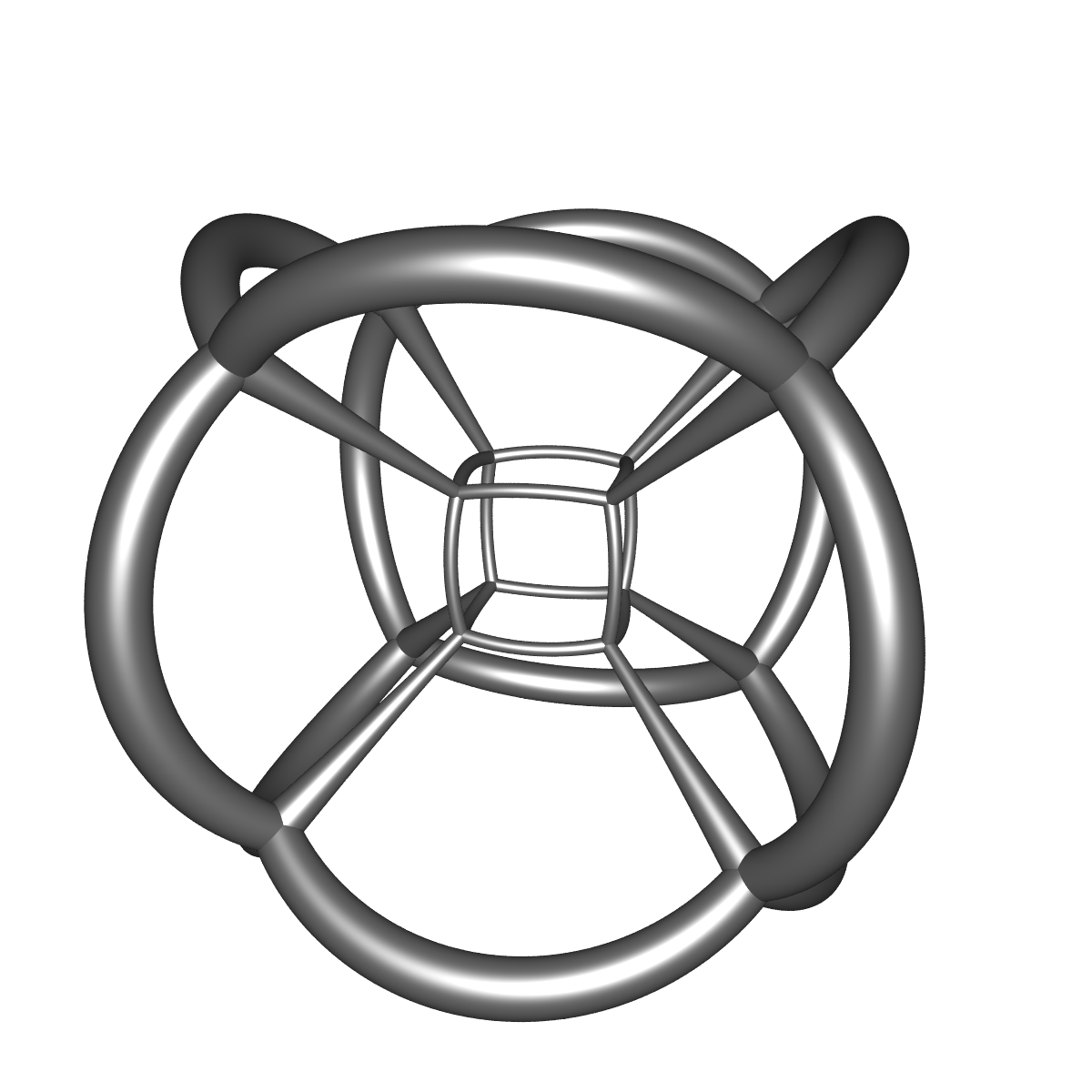}
\fi
\label{Fig:433}
}
\hspace{.1cm}
\subfloat[$\{4,3,4\}$]
{
\ifimages
\includegraphics[width=0.3\textwidth]{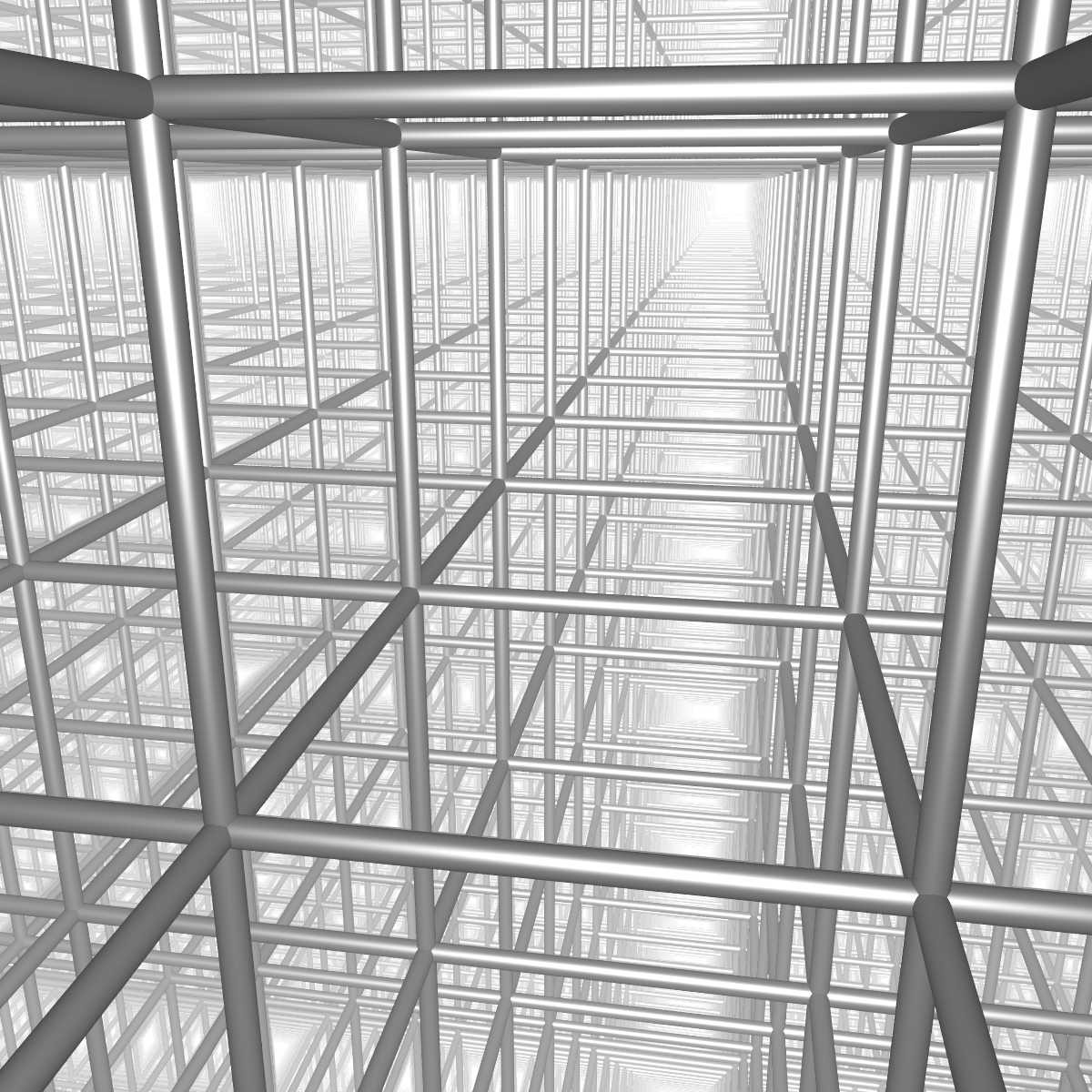}
\fi
\label{Fig:434}
}
\hspace{.5cm}
\subfloat[$\{4,3,5\}$]
{
\ifimages
\includegraphics[width=0.3\textwidth]{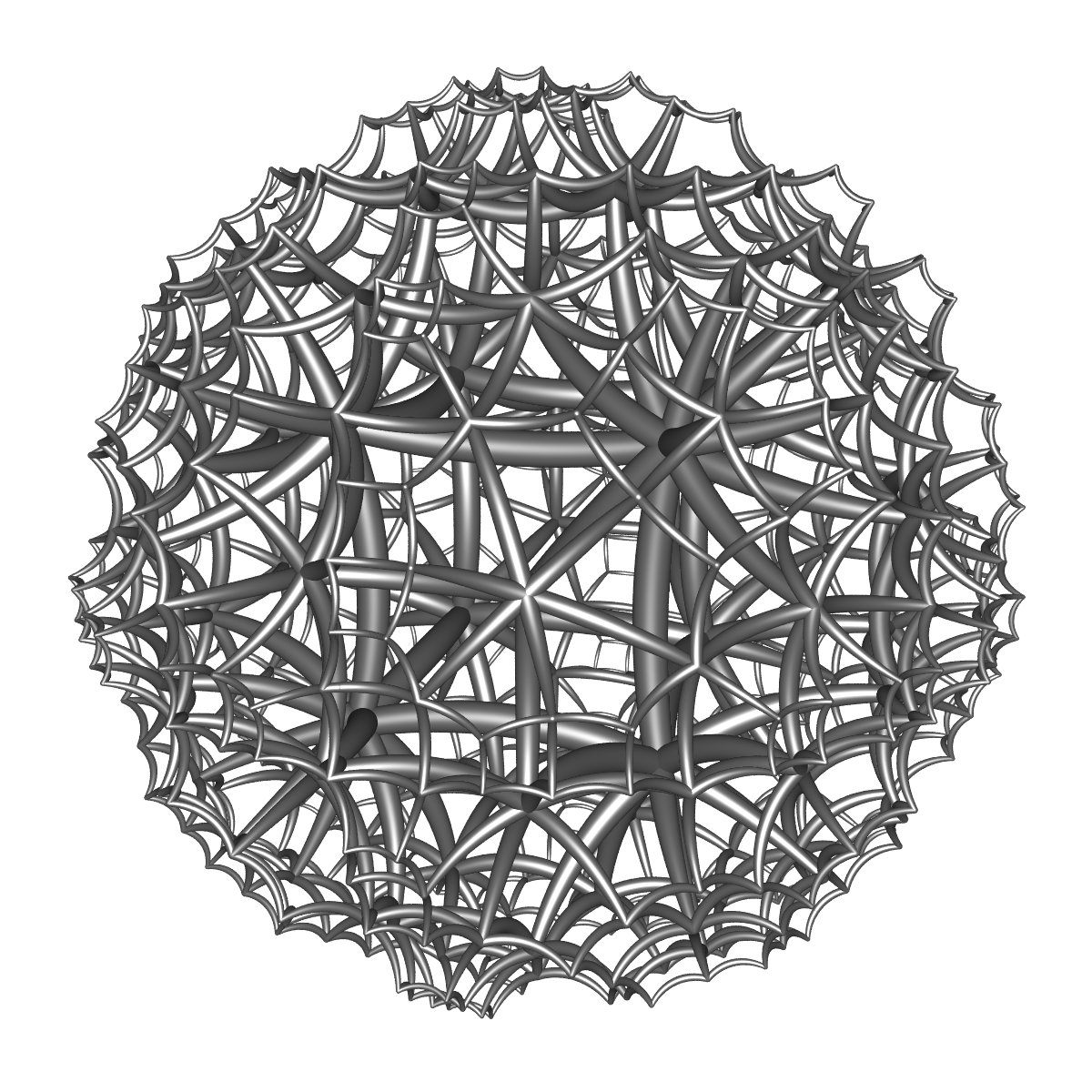}
\fi
\label{Fig:435_render}
}
\caption{Three-dimensional honeycombs with cubical cells.}
\label{Fig:433_434_435}
\end{figure}

Three-dimensional Schl\"afli symbols also represent tilings of various isotropic spaces. Here, we tile the three-sphere, $S^3$, three-dimensional euclidean space, $\EE^3$, and three-dimensional hyperbolic space, $\HH^3$. In three dimensions, we also call these tilings \emph{honeycombs}. Figure \ref{Fig:433_434_435} shows spherical, euclidean and hyperbolic honeycombs made out of cubes.

\begin{figure}[htbp]
\centering 
\subfloat[The central cube.]
{
\ifimages
\includegraphics[width=0.3\textwidth]{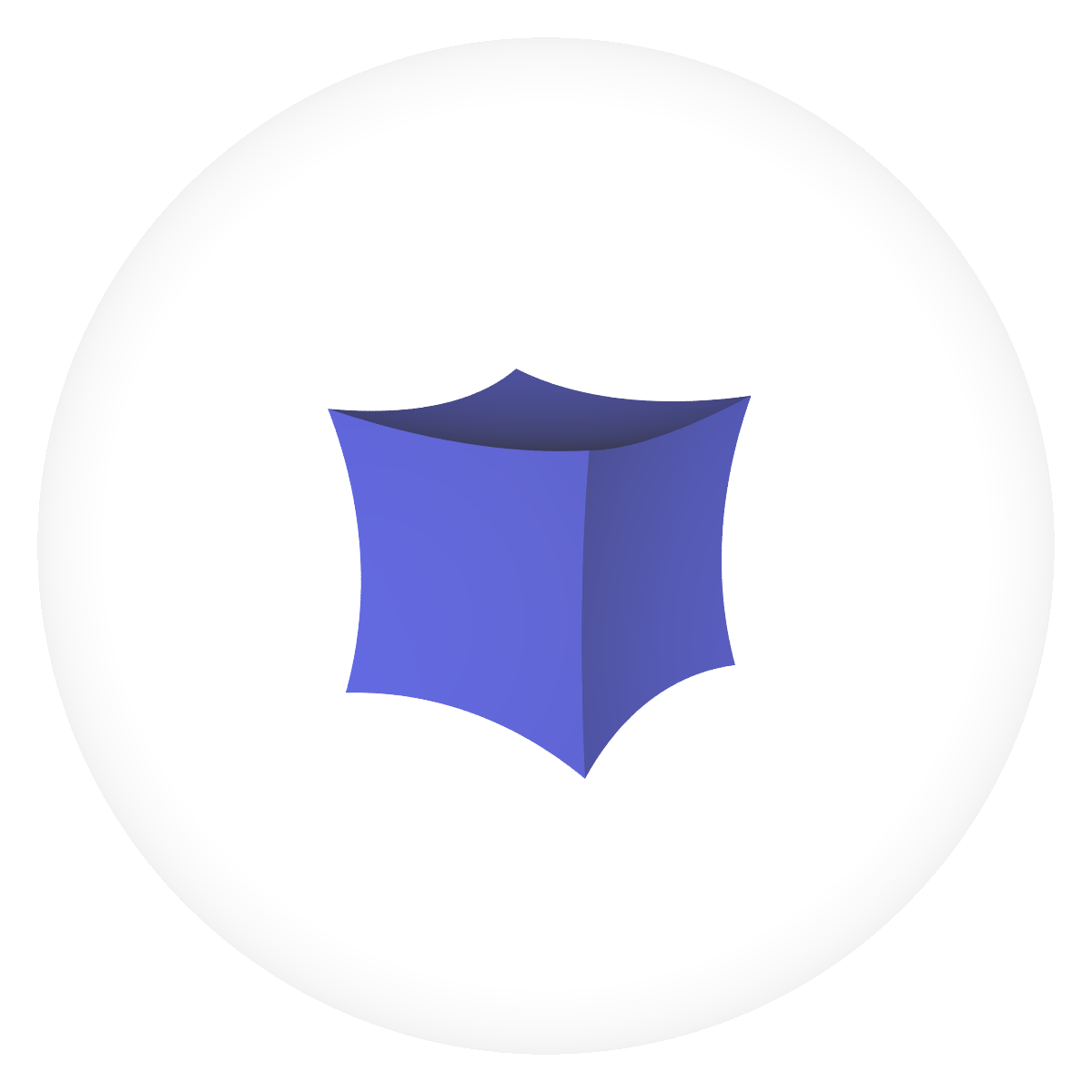}
\fi
\label{Fig:435_1}
}
\subfloat[Six more cubes.]
{
\ifimages
\includegraphics[width=0.3\textwidth]{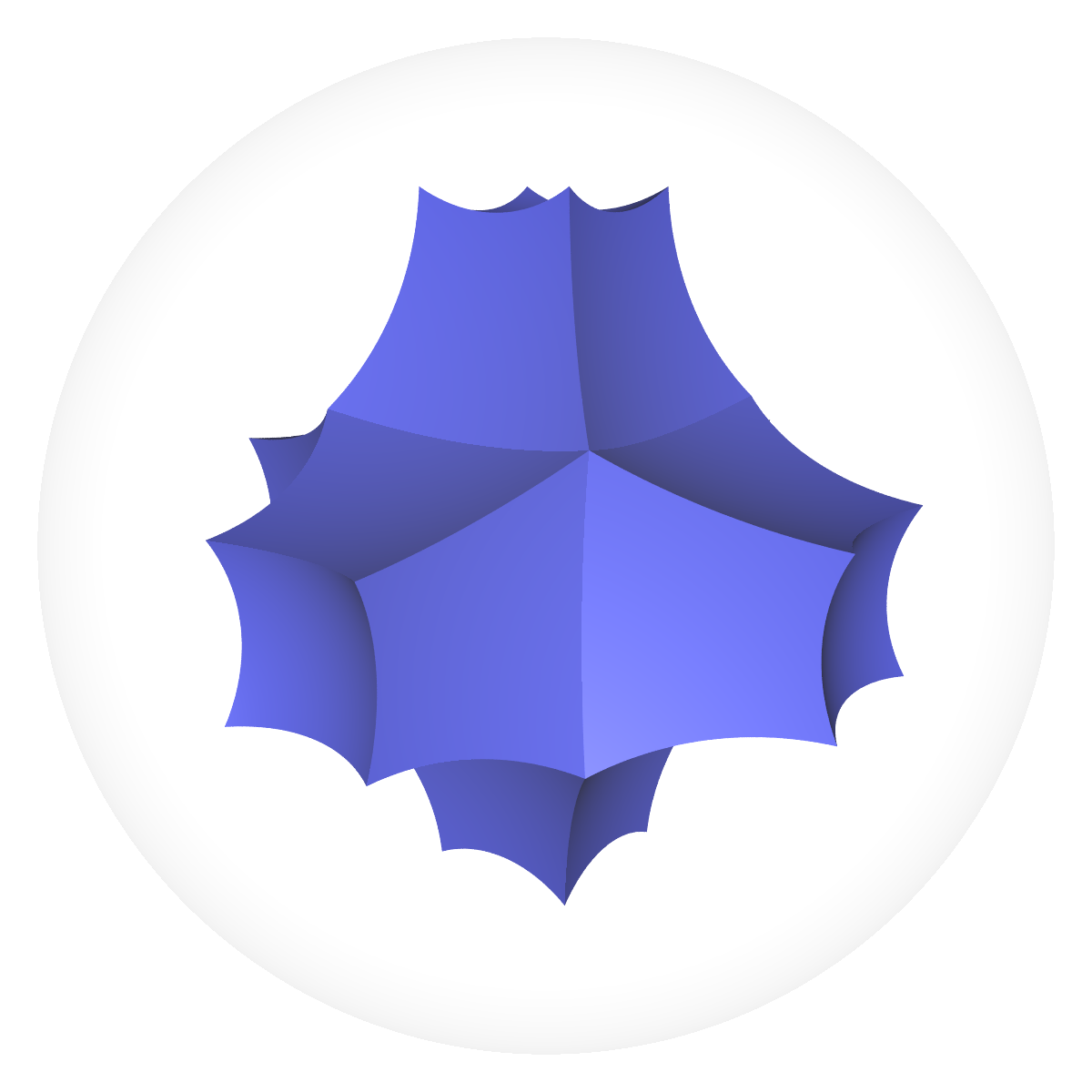}
\fi
\label{Fig:435_2}
}
\subfloat[Thirty more cubes.]
{
\ifimages
\includegraphics[width=0.3\textwidth]{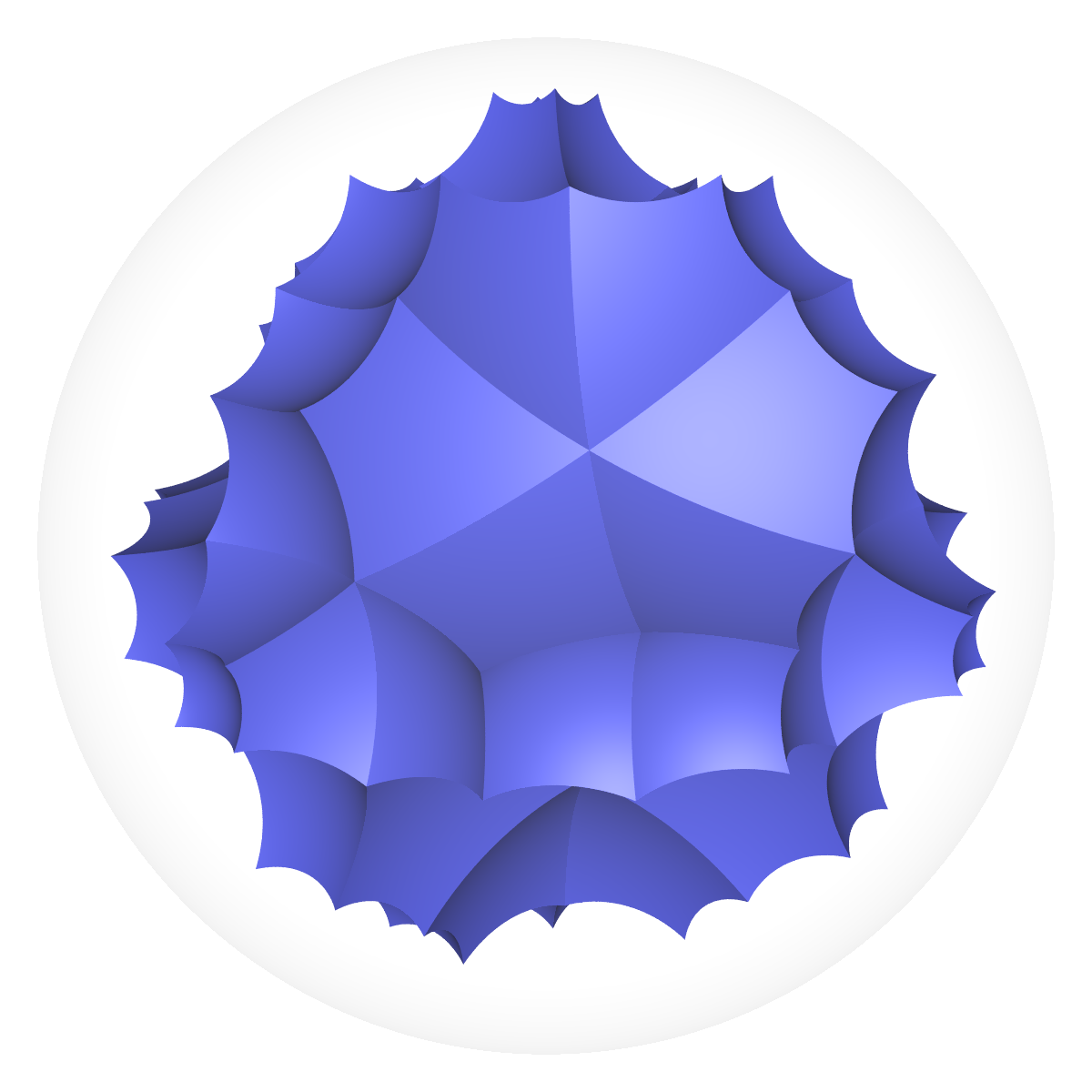}
\fi
\label{Fig:435_3}
}
\caption{Three layers of cubes in the $\{4,3,5\}$ honeycomb.}
\label{Fig:435_build_layers}
\end{figure}

Figure \ref{Fig:433} shows (the edges of) the tiling of $S^3$ by eight cubes. As a four-dimensional polytope, this is the hypercube -- three cubes meet around each edge, corresponding to the last term of its Schl\"afli symbol $\{4,3,3\}$.  As with the two-dimensional tilings of $S^2$ we can get from the four-dimensional polytope to the tiling (honeycomb) of $S^3$ by radially projecting the polytope onto $S^3$.  We draw the hypercube by first doing this radial projection, then stereographically projecting to three-dimensional euclidean space~\cite{schleimer2012sculptures}. In three dimensions, stereographic projection is a one-to-one map from $S^3$ to euclidean space with one point added at infinity. Here, seven of the eight cells of the hypercube are mapped inside of the shape shown in the figure, and the eighth is mapped outside, covering the entire rest of euclidean space and the point at infinity. Figure \ref{Fig:434} shows the familiar tiling of $\EE^3$ by cubes. Figure \ref{Fig:435_render} shows the $\{4,3,5\}$ tiling of the Poincar\'e ball model of $\HH^3$. Or rather, this shows a small part of the tiling -- there is a central cube, surrounded by a number of layers of cubes. The full tiling fills up all of $\HH^3$, so it would be impossible to see into the interior of the structure if we showed all of it. Figure \ref{Fig:435_build_layers} shows three layers of the $\{4,3,5\}$ tiling: the central cube, the cubes adjacent to the central cube, and the cubes adjacent to those.

\subsection{Duality}

\begin{figure}[htbp]
\centering 
\subfloat[$\{3,3,4\}$]
{
\ifimages
\includegraphics[width=0.3\textwidth]{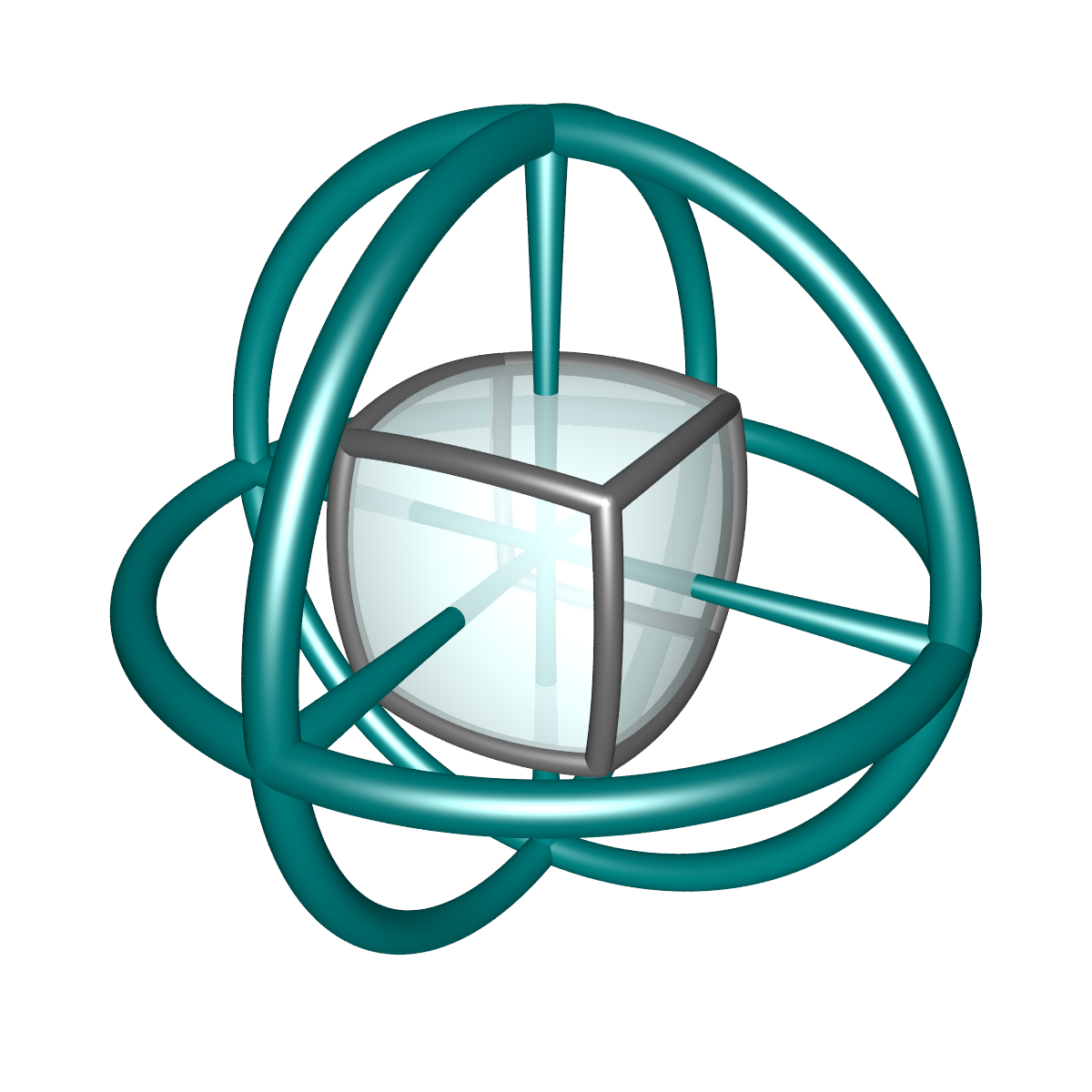}
\fi
\label{Fig:433_dual}
}
\subfloat[$\{4,3,4\}$]
{
\ifimages
\includegraphics[width=0.3\textwidth]{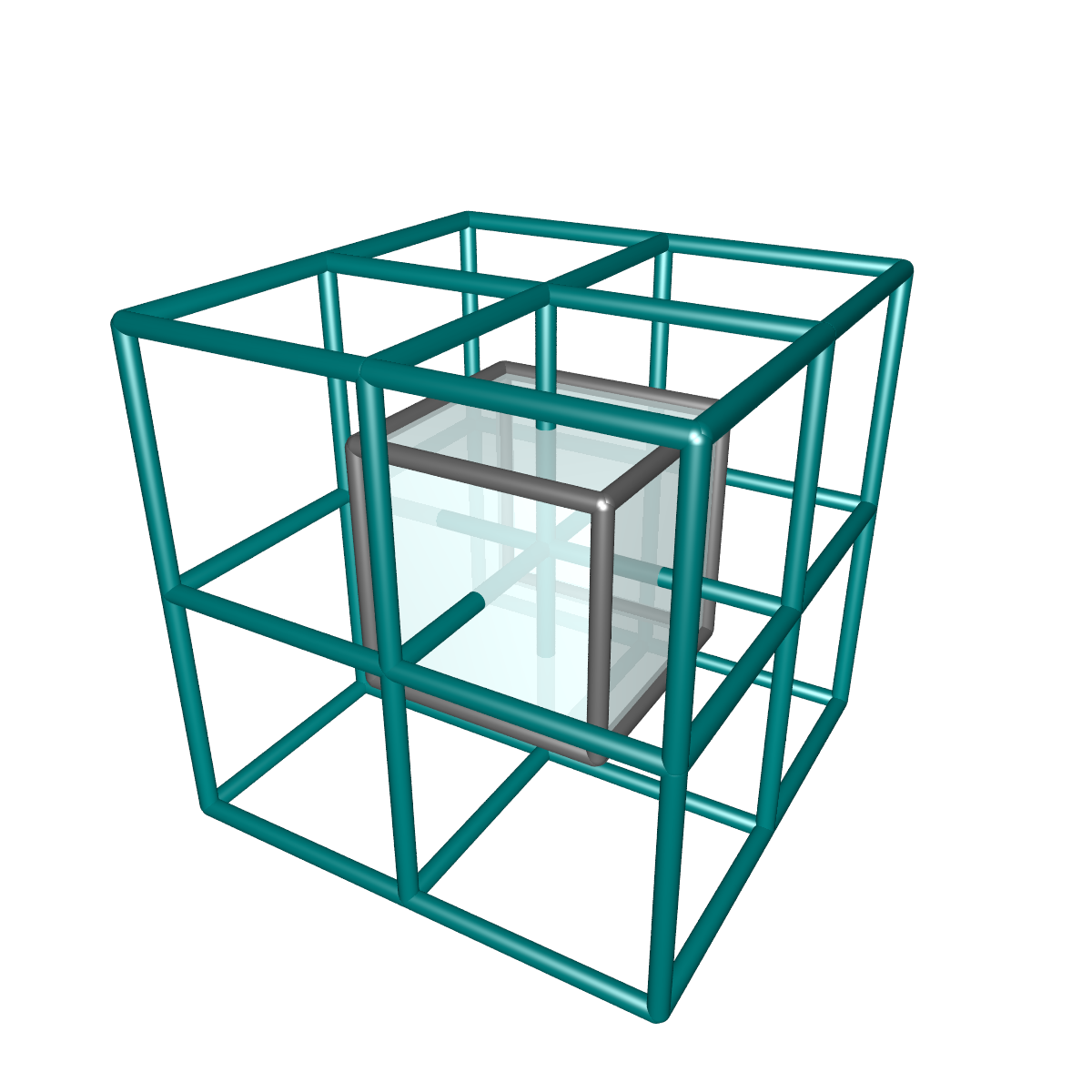}
\fi
\label{Fig:434_dual}
}
\subfloat[$\{5,3,4\}$]
{
\ifimages
\includegraphics[width=0.3\textwidth]{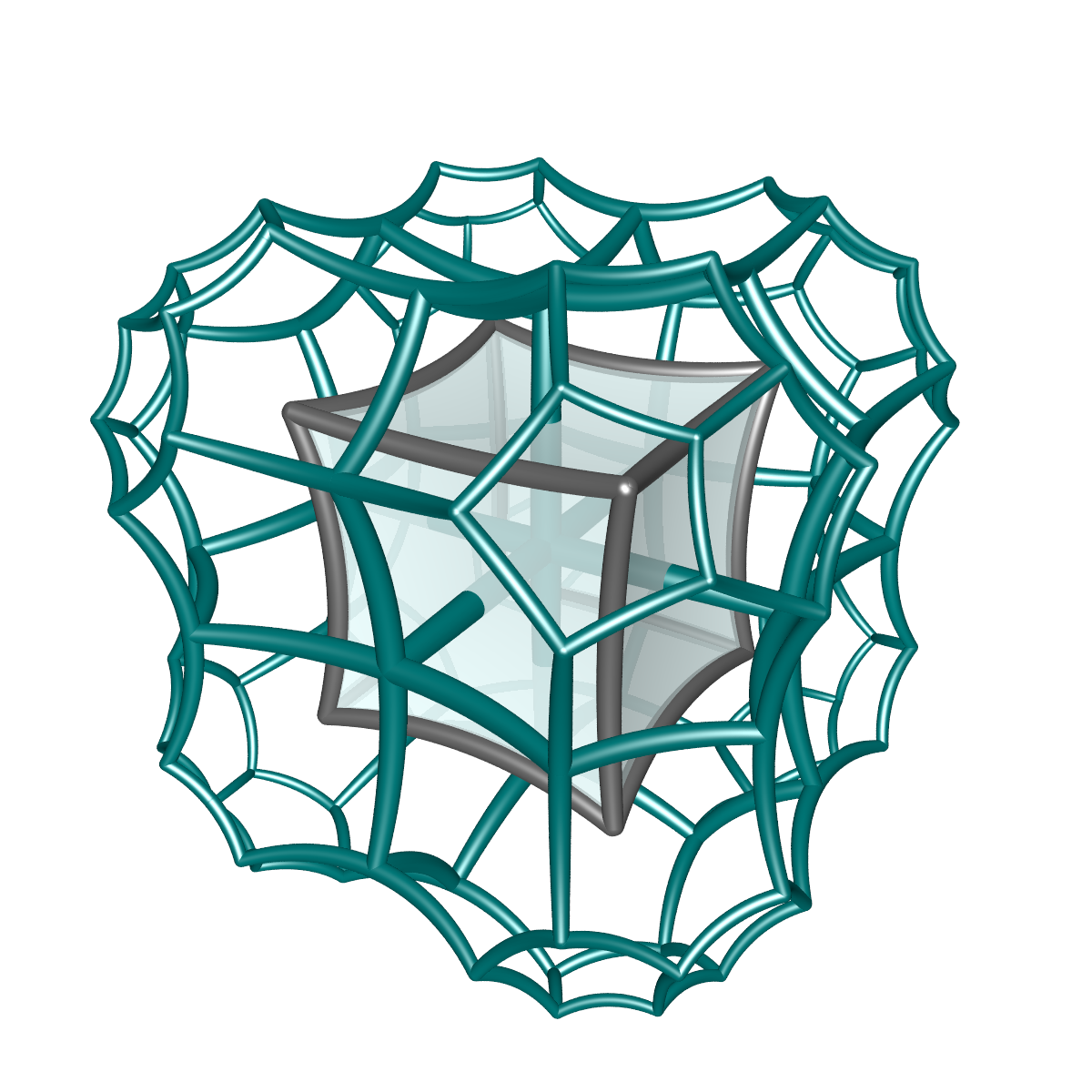}
\fi
\label{Fig:435_dual}
}
\caption{Three-dimensional honeycombs dual to the honeycombs in Figure \ref{Fig:433_434_435}. Each image shows eight cells from the honeycomb, and one cubical cell from the corresponding dual honeycomb.}
\label{Fig:3D_duals}
\end{figure}

As in two dimensions, the three-dimensional Schl\"afli symbols also have a notion of duality, and again we can get from a tiling to its dual by reversing the Schl\"afli symbol. See Figure \ref{Fig:3D_duals}. Geometrically, the process to get the dual of a honeycomb goes as follows: First we put a new vertex at the center of each three-dimensional cell of the intial honeycomb. Next, we join two new vertices by a new edge if the corresponding cells of the old honeycomb share a face. We add a new face cutting through each edge of the old honeycomb, and finally we add a new cell corresponding to each vertex of the old honeycomb.  We discuss our construction in detail in Appendix \ref{Sec:simplex_construction},
and there see that a honeycomb and its dual share the same \emph{fundamental simplex}.

\subsection{Which three-dimensional Schl\"afli symbols produce honeycombs of which geometries?}

As in two dimensions, there is a formula that tells us whether the honeycomb corresponding to the Schl\"afli symbol $\{p,q,r\}$ tiles $S^3, \EE^3$ or $\HH^3$. Instead of adding up the internal angles of polygons around a vertex, we add up the dihedral angles of polyhedra around an edge. For a euclidean polyhedron given by the Schl\"afli symbol $\{p,q\}$, the dihedral angle is $ 2 \arcsin\left(  \cos\frac{\pi}{q} / \sin\frac{\pi}{p} \right),$ for a total angle of \\
$ 2 r \arcsin\left(  \cos\frac{\pi}{q} / \sin\frac{\pi}{p} \right)$~\cite[p. 157]{Coxeter1954}.  We will be able to tile $\EE^3$ if this value is equal to $2\pi$, or equivalently if 

\begin{equation}
\cos \frac{\pi}{q} = \sin \frac{\pi}{p} \sin \frac{\pi}{r}.
\label{3d_equation}
\end{equation}

If the left hand side is less than the right, then we instead get a four-dimensional polytope, which as before we view as a tiling of $S^3$. If the left hand side is greater than the right, then we can tile $\HH^3$~\cite[p. 157]{Coxeter1954}.  Note that the dual of a honeycomb will tile the same space because swapping $p$ with $r$ doesn't change the equation.

Other than the complexity of drawing the figures, and the complexity of the formula, the story in three dimensions is so far very similar to the story in two dimensions. However, some new features arise in three dimensions, as we will see in the next sections. 

\section{Idealization}
\label{Sec:ideal}
\subsection{Ideal vertices}
\label{Sec:idealvertices}

\begin{figure}[htbp]
\centering 
\subfloat[$\{4,6\}$]
{
\ifimages
\includegraphics[width=0.23\textwidth]{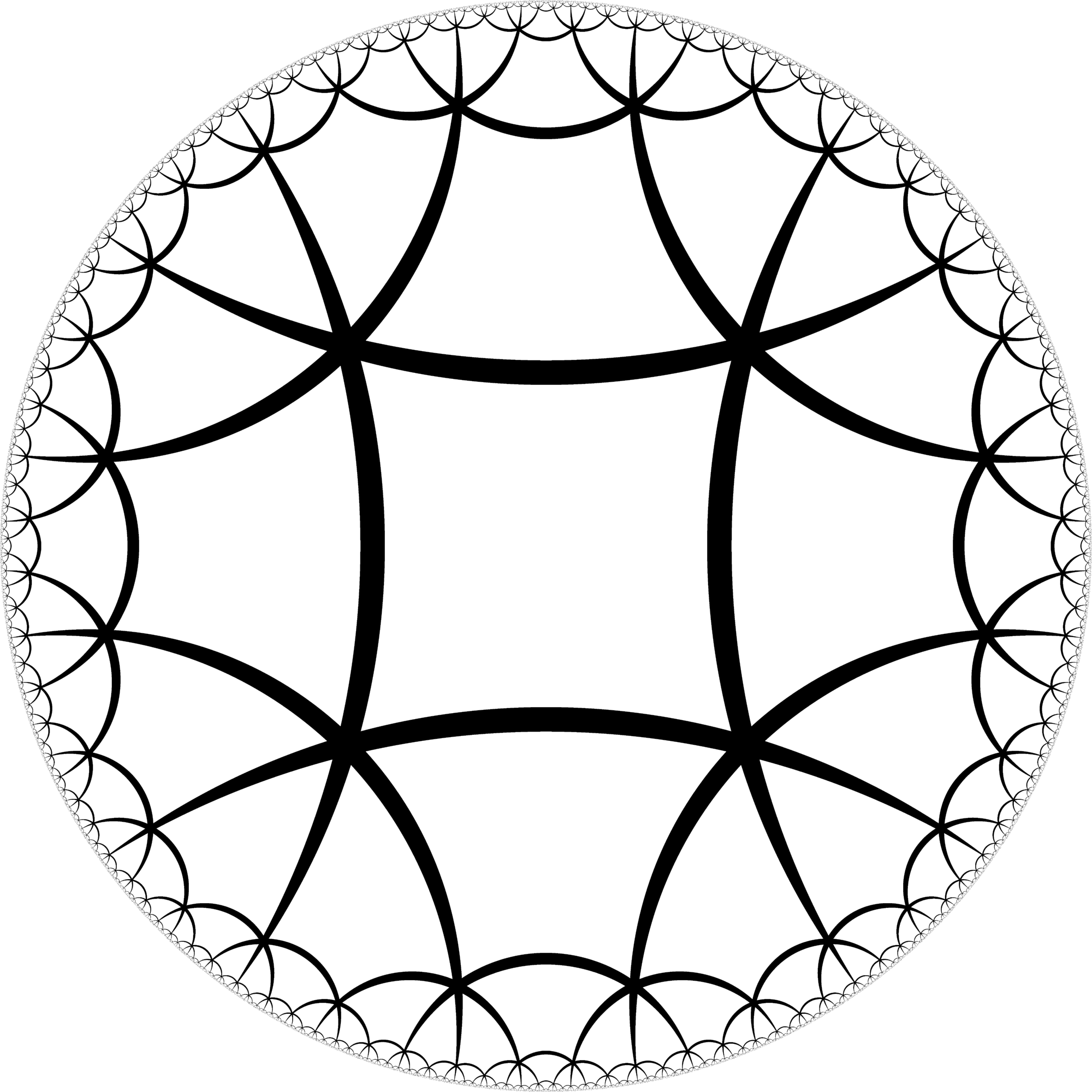}
\fi
\label{Fig:46}
}
\subfloat[$\{4,10\}$]
{
\ifimages
\includegraphics[width=0.23\textwidth]{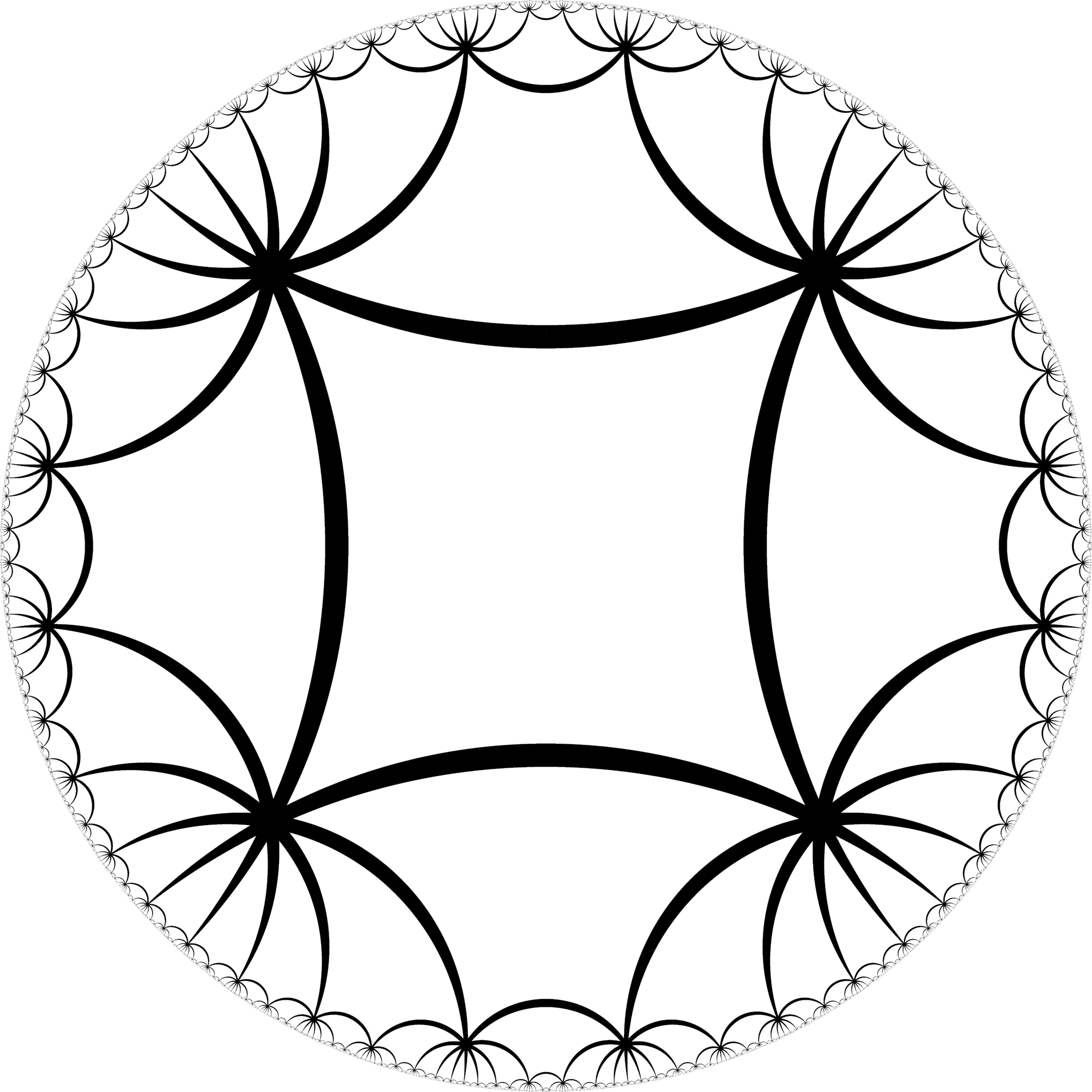}
\fi
\label{Fig:4_10}
}
\subfloat[$\{4,20\}$]
{
\ifimages
\includegraphics[width=0.23\textwidth]{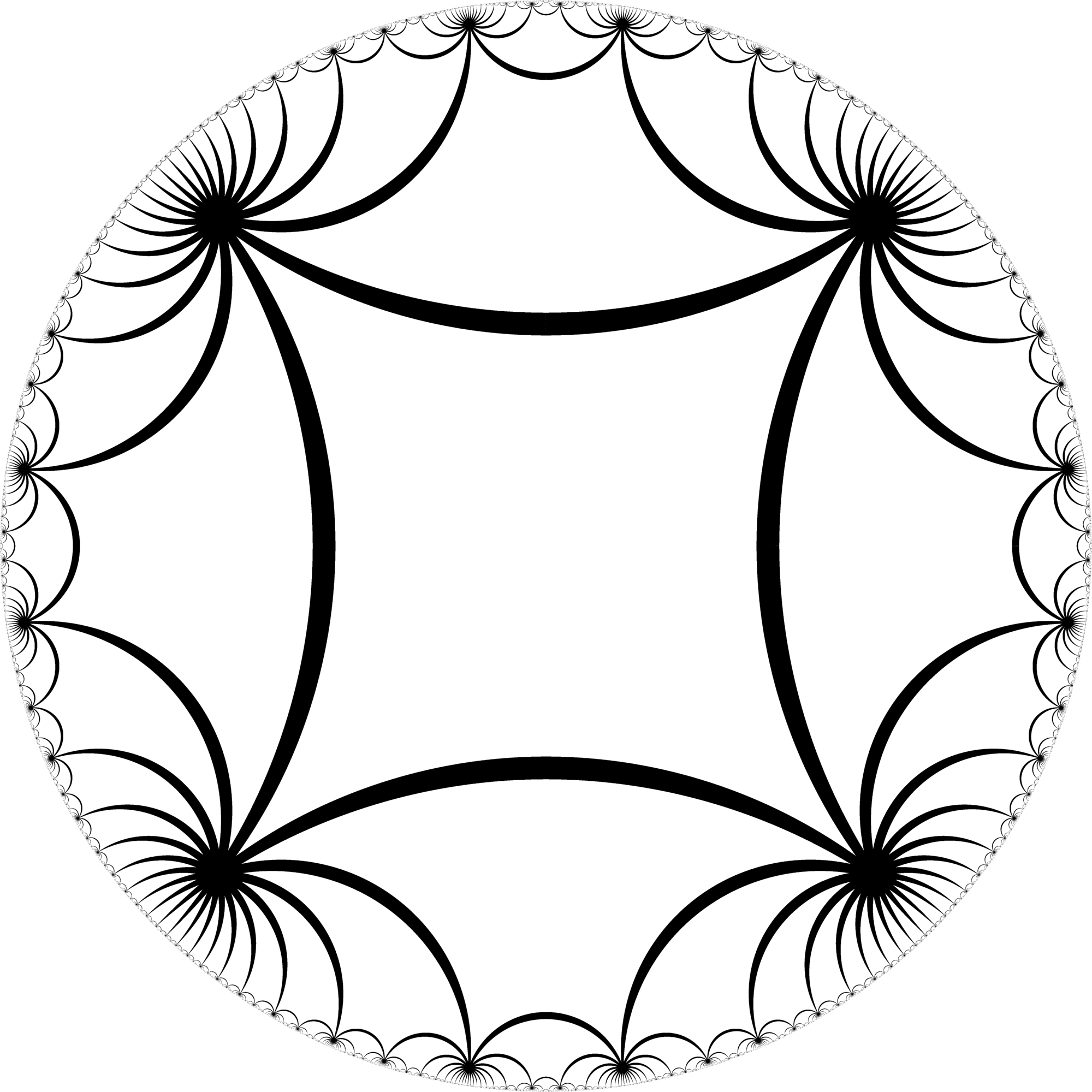}
\fi
\label{Fig:4_20}
}
\subfloat[$\{4,\infty\}$]
{
\ifimages
\includegraphics[width=0.23\textwidth]{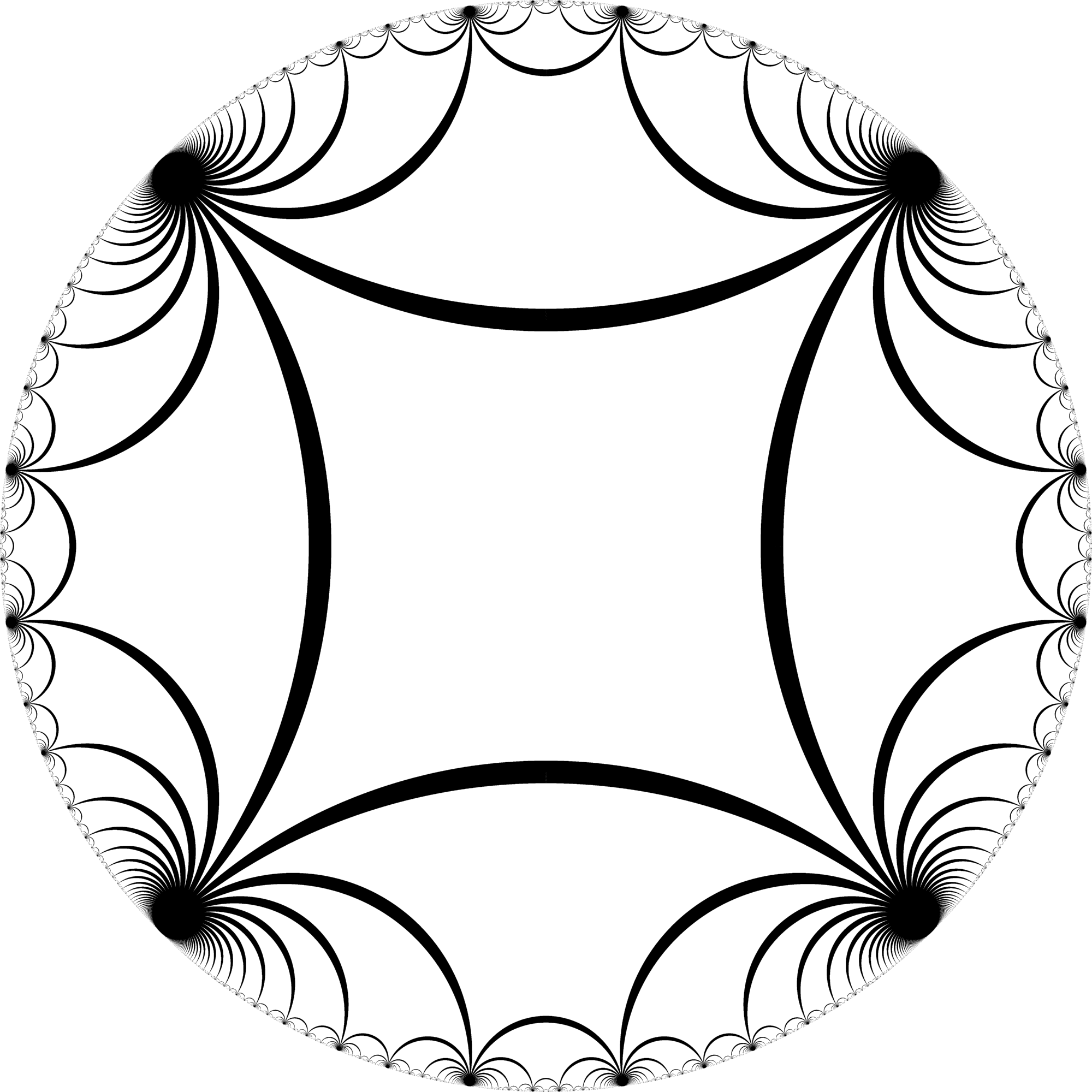}
\fi
\label{Fig:4i}
}
\caption{Tilings given by Schl\"afli symbols $\{4,q\}$.}
\label{Fig:4q}
\end{figure}

If we continue the sequence of tilings shown in Figure \ref{Fig:43_44_45}, all Schl\"afli symbols $\{4,q\}, q \geq 5$, are tilings of $\HH^2$. See Figure \ref{Fig:4q}. For each such tiling, the internal angle of each square must be $2\pi/q$. So as $q$ increases, the square gets spikier, and the vertices move farther and farther away from the center of the square. In the limit, the Schl\"afli symbol $\{4,\infty\}$ corresponds to a tiling with infinitely many squares around each vertex, each of which has internal angle zero, and the vertices are infinitely far away, on the boundary of $\HH^2$. See Figure \ref{Fig:4i}. 

\begin{figure}[htbp]
\centering 
\subfloat[The central cube.]
{
\ifimages
\includegraphics[width=0.3\textwidth]{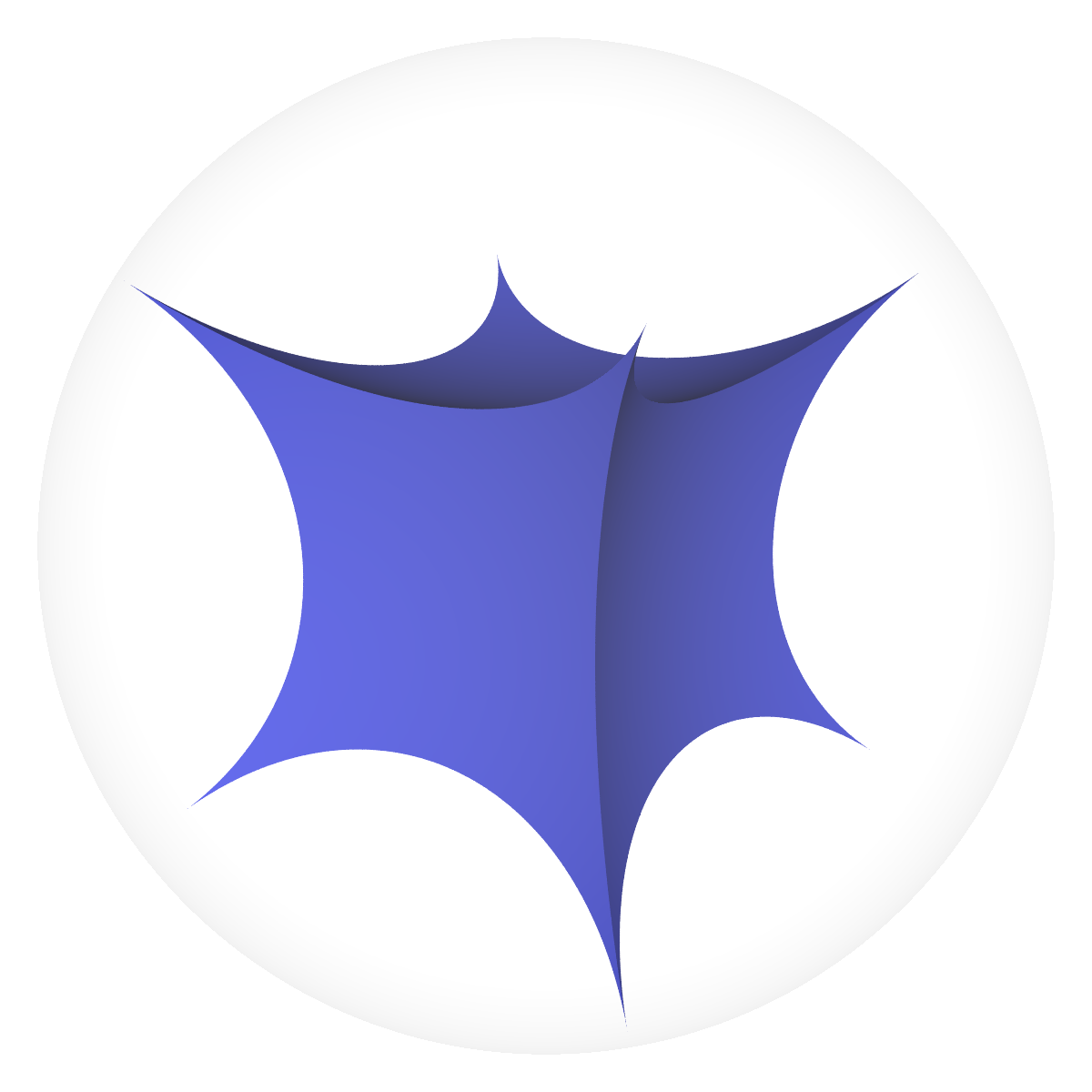}
\fi
\label{Fig:436_1}
}
\subfloat[Six more cubes.]
{
\ifimages
\includegraphics[width=0.3\textwidth]{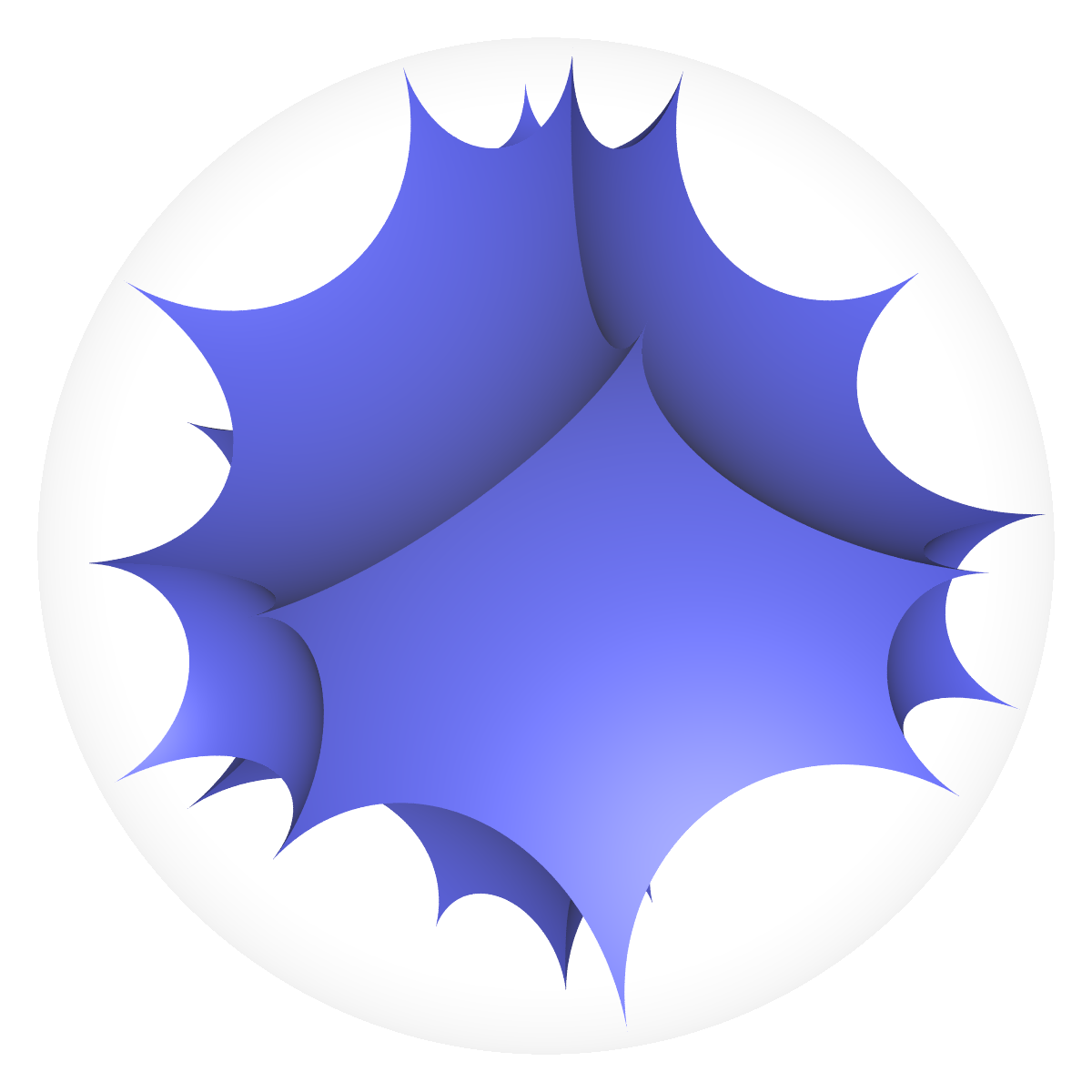}
\fi
\label{Fig:436_2}
}
\subfloat[Thirty more cubes.]
{
\ifimages
\includegraphics[width=0.3\textwidth]{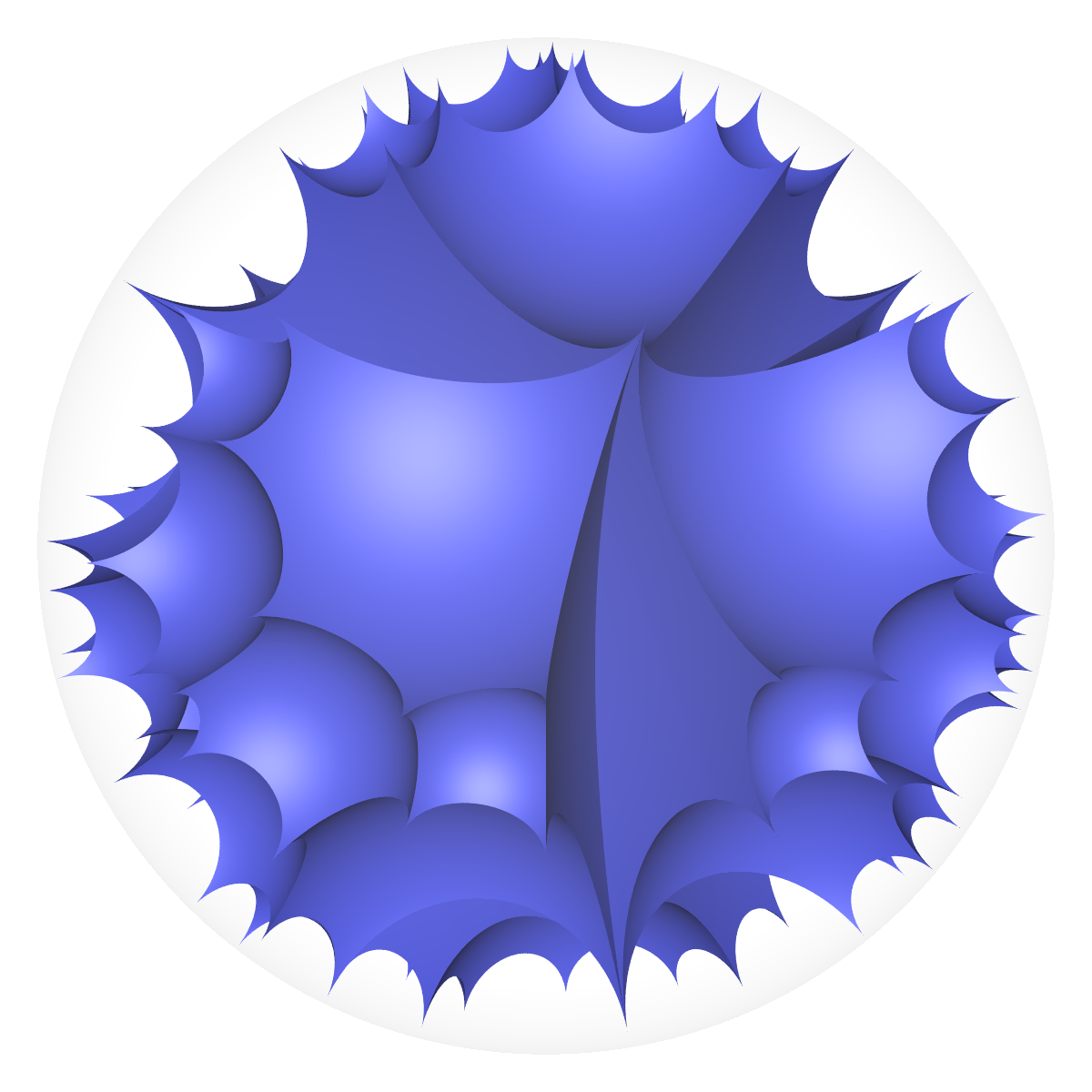}
\fi
\label{Fig:436_3}
}
\caption{Three layers of cubes in the $\{4,3,6\}$ honeycomb.}
\label{Fig:436_build_layers}
\end{figure}

A vertex that is on the boundary of hyperbolic space is called an \emph{ideal vertex}. We call a vertex that is inside hyperbolic space a \emph{material vertex}. In three dimensions, we don't need to go to Schl\"afli symbols with infinite terms to find ideal vertices. In fact, the next Schl\"afli symbol in the sequence shown in Figure \ref{Fig:433_434_435}, $\{4,3,6\}$, corresponds to a honeycomb of $\HH^3$ by cubes with ideal vertices. See Figure \ref{Fig:436_build_layers}. At first, it is a little unclear why $\{4,3,6\}$ should have ideal vertices while $\{4,3,5\}$ does not. This is easiest to see with the dual honeycomb, which we investigate in the next section.  Later, in Section \ref{Sec:hyperidealvertices}, we will see what happens when $r\rightarrow\infty$, which provides further intuition.


\subsection{Ideal tiles}

\begin{figure}[htbp]
\centering 
\subfloat[$\{6,4\}$]
{
\ifimages
\includegraphics[width=0.23\textwidth]{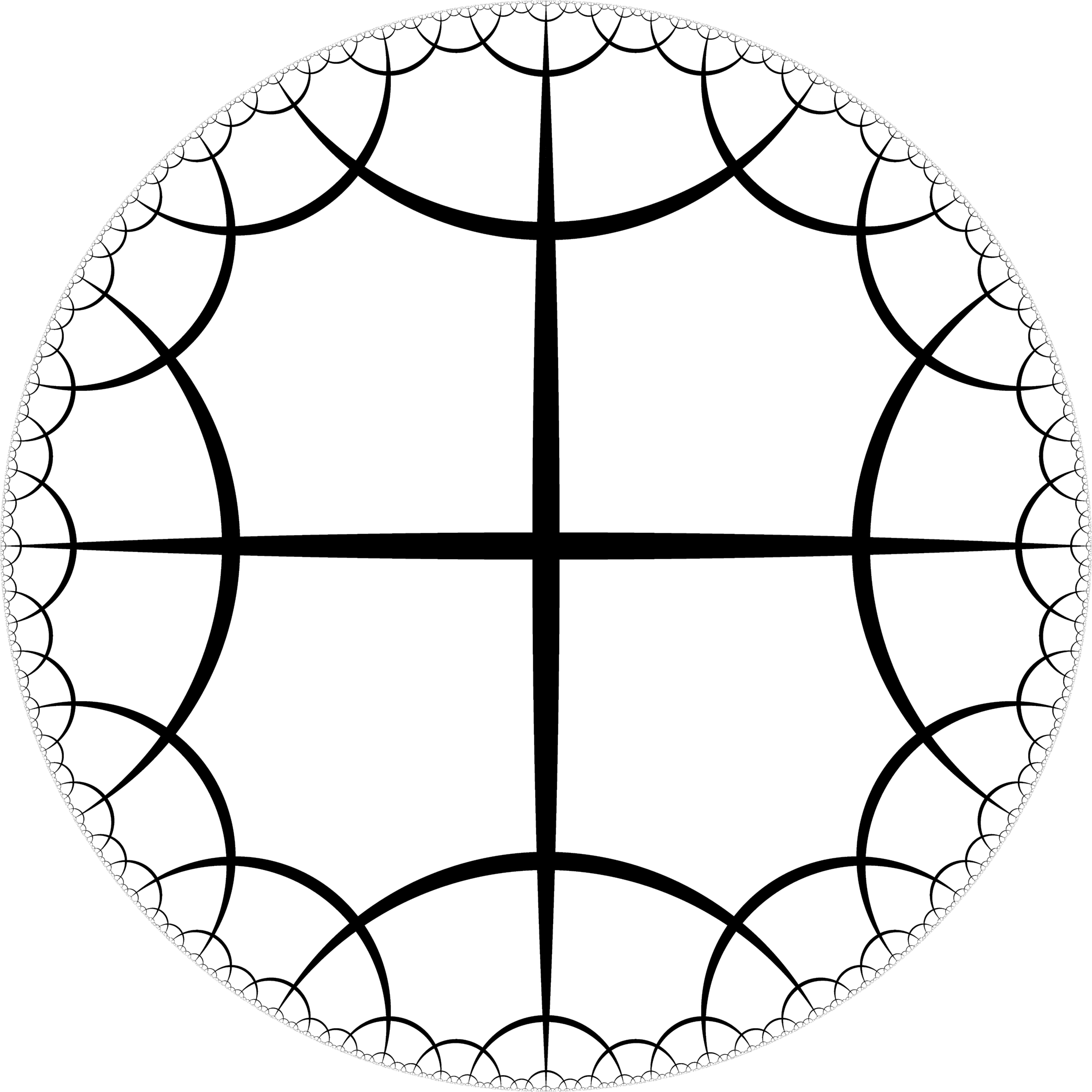}
\fi
\label{Fig:64}
}
\subfloat[$\{10,4\}$]
{
\ifimages
\includegraphics[width=0.23\textwidth]{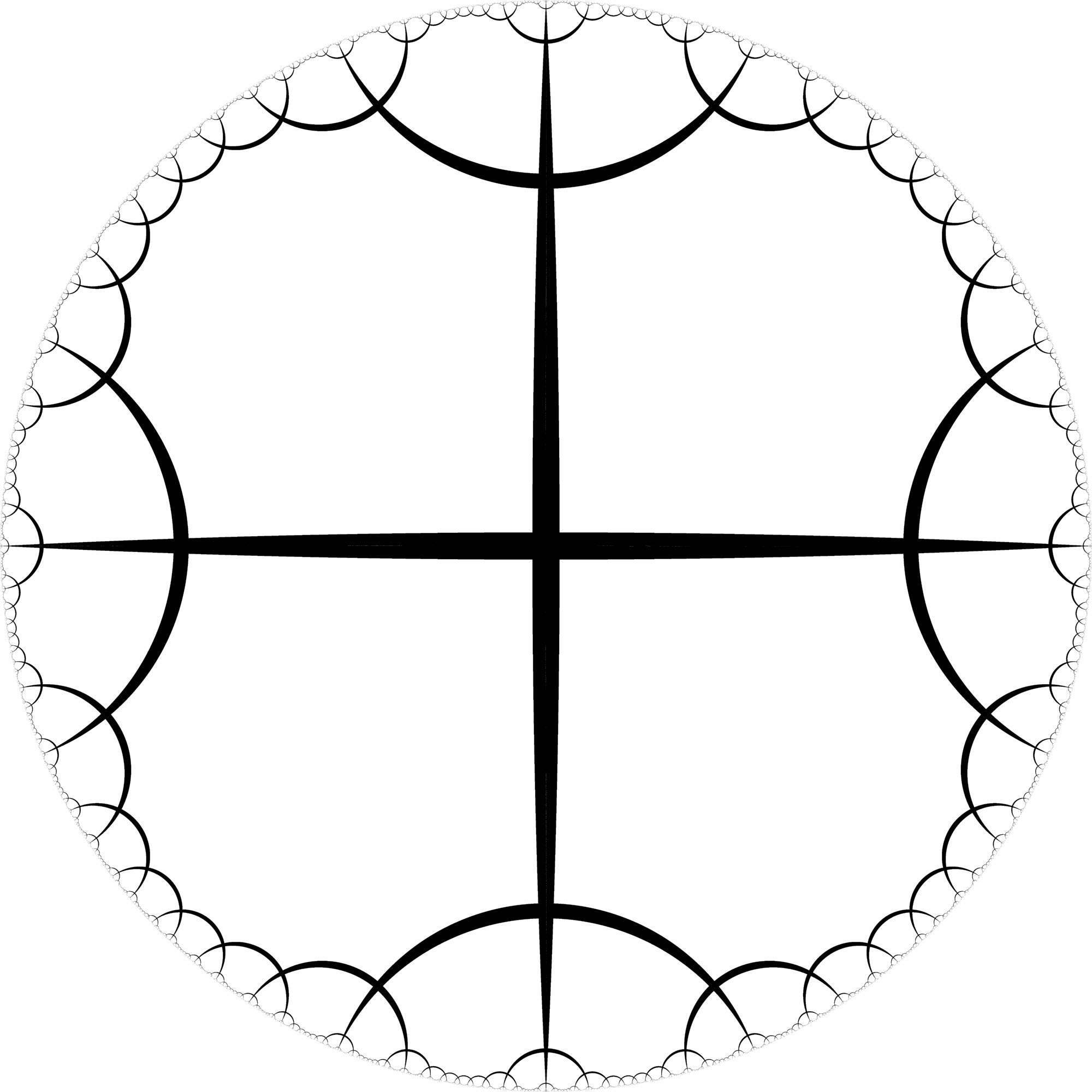}
\fi
\label{Fig:10_4}
}
\subfloat[$\{20,4\}$]
{
\ifimages
\includegraphics[width=0.23\textwidth]{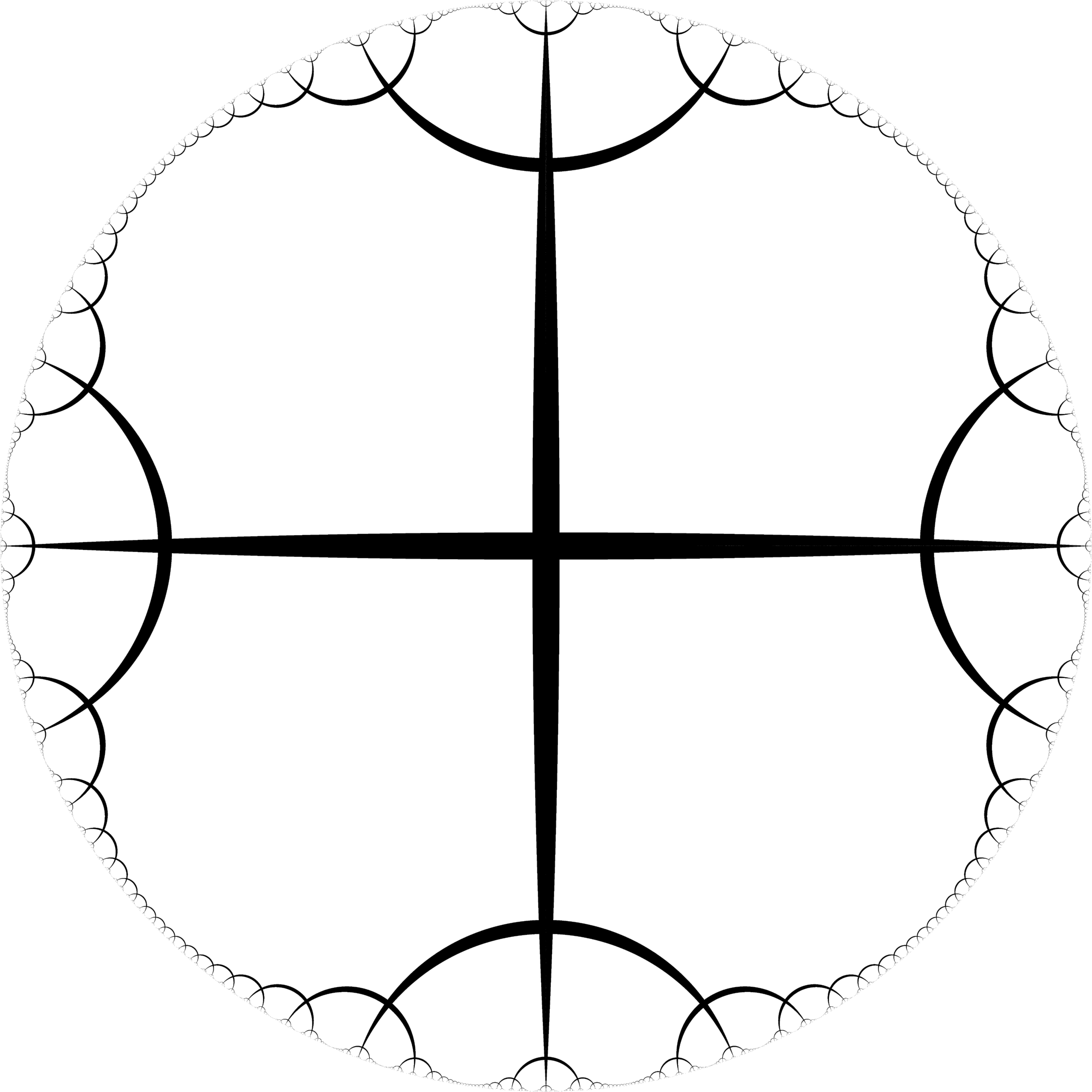}
\fi
\label{Fig:20_4}
}
\subfloat[$\{\infty,4\}$]
{
\ifimages
\includegraphics[width=0.23\textwidth]{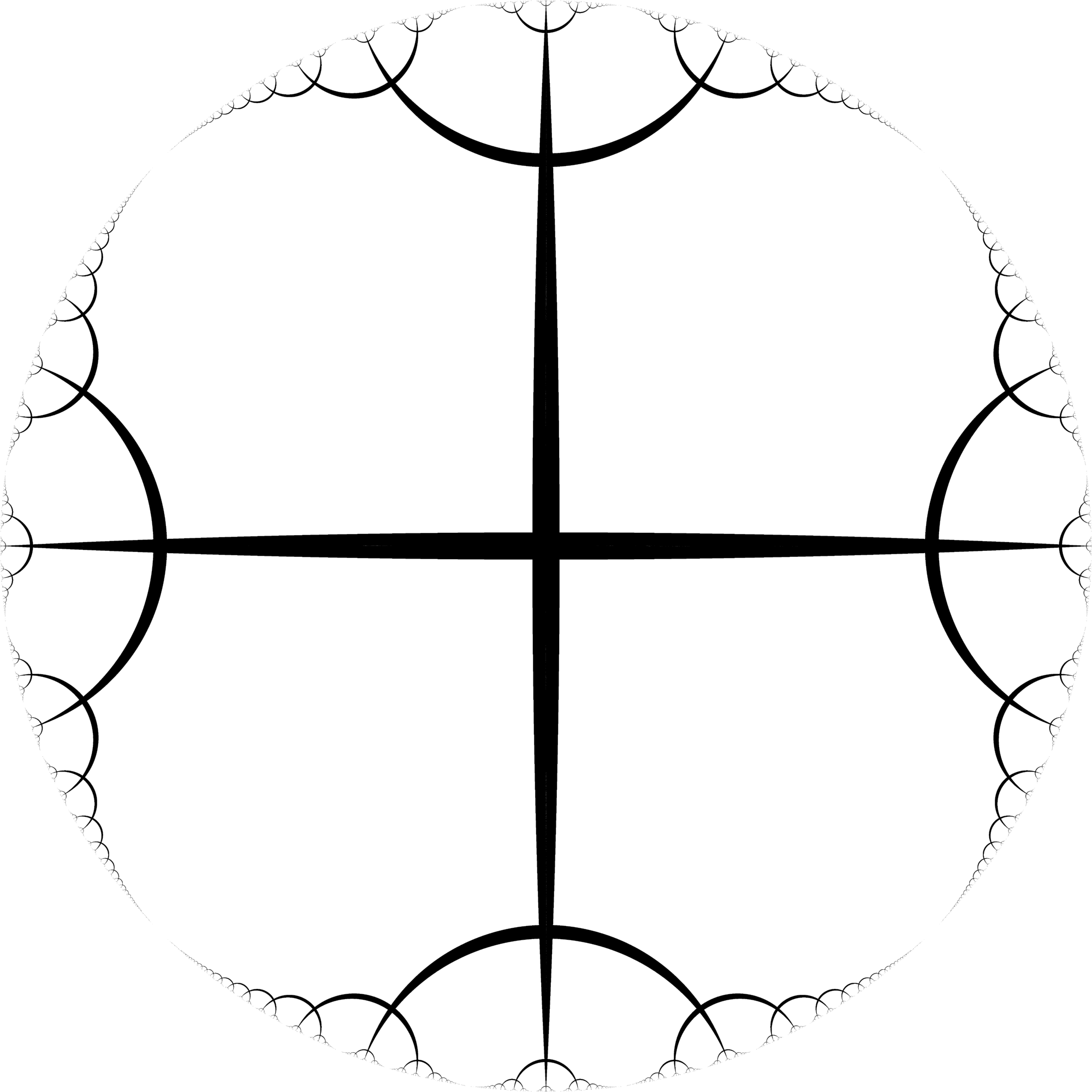}
\fi
\label{Fig:i4}
}
\caption{Tilings given by Schl\"afli symbols $\{p,4\}$.}
\label{Fig:p4}
\end{figure}

If a tiling has ideal vertices, then its dual will have ideal tiles (i.e. ideal faces for two-dimensional tilings, and ideal cells for three-dimensional honeycombs). We can see this in the two-dimensional case by looking at the dual sequence of tilings to the sequence shown in Figure \ref{Fig:4q}. See Figure \ref{Fig:p4}.

\begin{figure}[htbp]
\centering 
\includegraphics[width=0.5\textwidth]{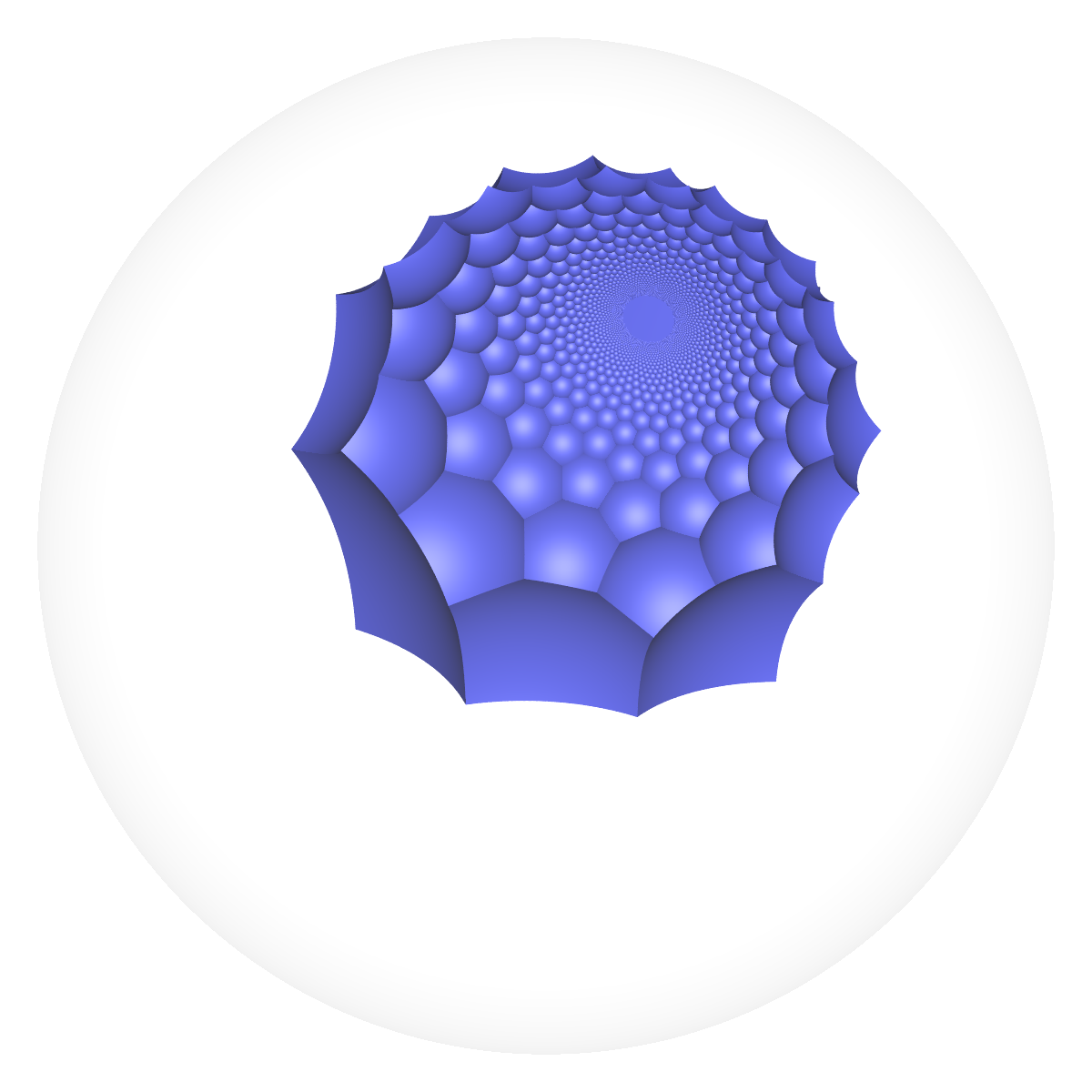}
\caption{A single cell from the honeycomb given by the Schl\"afli symbol $\{6,3,4\}$.}
\label{Fig:634_cell}
\end{figure}

We cannot draw a vertex-centered picture of a tiling with ideal vertices, since the vertices cannot be put at the center of the Poincar\'e disk. Similarly, it is impossible to draw a tile-centered picture of a tiling with ideal tiles. Figure \ref{Fig:634_cell} shows one ideal cell from the $\{6,3,4\}$ honeycomb. The cell has infinitely many hexagonal faces, as you would expect since its Schl\"afli symbol is $\{6,3\}$. Geometrically, the vertices of the cell sit on a \emph{horosphere}~\cite[p. 326]{Needham199902} of $\HH^3$. The Poincar\'e ball model represents $\HH^3$ as a euclidean ball in $\EE^3$. A horosphere is represented in the model as a euclidean sphere which is contained in and touches the boundary of the Poincar\'e ball. As viewed in the negatively curved metric in $\HH^3$ however, a horosphere is isometric to a copy of the euclidean plane $\EE^2$. This makes sense, since $\{6,3\}$ is a tiling of the euclidean plane. 

In general, an \emph{ideal cell} is a three-dimensional cell whose vertices lie on a horosphere. We call a cell whose vertices lie on a sphere that is strictly inside of $\HH^3$ a \emph{material cell}.

From this perspective we can see why $\{4,3,6\}$ has ideal vertices: the dual honeycomb has cells with Schl\"afli symbol $\{6,3\}$. Since this is a tiling of the euclidean plane rather than a sphere, the cell vertices must live on a horosphere. The original honeycomb $\{4,3,6\}$ has a vertex at the center of this dual cell, and the center of a horosphere is on the boundary of $\HH^3$. 

Note that only the vertices of an ideal cell live on the horosphere, not the edges or faces. In hyperbolic space, the euclidean geometry of the horosphere curves away from the hyperbolic geodesics of the edges and faces. 

\section{The three-dimensional Schl\"afli symbol map} 
\label{Sec:3d_map}

\begin{figure}[htbp]
\centering 

\labellist
\small\hair 2pt
\pinlabel 3 at 19 78
\pinlabel 4 at 89 57
\pinlabel 5 at 165 36
\pinlabel 6 at 246 13
\pinlabel 7 at 333 -9
\pinlabel $p$ at 160 20

\pinlabel 3 at 346 0
\pinlabel 4 at 380 43
\pinlabel 5 at 411 80
\pinlabel 6 at 435 115
\pinlabel 7 at 458 145
\pinlabel $q$ at 425 68

\pinlabel 3 at 11 90
\pinlabel 4 at 7 167
\pinlabel 5 at 3.5 248
\pinlabel 6 at -1 334
\pinlabel 7 at -6 424
\pinlabel $r$ at -10 240

\pinlabel $B$ at 240 381
\pinlabel $C$ at 332 263
\pinlabel $A$ at 130 84

\endlabellist
\includegraphics[width=0.8\textwidth]{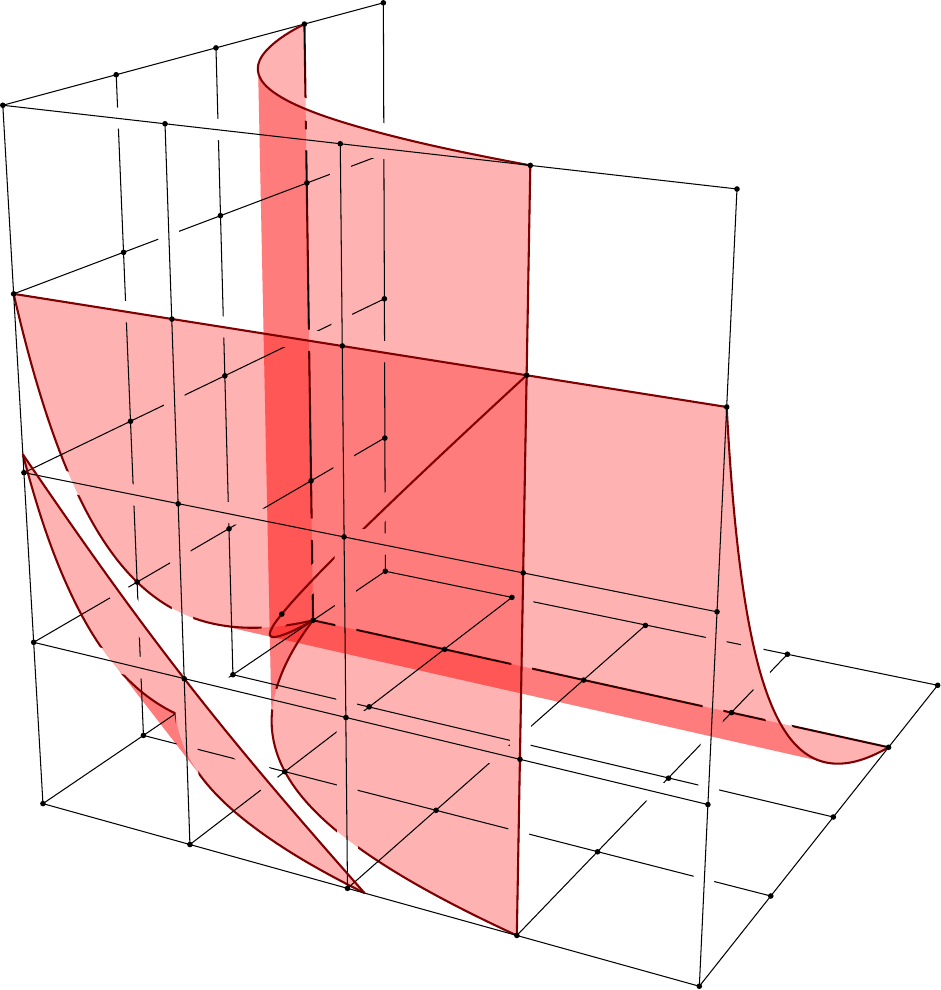}

\caption{Geometries and cell/vertex types associated to the three-dimensional Schl\"afli symbols. Surface $A$ is given by $\cos\frac{\pi}{q}=\sin\frac{\pi}{p}\sin\frac{\pi}{r}$, surface $B$ by $(p-2)(q-2)=4$, and surface $C$ by $(r-2)(q-2)=4$.}
\label{Fig:3D_Schlafli_map}
\end{figure}
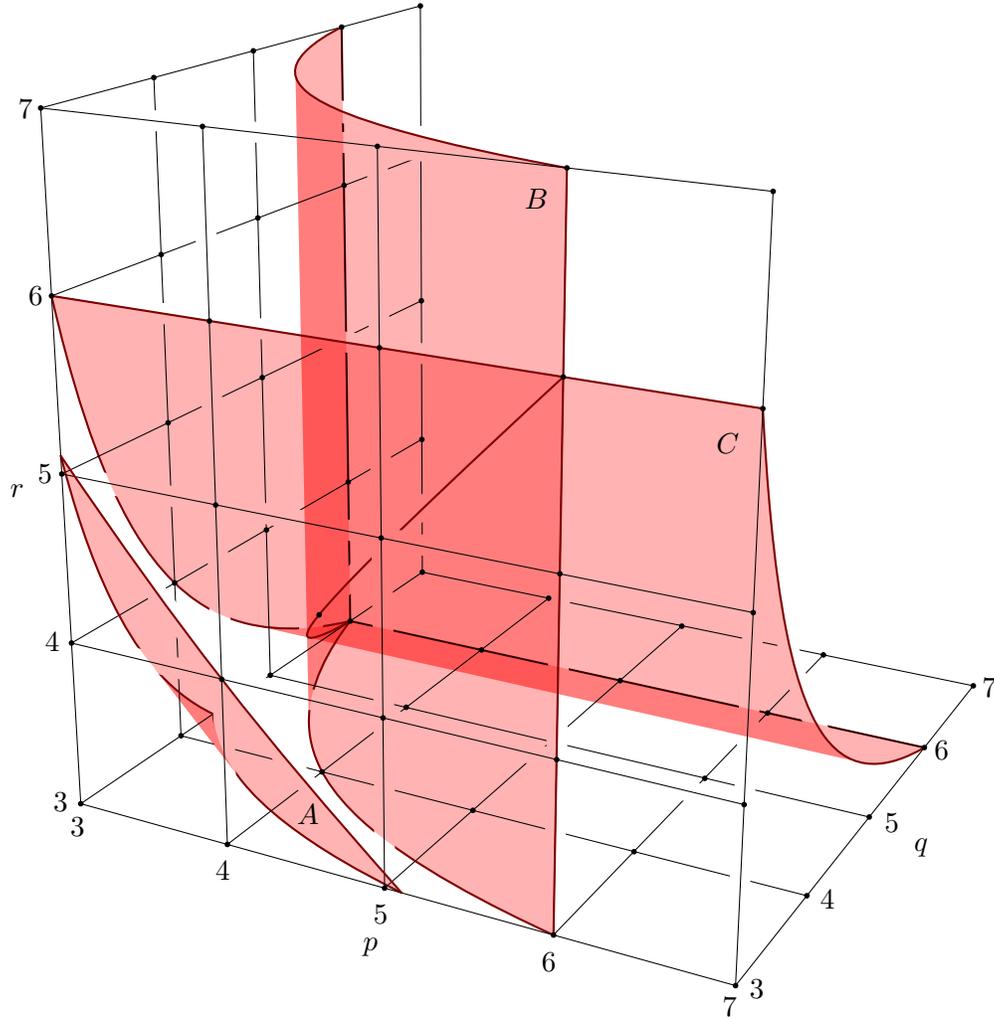

We can now describe the map of the kinds of geometry, vertices and cells we see for three-dimensional Schl\"afli symbols. See Figure \ref{Fig:3D_Schlafli_map}. Three surfaces are shown. In the bottom left is the surface associated to Equation \ref{3d_equation}. Below the surface the honeycomb is spherical, on the surface the honeycomb is euclidean, and above the surface the honeycomb is hyperbolic. The other two surfaces are associated to Equation \ref{2d_equation}, the first using $p$ and $q$ from the three-dimensional Schl\"afli symbol $\{p,q,r\}$, and the second using $r$ and $q$. Honeycombs on the first surface have ideal cells, while those closer to the origin have material cells. Honeycombs on the second surface have ideal vertices, while those closer to the origin have material vertices.  The three honeycombs lying on both of these surfaces have ideal cells and ideal vertices (Section \ref{Sec:ideal_verts_and_cells}).

In the next section, we explore what happens when we move to the far side of these two surfaces.

\section{Hyperidealization}
\label{Sec:hyperidealization}
\subsection{Hyperideal vertices}
\label{Sec:hyperidealvertices} 

\begin{figure}[htbp]
\centering 
\subfloat[Material, ideal, and hyperideal squares in the Poincar\'e ball model.]
{
\includegraphics[width=0.35\textwidth]{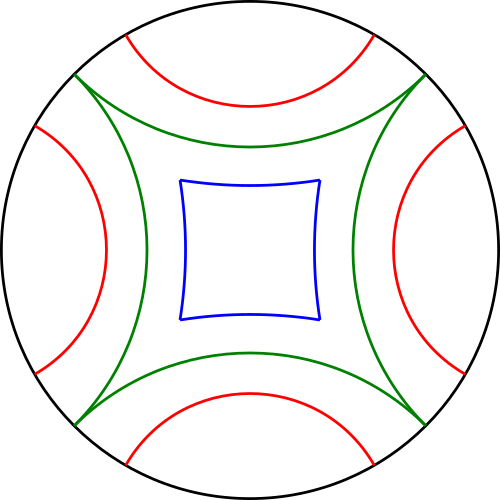}
\label{Fig:hyperideal_ball}
}
\hspace{1.0cm}
\subfloat[The same squares in the Klein model.]
{
\includegraphics[width=0.35\textwidth]{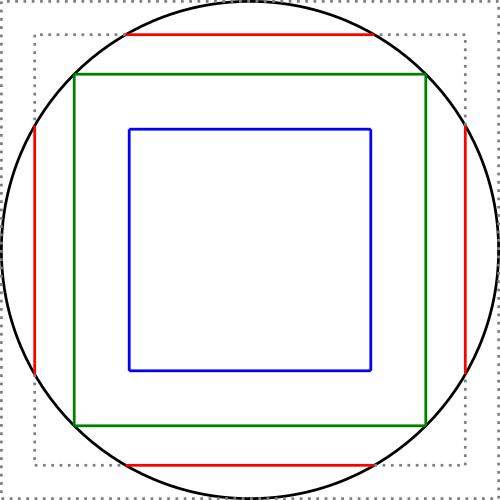}
\label{Fig:hyperideal_klein}
}
\caption{Hyperidealization.  Hyperideal vertices of the squares can be seen in the Klein model.}
\label{Fig:hyperideal}
\end{figure}

Let's look at what happens when we make the $\{4,\infty\}$ squares of Figure \ref{Fig:4q} even bigger.  Figure \ref{Fig:hyperideal} shows four squares in both the Poincar\'e and Klein models~\cite[pp. 322-323]{Needham199902}, one with material vertices, one with ideal vertices, and two with \emph{hyperideal} vertices, that is vertices living beyond the boundary of $\HH^2$. (Hyperideal points can be defined in terms of their dual relationship to planes within hyperbolic space~\cite[p. 71]{Levy199701}.) The largest hyperideal square here has no material points and only four ideal points on the boundary of $\HH^2$. We will shortly see similar degenerate edges in some of our three-dimensional honeycombs.  

\begin{figure}[htbp]
\centering 
\subfloat[The central cube.]
{
\ifimages
\includegraphics[width=0.3\textwidth]{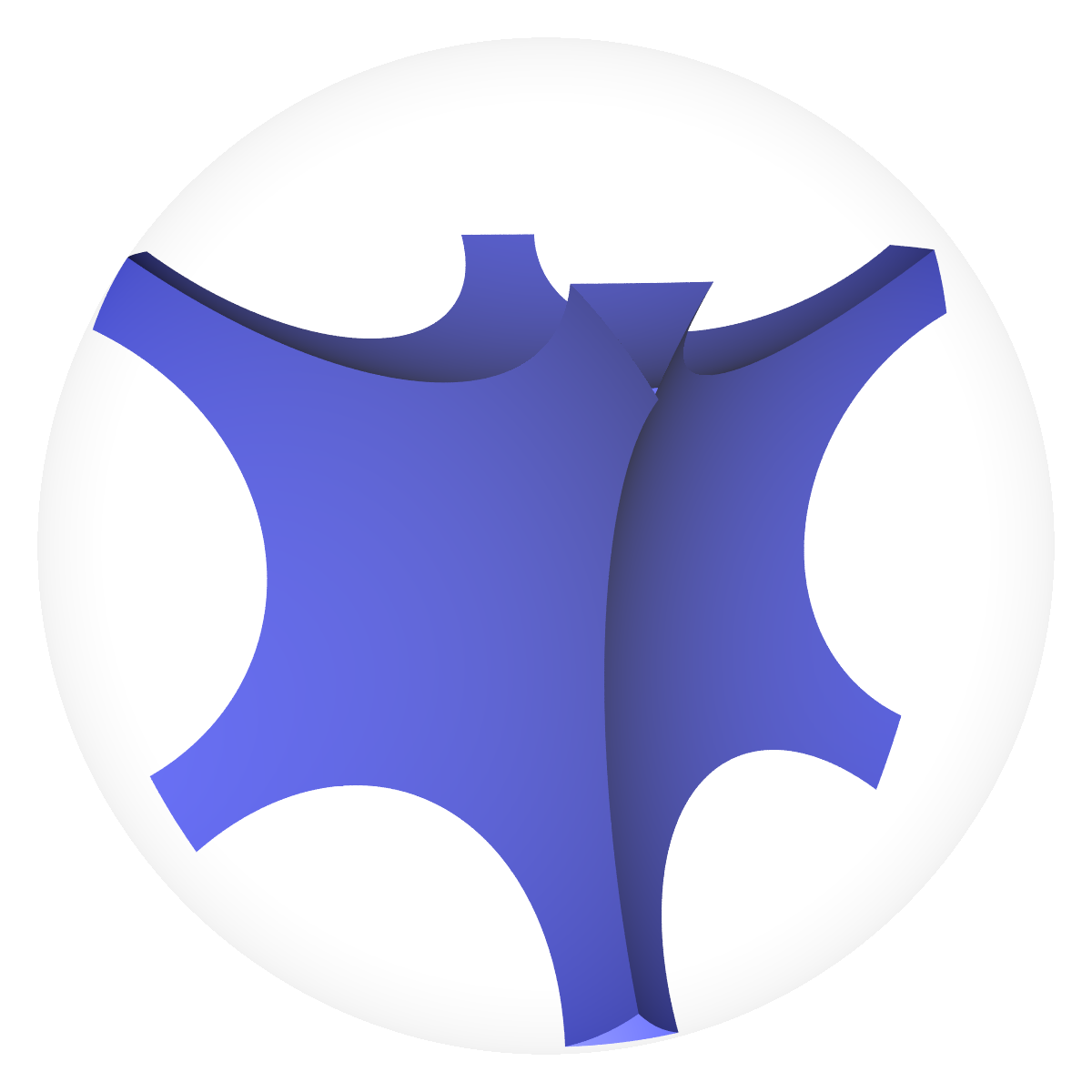}
\fi
\label{Fig:437_1}
}
\subfloat[Six more cubes.]
{
\ifimages
\includegraphics[width=0.3\textwidth]{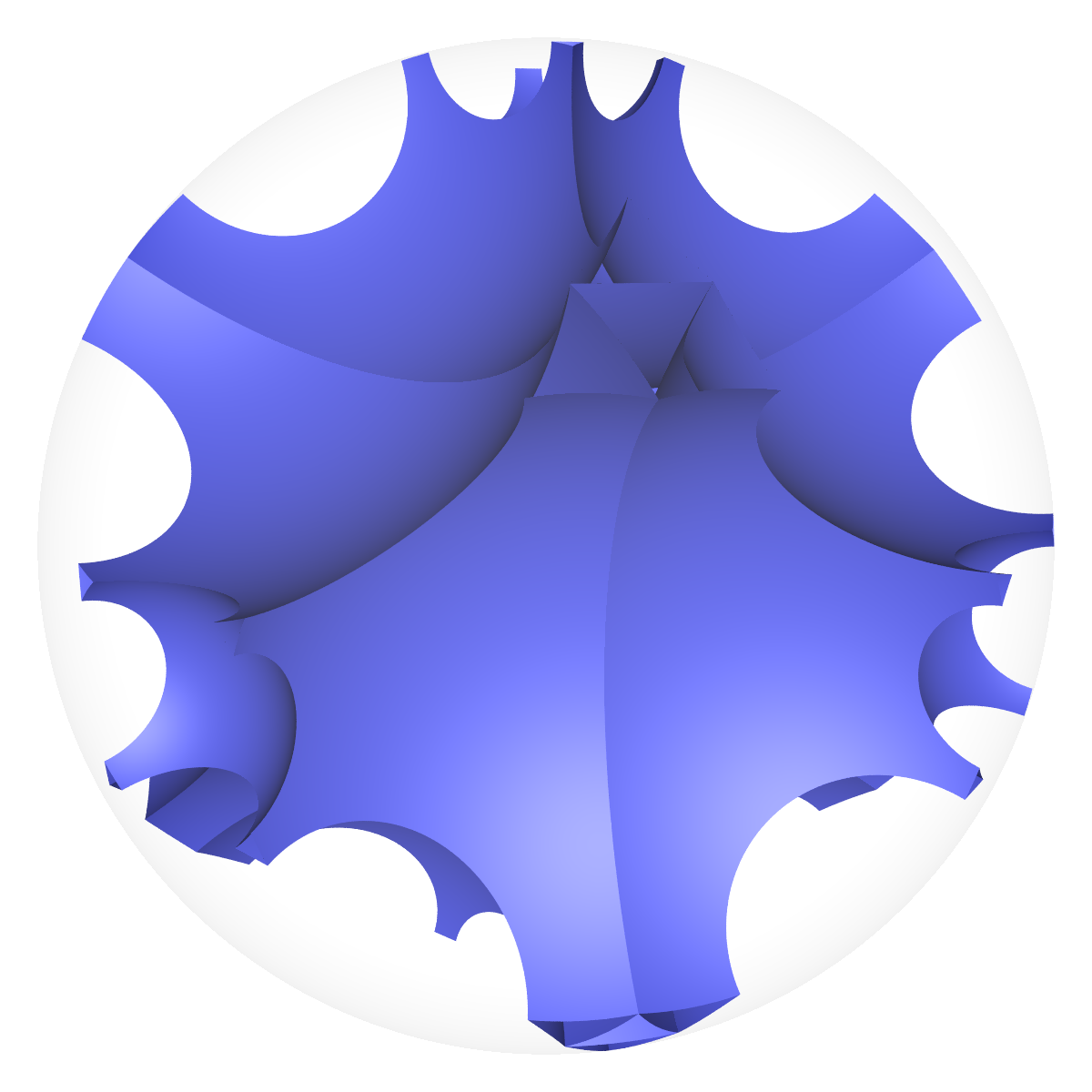}
\fi
\label{Fig:437_2}
}
\subfloat[Thirty more cubes.]
{
\ifimages
\includegraphics[width=0.3\textwidth]{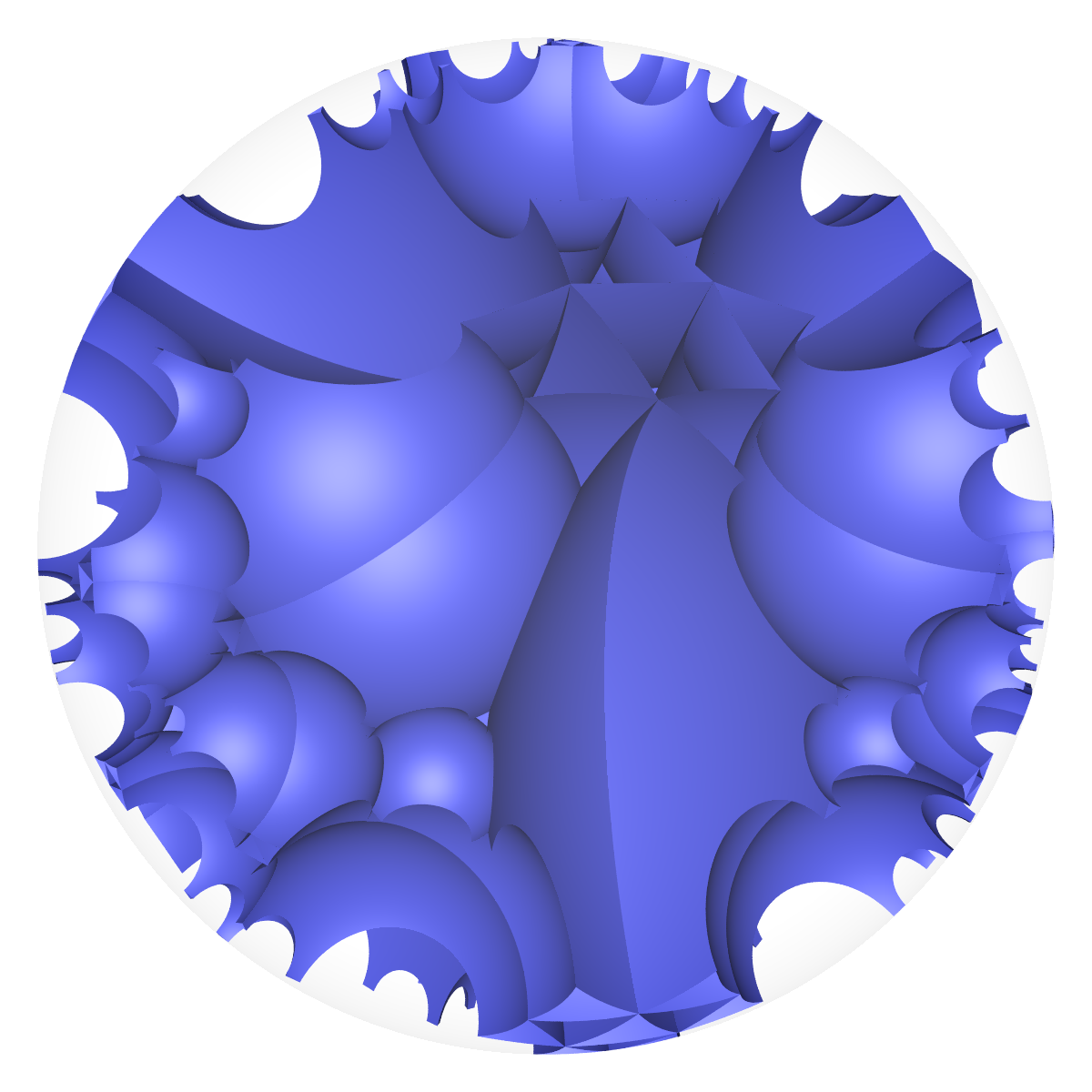}
\fi
\label{Fig:437_3}
}
\caption{Three layers of cubes in the $\{4,3,7\}$ honeycomb. In these images, the corners of the cubes are cut off by the boundary of $\HH^3$, and you can see into the interiors of the cubes through the resulting holes.}
\label{Fig:437_build_layers}
\end{figure}

Continuing our $\{4,3,r\}$ sequence in three dimensions, the $\{4,3,7\}$ honeycomb has cubical cells with hyperideal vertices. As in two dimensions, the vertices no longer exist, even on the boundary of $\HH^3$. Instead the cell has eight \emph{legs} that end in triangular patches of the boundary of $\HH^3$. See Figure \ref{Fig:437_build_layers}. In general, following Equation \ref{2d_equation}, a honeycomb $\{p,q,r\}$ has hyperideal vertices when $(r-2)(q-2)>4$. 

\begin{figure}[htbp]
\centering 
\subfloat[One cube with the induced pattern on the boundary of $\HH^3$.]
{
\includegraphics[width=0.4\textwidth]{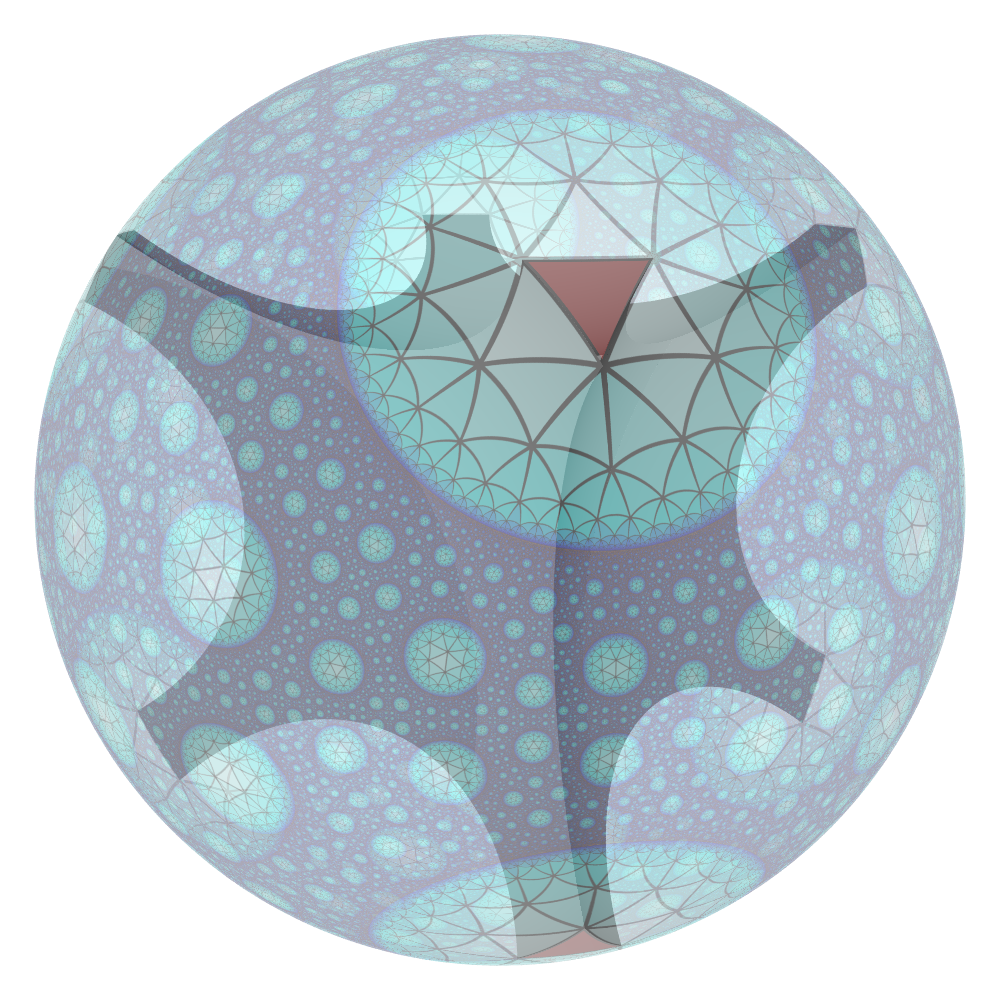}
\label{Fig:437_ball}
}
\hspace{1.0cm}
\subfloat[Stereographic projection of the boundary pattern.]
{
\includegraphics[width=0.4\textwidth]{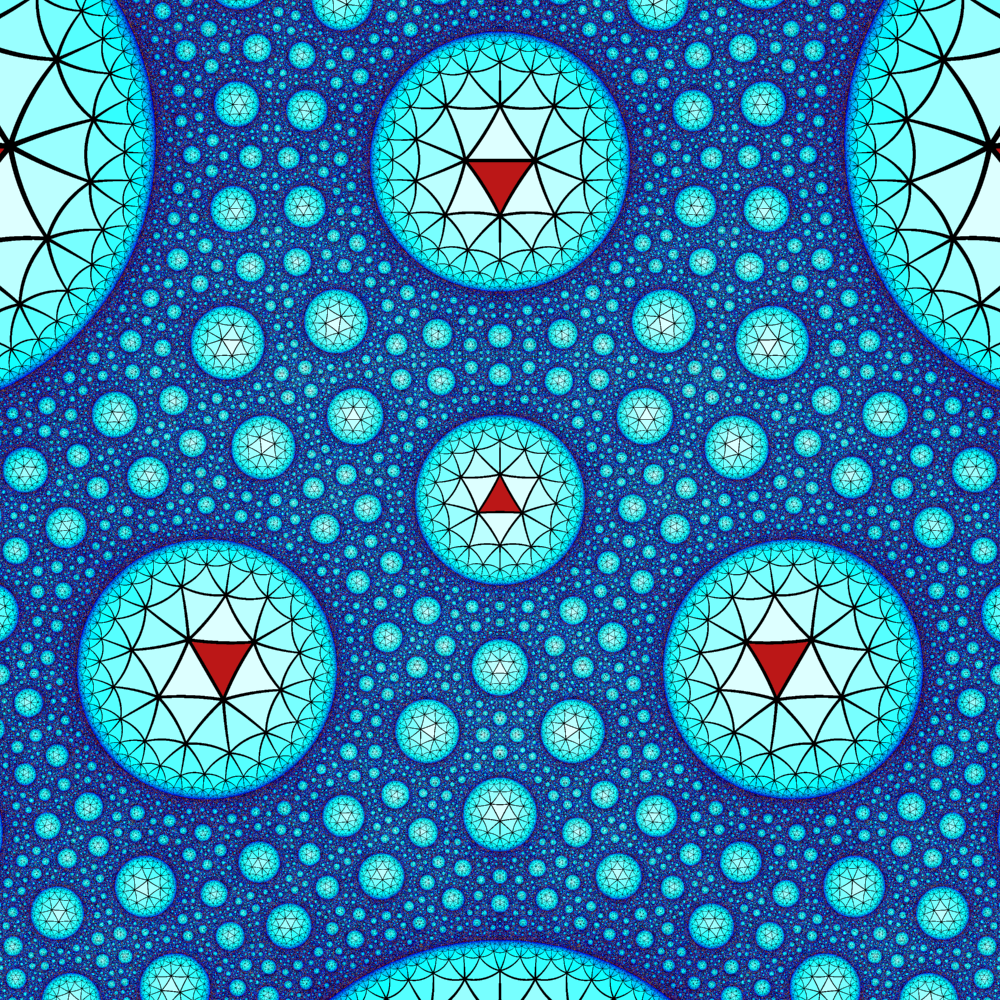}
\label{Fig:437_boundary}
}
\caption{The boundary pattern for the $\{4,3,7\}$ honeycomb. The eight \emph{legs} of one cell are highlighted red.} 
\label{Fig:437_boundary_pattern}
\end{figure}

When our honeycomb has hyperideal vertices, the geometry inside $\HH^3$ becomes less interesting in a sense, because the edges never meet each other. It becomes harder to understand what is going on by looking only at the edges going through the space, as we did in Figure \ref{Fig:433_434_435}. Instead, we can look at the pattern we get on the boundary of $\HH^3$ from the hyperideal cells.  Figure \ref{Fig:437_ball} shows one cube of the $\{4,3,7\}$ honeycomb. It intersects the boundary in eight triangles at the ends of its legs, which are drawn on the surface of the sphere in red. Also drawn are the triangular intersections with the boundary of all other cubes in the honeycomb. Figure \ref{Fig:437_boundary} shows this boundary pattern, stereographically projected from the sphere to the plane, with the projection point at the center of a leg. Equivalently, we can think of this as being the boundary pattern of the $\{4,3,7\}$ honeycomb as drawn in the upper half space model of $\HH^3$. The triangular ends of the central cube are again shaded red, with progressively darker shades of blue for triangular ends of cubes farther away from the central cube. Roughly one hemisphere is visible, and so only four of the eight legs are visible in this image.  For more details on the colouring scheme used here, see Appendix \ref{Sec:colouring}.

\begin{figure}[htbp]
\centering 

\subfloat[$\{4,3,10\}$]
{
\includegraphics[width=0.23\textwidth]{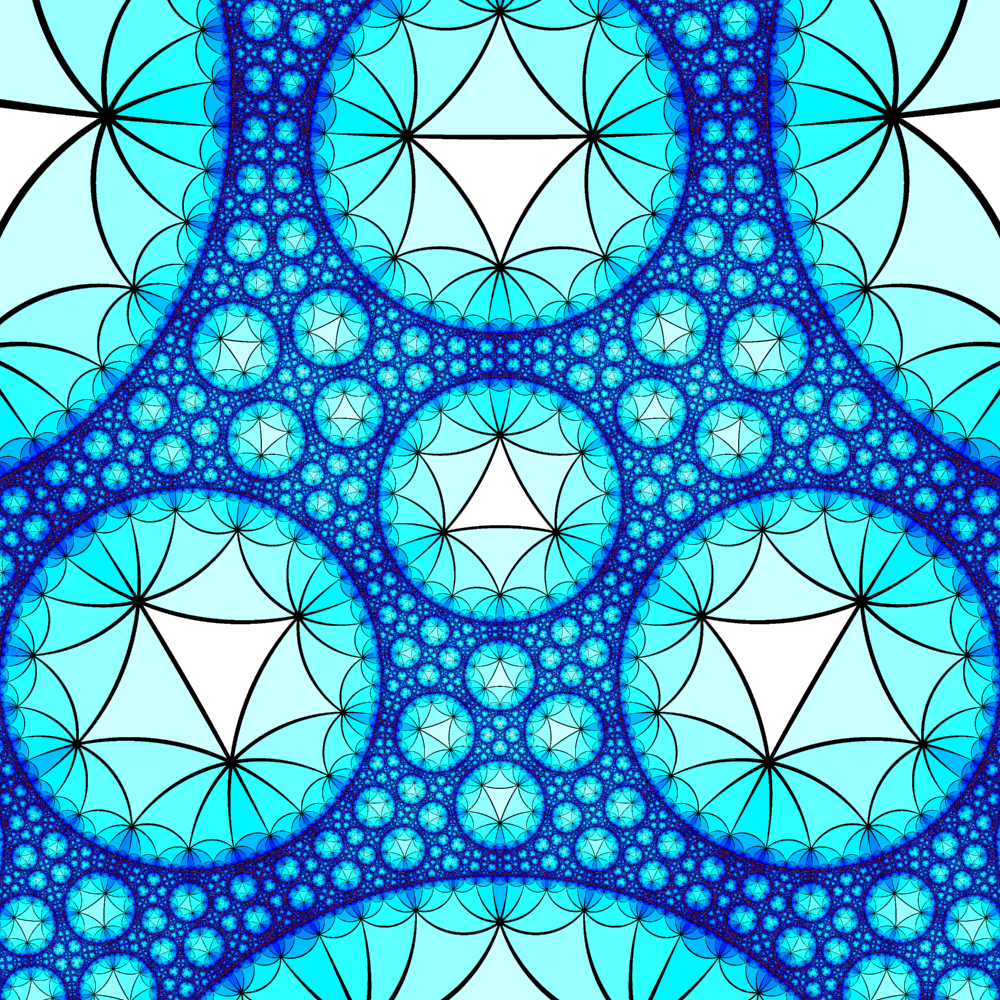}
\label{Fig:43_10}
}
\subfloat[$\{4,3,15\}$]
{
\includegraphics[width=0.23\textwidth]{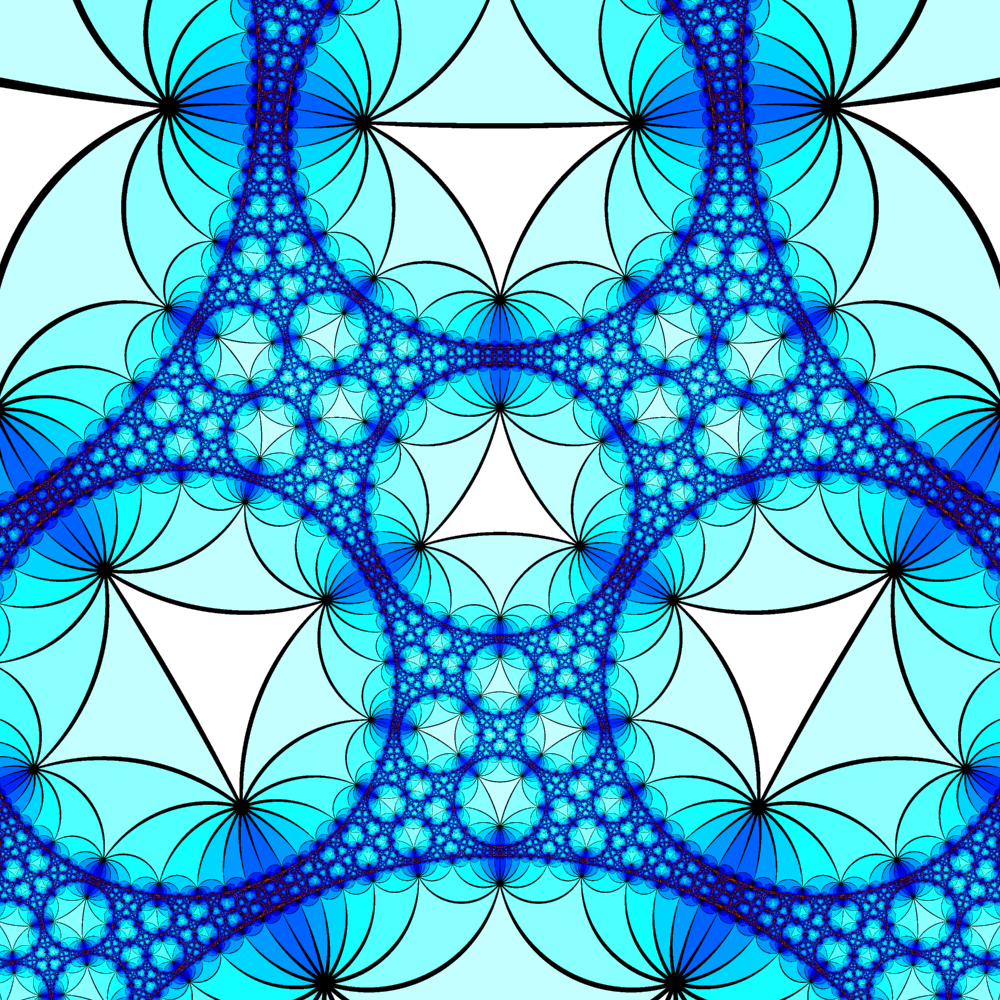}
\label{Fig:43_15}
}
\subfloat[$\{4,3,30\}$]
{
\includegraphics[width=0.23\textwidth]{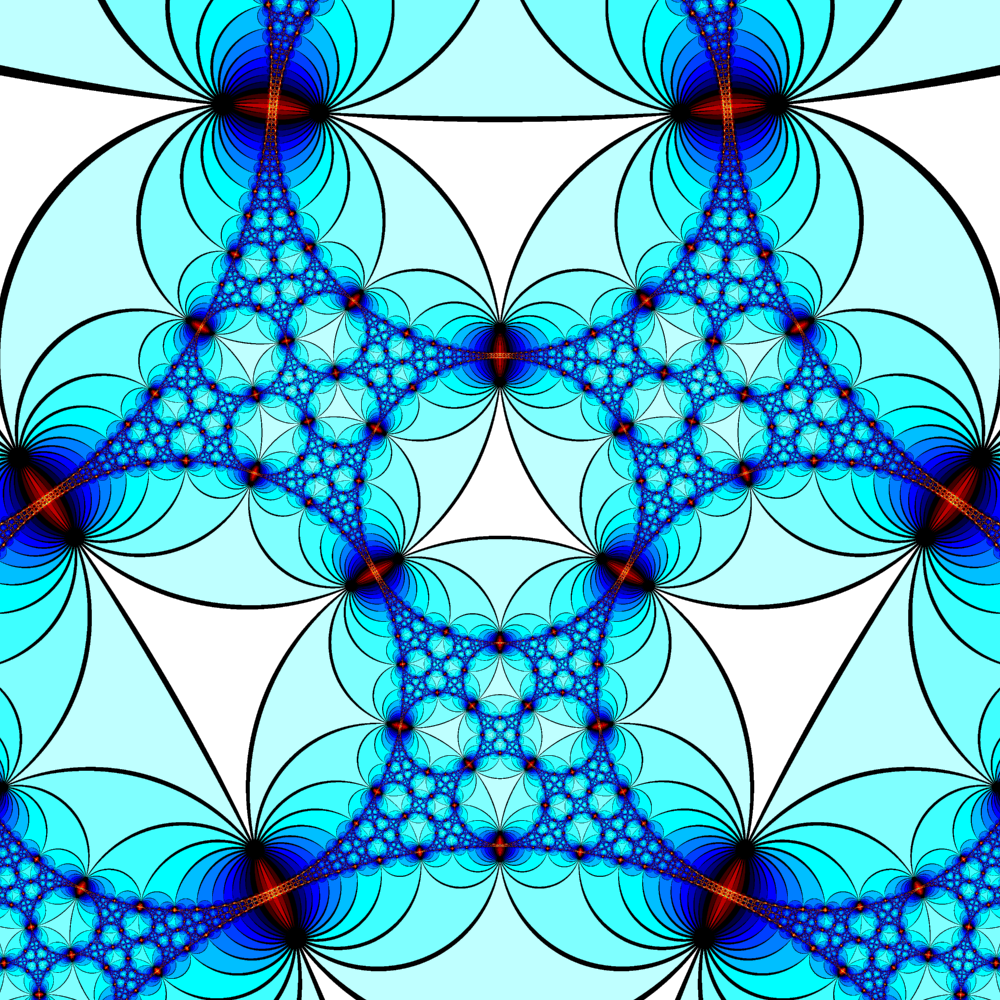}
\label{Fig:43_30}
}
\subfloat[$\{4,3,\infty\}$]
{
\includegraphics[width=0.23\textwidth]{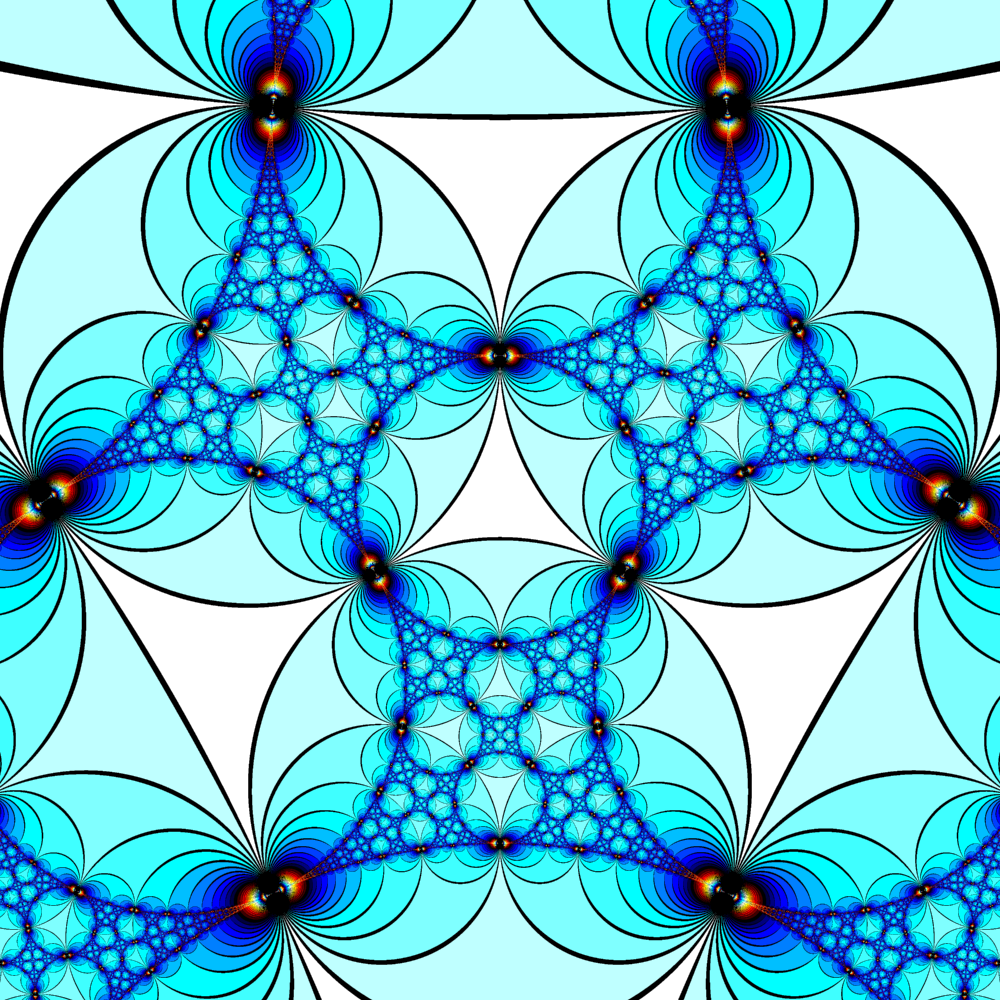}
\label{Fig:43i}
}

\caption{$\{4,3,r\}$ honeycombs as $r\rightarrow\infty$. Compare with the $\{4,3,7\}$ honeycomb in Figure \ref{Fig:437_boundary}. 
In the limit, the edges become hyperideal except for a single ideal point, located where the apparent $\{3,\infty\}$ tilings kiss.}
\label{Fig:43r_progression}
\end{figure}

Edges with ideal or hyperideal vertices usually intersect the boundary at two ideal points.  However, as $r\rightarrow\infty$, less and less of each edge intersects $\HH^3$.  Figure \ref{Fig:43r_progression} shows a progression of $\{4,3,r\}$ honeycombs as r becomes large. We can see from this picture where the endpoints of the edges are. For example, two nearby white triangles have an edge between their closest vertices (compare with Figure \ref{Fig:437_boundary_pattern}).  In the limiting $\{p,q,\infty\}$ honeycombs, each edge has only one ideal point, the rest of the edge having become hyperideal.  In this sense, the edges of a $\{p,q,\infty\}$ honeycomb are like the edges of the largest square of Figure \ref{Fig:hyperideal}.  
 
In the two-dimensional case, when $q$ is infinite, the vertices become ideal.
In the three-dimensional case, when $r$ is infinite, the edges rather than the vertices have a single remaining ideal point.  In the four-dimensional case, the faces have a single remaining ideal point when the last term in the Schl\"afli symbol becomes infinite. In general, each facet of codimension two in the polytope has only a single ideal point when the last term of the Schl\"afli symbol becomes infinite.


As the last term of the Schl\"afli symbol is increased, vertices become ideal, then hyperideal, and sooner in higher dimensions.
This causes honeycombs with ideal vertices to become rarer as we move up dimensions, and in fact they do not exist in dimensions six and above~\cite{Coxeter1954}.

\subsection{Hyperideal cells}

As with ideal vertices, we can tell if a tiling will have hyperideal vertices by looking at the dual tiling and asking what kind of geometry the cells have. In the case of the $\{4,3,7\}$ honeycomb, the dual honeycomb has cells with Schl\"afli symbol $\{7,3\}$, which is a tiling of $\HH^2$ with three heptagons around each vertex. We can see many copies of the dual $\{3,7\}$ tiling in Figure \ref{Fig:437_boundary}, one for each hyperideal vertex: the hyperideal cells of the $\{7,3,4\}$ honeycomb are centered on the hyperideal vertices of the $\{4,3,7\}$ honeycomb, with a face corresponding to each edge of the $\{4,3,7\}$ honeycomb.

\begin{figure}[htbp]
\centering 
\subfloat[One $\{7,3\}$ cell with the induced pattern on the boundary of $\HH^3$.]
{
\includegraphics[width=0.4\textwidth]{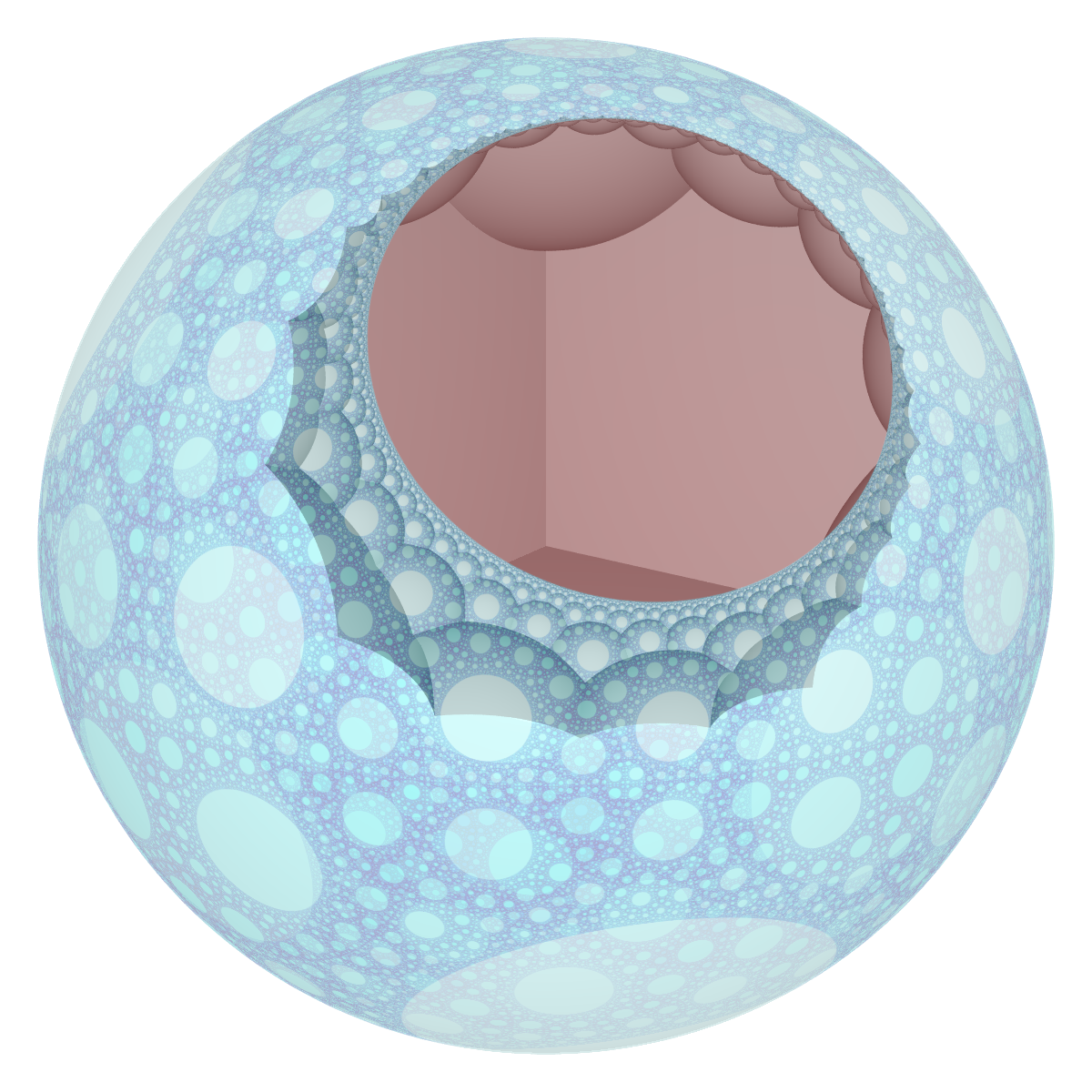}
\label{Fig:734_ball}
}
\hspace{1.0cm}
\subfloat[Stereographic projection of the boundary pattern.]
{
\includegraphics[width=0.4\textwidth]{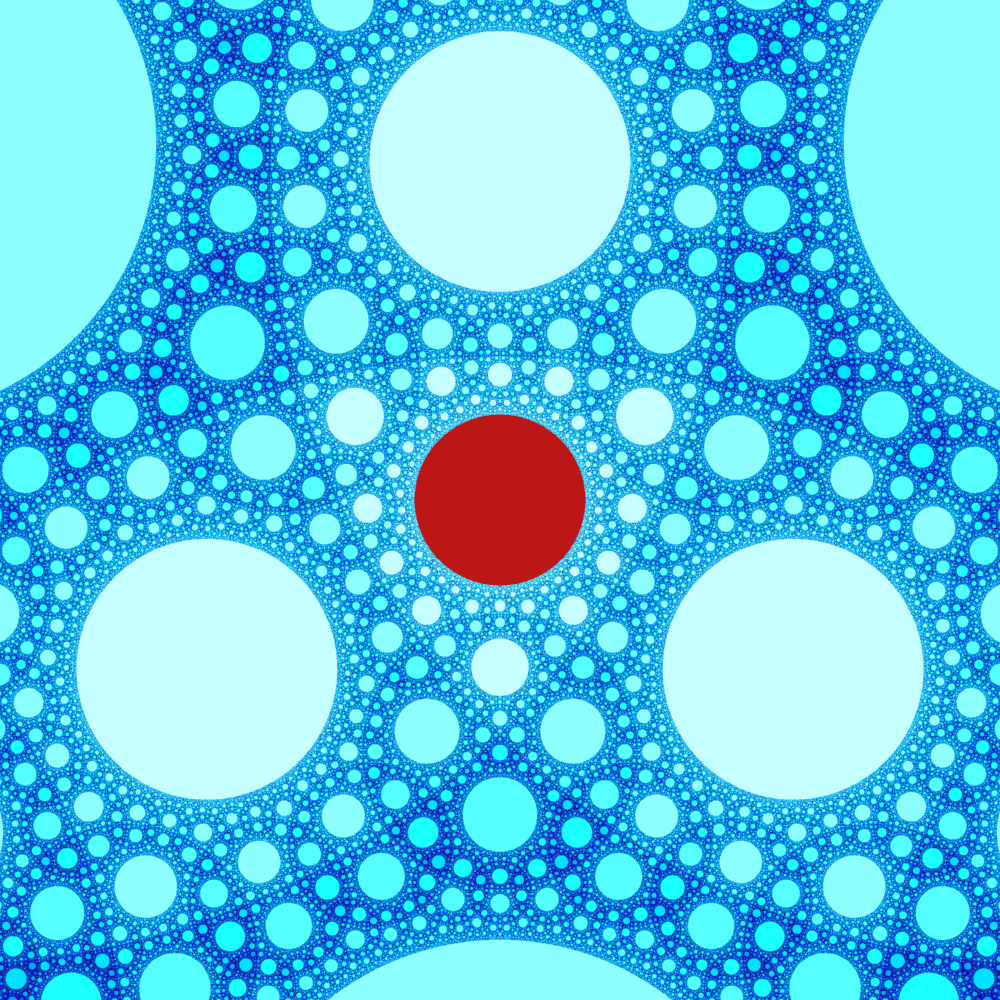}
\label{Fig:734_boundary}
}
\caption{The boundary pattern for the $\{7,3,4\}$ honeycomb.  The single \emph{head} of one cell is highlighted red.  Compare with Figure \ref{Fig:437_boundary_pattern}.}
\label{Fig:734_boundary_pattern}
\end{figure}

Hyperideal cells have a new feature: each such cell meets the boundary of $\HH^3$ in a circular region we call the \emph{head} of the cell, as shown in Figure \ref{Fig:734_boundary_pattern}.  In general, a honeycomb $\{p,q,r\}$ has hyperideal cells when $(p-2)(q-2)>4$.

\subsection{Ideal vertices and ideal cells}
\label{Sec:ideal_verts_and_cells}

\begin{figure}[htbp]
\centering 
\includegraphics[width=0.5\textwidth]{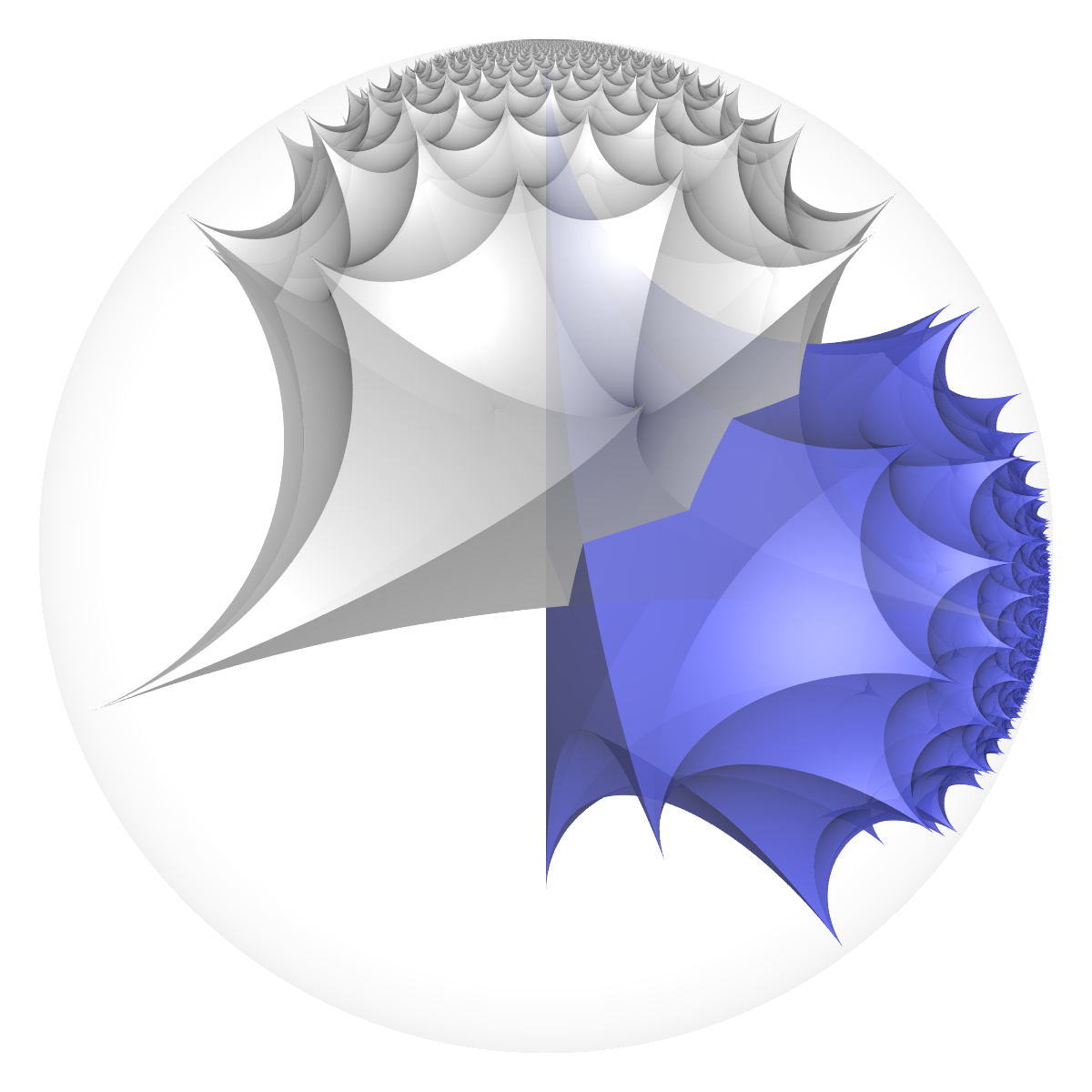}
\caption{A cell from the $\{3,6,3\}$ honeycomb, overlaid with a dual cell.}
\label{Fig:363_dual_cells}
\end{figure}

A honeycomb that has both ideal vertices and ideal cells cannot be represented in either a vertex-centered or a cell-centered picture. Instead, the most symmetrical pictures we can draw in the Poincar\'e ball model are edge-centered and face-centered. See for example Figure \ref{Fig:363_dual_cells}, showing a cell from the $\{3,6,3\}$ honeycomb and a cell from the dual honeycomb, another copy of $\{3,6,3\}$.  The blue cell is edge centered, with an edge of the cell going vertically through the center of the ball.  The white cell is face centered, with a triangular face of the cell on the horizontal plane through the center of the ball.  One ideal vertex of each cell is incident with the ideal center of its dual cell.  

\subsection{Hyperideal vertices and hyperideal cells}

\begin{figure}[htbp]
\centering 
\includegraphics[width=0.5\textwidth]{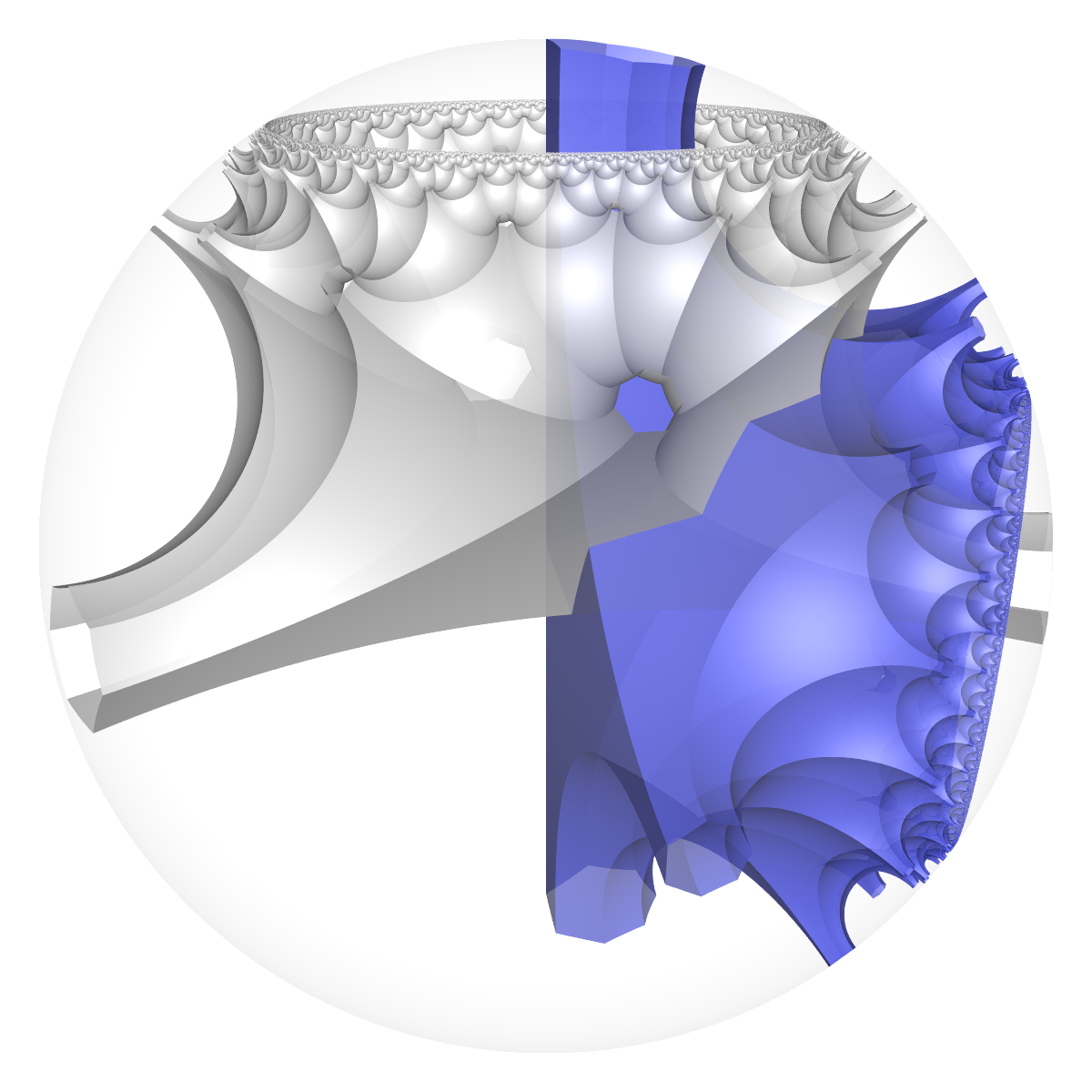}
\caption{A cell from the $\{3,7,3\}$, overlaid with a dual cell.}
\label{Fig:373_dual_cells}
\end{figure}

Similarly, honeycombs with hyperideal vertices and hyperideal cells cannot be represented in either vertex-centered or cell-centered pictures. For example, Figure \ref{Fig:373_dual_cells} shows a cell from the $\{3,7,3\}$ honeycomb and a cell from the dual honeycomb, another copy of $\{3,7,3\}$.  As in Figure \ref{Fig:363_dual_cells}, the blue cell is edge centered, with an edge of the cell going vertically through the center of the ball.  The white cell is face centered, with a hyperideal triangular face of the cell on the horizontal plane through the center of the ball.  One hyperideal vertex of each cell is incident with the hyperideal center of the dual cell, visible as a leg of the cell protruding through the dual cell's head.

\section{The Schl\"afli cube}
\label{Sec:schlafli_cube}

We can now draw pictures of every Schl\"afli honeycomb. Each Schl\"afli honeycomb corresponds to an integral point of the infinite cube $\{(p,q,r) \mid 3\leq p,q,r \leq\infty\}$. In this section we draw tables showing the honeycombs on the six faces of this cube, and on the two surfaces in Figure \ref{Fig:3D_Schlafli_map} corresponding to the honeycombs with ideal vertices and ideal cells.  High resolution images of a selection of honeycombs throughout the cube are available at \htmladdnormallink{hyperbolichoneycombs.org}{http://hyperbolichoneycombs.org}. Appendix \ref{Sec:simplex_construction} describes the construction in more detail, and provides some general observations about the pattern on the boundary for arbitrary $p$, $q$, and $r$.

For honeycombs with neither hyperideal vertices nor hyperideal cells, we draw (using POV-Ray~\cite{Povray}) the in-space view -- a view of the edges of the honeycomb as seen from within the space, assuming that light rays travel along geodesics. These are the fifteen hyperbolic honeycombs listed by Schlegel and Coxeter~\cite{Coxeter1954}, the one euclidean honeycomb $\{4,3,4\}$, and the six spherical honeycombs corresponding to the regular four-dimensional polychora. For all other honeycombs, we draw the boundary picture, as described in Section \ref{Sec:hyperidealization}.

A number of interesting patterns appear as the terms of the Schl\"afli symbol approach infinity. As $p\rightarrow\infty$, cell heads become tangent to each other.  As $q\rightarrow\infty$, cell heads become tangent to hyperideal vertices.  As $r\rightarrow\infty$, hyperideal vertices become tangent to each other.  Only for $\{\infty,\infty,\infty\}$ do all of these tangencies exist simultaneously.  

Once we started drawing pictures of hyperideal honeycombs, we began to see them in a number of places, especially honeycombs with hyperideal cells and material vertices.  Many of these correspondences relate to Kleinian groups, fractals, and Apollonian gaskets.  \emph{Indra's Pearls}, by David Mumford, Caroline Series, and David Wright, shows an image of a Kleinian group which is visually the same as our boundary image of the $\{3,3,\infty\}$ honeycomb~\cite[p. 200]{Mumford2002}.  The cover of \emph{Multifractals and 1/f Noise: Wild Self-Affinity in Physics} by Benoit B. Mandelbrot shows a fractal known as the ``The Pharaoh's Breastplate''~\cite{benoit1998multifractals}. This is also a boundary image of the $\{\infty,3,4\}$ honeycomb.  Beautiful images of Schmidt arrangements in the paper \emph{Visualising the arithmetic of imaginary quadratic fields} by Katherine E. Stange are intimately connected to ideal honeycombs such as the $\{4,4,4\}$~\cite[p. 8]{stange2014visualising}.  Visualizations of ``limit root sets" may end up looking like our boundary images, for example the $\{7,3,3\}$~\cite{chen2015lorentzian}.  Random groups are related to Coxeter groups of honeycombs of the form $\{p,3,3\}$ ~\cite{calegari20143}.

\iftables
\def\w{0.14}
\setlength\tabcolsep{2pt}
\newcolumntype{V}[1]{>{\centering\arraybackslash} m{#1\textwidth} <{\vspace{-1pt}}}

\begin{table}[htbp]
\centering
\begin{tabular}{V{0.05}|V{\w}V{\w}V{\w}V{\w}V{\w}V{0.06}V{\w}}

  $r\backslash{q}$ & 3 & 4 & 5 & 6 & 7 & $\cdots$ & $\infty$ \\	\hline \\ [-12pt]
    3	& \includegraphics[width=\w\textwidth]{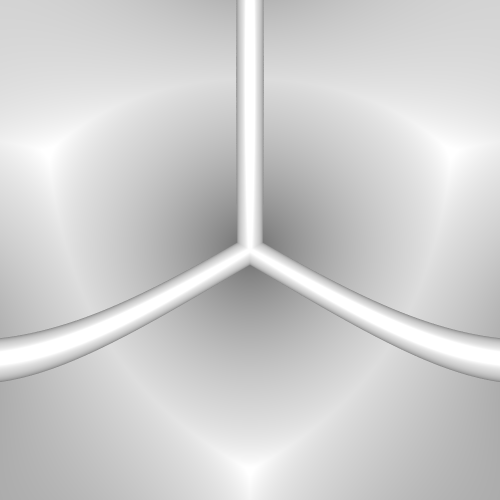} 
	& \includegraphics[width=\w\textwidth]{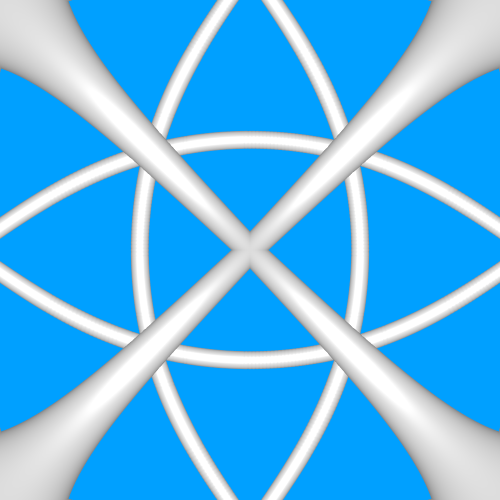} 
	& \includegraphics[width=\w\textwidth]{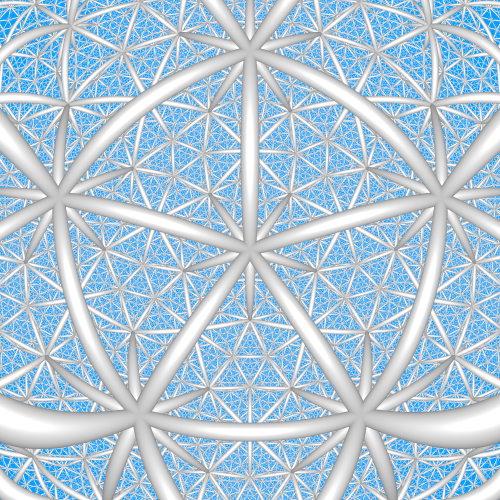} 
 	& \includegraphics[width=\w\textwidth]{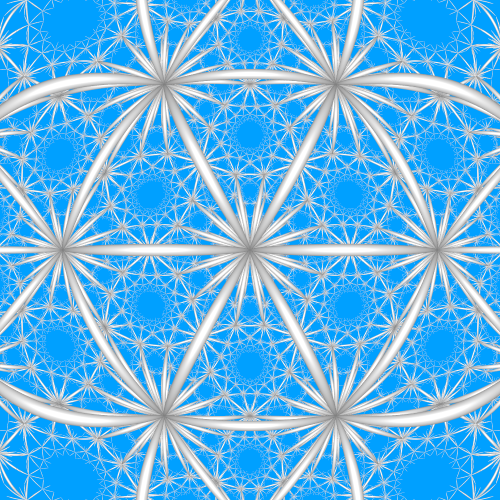} 
	& \includegraphics[width=\w\textwidth]{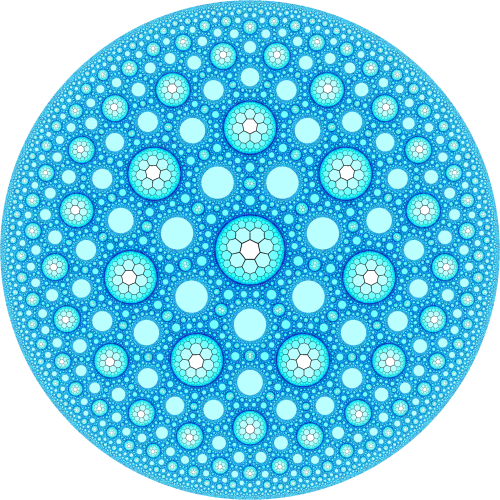} 
	& $\cdots$ 
	& \includegraphics[width=\w\textwidth]{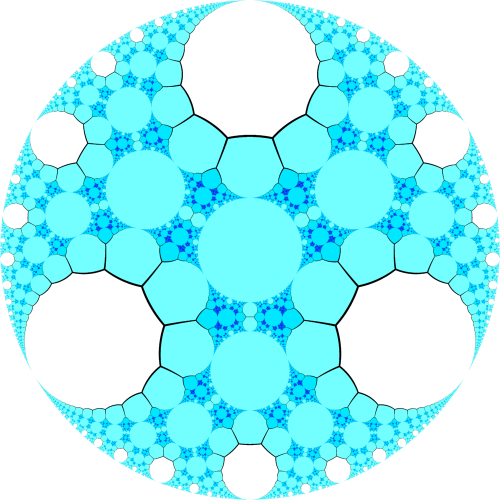}  \\
    4	& \includegraphics[width=\w\textwidth]{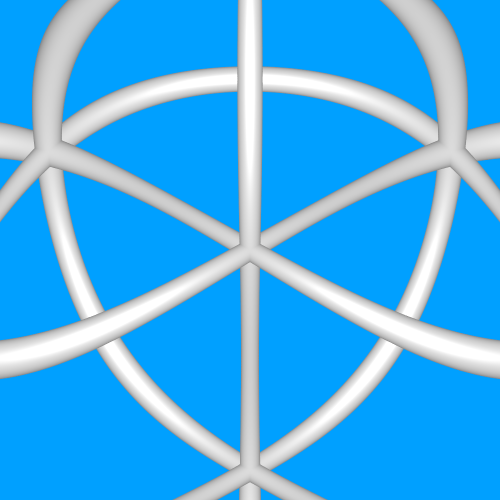} 
	& \includegraphics[width=\w\textwidth]{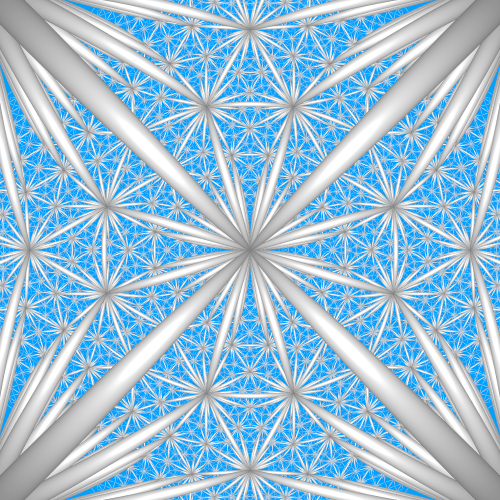} 
	& \includegraphics[width=\w\textwidth]{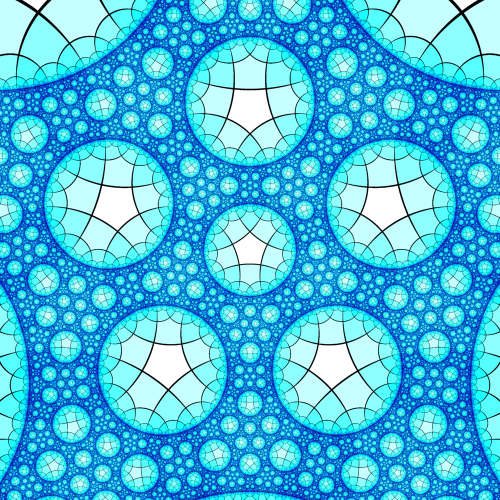} 
 	& \includegraphics[width=\w\textwidth]{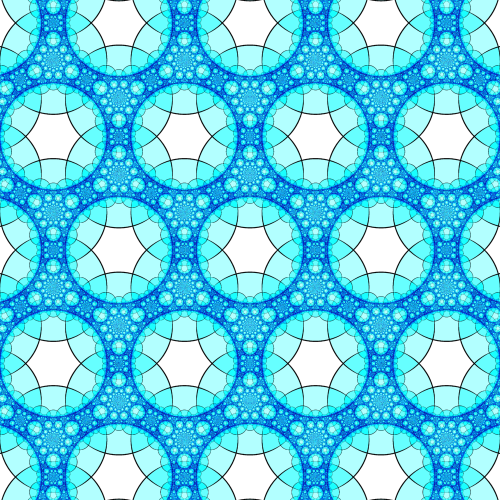} 
	& \includegraphics[width=\w\textwidth]{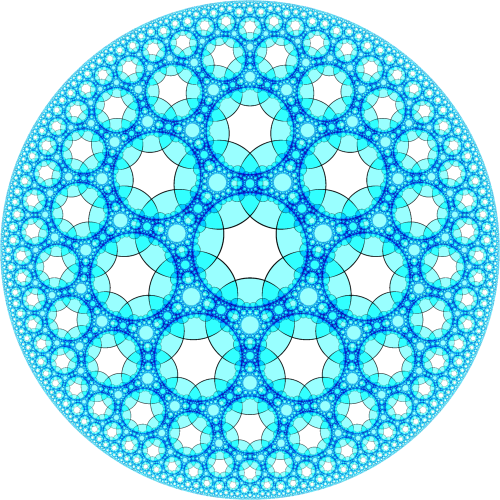} 
	& $\cdots$ 
	& \includegraphics[width=\w\textwidth]{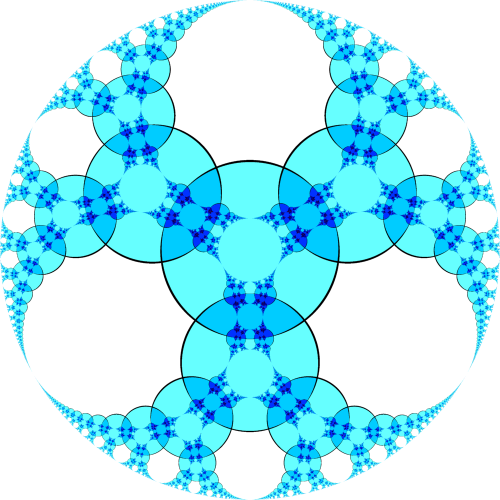}  \\
    5	& \includegraphics[width=\w\textwidth]{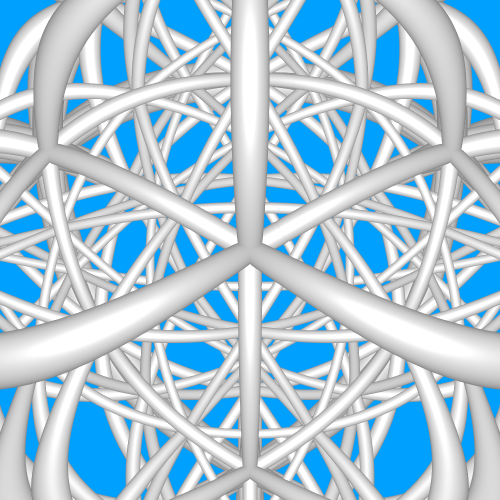} 
	& \includegraphics[width=\w\textwidth]{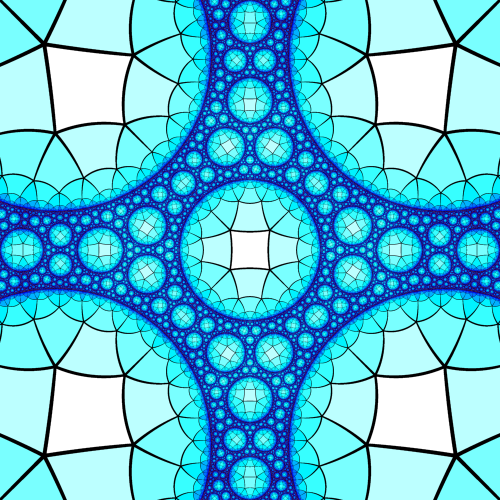} 
	& \includegraphics[width=\w\textwidth]{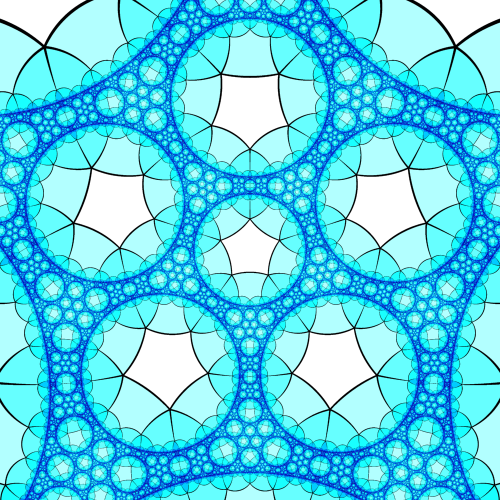} 
 	& \includegraphics[width=\w\textwidth]{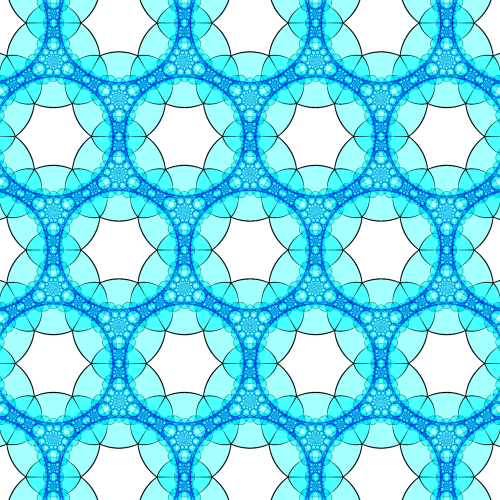} 
	& \includegraphics[width=\w\textwidth]{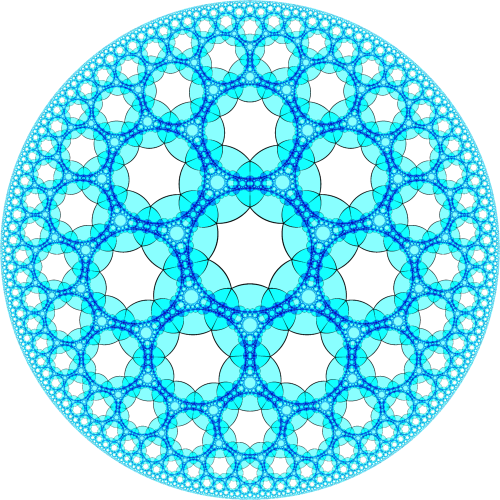} 
	& $\cdots$ 
	& \includegraphics[width=\w\textwidth]{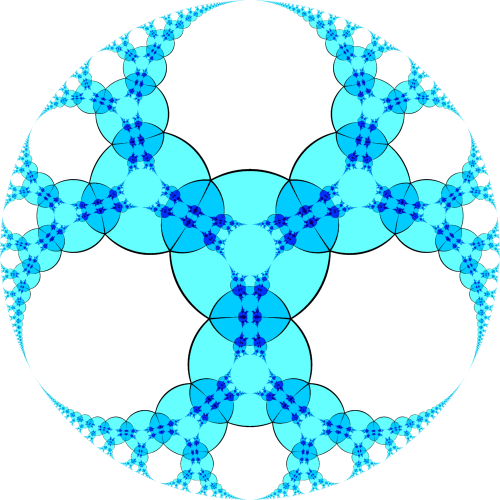}  \\
    6	& \includegraphics[width=\w\textwidth]{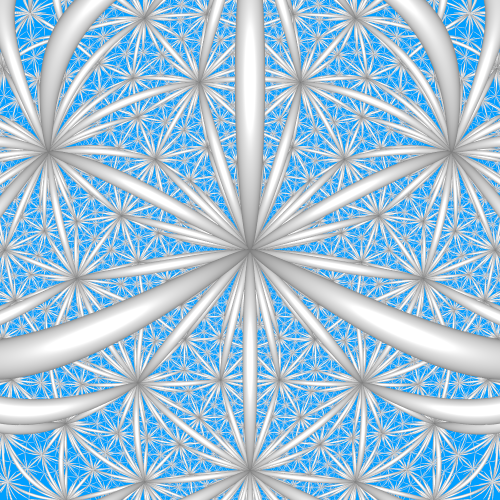} 
	& \includegraphics[width=\w\textwidth]{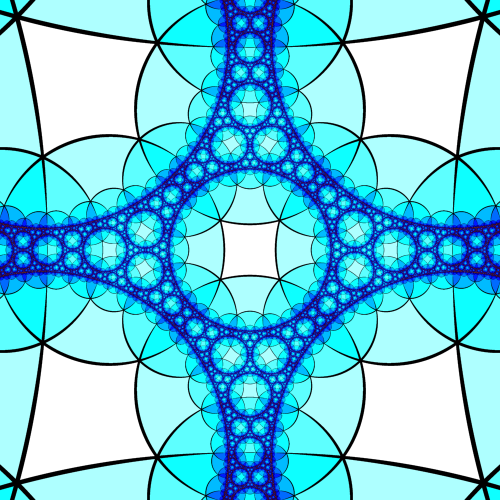} 
	& \includegraphics[width=\w\textwidth]{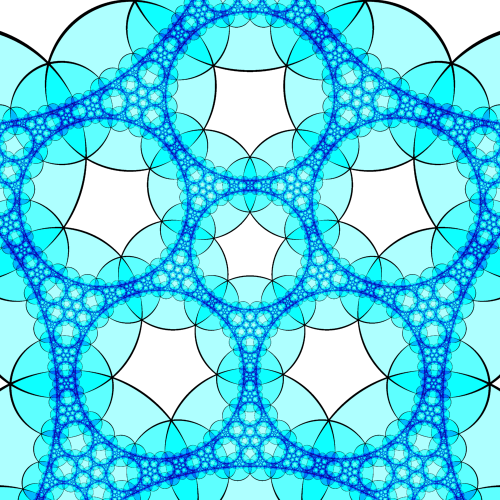} 
 	& \includegraphics[width=\w\textwidth]{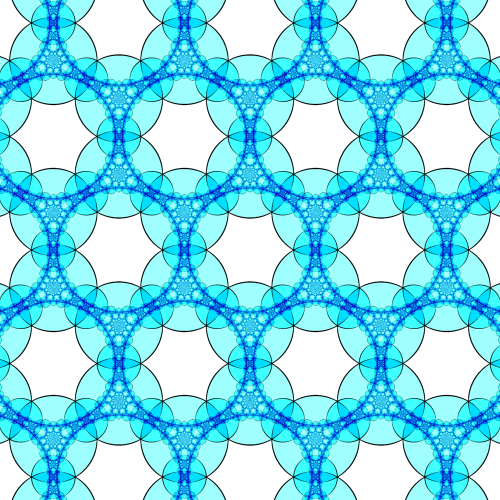} 
	& \includegraphics[width=\w\textwidth]{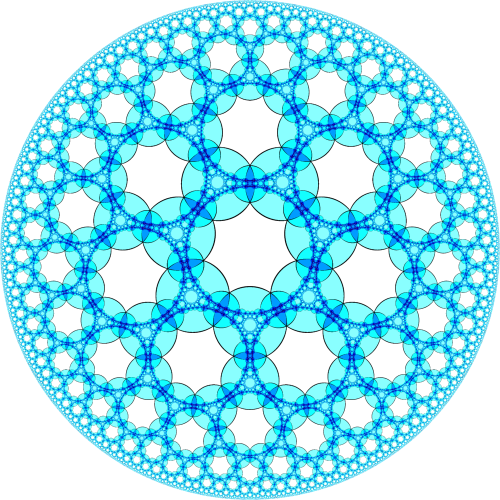} 
	& $\cdots$ 
	& \includegraphics[width=\w\textwidth]{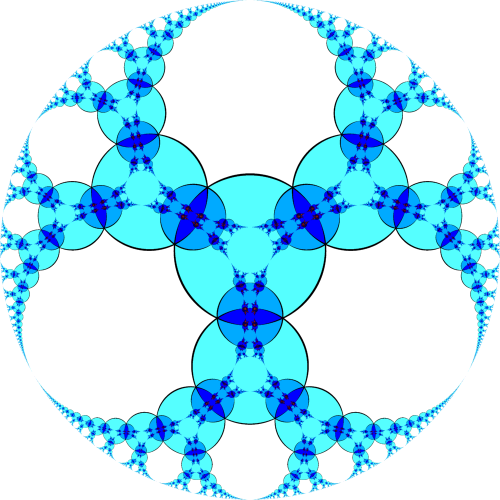}  \\
    7	& \includegraphics[width=\w\textwidth]{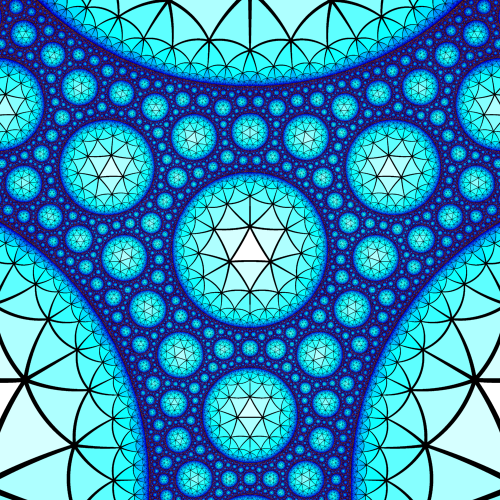} 
	& \includegraphics[width=\w\textwidth]{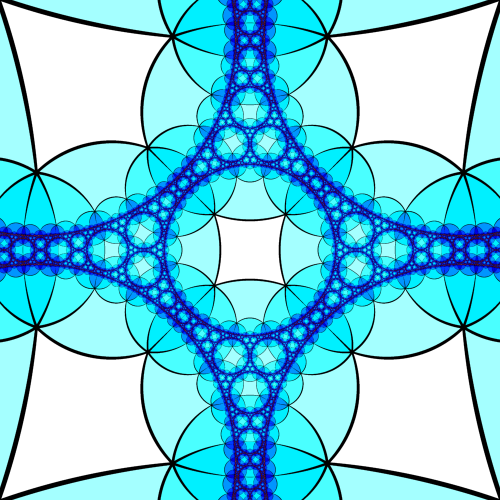} 
	& \includegraphics[width=\w\textwidth]{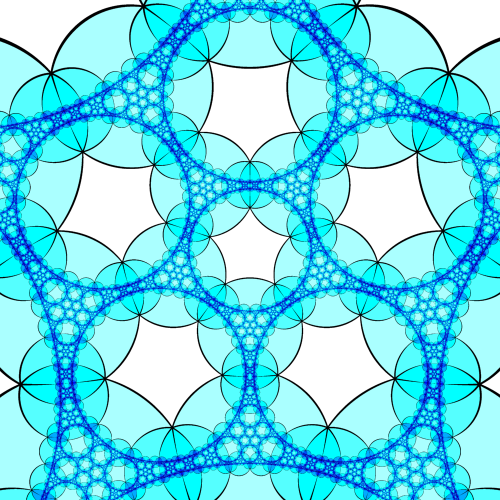} 
 	& \includegraphics[width=\w\textwidth]{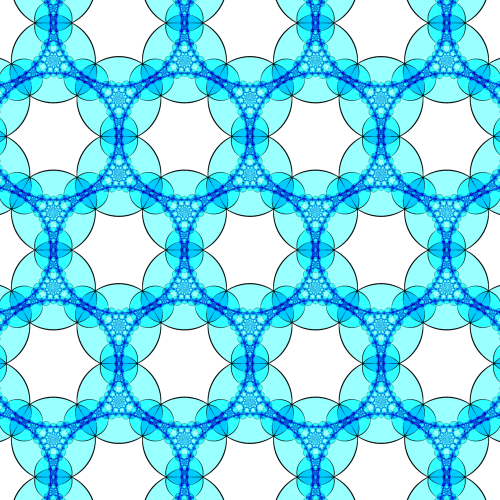} 
	& \includegraphics[width=\w\textwidth]{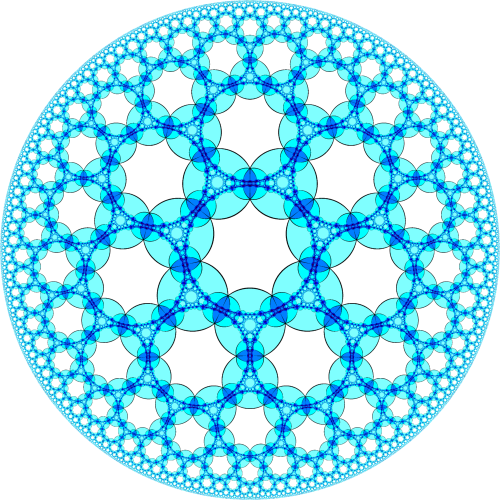} 
	& $\cdots$ 
	& \includegraphics[width=\w\textwidth]{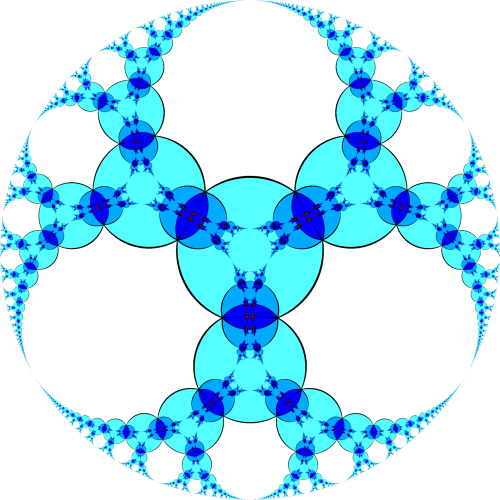}  \\ [-10pt] \\
  $\vdots$ &   $\vdots$  &   $\vdots$ &  $\vdots$ &   $\vdots$&   $\vdots$& $\ddots$ & $\vdots$  \\ \\
  $\infty$ 
	& \includegraphics[width=\w\textwidth]{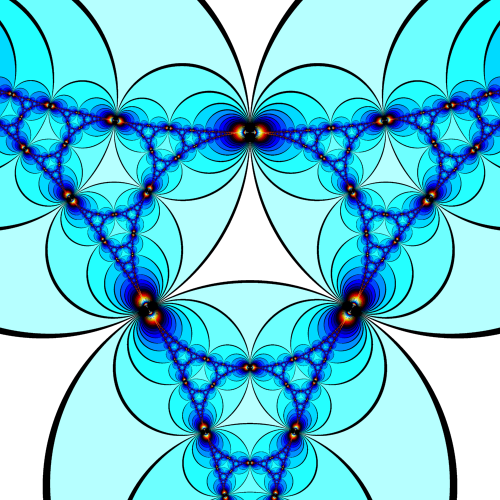} 
	& \includegraphics[width=\w\textwidth]{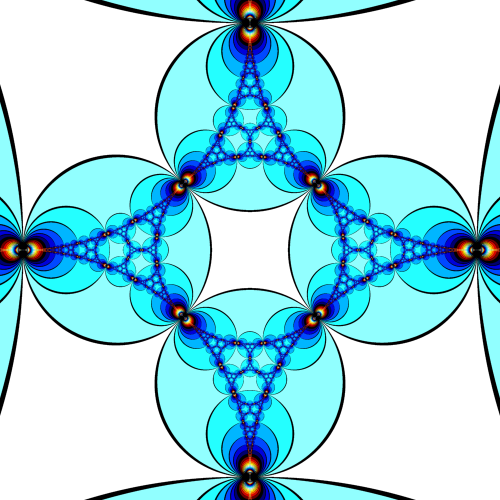} 
	& \includegraphics[width=\w\textwidth]{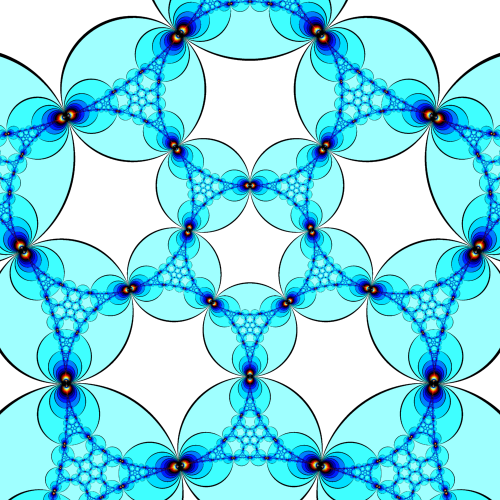} 
 	& \includegraphics[width=\w\textwidth]{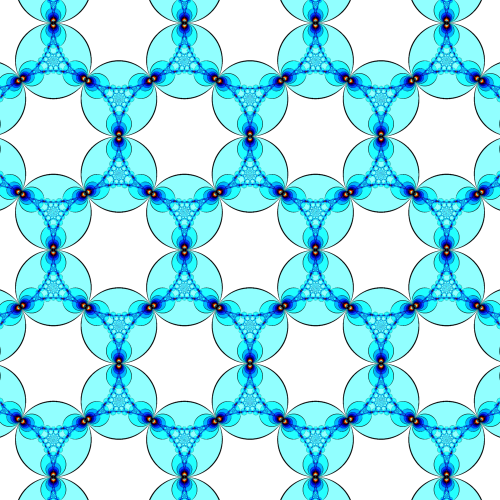} 
	& \includegraphics[width=\w\textwidth]{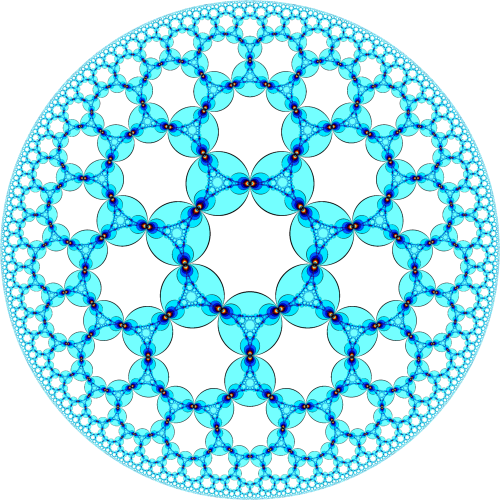} 
	& $\cdots$ 
	& \includegraphics[width=\w\textwidth]{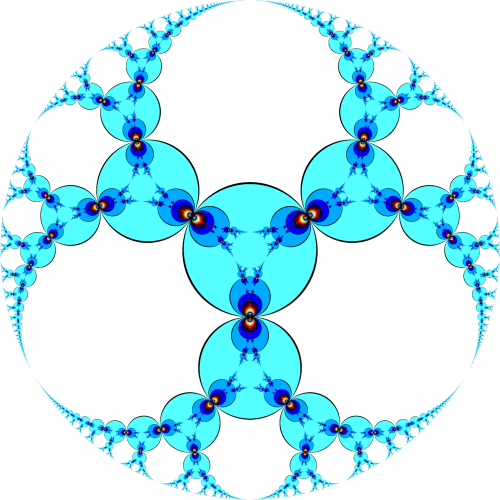}  \\

\end{tabular}
\vspace{15pt}
\caption{The $\{3,q,r\}$ honeycombs.}
\label{Table:3qr}
\end{table}

\begin{table}[htbp]
\centering
\begin{tabular}{V{0.05}|V{\w}V{\w}V{\w}V{\w}V{\w}V{0.06}V{\w}}

  $r\backslash{p}$ & 3 & 4 & 5 & 6 & 7 & $\cdots$ & $\infty$ \\	\hline \\ [-12pt]
    3	& \includegraphics[width=\w\textwidth]{Figures/compact/333.png} 
	& \includegraphics[width=\w\textwidth]{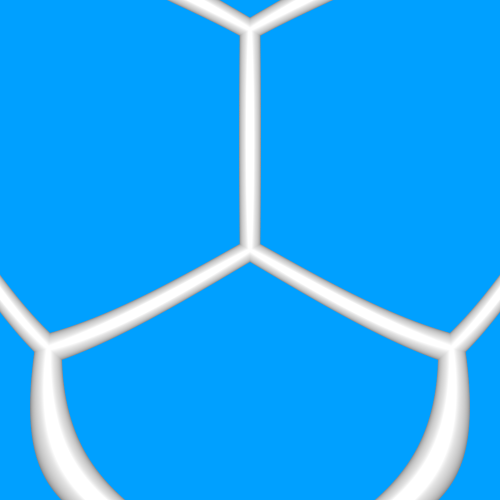} 
	& \includegraphics[width=\w\textwidth]{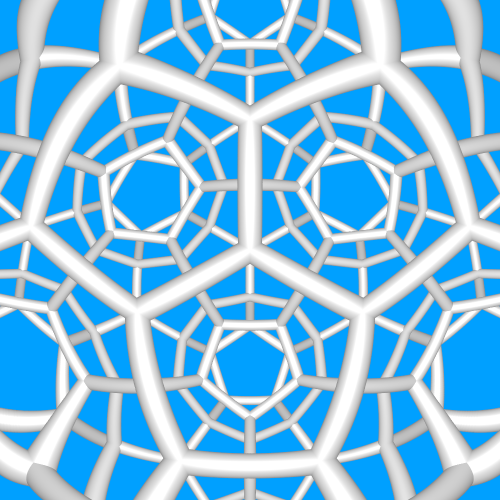} 
 	& \includegraphics[width=\w\textwidth]{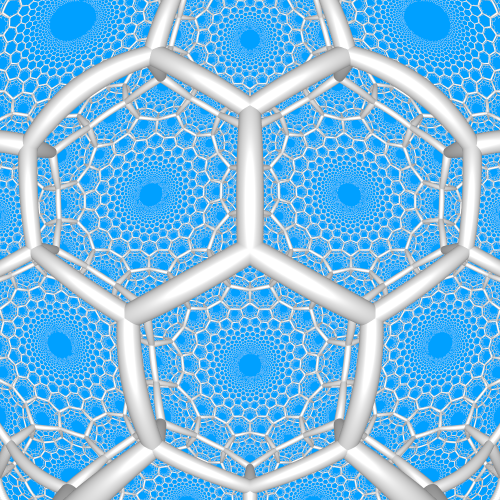} 
	& \includegraphics[width=\w\textwidth]{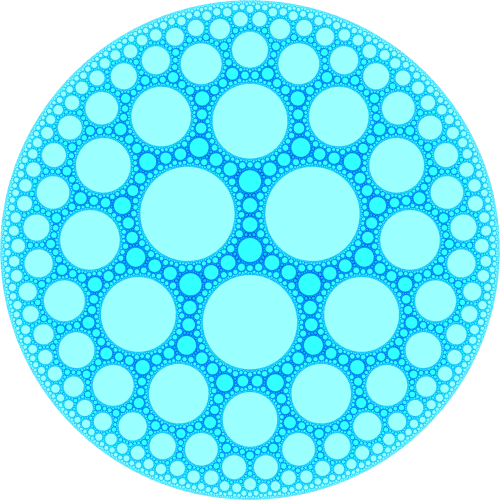} 
	& $\cdots$ 
	& \includegraphics[width=\w\textwidth]{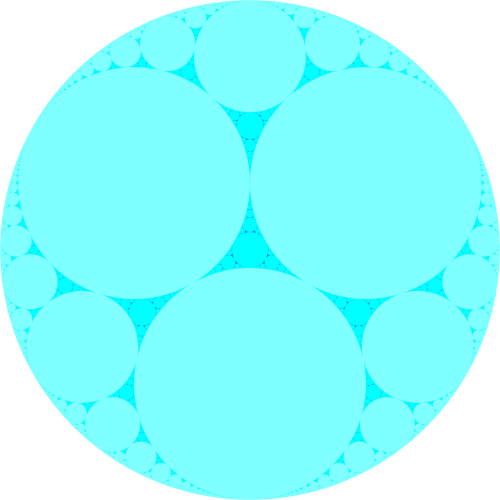}  \\
    4	& \includegraphics[width=\w\textwidth]{Figures/compact/334.png} 
	& \includegraphics[width=\w\textwidth]{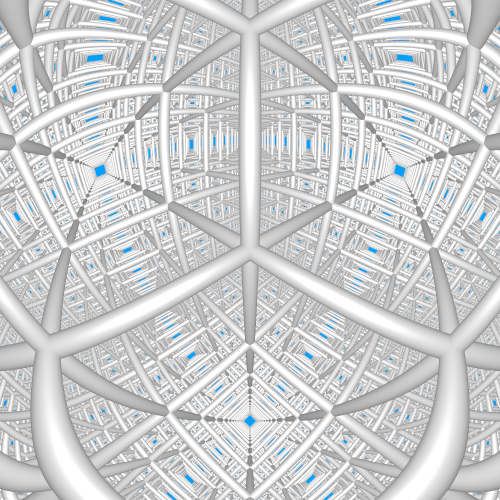} 
	& \includegraphics[width=\w\textwidth]{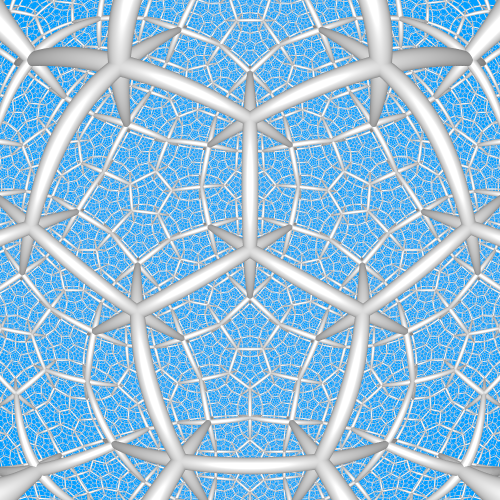} 
 	& \includegraphics[width=\w\textwidth]{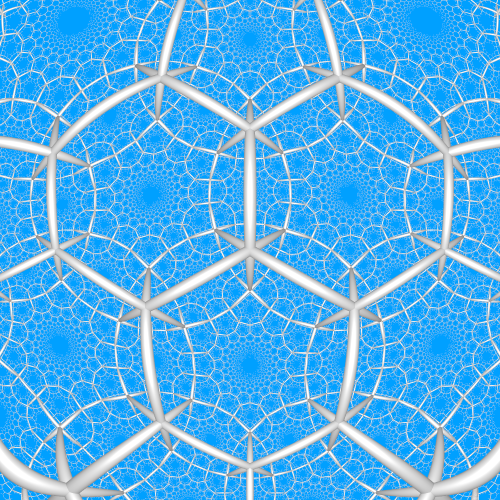} 
	& \includegraphics[width=\w\textwidth]{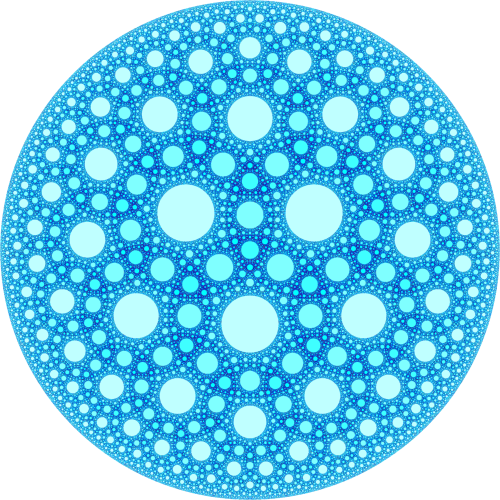} 
	& $\cdots$ 
	& \includegraphics[width=\w\textwidth]{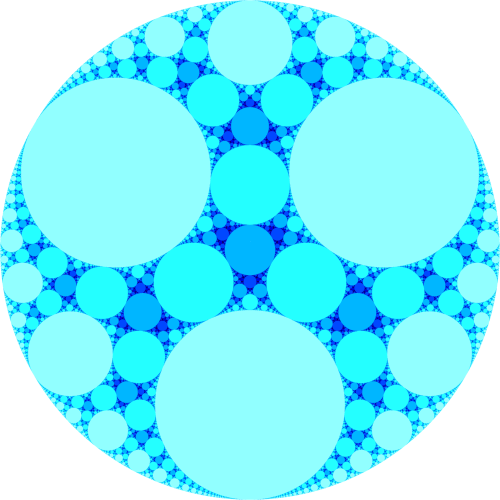}  \\
    5	& \includegraphics[width=\w\textwidth]{Figures/compact/335.png} 
	& \includegraphics[width=\w\textwidth]{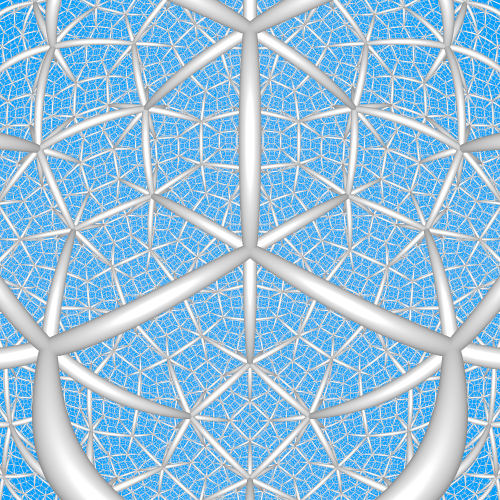} 
	& \includegraphics[width=\w\textwidth]{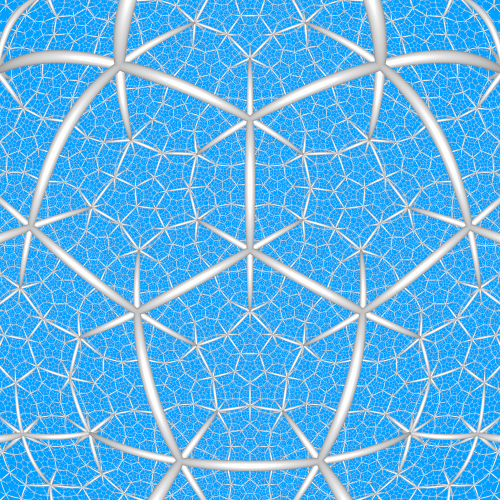} 
 	& \includegraphics[width=\w\textwidth]{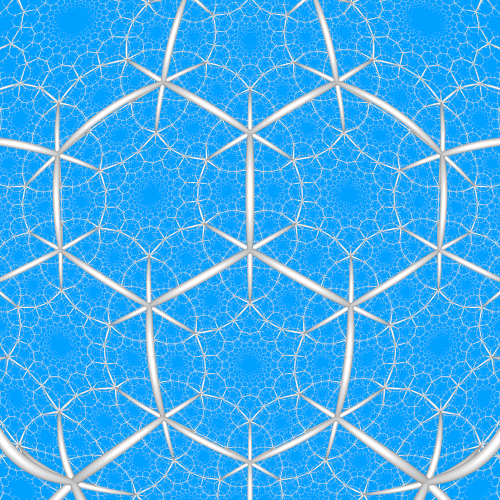} 
	& \includegraphics[width=\w\textwidth]{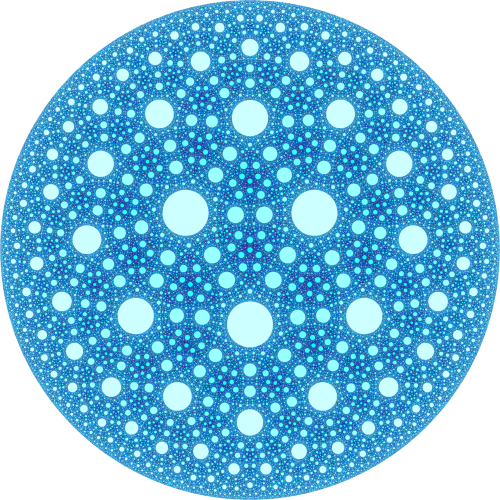} 
	& $\cdots$ 
	& \includegraphics[width=\w\textwidth]{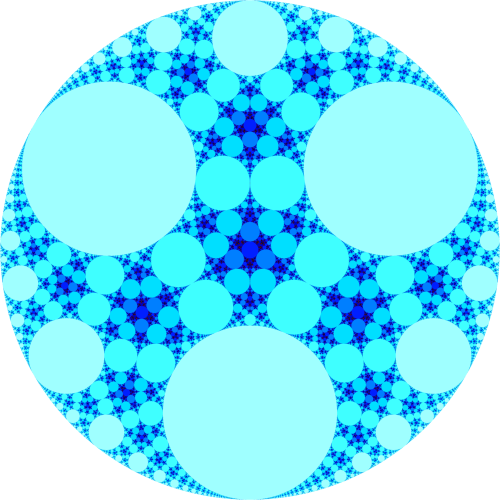}  \\
    6	& \includegraphics[width=\w\textwidth]{Figures/compact/336.png} 
	& \includegraphics[width=\w\textwidth]{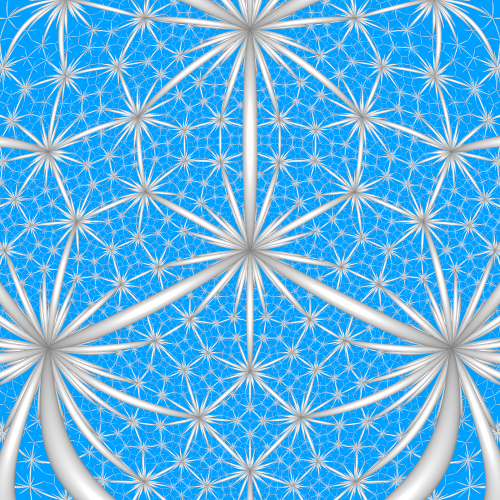} 
	& \includegraphics[width=\w\textwidth]{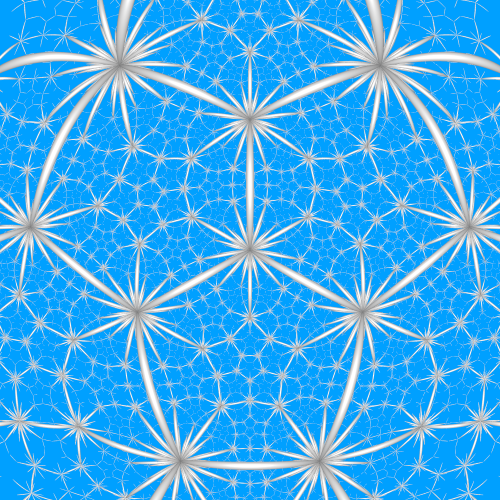} 
 	& \includegraphics[width=\w\textwidth]{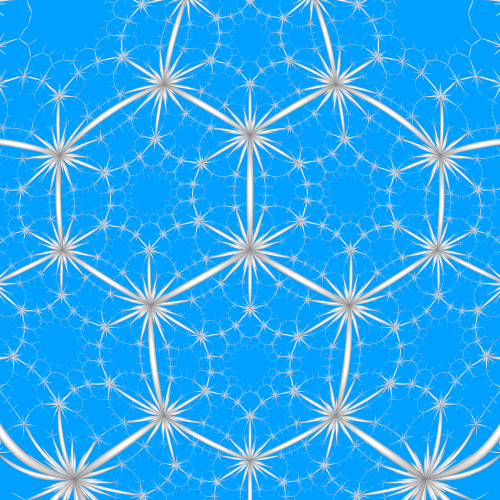} 
	& \includegraphics[width=\w\textwidth]{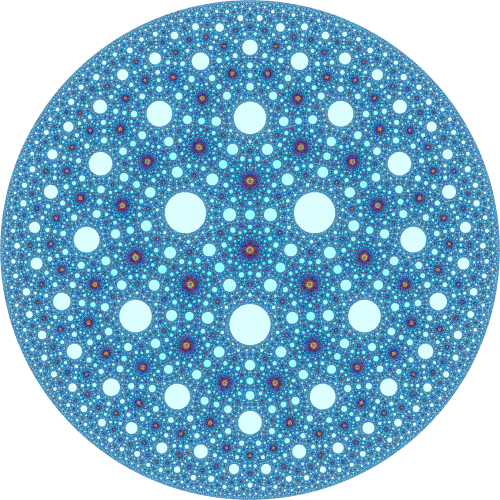} 
	& $\cdots$ 
	& \includegraphics[width=\w\textwidth]{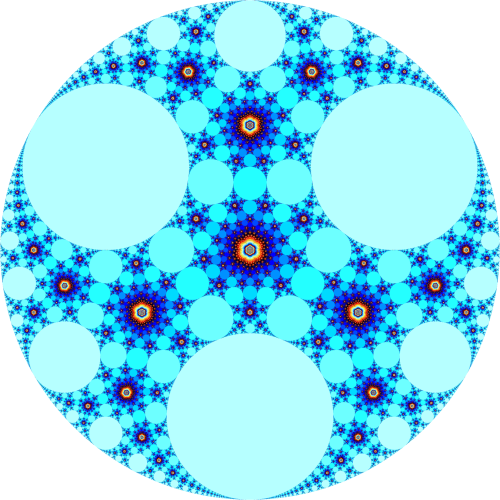}  \\
    7	& \includegraphics[width=\w\textwidth]{Figures/noncompact/337.png} 
	& \includegraphics[width=\w\textwidth]{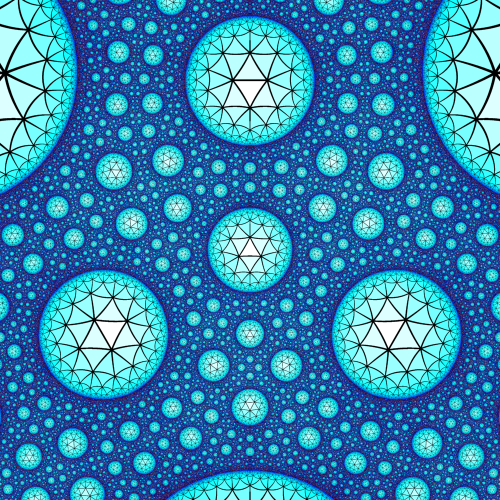} 
	& \includegraphics[width=\w\textwidth]{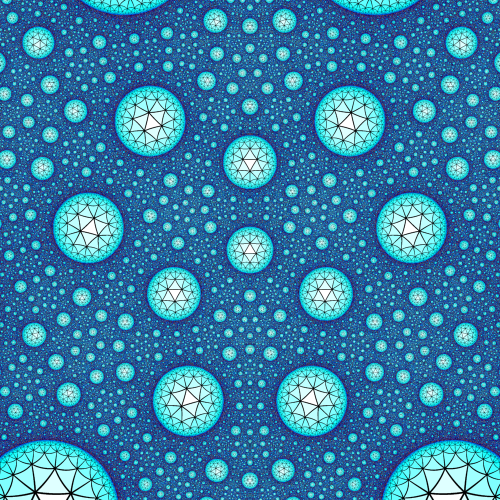} 
 	& \includegraphics[width=\w\textwidth]{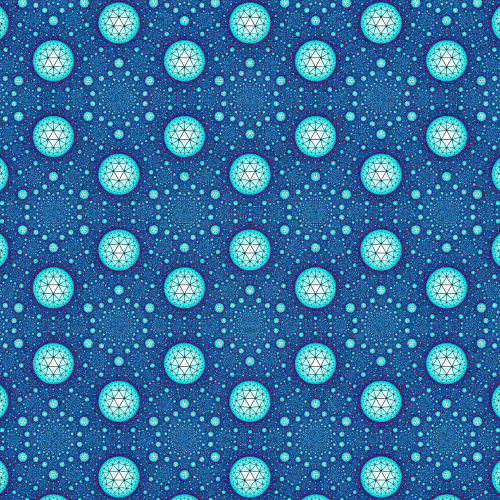} 
	& \includegraphics[width=\w\textwidth]{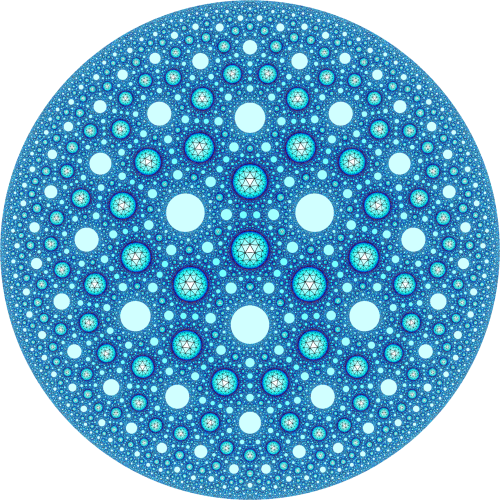} 
	& $\cdots$ 
	& \includegraphics[width=\w\textwidth]{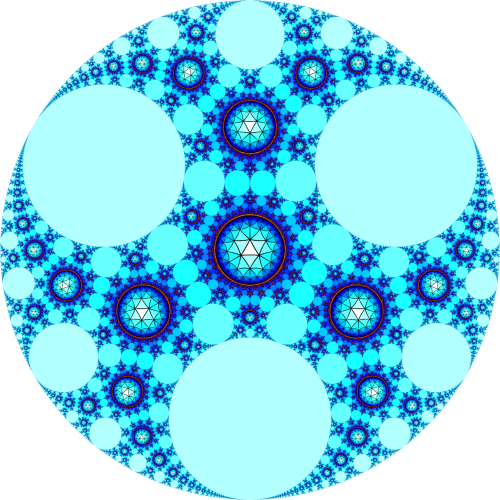}  \\ [-10pt] \\
  $\vdots$ &   $\vdots$  &   $\vdots$ &  $\vdots$ &   $\vdots$&   $\vdots$& $\ddots$ & $\vdots$  \\ \\
  $\infty$ 
	& \includegraphics[width=\w\textwidth]{Figures/noncompact/33i.png} 
	& \includegraphics[width=\w\textwidth]{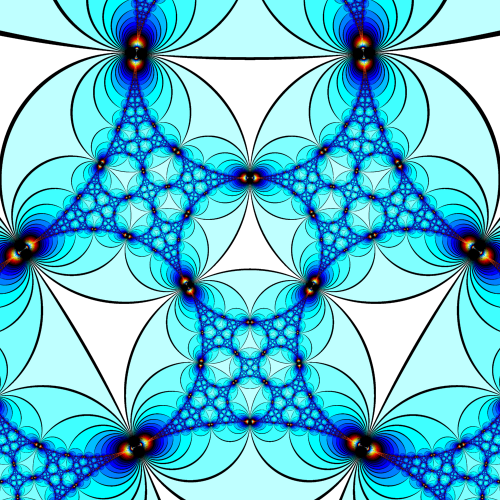} 
	& \includegraphics[width=\w\textwidth]{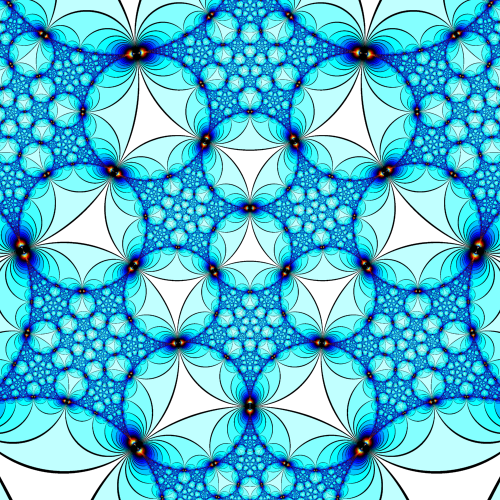} 
 	& \includegraphics[width=\w\textwidth]{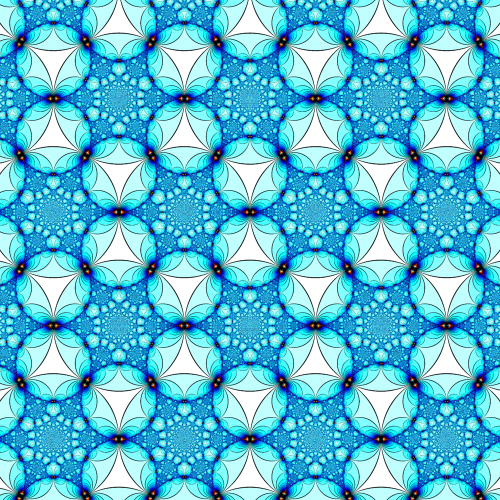} 
	& \includegraphics[width=\w\textwidth]{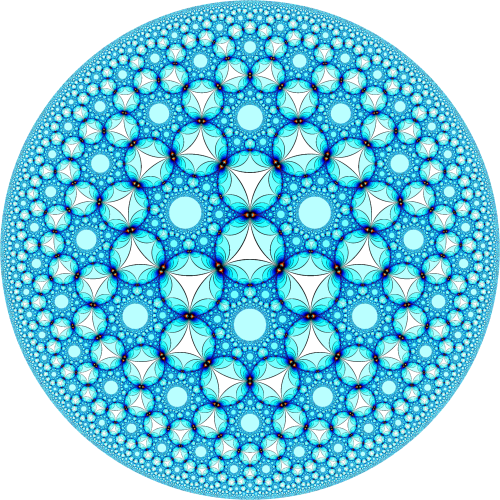} 
	& $\cdots$ 
	& \includegraphics[width=\w\textwidth]{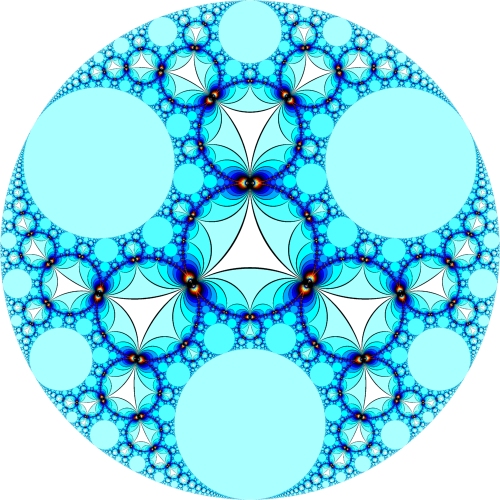}  \\

\end{tabular}
\vspace{15pt}
\caption{The $\{p,3,r\}$ honeycombs.}
\label{Table:p3r}
\end{table}

\begin{table}[htbp]
\centering
\begin{tabular}{V{0.05}|V{\w}V{\w}V{\w}V{\w}V{\w}V{0.06}V{\w}}

  $q\backslash{p}$ & 3 & 4 & 5 & 6 & 7 & $\cdots$ & $\infty$ \\	\hline \\ [-12pt]
    3	& \includegraphics[width=\w\textwidth]{Figures/compact/333.png} 
	& \includegraphics[width=\w\textwidth]{Figures/compact/433.png} 
	& \includegraphics[width=\w\textwidth]{Figures/compact/533.png} 
 	& \includegraphics[width=\w\textwidth]{Figures/compact/633.png} 
	& \includegraphics[width=\w\textwidth]{Figures/noncompact/733.png} 
	& $\cdots$ 
	& \includegraphics[width=\w\textwidth]{Figures/noncompact/i33.png}  \\
    4	& \includegraphics[width=\w\textwidth]{Figures/compact/343.png} 
	& \includegraphics[width=\w\textwidth]{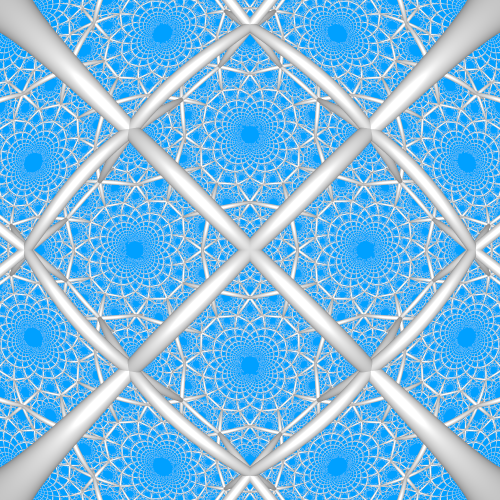} 
	& \includegraphics[width=\w\textwidth]{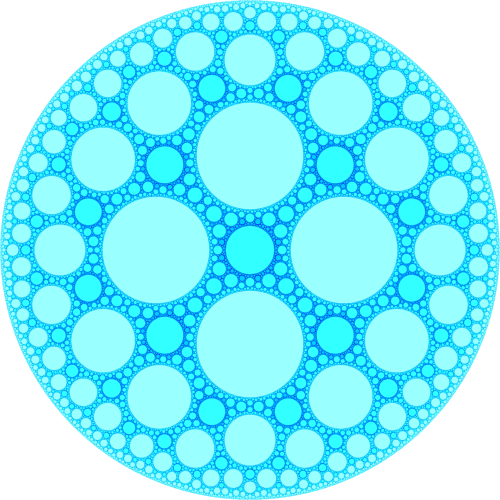} 
 	& \includegraphics[width=\w\textwidth]{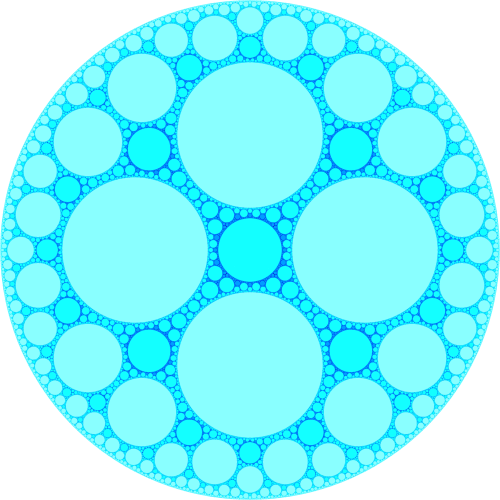} 
	& \includegraphics[width=\w\textwidth]{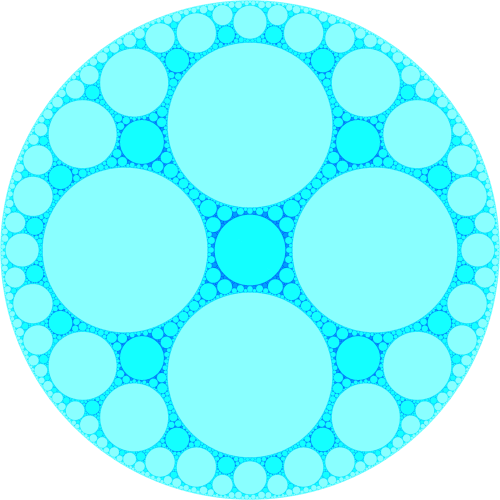} 
	& $\cdots$ 
	& \includegraphics[width=\w\textwidth]{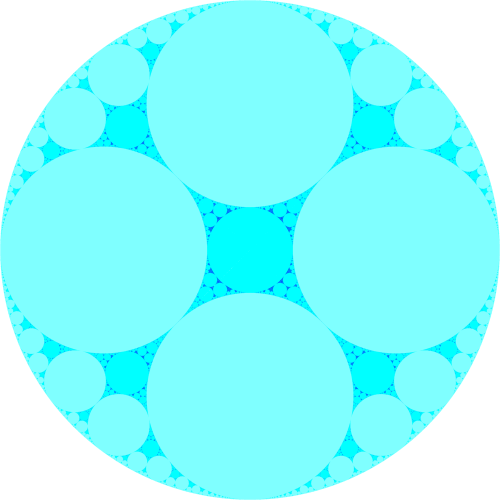}  \\
    5	& \includegraphics[width=\w\textwidth]{Figures/compact/353.png} 
	& \includegraphics[width=\w\textwidth]{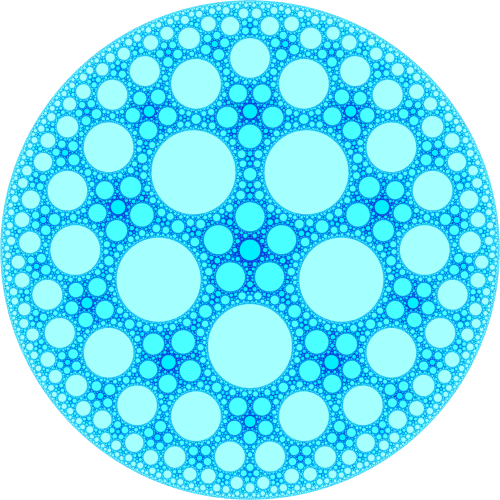} 
	& \includegraphics[width=\w\textwidth]{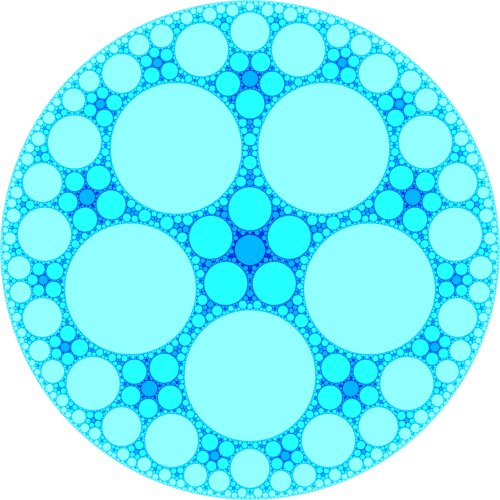} 
 	& \includegraphics[width=\w\textwidth]{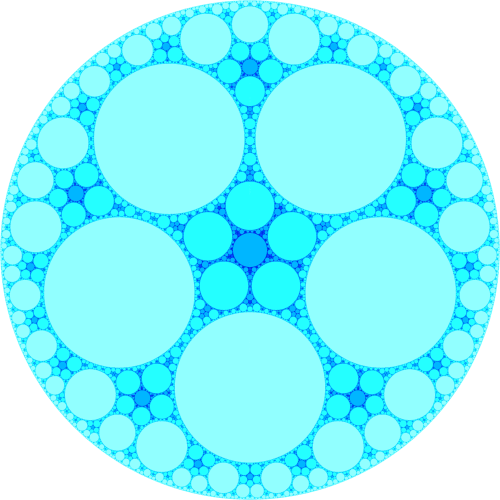} 
	& \includegraphics[width=\w\textwidth]{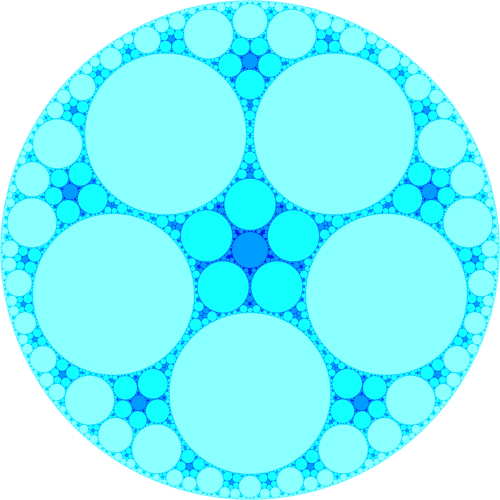} 
	& $\cdots$ 
	& \includegraphics[width=\w\textwidth]{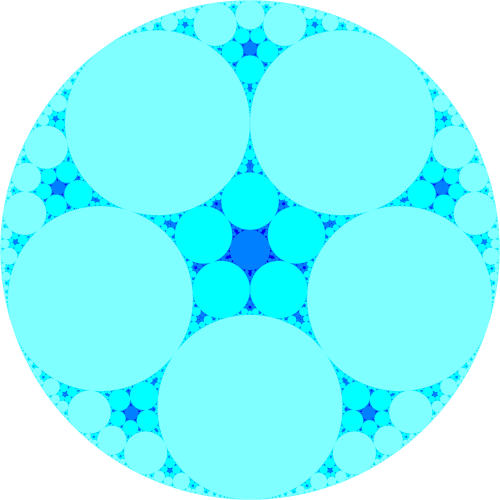}  \\
    6	& \includegraphics[width=\w\textwidth]{Figures/compact/363.png} 
	& \includegraphics[width=\w\textwidth]{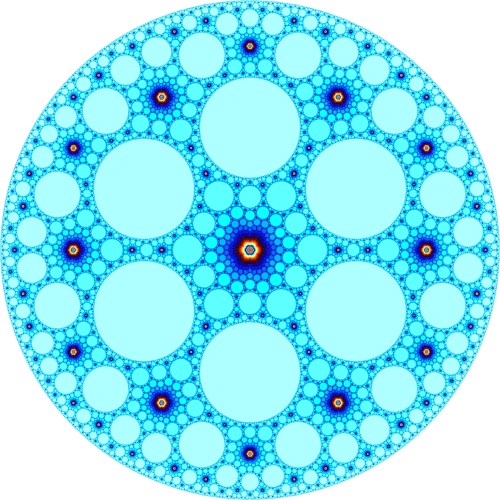} 
	& \includegraphics[width=\w\textwidth]{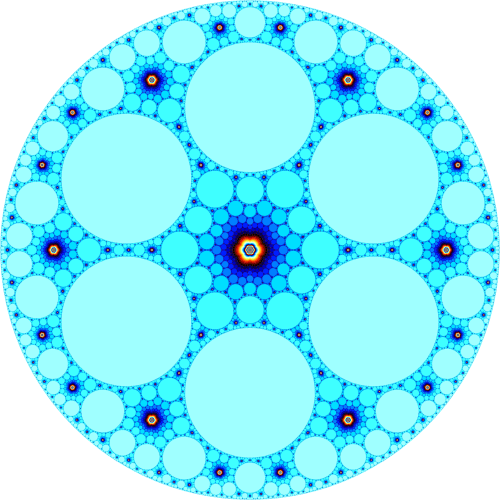} 
 	& \includegraphics[width=\w\textwidth]{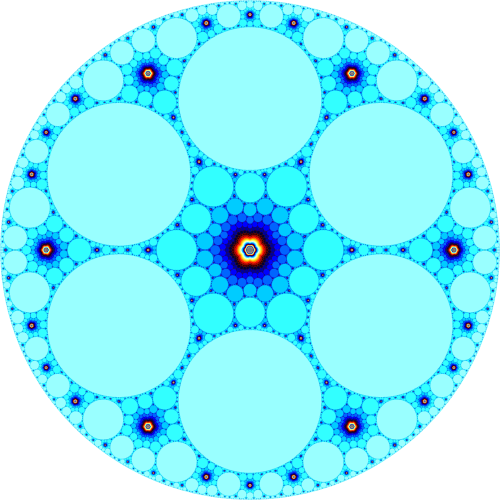} 
	& \includegraphics[width=\w\textwidth]{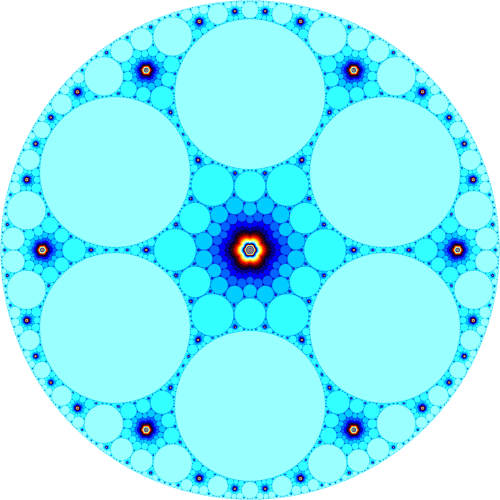} 
	& $\cdots$ 
	& \includegraphics[width=\w\textwidth]{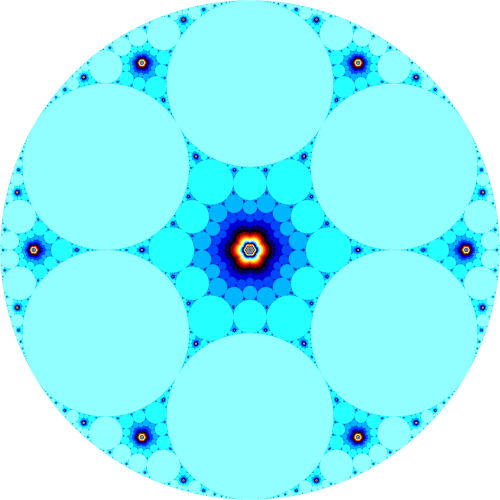}  \\
    7	& \includegraphics[width=\w\textwidth]{Figures/noncompact/373.png} 
	& \includegraphics[width=\w\textwidth]{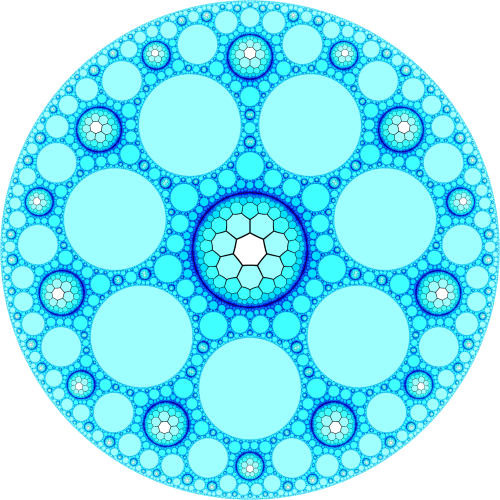} 
	& \includegraphics[width=\w\textwidth]{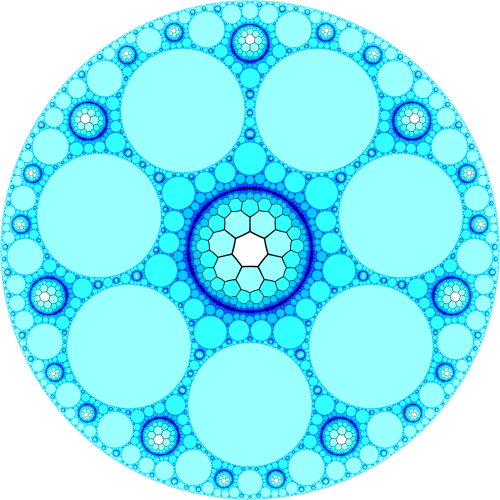} 
 	& \includegraphics[width=\w\textwidth]{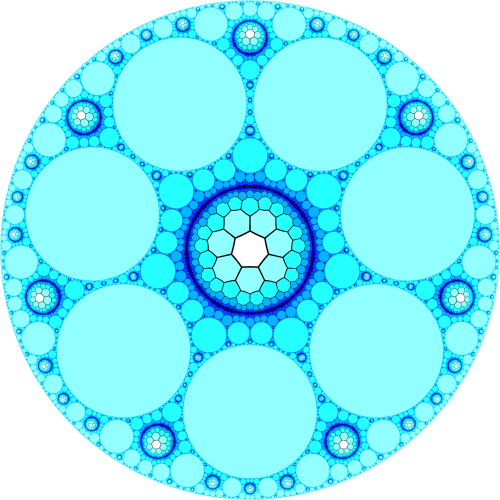} 
	& \includegraphics[width=\w\textwidth]{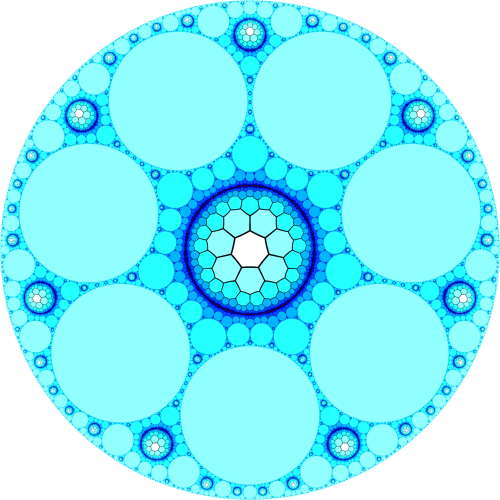} 
	& $\cdots$ 
	& \includegraphics[width=\w\textwidth]{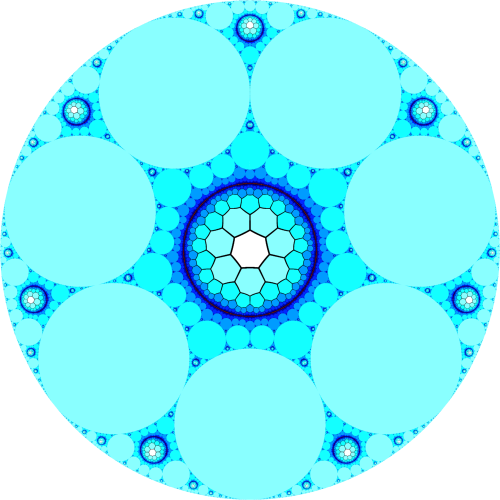}  \\ [-10pt] \\
  $\vdots$ &   $\vdots$  &   $\vdots$ &  $\vdots$ &   $\vdots$&   $\vdots$& $\ddots$ & $\vdots$  \\ \\
  $\infty$ 
	& \includegraphics[width=\w\textwidth]{Figures/noncompact/3i3.png} 
	& \includegraphics[width=\w\textwidth]{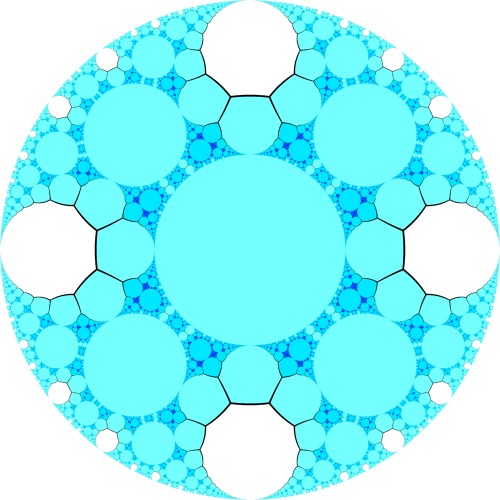} 
	& \includegraphics[width=\w\textwidth]{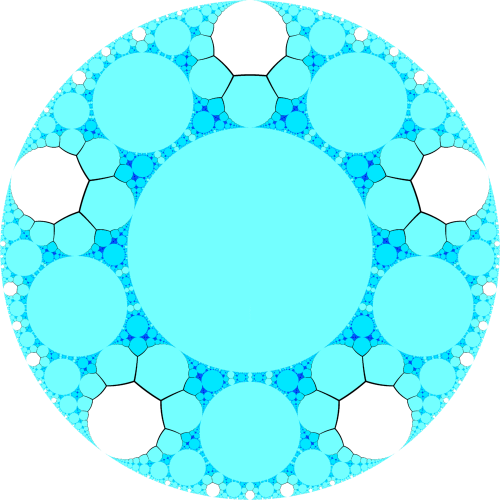} 
 	& \includegraphics[width=\w\textwidth]{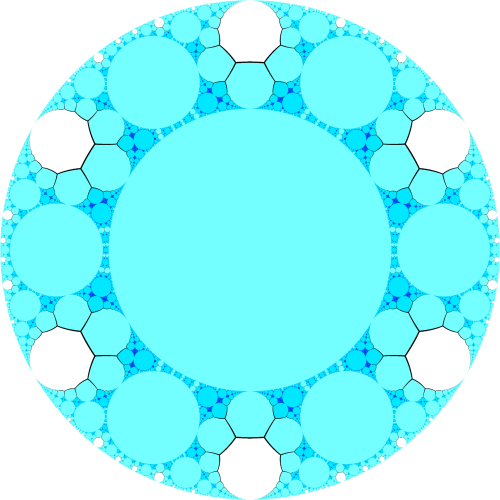} 
	& \includegraphics[width=\w\textwidth]{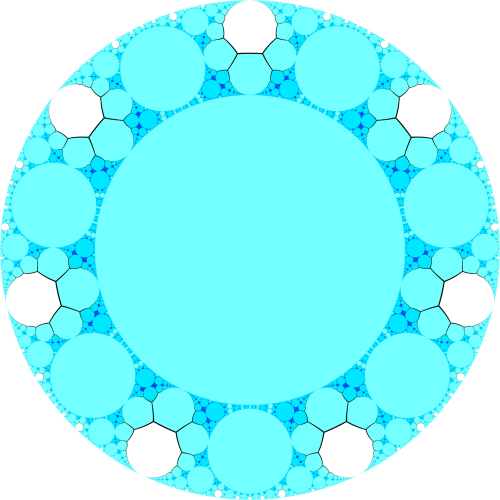} 
	& $\cdots$ 
	& \includegraphics[width=\w\textwidth]{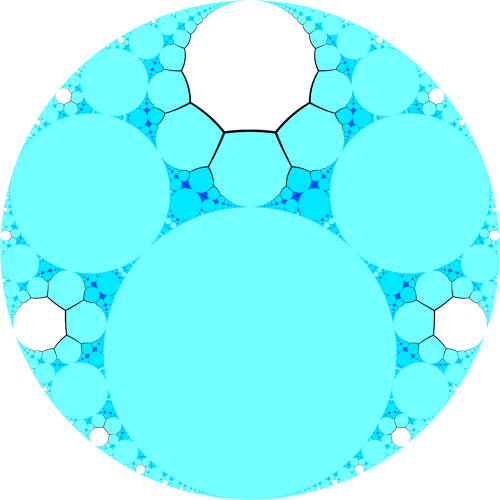}  \\

\end{tabular}
\vspace{15pt}
\caption{The $\{p,q,3\}$ honeycombs.}
\label{Table:pq3}
\end{table}

\begin{table}[htbp]
\centering
\begin{tabular}{V{0.05}|V{\w}V{\w}V{\w}V{\w}V{\w}V{0.06}V{\w}}

  $r\backslash{q}$ & 3 & 4 & 5 & 6 & 7 & $\cdots$ & $\infty$ \\	\hline \\ [-12pt]
    3	& \includegraphics[width=\w\textwidth]{Figures/noncompact/i33.png} 
	& \includegraphics[width=\w\textwidth]{Figures/noncompact/i43.png} 
	& \includegraphics[width=\w\textwidth]{Figures/noncompact/i53.png} 
 	& \includegraphics[width=\w\textwidth]{Figures/noncompact/i63.png} 
	& \includegraphics[width=\w\textwidth]{Figures/noncompact/i73.png} 
	& $\cdots$ 
	& \includegraphics[width=\w\textwidth]{Figures/noncompact/ii3.png}  \\
    4	& \includegraphics[width=\w\textwidth]{Figures/noncompact/i34.png} 
	& \includegraphics[width=\w\textwidth]{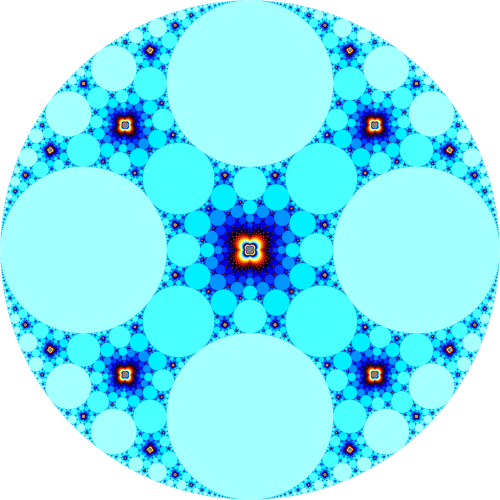} 
	& \includegraphics[width=\w\textwidth]{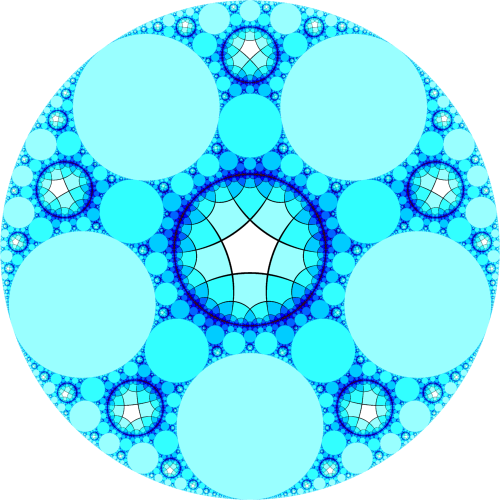} 
 	& \includegraphics[width=\w\textwidth]{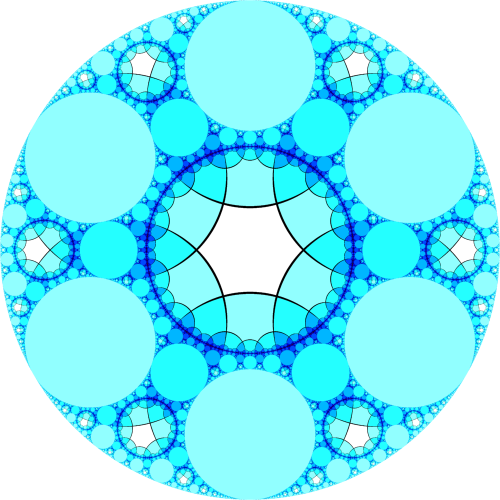} 
	& \includegraphics[width=\w\textwidth]{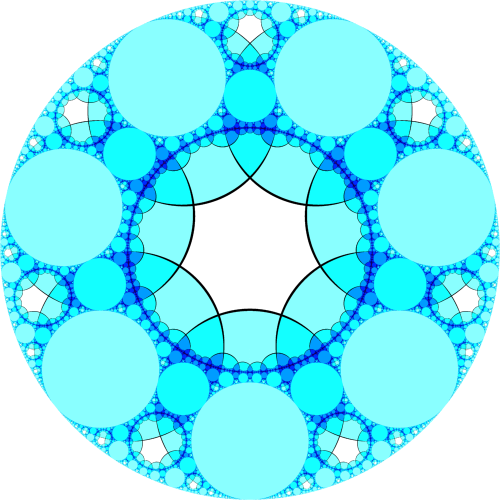} 
	& $\cdots$ 
	& \includegraphics[width=\w\textwidth]{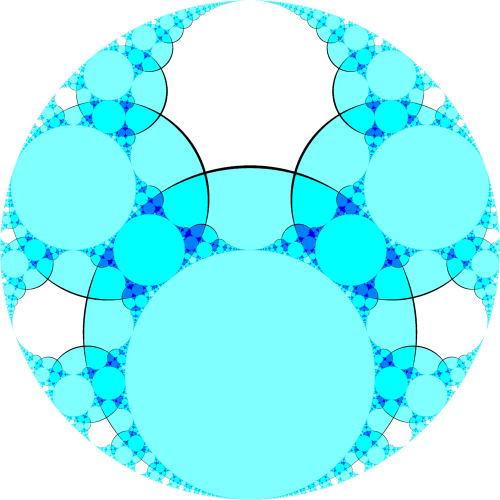}  \\
    5	& \includegraphics[width=\w\textwidth]{Figures/noncompact/i35.png} 
	& \includegraphics[width=\w\textwidth]{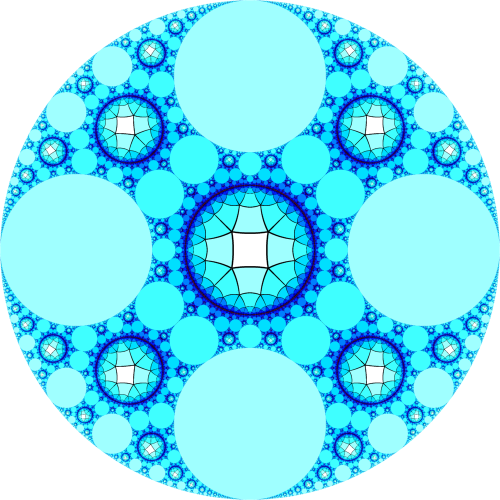} 
	& \includegraphics[width=\w\textwidth]{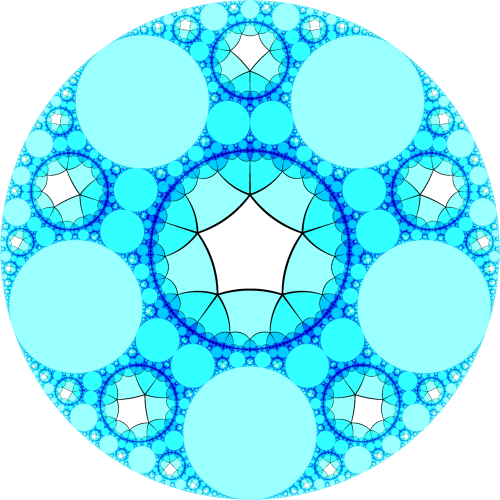} 
 	& \includegraphics[width=\w\textwidth]{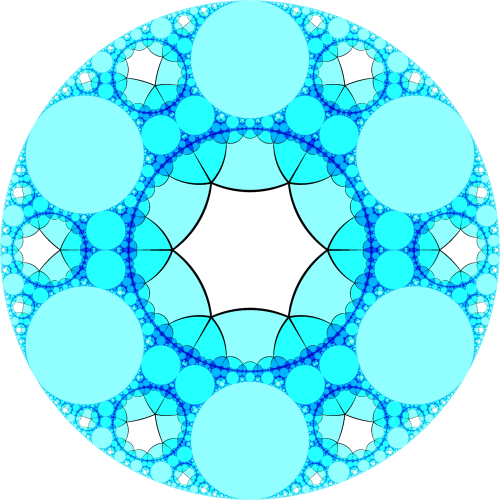} 
	& \includegraphics[width=\w\textwidth]{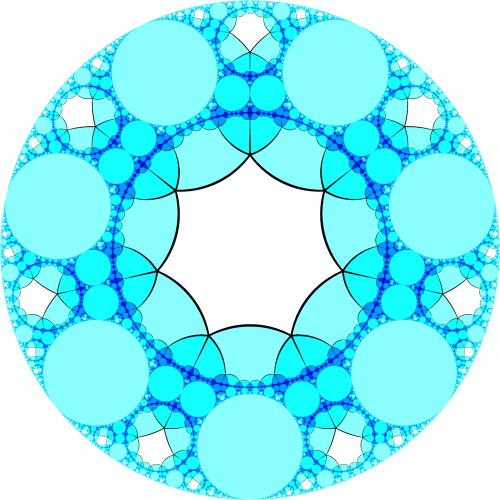} 
	& $\cdots$ 
	& \includegraphics[width=\w\textwidth]{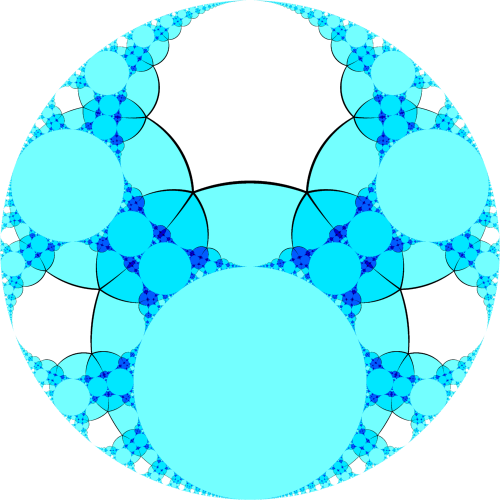}  \\
    6	& \includegraphics[width=\w\textwidth]{Figures/noncompact/i36.png} 
	& \includegraphics[width=\w\textwidth]{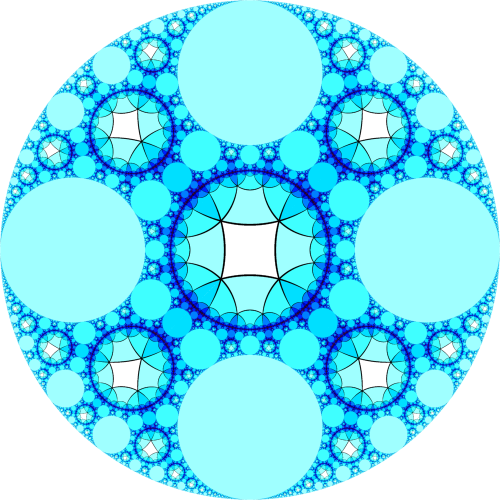} 
	& \includegraphics[width=\w\textwidth]{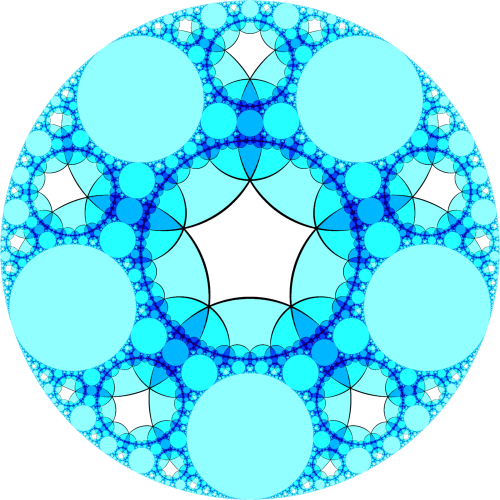} 
 	& \includegraphics[width=\w\textwidth]{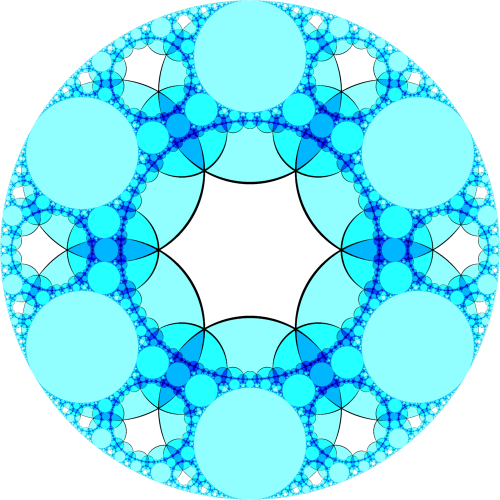} 
	& \includegraphics[width=\w\textwidth]{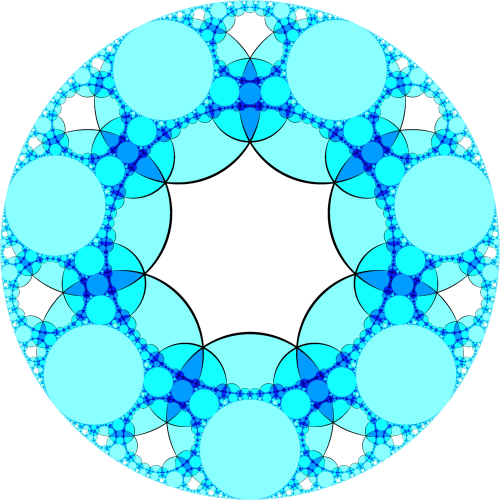} 
	& $\cdots$ 
	& \includegraphics[width=\w\textwidth]{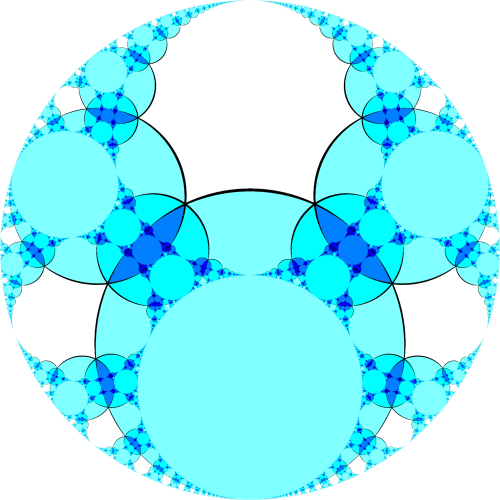}  \\
    7	& \includegraphics[width=\w\textwidth]{Figures/noncompact/i37.png} 
	& \includegraphics[width=\w\textwidth]{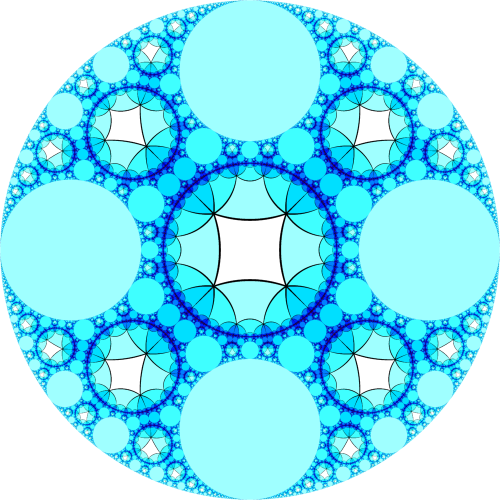} 
	& \includegraphics[width=\w\textwidth]{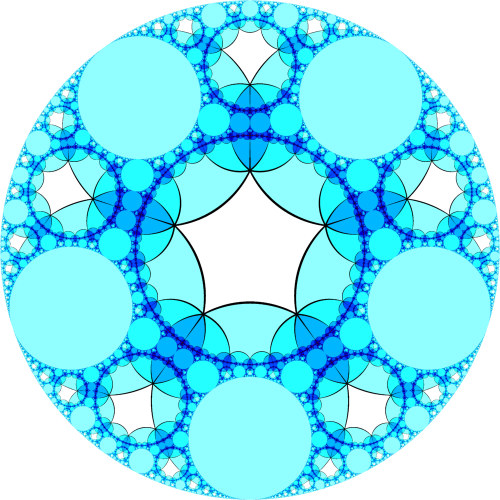} 
 	& \includegraphics[width=\w\textwidth]{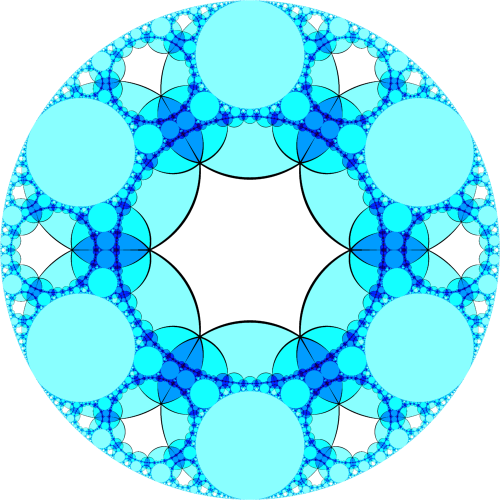} 
	& \includegraphics[width=\w\textwidth]{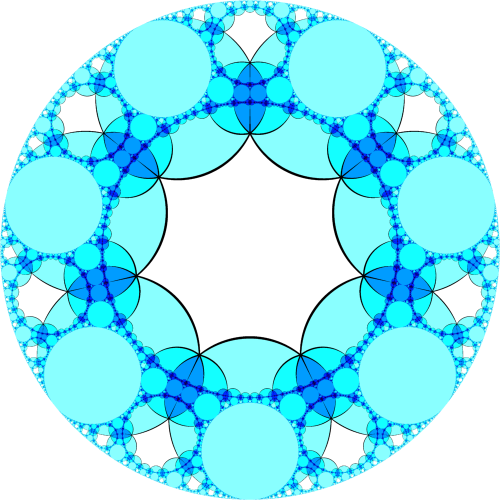} 
	& $\cdots$ 
	& \includegraphics[width=\w\textwidth]{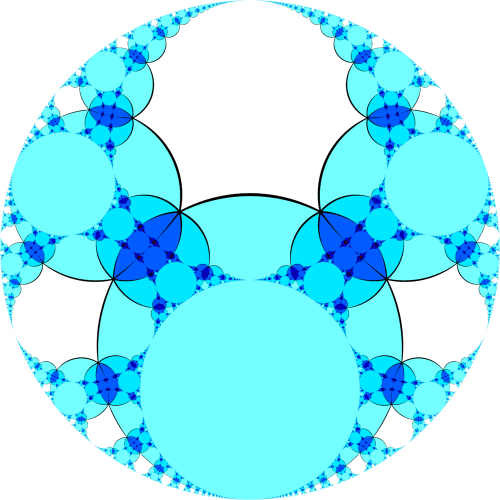}  \\ [-10pt] \\
  $\vdots$ &   $\vdots$  &   $\vdots$ &  $\vdots$ &   $\vdots$&   $\vdots$& $\ddots$ & $\vdots$  \\ \\
  $\infty$ 
	& \includegraphics[width=\w\textwidth]{Figures/noncompact/i3i.png} 
	& \includegraphics[width=\w\textwidth]{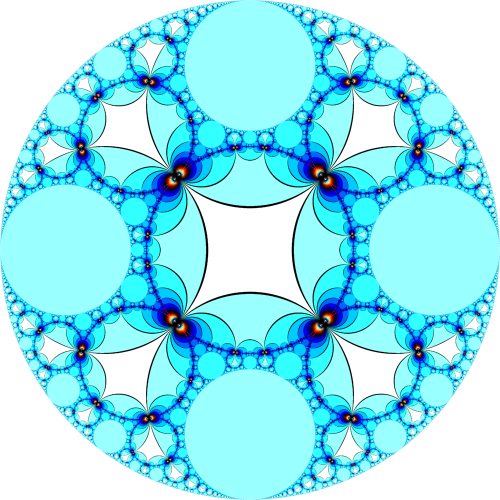} 
	& \includegraphics[width=\w\textwidth]{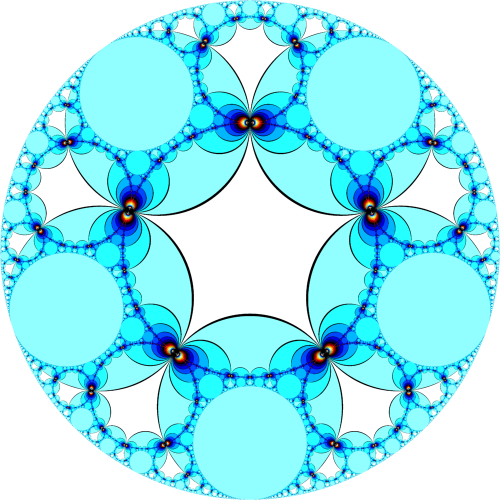} 
 	& \includegraphics[width=\w\textwidth]{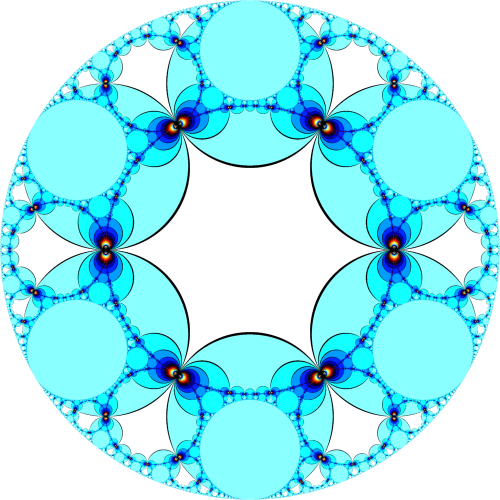} 
	& \includegraphics[width=\w\textwidth]{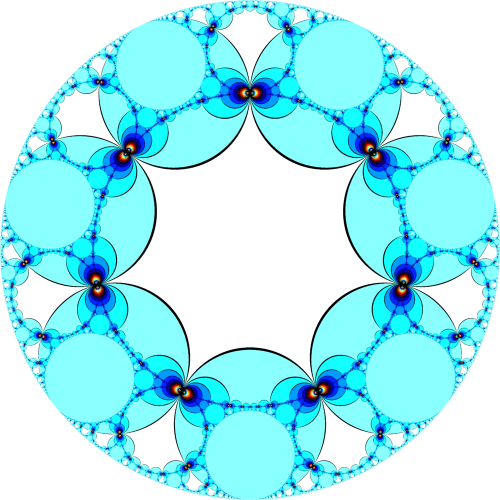} 
	& $\cdots$ 
	& \includegraphics[width=\w\textwidth]{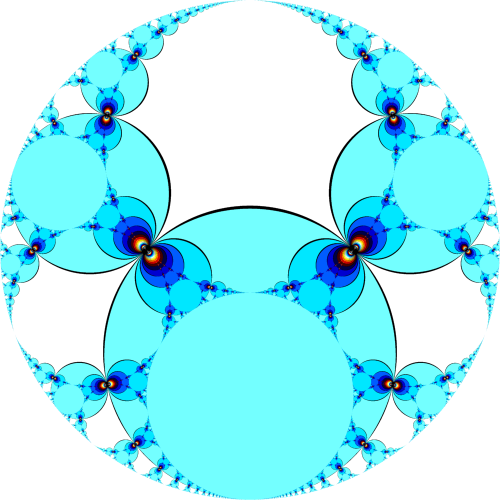}  \\

\end{tabular}
\vspace{15pt}
\caption{The $\{\infty,q,r\}$ honeycombs.}
\label{Table:iqr}
\end{table}

\begin{table}[htbp]
\centering
\begin{tabular}{V{0.05}|V{\w}V{\w}V{\w}V{\w}V{\w}V{0.06}V{\w}}

  $r\backslash{p}$ & 3 & 4 & 5 & 6 & 7 & $\cdots$ & $\infty$ \\	\hline \\ [-12pt]
    3	& \includegraphics[width=\w\textwidth]{Figures/noncompact/3i3.png} 
	& \includegraphics[width=\w\textwidth]{Figures/noncompact/4i3.png} 
	& \includegraphics[width=\w\textwidth]{Figures/noncompact/5i3.png} 
 	& \includegraphics[width=\w\textwidth]{Figures/noncompact/6i3.png} 
	& \includegraphics[width=\w\textwidth]{Figures/noncompact/7i3.png} 
	& $\cdots$ 
	& \includegraphics[width=\w\textwidth]{Figures/noncompact/ii3.png}  \\
    4	& \includegraphics[width=\w\textwidth]{Figures/noncompact/3i4.png} 
	& \includegraphics[width=\w\textwidth]{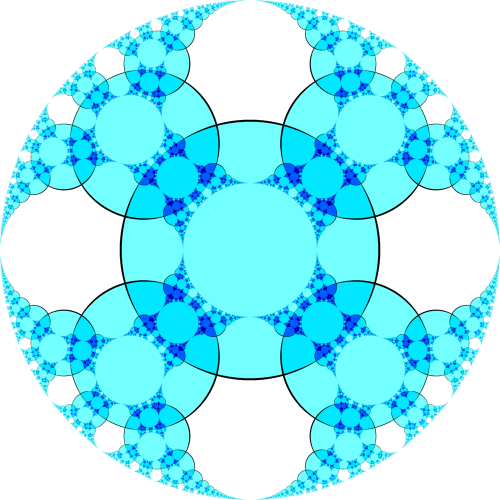} 
	& \includegraphics[width=\w\textwidth]{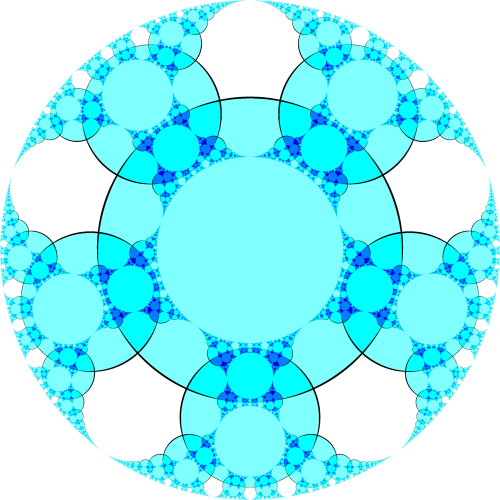} 
 	& \includegraphics[width=\w\textwidth]{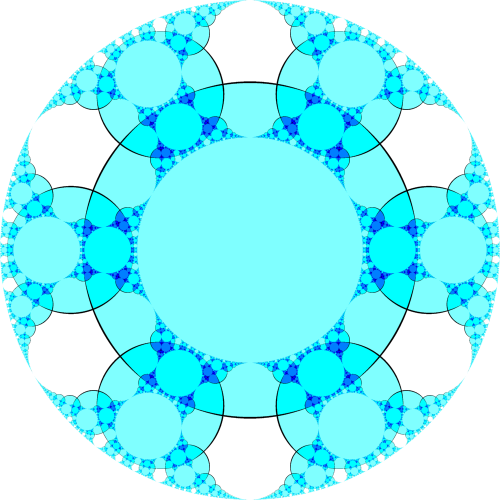} 
	& \includegraphics[width=\w\textwidth]{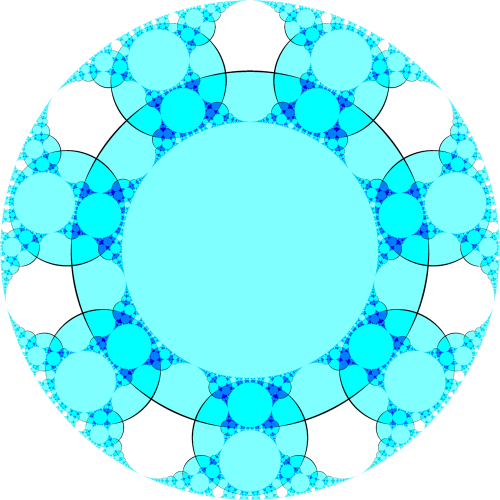} 
	& $\cdots$ 
	& \includegraphics[width=\w\textwidth]{Figures/noncompact/ii4.png}  \\
    5	& \includegraphics[width=\w\textwidth]{Figures/noncompact/3i5.png} 
	& \includegraphics[width=\w\textwidth]{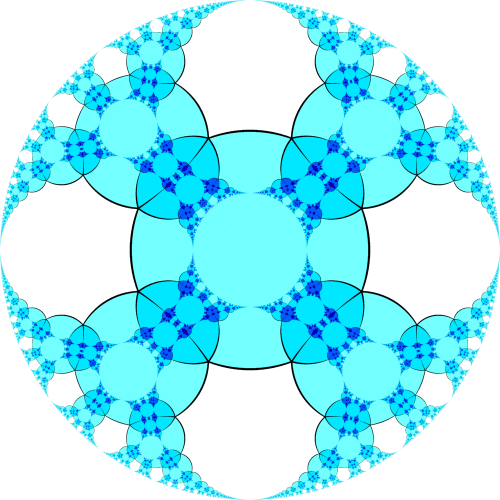} 
	& \includegraphics[width=\w\textwidth]{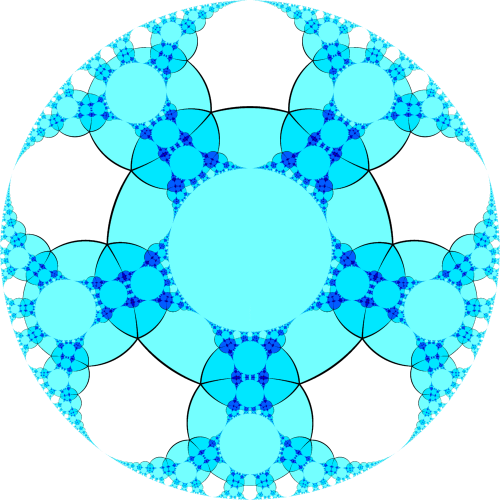} 
 	& \includegraphics[width=\w\textwidth]{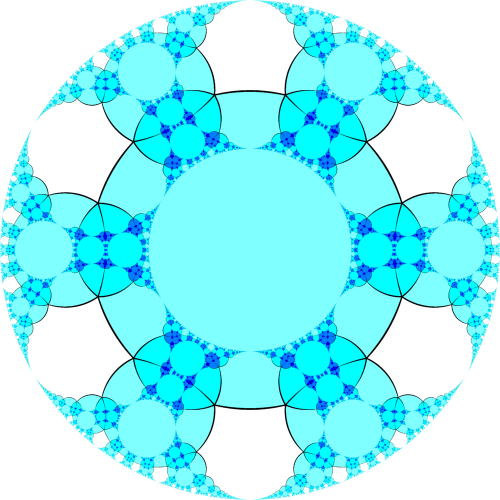} 
	& \includegraphics[width=\w\textwidth]{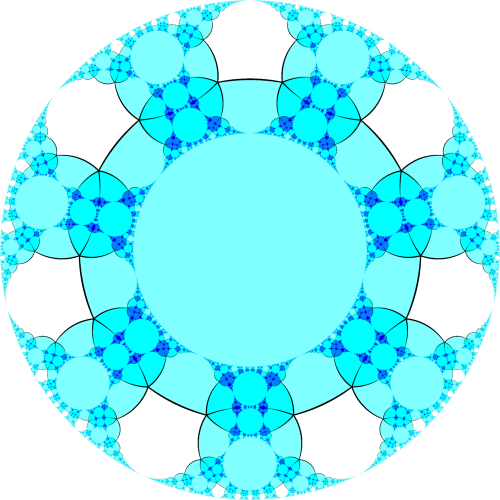} 
	& $\cdots$ 
	& \includegraphics[width=\w\textwidth]{Figures/noncompact/ii5.png}  \\
    6	& \includegraphics[width=\w\textwidth]{Figures/noncompact/3i6.png} 
	& \includegraphics[width=\w\textwidth]{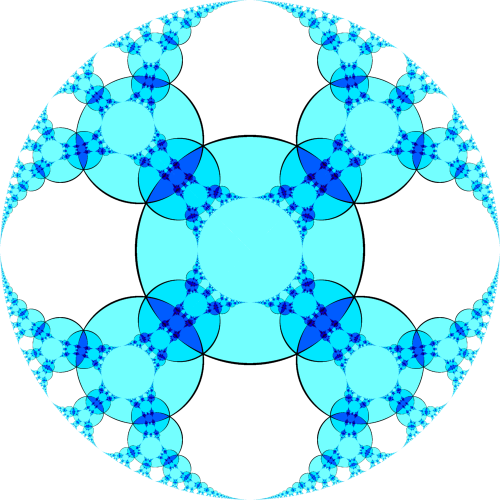} 
	& \includegraphics[width=\w\textwidth]{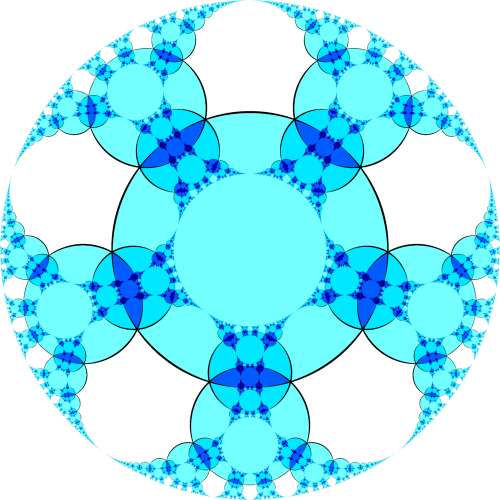} 
 	& \includegraphics[width=\w\textwidth]{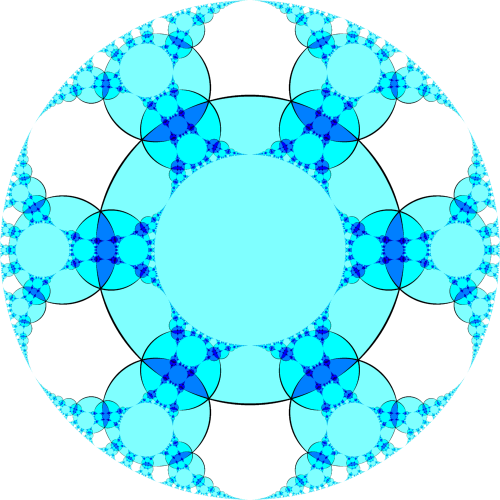} 
	& \includegraphics[width=\w\textwidth]{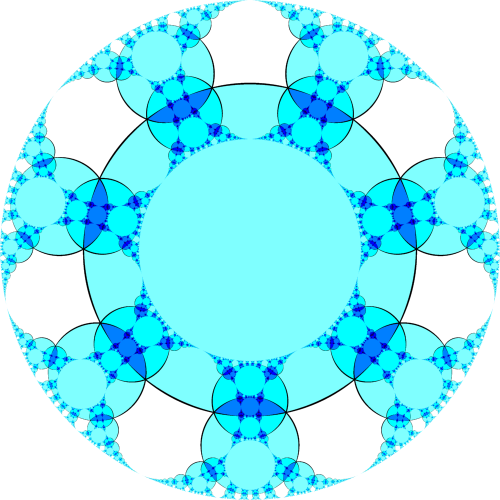} 
	& $\cdots$ 
	& \includegraphics[width=\w\textwidth]{Figures/noncompact/ii6.png}  \\
    7	& \includegraphics[width=\w\textwidth]{Figures/noncompact/3i7.png} 
	& \includegraphics[width=\w\textwidth]{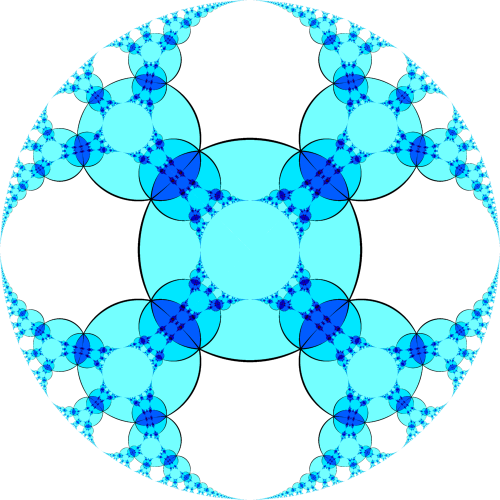} 
	& \includegraphics[width=\w\textwidth]{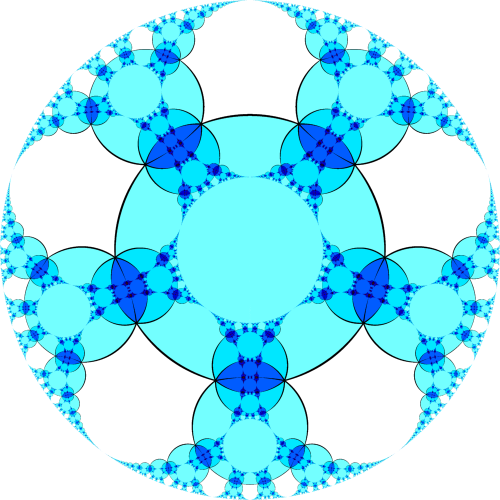} 
 	& \includegraphics[width=\w\textwidth]{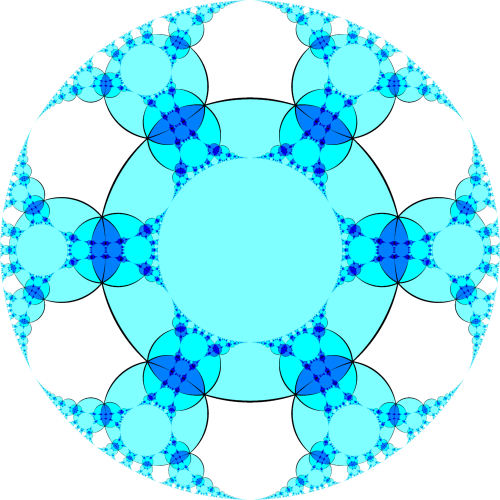} 
	& \includegraphics[width=\w\textwidth]{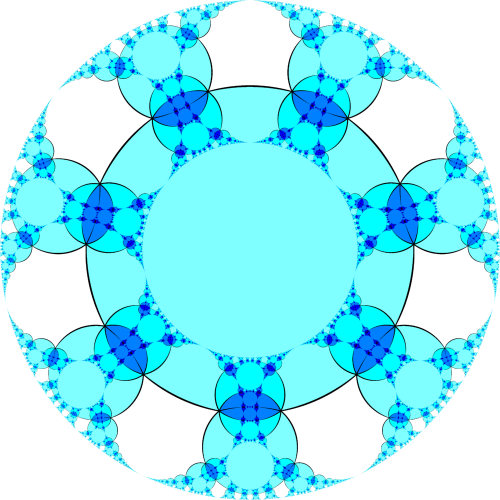} 
	& $\cdots$ 
	& \includegraphics[width=\w\textwidth]{Figures/noncompact/ii7.png}  \\ [-10pt] \\
  $\vdots$ &   $\vdots$  &   $\vdots$ &  $\vdots$ &   $\vdots$&   $\vdots$& $\ddots$ & $\vdots$  \\ \\
  $\infty$ 
	& \includegraphics[width=\w\textwidth]{Figures/noncompact/3ii.png} 
	& \includegraphics[width=\w\textwidth]{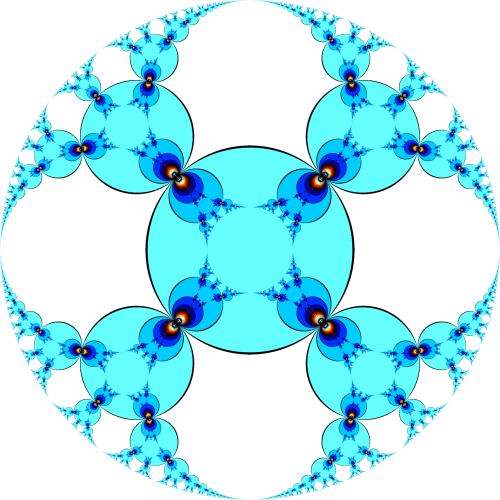} 
	& \includegraphics[width=\w\textwidth]{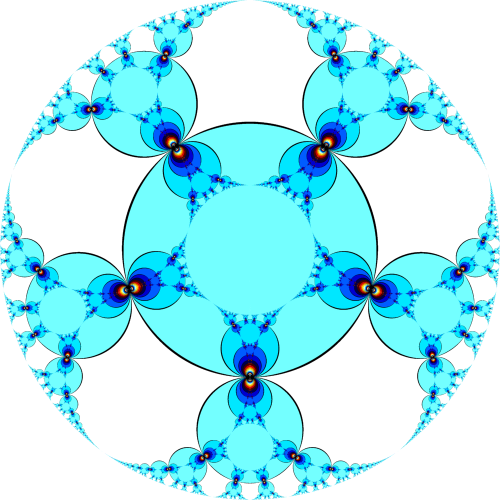} 
 	& \includegraphics[width=\w\textwidth]{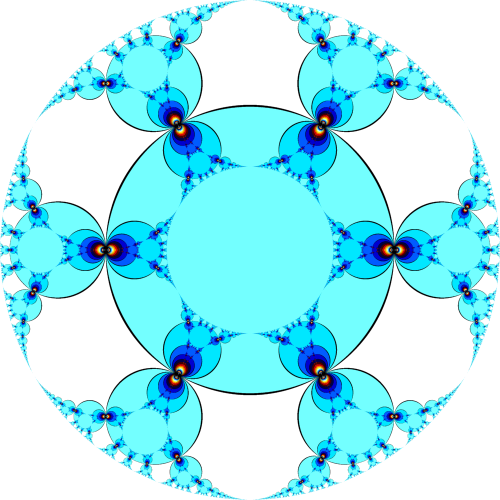} 
	& \includegraphics[width=\w\textwidth]{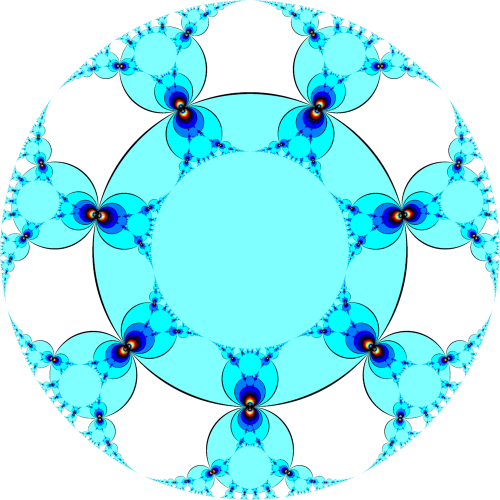} 
	& $\cdots$ 
	& \includegraphics[width=\w\textwidth]{Figures/noncompact/iii.png}  \\

\end{tabular}
\vspace{15pt}
\caption{The $\{p,\infty,r\}$ honeycombs.}
\label{Table:pir}
\end{table}

\begin{table}[htbp]
\centering
\begin{tabular}{V{0.05}|V{\w}V{\w}V{\w}V{\w}V{\w}V{0.06}V{\w}}

  $q\backslash{p}$ & 3 & 4 & 5 & 6 & 7 & $\cdots$ & $\infty$ \\	\hline \\ [-12pt]
    3	& \includegraphics[width=\w\textwidth]{Figures/noncompact/33i.png} 
	& \includegraphics[width=\w\textwidth]{Figures/noncompact/43i.png} 
	& \includegraphics[width=\w\textwidth]{Figures/noncompact/53i.png} 
 	& \includegraphics[width=\w\textwidth]{Figures/noncompact/63i.png} 
	& \includegraphics[width=\w\textwidth]{Figures/noncompact/73i.png} 
	& $\cdots$ 
	& \includegraphics[width=\w\textwidth]{Figures/noncompact/i3i.png}  \\
    4	& \includegraphics[width=\w\textwidth]{Figures/noncompact/34i.png} 
	& \includegraphics[width=\w\textwidth]{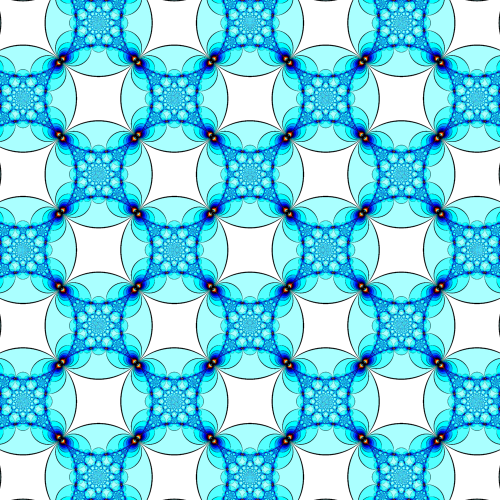} 
	& \includegraphics[width=\w\textwidth]{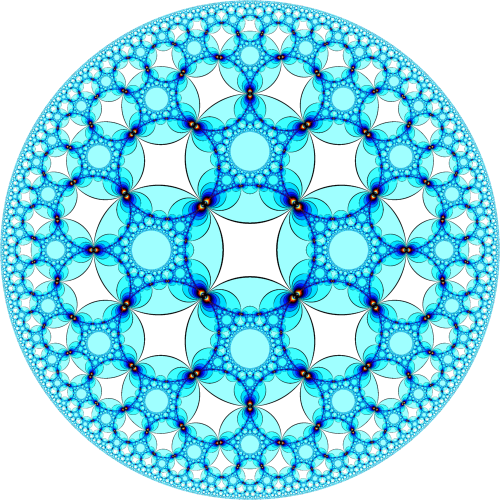} 
 	& \includegraphics[width=\w\textwidth]{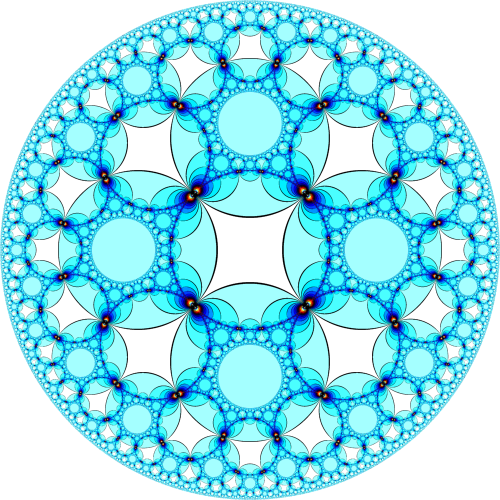} 
	& \includegraphics[width=\w\textwidth]{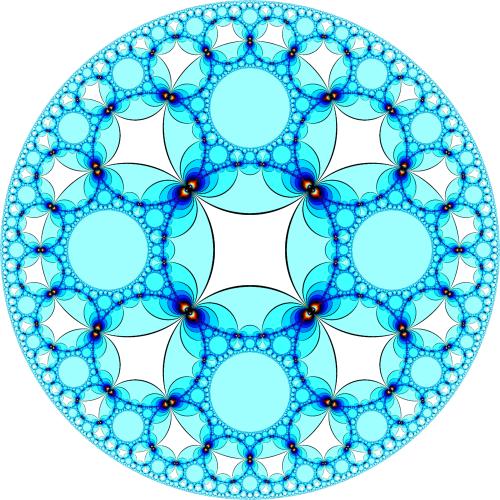} 
	& $\cdots$ 
	& \includegraphics[width=\w\textwidth]{Figures/noncompact/i4i.png}  \\
    5	& \includegraphics[width=\w\textwidth]{Figures/noncompact/35i.png} 
	& \includegraphics[width=\w\textwidth]{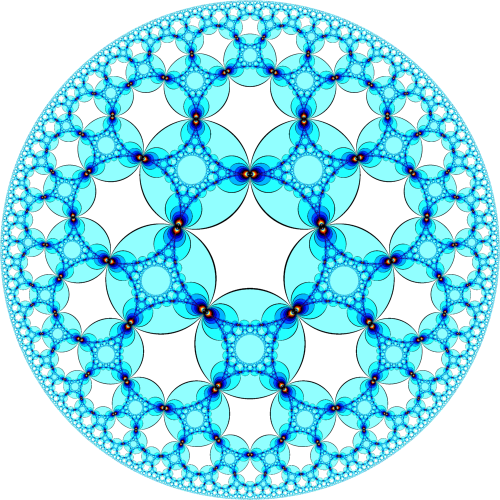} 
	& \includegraphics[width=\w\textwidth]{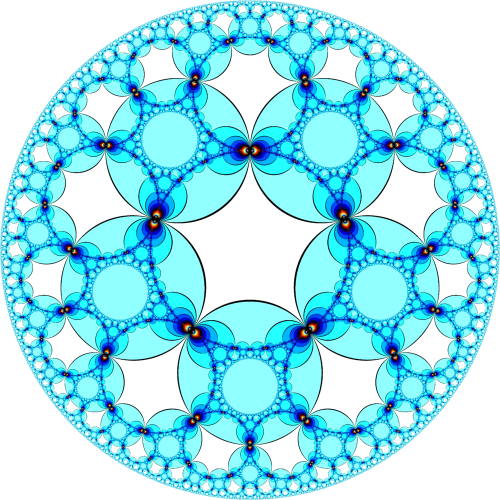} 
 	& \includegraphics[width=\w\textwidth]{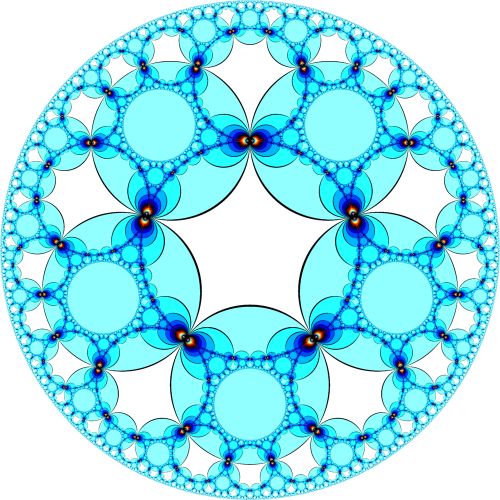} 
	& \includegraphics[width=\w\textwidth]{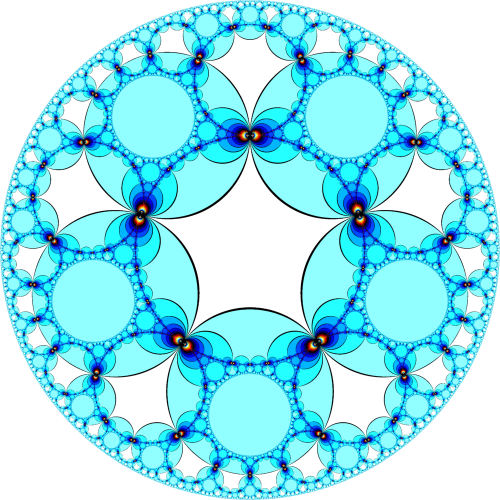} 
	& $\cdots$ 
	& \includegraphics[width=\w\textwidth]{Figures/noncompact/i5i.png}  \\
    6	& \includegraphics[width=\w\textwidth]{Figures/noncompact/36i.png} 
	& \includegraphics[width=\w\textwidth]{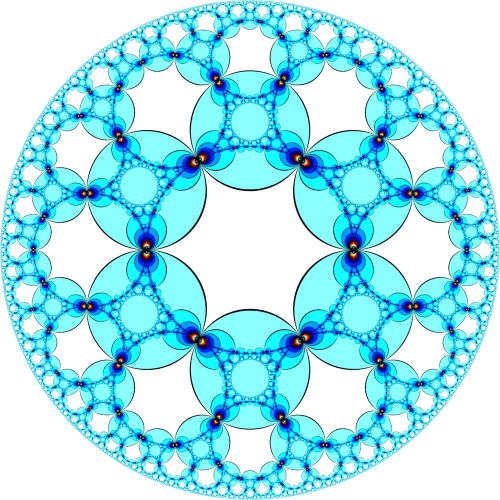} 
	& \includegraphics[width=\w\textwidth]{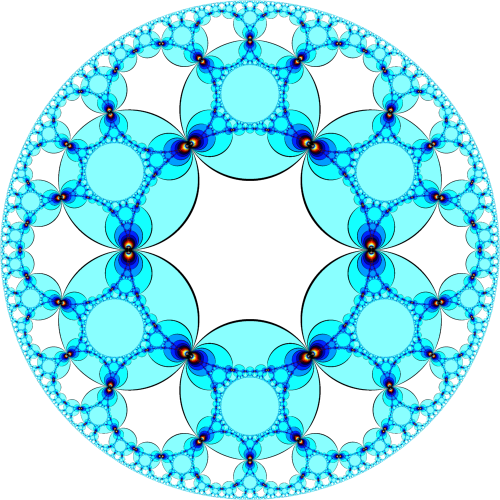} 
 	& \includegraphics[width=\w\textwidth]{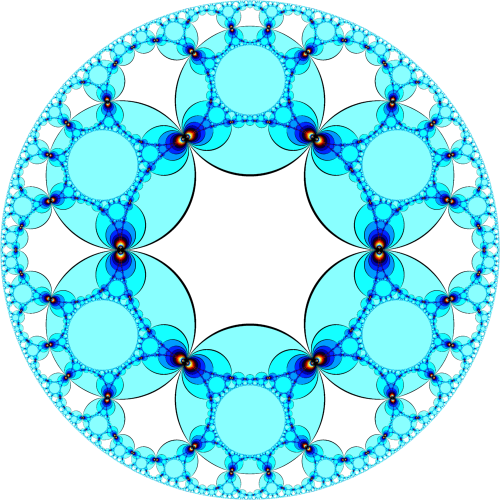} 
	& \includegraphics[width=\w\textwidth]{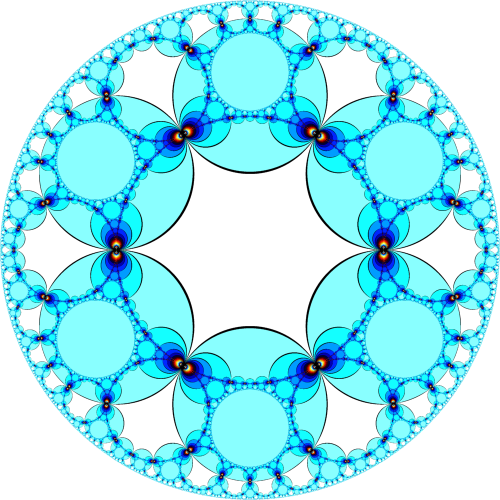} 
	& $\cdots$ 
	& \includegraphics[width=\w\textwidth]{Figures/noncompact/i6i.png}  \\
    7	& \includegraphics[width=\w\textwidth]{Figures/noncompact/37i.png} 
	& \includegraphics[width=\w\textwidth]{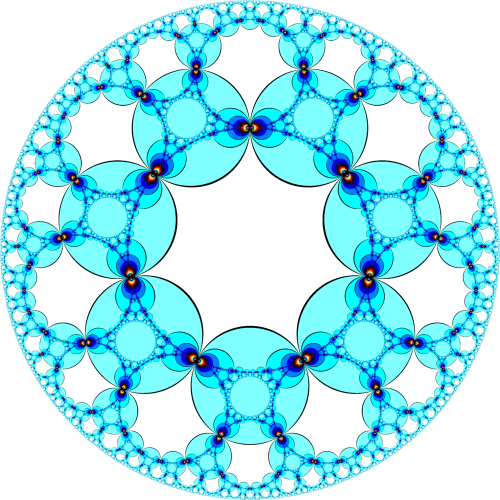} 
	& \includegraphics[width=\w\textwidth]{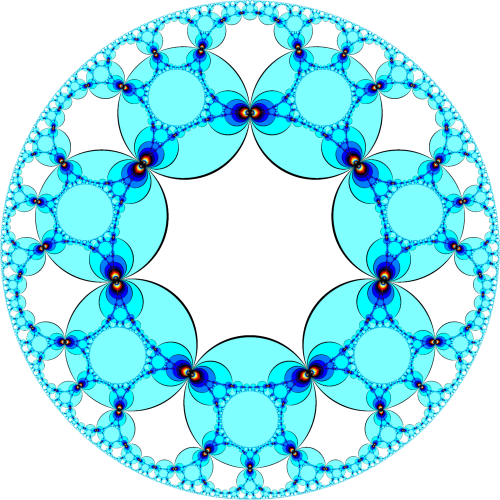} 
 	& \includegraphics[width=\w\textwidth]{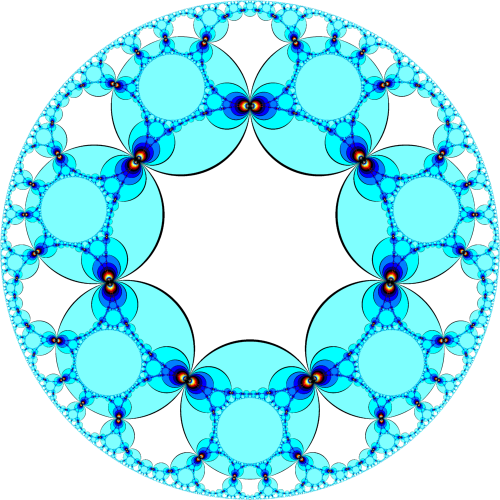} 
	& \includegraphics[width=\w\textwidth]{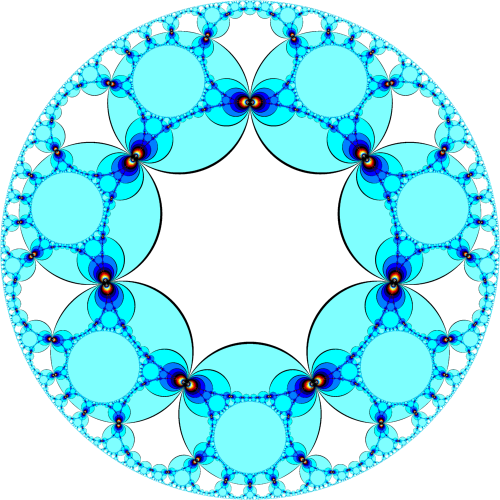} 
	& $\cdots$ 
	& \includegraphics[width=\w\textwidth]{Figures/noncompact/i7i.png}  \\ [-10pt] \\
  $\vdots$ &   $\vdots$  &   $\vdots$ &  $\vdots$ &   $\vdots$&   $\vdots$& $\ddots$ & $\vdots$  \\ \\
  $\infty$ 
	& \includegraphics[width=\w\textwidth]{Figures/noncompact/3ii.png} 
	& \includegraphics[width=\w\textwidth]{Figures/noncompact/4ii.png} 
	& \includegraphics[width=\w\textwidth]{Figures/noncompact/5ii.png} 
 	& \includegraphics[width=\w\textwidth]{Figures/noncompact/6ii.png} 
	& \includegraphics[width=\w\textwidth]{Figures/noncompact/7ii.png} 
	& $\cdots$ 
	& \includegraphics[width=\w\textwidth]{Figures/noncompact/iii.png}  \\

\end{tabular}
\vspace{15pt}
\caption{The $\{p,q,\infty\}$ honeycombs.}
\label{Table:pqi}
\end{table}

\def\w{0.2}
\setlength\tabcolsep{3pt}

\begin{table}[htbp]
\centering
\begin{tabular}{V{0.1}|V{\w}V{\w}V{\w}V{\w}}

  $p\backslash{\{q,r\}}$ & $\{3,6\}$ & $\{4,4\}$ & $\{6,3\}$ \\	\hline \\ [-12pt]
    3	& \includegraphics[width=\w\textwidth]{Figures/compact/336.png} 
	& \includegraphics[width=\w\textwidth]{Figures/compact/344.png} 
	& \includegraphics[width=\w\textwidth]{Figures/compact/363.png}  
\\
    4	& \includegraphics[width=\w\textwidth]{Figures/compact/436.png} 
	& \includegraphics[width=\w\textwidth]{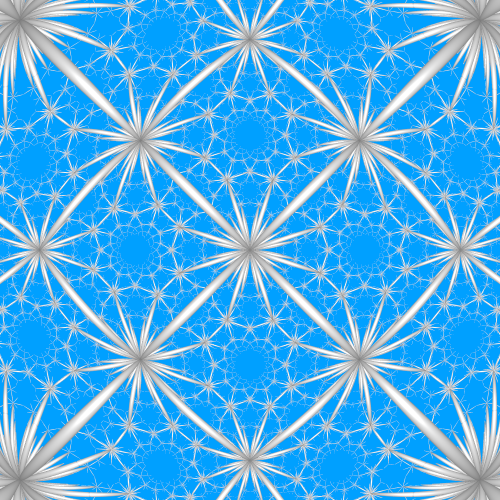} 
	& \includegraphics[width=\w\textwidth]{Figures/noncompact/463.png}  
\\
    5	& \includegraphics[width=\w\textwidth]{Figures/compact/536.png} 
	& \includegraphics[width=\w\textwidth]{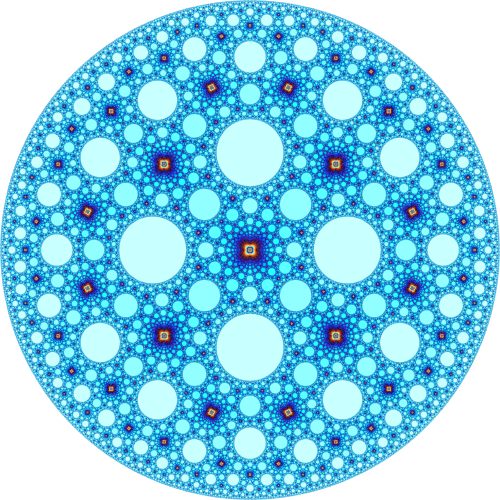} 
	& \includegraphics[width=\w\textwidth]{Figures/noncompact/563.png}  
\\
    6	& \includegraphics[width=\w\textwidth]{Figures/compact/636.png} 
	& \includegraphics[width=\w\textwidth]{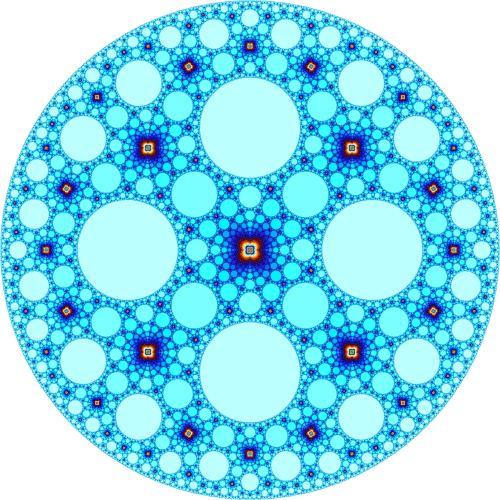} 
	& \includegraphics[width=\w\textwidth]{Figures/noncompact/663.png}  
\\
    7	& \includegraphics[width=\w\textwidth]{Figures/noncompact/736.png} 
	& \includegraphics[width=\w\textwidth]{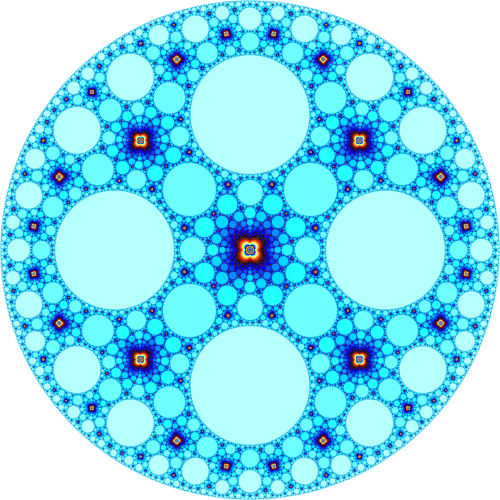} 
	& \includegraphics[width=\w\textwidth]{Figures/noncompact/763.png}  \\ [-10pt] \\
  $\vdots$ &   $\vdots$  &   $\vdots$ &  $\vdots$  \\ \\
  $\infty$ 
	& \includegraphics[width=\w\textwidth]{Figures/noncompact/i36.png} 
	& \includegraphics[width=\w\textwidth]{Figures/noncompact/i44.png} 
	& \includegraphics[width=\w\textwidth]{Figures/noncompact/i63.png}  \\

\end{tabular}
\vspace{15pt}
\caption{Honeycombs with ideal vertices. }
\label{Table:ideal_verts}
\end{table}

\begin{table}[htbp]
\centering
\begin{tabular}{V{0.1}|V{\w}V{\w}V{\w}V{\w}}

  $r\backslash{\{p,q\}}$ & $\{3,6\}$ & $\{4,4\}$ & $\{6,3\}$ \\	\hline \\ [-12pt]
    3	& \includegraphics[width=\w\textwidth]{Figures/compact/363.png} 
	& \includegraphics[width=\w\textwidth]{Figures/compact/443.png} 
	& \includegraphics[width=\w\textwidth]{Figures/compact/633.png}  
\\
    4	& \includegraphics[width=\w\textwidth]{Figures/noncompact/364.png} 
	& \includegraphics[width=\w\textwidth]{Figures/compact/444.png} 
	& \includegraphics[width=\w\textwidth]{Figures/compact/634.png}  
\\
    5	& \includegraphics[width=\w\textwidth]{Figures/noncompact/365.png} 
	& \includegraphics[width=\w\textwidth]{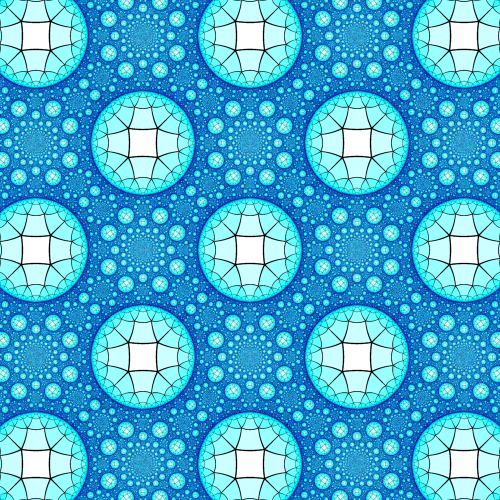} 
	& \includegraphics[width=\w\textwidth]{Figures/compact/635.png}  
\\
    6	& \includegraphics[width=\w\textwidth]{Figures/noncompact/366.png} 
	& \includegraphics[width=\w\textwidth]{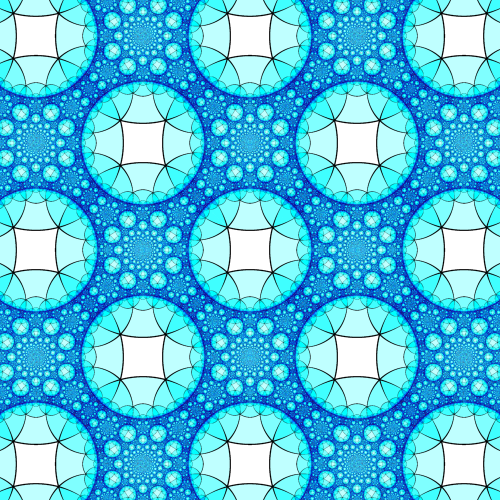} 
	& \includegraphics[width=\w\textwidth]{Figures/compact/636.png}  
\\
    7	& \includegraphics[width=\w\textwidth]{Figures/noncompact/367.png} 
	& \includegraphics[width=\w\textwidth]{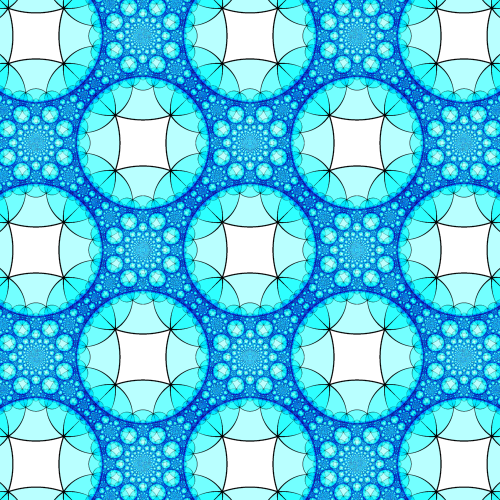} 
	& \includegraphics[width=\w\textwidth]{Figures/noncompact/637.png}  \\ [-10pt] \\
  $\vdots$ &   $\vdots$  &   $\vdots$ &  $\vdots$  \\ \\
  $\infty$ 
	& \includegraphics[width=\w\textwidth]{Figures/noncompact/36i.png} 
	& \includegraphics[width=\w\textwidth]{Figures/noncompact/44i.png} 
	& \includegraphics[width=\w\textwidth]{Figures/noncompact/63i.png}  \\

\end{tabular}
\vspace{15pt}
\caption{Honeycombs with ideal cells. }
\label{Table:ideal_cells}
\end{table}

\fi  

\FloatBarrier
\section{Artwork}
\label{Sec:artwork}
\subsection{3D printed sculptures}

A hyperbolic honeycomb represented using only edges as in Figure \ref{Fig:435_3} is difficult to interpret from a two-dimensional picture, but a 3D printed sculpture is much more useful. See Figures \ref{Fig:435_535} through \ref{Fig:444}. Also see Appendices \ref{Sec:dupin} and \ref{Sec:3d_printing} for technical details on designing these sculptures.

\begin{figure}[htbp]
\centering 
\hspace{-.7cm}
\subfloat[]
{
\ifimages
\includegraphics[width=0.43\textwidth]{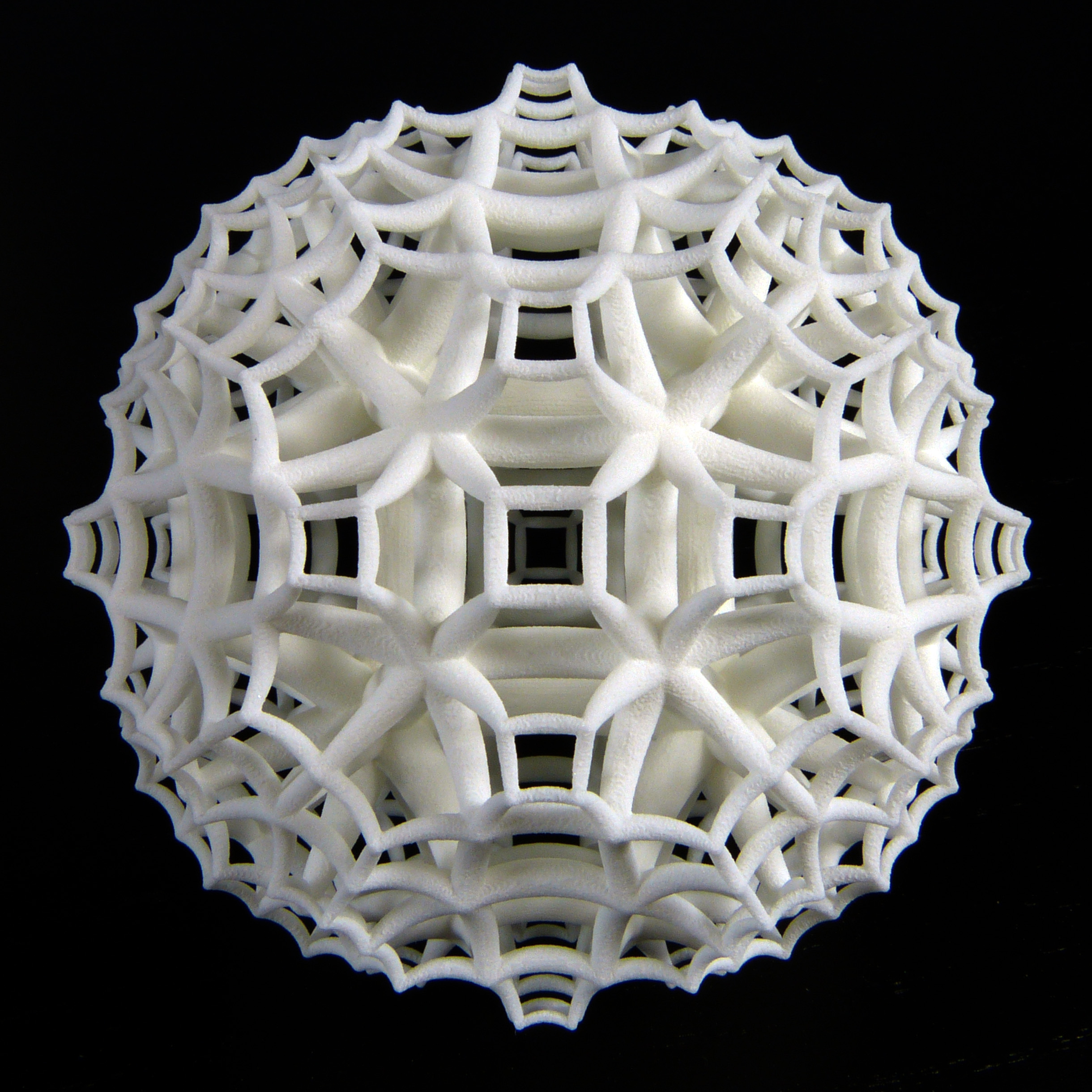}
\fi
\label{Fig:435}
} 
\hspace{.7cm}
\subfloat[]
{
\ifimages
\includegraphics[width=0.43\textwidth]{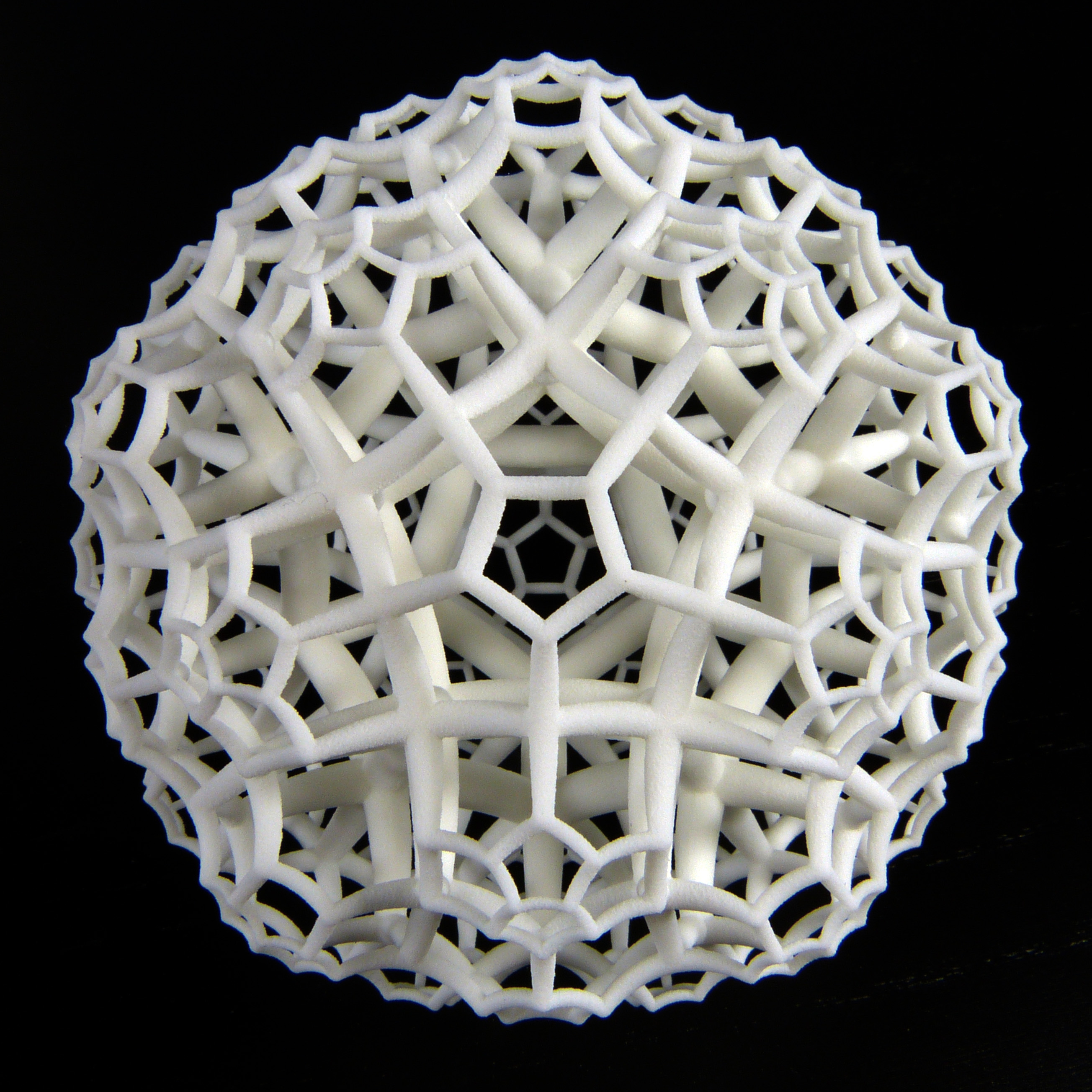}
\fi
\label{Fig:534}
}
\caption{The dual $\{4,3,5\}$ and $\{5,3,4\}$ honeycombs.}
\label{Fig:435_535}
\end{figure}

\begin{figure}[htbp]
\centering 
\hspace{-.7cm}
\subfloat[]
{
\ifimages
\includegraphics[width=0.43\textwidth]{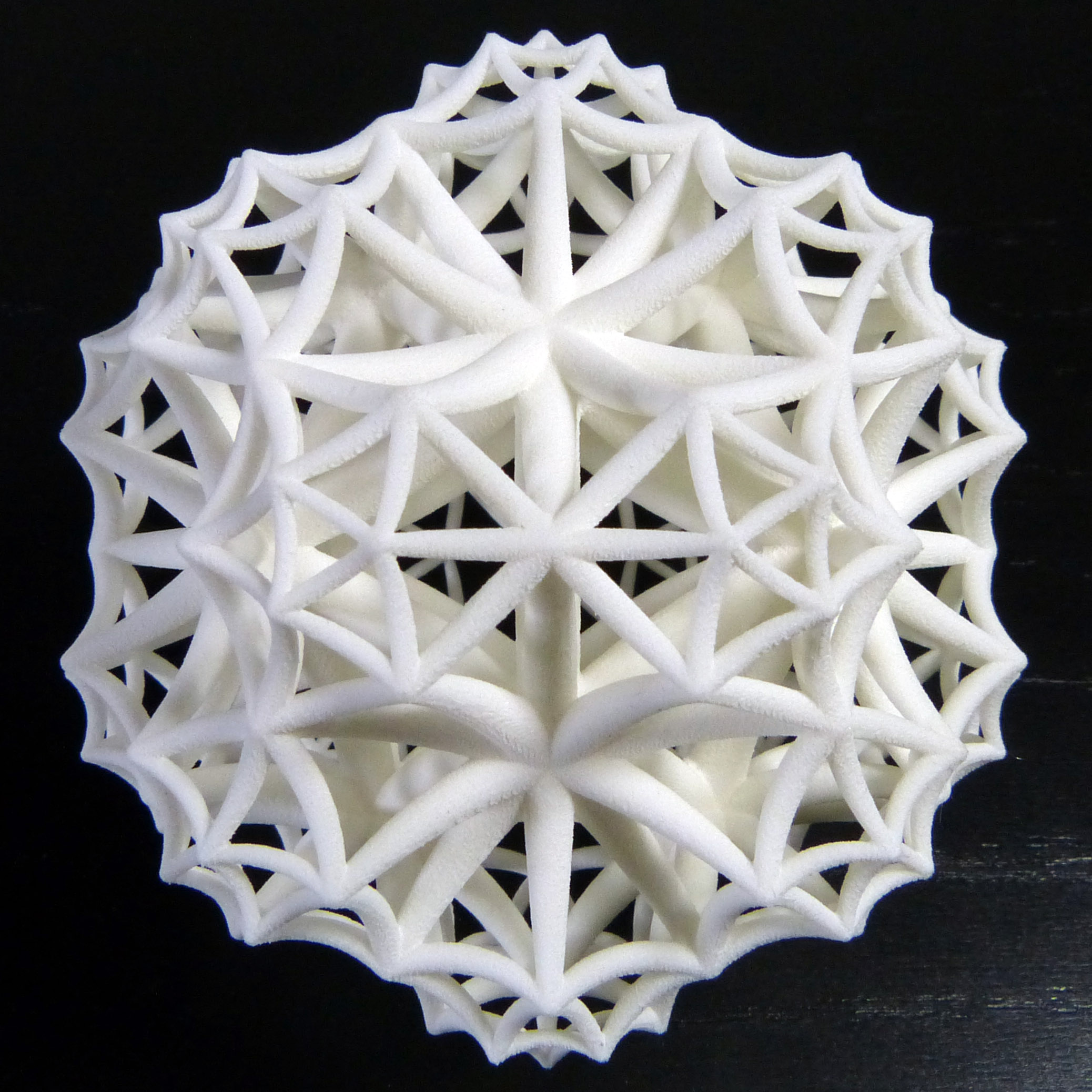}
\fi
\label{Fig:353_H3_honeycomb_sym22}
} 
\hspace{.7cm}
\subfloat[]
{
\ifimages
\includegraphics[width=0.43\textwidth]{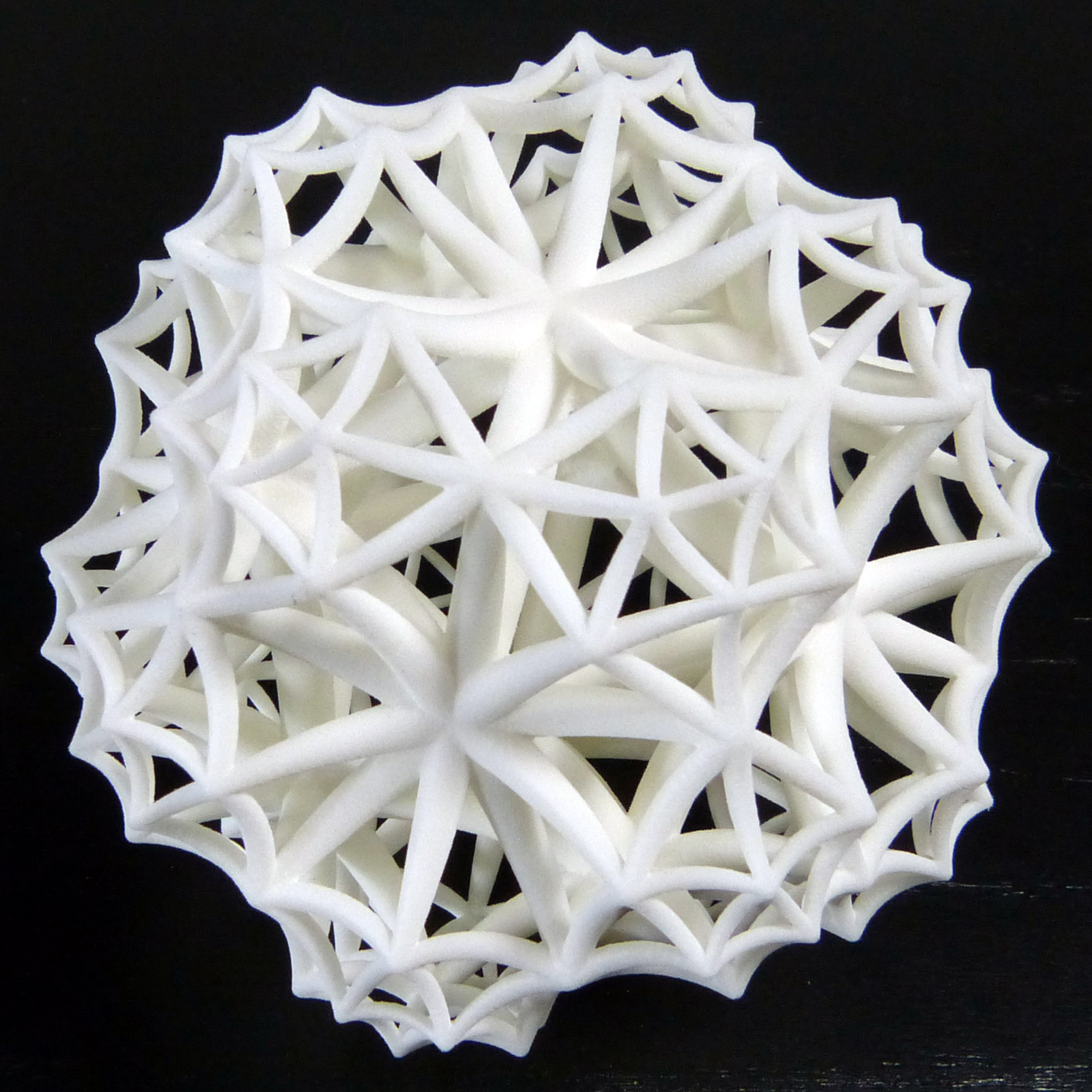}
\fi
\label{Fig:353_H3_honeycomb_1}
}
\caption{Two views of the self-dual $\{3,5,3\}$ honeycomb.}
\label{Fig:353}
\end{figure}

\begin{figure}[htbp]
\centering 
\hspace{-.7cm}
\subfloat[]
{
\ifimages
\includegraphics[width=0.43\textwidth]{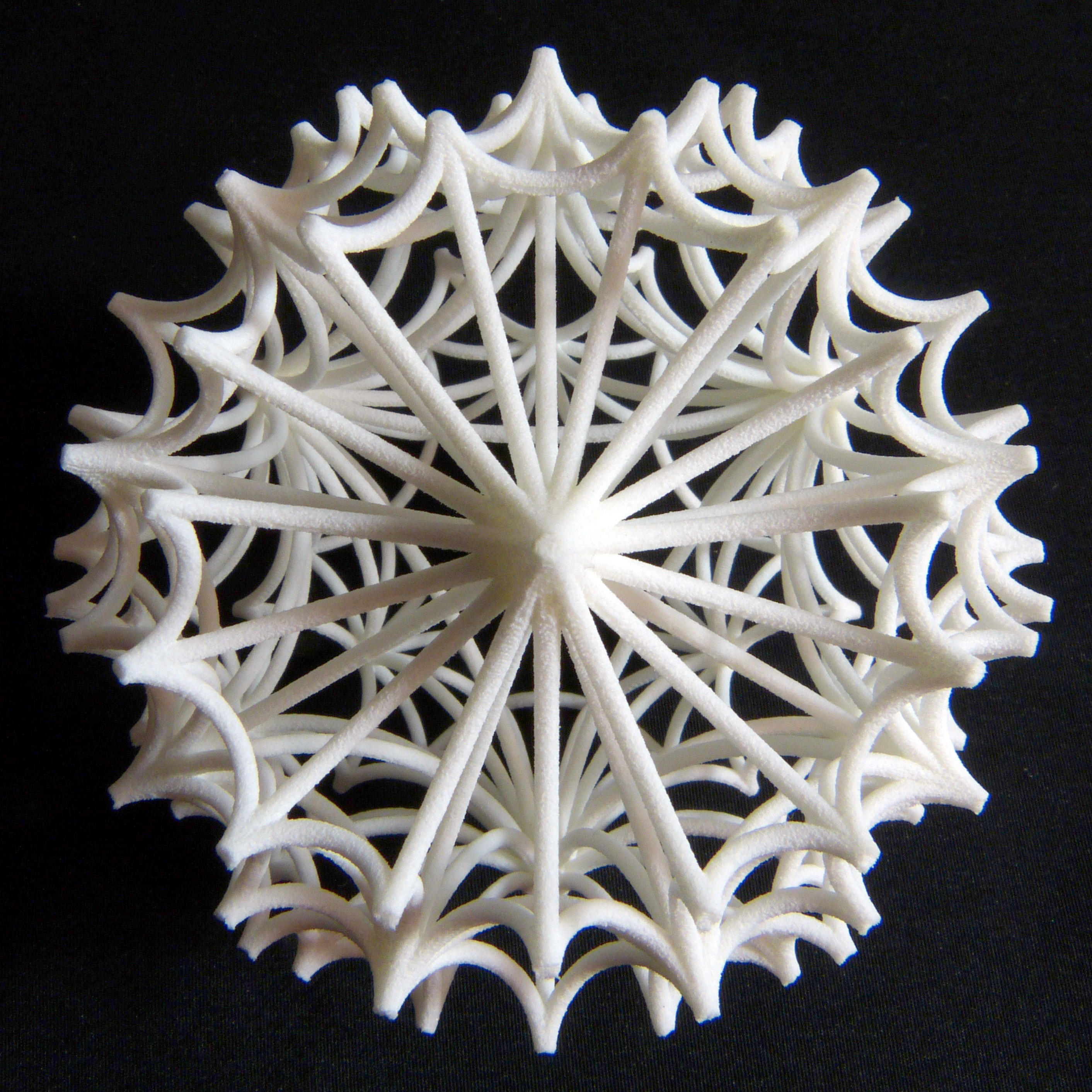}
\fi
\label{Fig:336}
} 
\hspace{.7cm}
\subfloat[]
{
\ifimages
\includegraphics[width=0.43\textwidth]{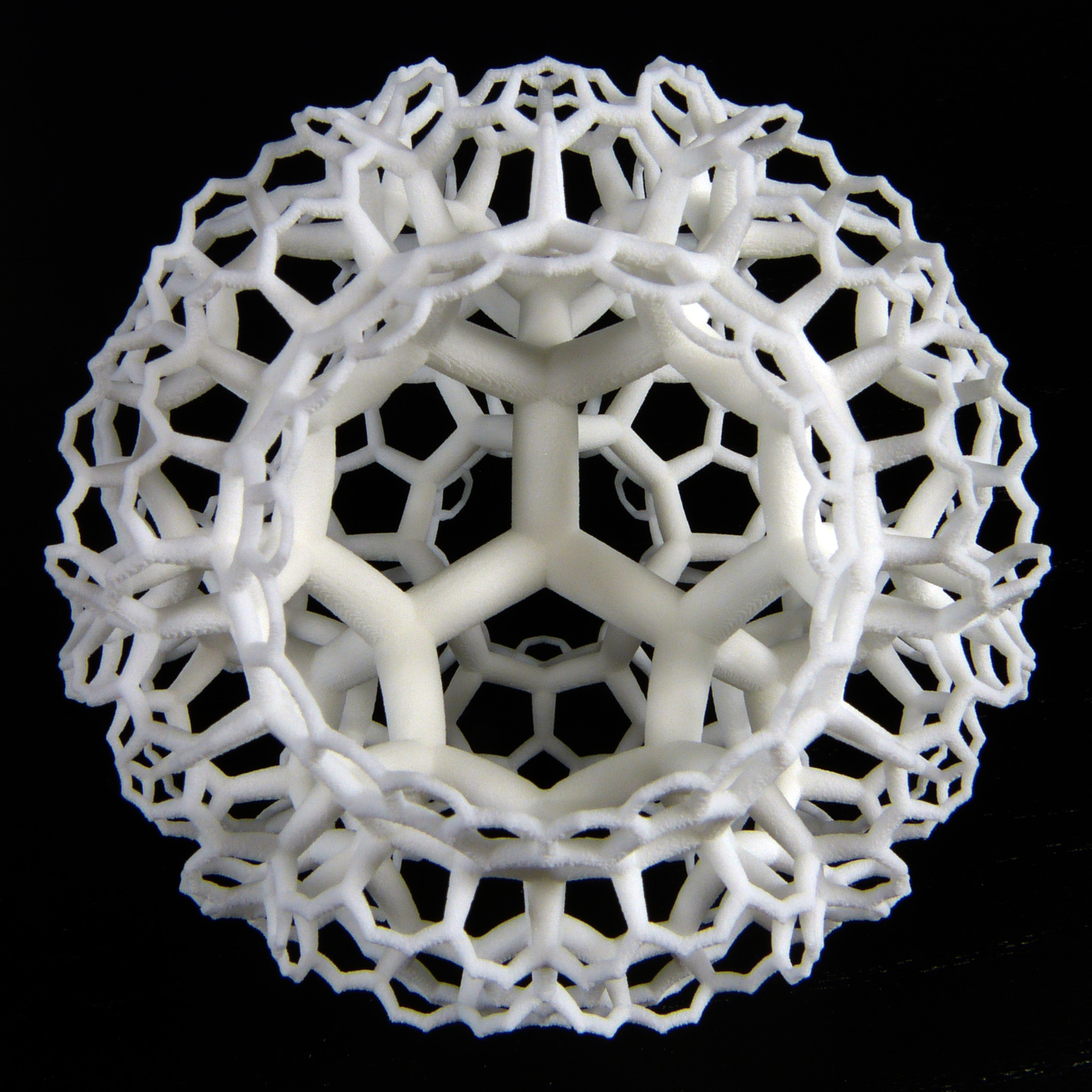}
\fi
\label{Fig:633}
}
\caption{The dual $\{3,3,6\}$ and $\{6,3,3\}$ honeycombs.}
\label{Fig:336_633}
\end{figure}

\begin{figure}[htbp]
\centering 
\hspace{-.7cm}
\subfloat[]
{
\ifimages
\includegraphics[width=0.43\textwidth]{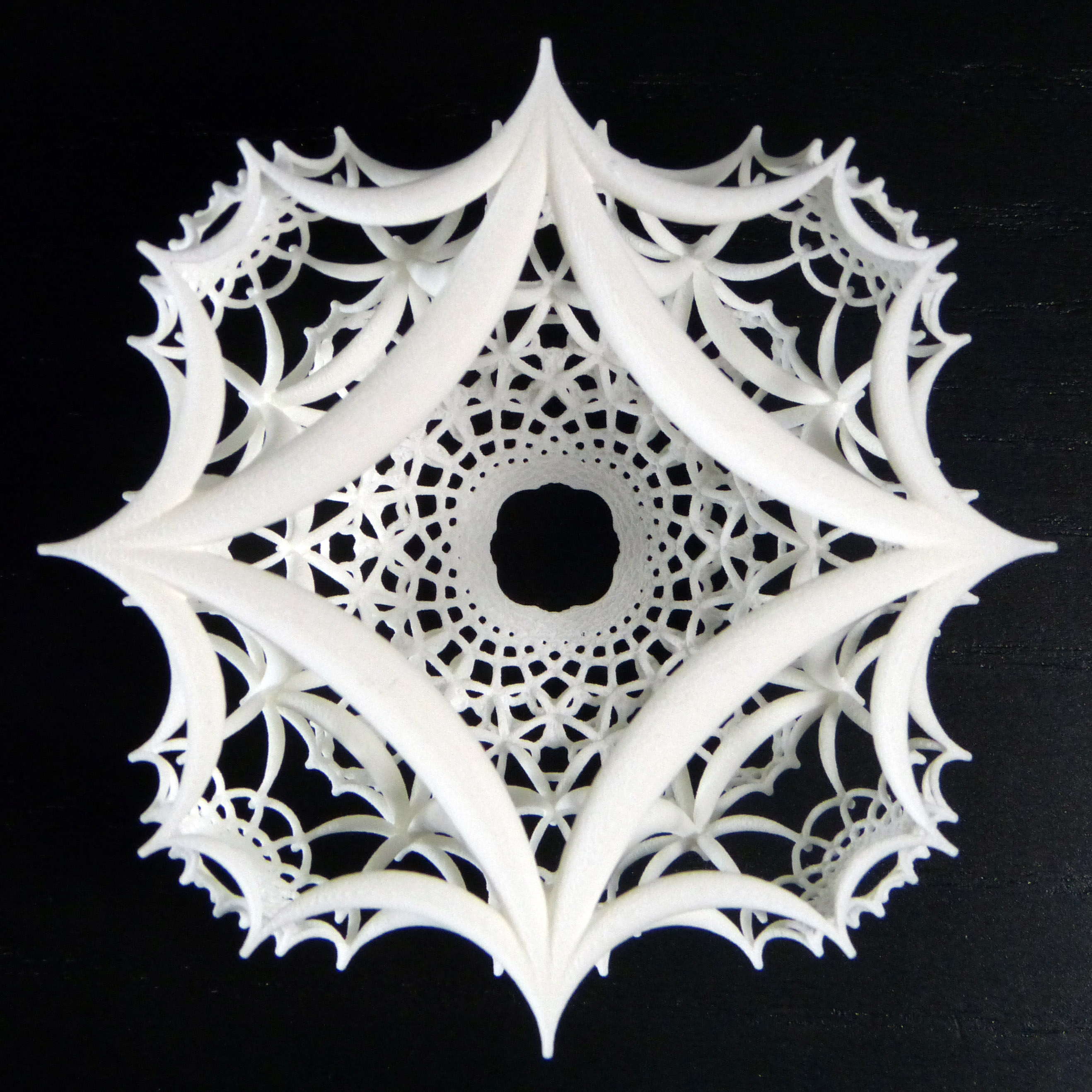}
\fi
\label{Fig:444_honeycomb_sym4}
} 
\hspace{.7cm}
\subfloat[]
{
\ifimages
\includegraphics[width=0.43\textwidth]{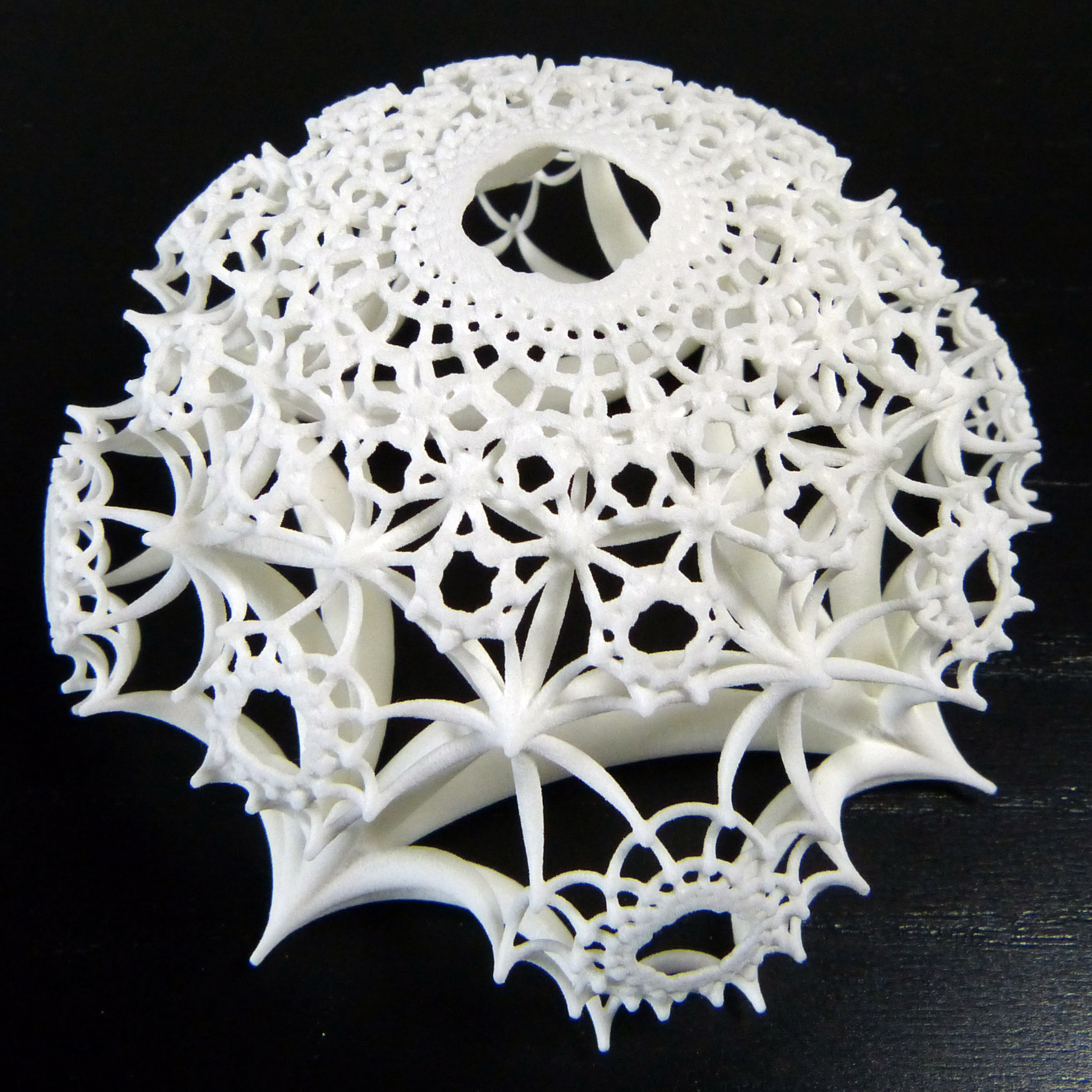}
\fi
\label{Fig:half_444_H3_honeycomb_1}
}
\caption{Two views of the self-dual $\{4,4,4\}$ honeycomb.  Only half of the Poincar\'e ball is rendered in this print.}
\label{Fig:444}
\end{figure}

\subsection{Two-dimensional images}

As we saw in Section \ref{Sec:schlafli_cube}, boundary patterns of hyperbolic honeycombs provide a huge collection of strikingly varied imagery. There are many other ways to show these honeycombs, giving yet more images.

\begin{figure}[htbp]
\centering 
\includegraphics[width=0.8\textwidth]{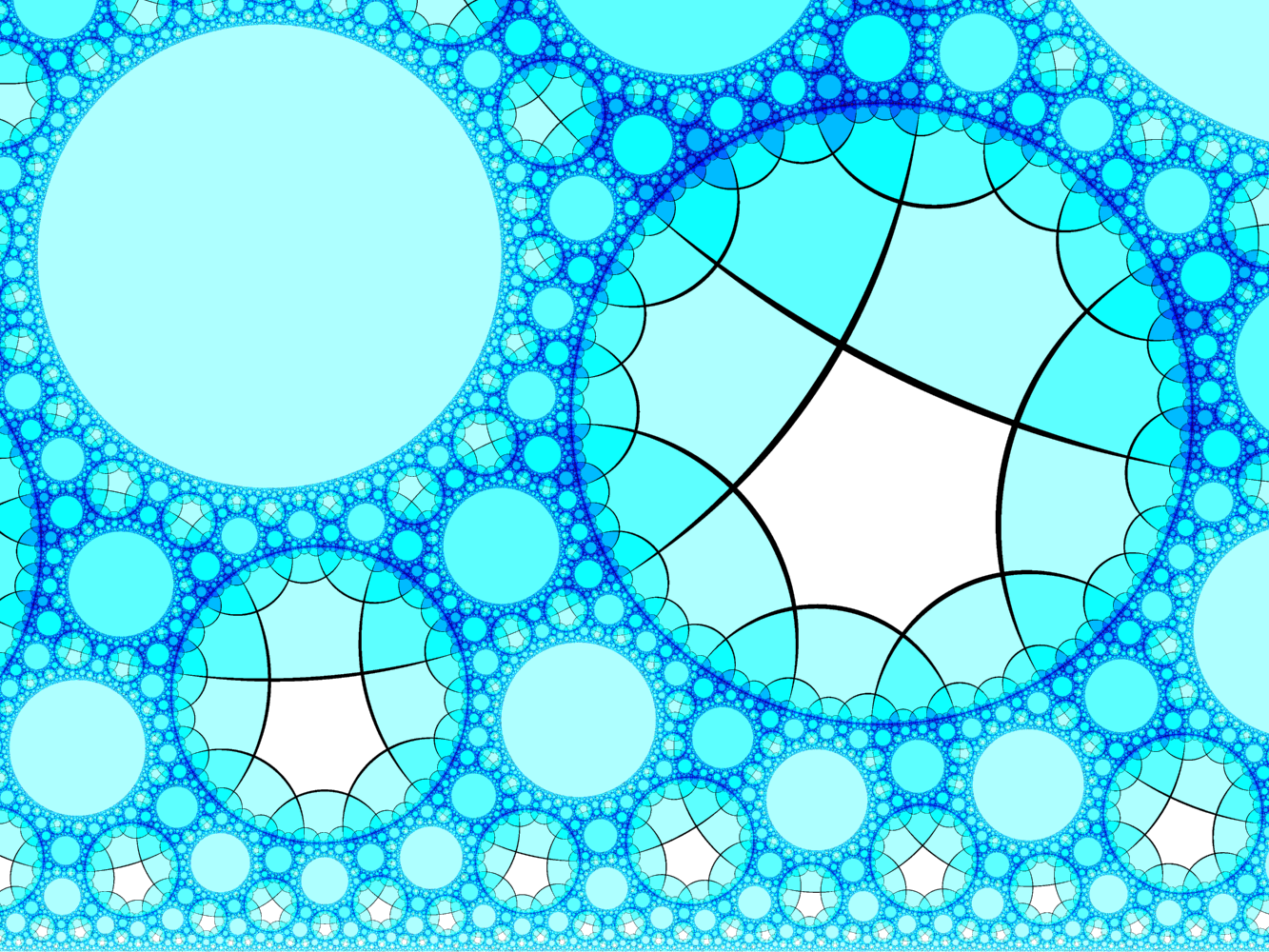}
\caption{Part of the boundary of the $\{4,6,4\}$ honeycomb.}
\label{Fig:464}
\end{figure}

\begin{figure}[htbp]
\centering 
\hspace{-.7cm}
\subfloat[The $\{7,3,\infty\}$ honeycomb, with a point at which two legs of a hyperideal cell meet sent to infinity.  Compare with the corresponding image in Table \ref{Table:pqi}.]
{
\includegraphics[width=0.43\textwidth]{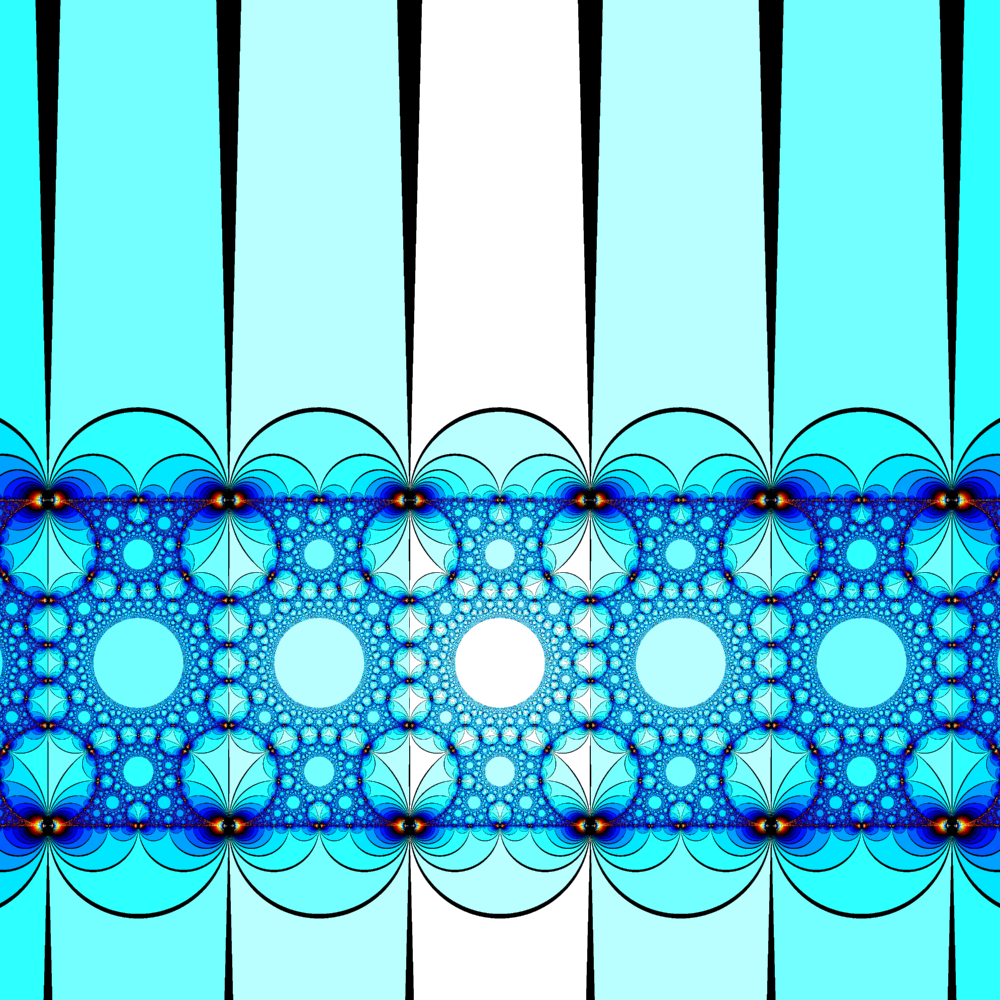}
\label{Fig:73i}
} 
\hspace{.7cm}
\subfloat[The $\{4,3,7\}$ honeycomb, oriented so that the image is self-similar.]
{
\includegraphics[width=0.43\textwidth]{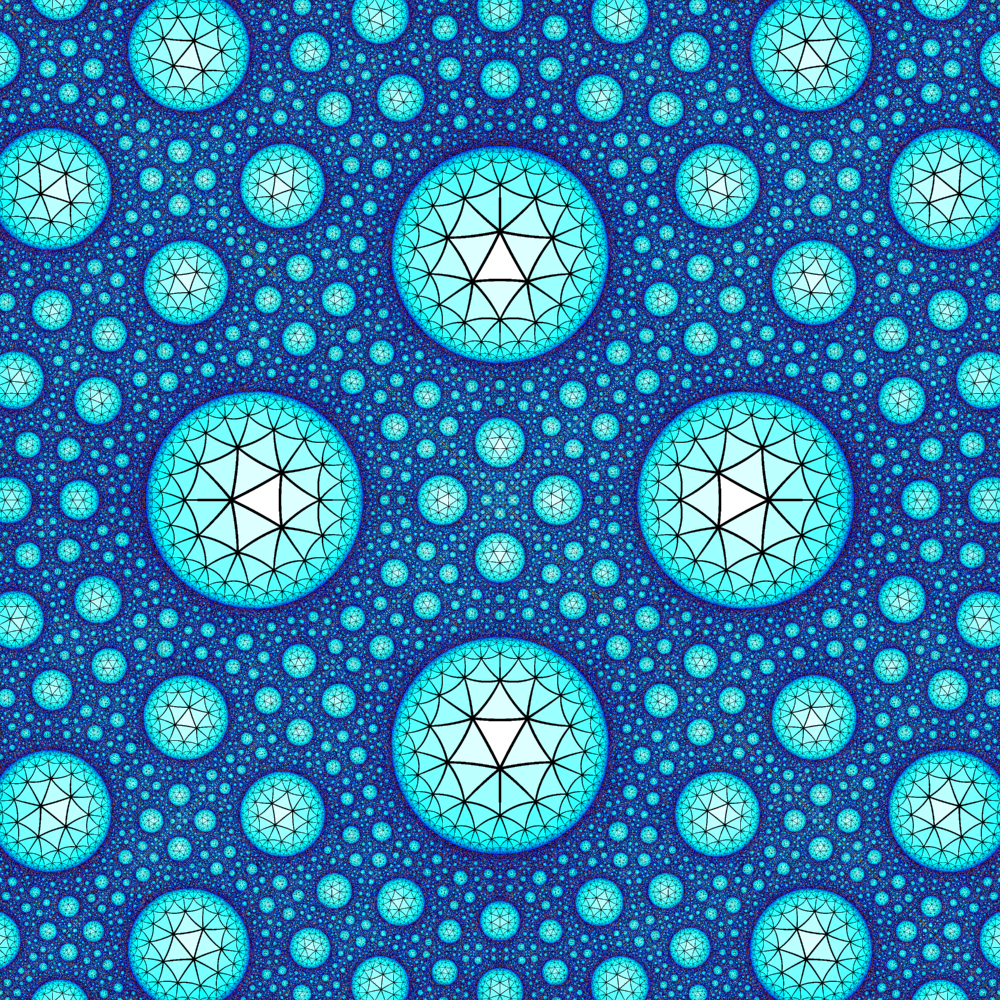}
\label{Fig:437_self_similar}
}
\caption{Further variations, applying isometries to the boundary image.}
\label{Fig:isometries}
\end{figure}

One way we can alter the images is by applying hyperbolic isometries.  Figure \ref{Fig:464} shows part of the $\{4,6,4\}$ honeycomb, drawn on the boundary of the upper half space model of $\HH^3$, as in Section \ref{Sec:schlafli_cube}. Here however, we apply an isometry so that the head of the white cell takes up the whole of the lower half plane.  Figure \ref{Fig:73i} shows a variation that results in a pattern with translational symmetry.  Figure \ref{Fig:437_self_similar} shows the $\{4,3,7\}$ honeycomb arranged so that the image is self-similar after a uniform scaling of 
$$e^{2 \thinspace\text{arccosh} \left( \sqrt{2}\cos \frac{\pi}{7} \right)} \approx  4.25917.$$
If we look down on the upper half space model, with the boundary of $\HH^3$ being the page, then for this honeycomb we are looking down an infinite tower of cubes. The isometry that moves one cube in the tower to the next corresponds to scaling the boundary pattern, and so the boundary pattern is self-similar.
 The scaling factor may be calculated from the distance between faces in this tower of cubes, which is twice the inradius of the honeycomb. The inradius can be calculated using a formula due to Coxeter~\cite[p. 158]{Coxeter1954}. 

\begin{figure}[htbp]
\centering 
\includegraphics[width=0.8\textwidth]{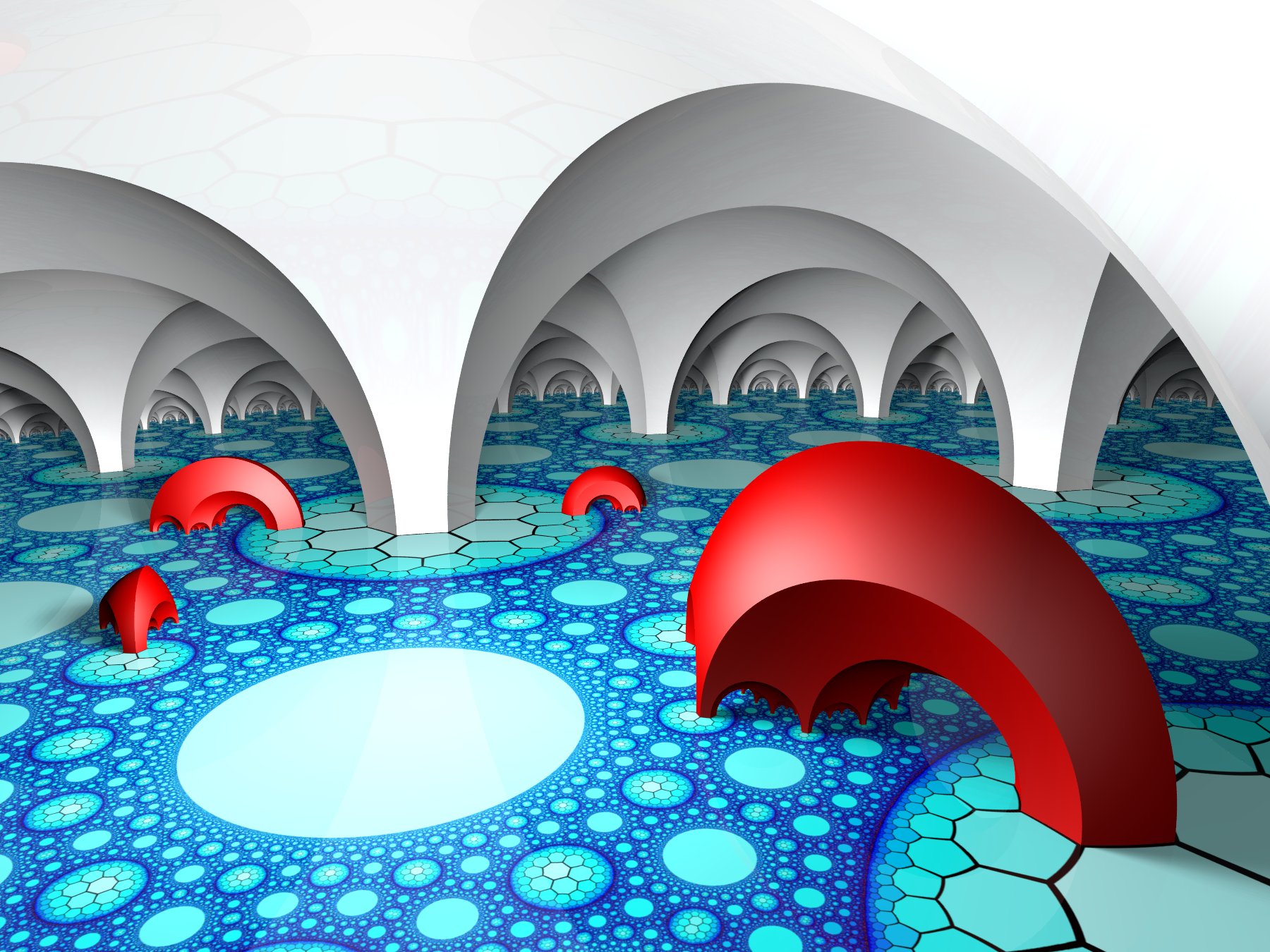}
\caption{Hyperbolic Catacombs.}
\label{Fig:catacombs}
\end{figure}

Figure \ref{Fig:catacombs} is based on the $\{3,7,3\}$ honeycomb. The floor of the ``catacombs'' is the boundary of the $\{3,7,3\}$ honeycomb as shown in Tables \ref{Table:3qr} and \ref{Table:pq3}. The image as a whole shows a view of the upper half space model of $\HH^3$, with only five cells of the honeycomb drawn: The ceiling of the catacombs is one cell, the other four are the family of red ``creatures''.

\section{Future directions}
\label{Sec:future}

In this paper we considered only regular honeycombs.  We have rendered a number of in-space views of archimedean honeycombs for Wikipedia, but have yet to explore hyperideal archimedean honeycombs.  By \emph{archimedean}, we mean honeycombs not transitive on all flags: for example, honeycombs with archimedean cells or honeycombs with more than than one cell type.  By relaxing the regularity constraint, there are many further possibilities for 3D prints, as well as new upper half space boundary images.

Schl\"afli symbols can be generalized further, to four-dimensional tilings and higher. Four-dimensional tilings are difficult to visualize. For example we have to cheat in Figure \ref{Fig:433} to visualize a single hypercube -- really we are only seeing the edges of the hypercube. The four-dimensional tiling $\{4,3,3,4\}$, consisting of infinitely many hypercubes, seems to be out of reach of our techniques. Paradoxically, we may have better luck visualizing four-dimensional tilings with hyperideal vertices of cells, because just as in three dimensions, we can see what is happening on the boundary of $\HH^4$, which is $S^3$. Using stereographic projection again, we can visualize $S^3$.

\section*{Acknowledgments}
We would like to thank Don Hatch and Nan Ma for early discussion about visualizing hyperideal honeycombs.  The exciting thread that seeded our approach of boundary images is available online~\cite{Cubing2012}.  We thank Tom Ruen for useful discussions, and for suggesting the possibility of rendering $\{\infty,\infty,\infty\}$. We thank Saul Schleimer and Christopher Tuffley for discussions on calculating distance in the cell adjacency graph of a honeycomb. We also thank the reviewers for their valuable feedback.

\appendix

\section{Implementation details}
\label{Sec:implementation}

Code for generating the images and models in this paper is available online \cite{hyp_honeycombs_code}.

\subsection{Honeycomb construction}
\label{Sec:honeycomb_construction}

As mentioned in the introduction, our honeycombs and two-dimensional tilings are regular tilings of space, meaning that the entire tiling is determined by some fundamental domain, under the action of some symmetry group. We construct our honeycombs in the conformal models: stereographic projection for spherical space, ordinary euclidean space, and the upper half space model for hyperbolic space.  In these models, geodesic surfaces are spheres (or planes, i.e. spheres of infinite radius), and we can take the fundamental domain of a tiling to be a simplex (i.e. a triangle in two dimensions or a tetrahedron in three dimensions) bounded by such geodesic surfaces. Moreover, the symmetry group is generated by reflections in these geodesic surfaces. In our constructions, we calculate the surfaces for the fundamental simplex of the honeycomb according to the dihedral angles of the simplex (see Appendix \ref{Sec:simplex_construction}). The fundamental simplex then determines the honeycomb: we build the honeycomb by recursive reflections (sphere inversions~\cite[pp. 124-126]{Needham199902}) in the surfaces of the simplex.  

In terms of implementation, we use the faces of the simplex differently depending on the context.  For the 3D prints and POV-Ray renderings, we start with some honeycomb element, say an edge or a face, and reflect it outwards to build up the honeycomb, using some size threshold to stop the recursion. We used a hash map to avoid constructing multiple identical copies of each element. Gunn~\cite{gunn1993discrete} instead used the theory of automatic groups to address this problem. 

For the upper half space boundary images, we work in reverse.  For every pixel in the image, we move inwards, reflecting its position in the simplex faces until it is within the fundamental simplex.  We then colour the pixel based on where it ends up in the simplex (see Appendix \ref{Sec:drawing_lines}), and the distance of the cell the pixel is in from the cell containing the fundamental simplex, where distance is measured in the cell adjacency graph of the honeycomb (also see Appendix \ref{Sec:colouring} for details on the colour choices).  

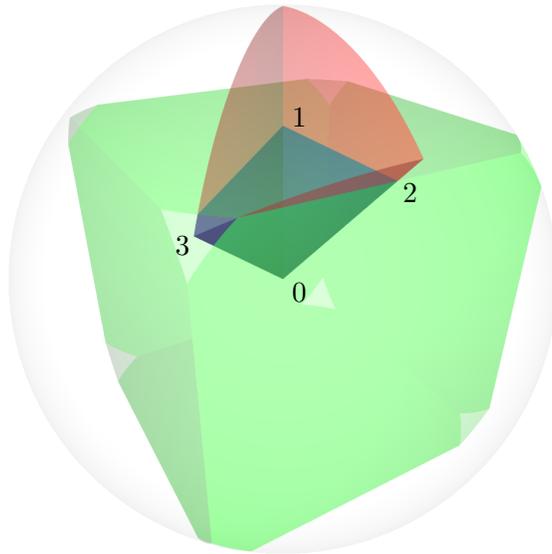
\begin{figure}[htbp]
\centering 

\labellist
\small\hair 2pt
\pinlabel 0 at 530 480
\pinlabel 1 at 530 780
\pinlabel 2 at 720 650
\pinlabel 3 at 330 560
\endlabellist

\includegraphics[width=0.5\textwidth]{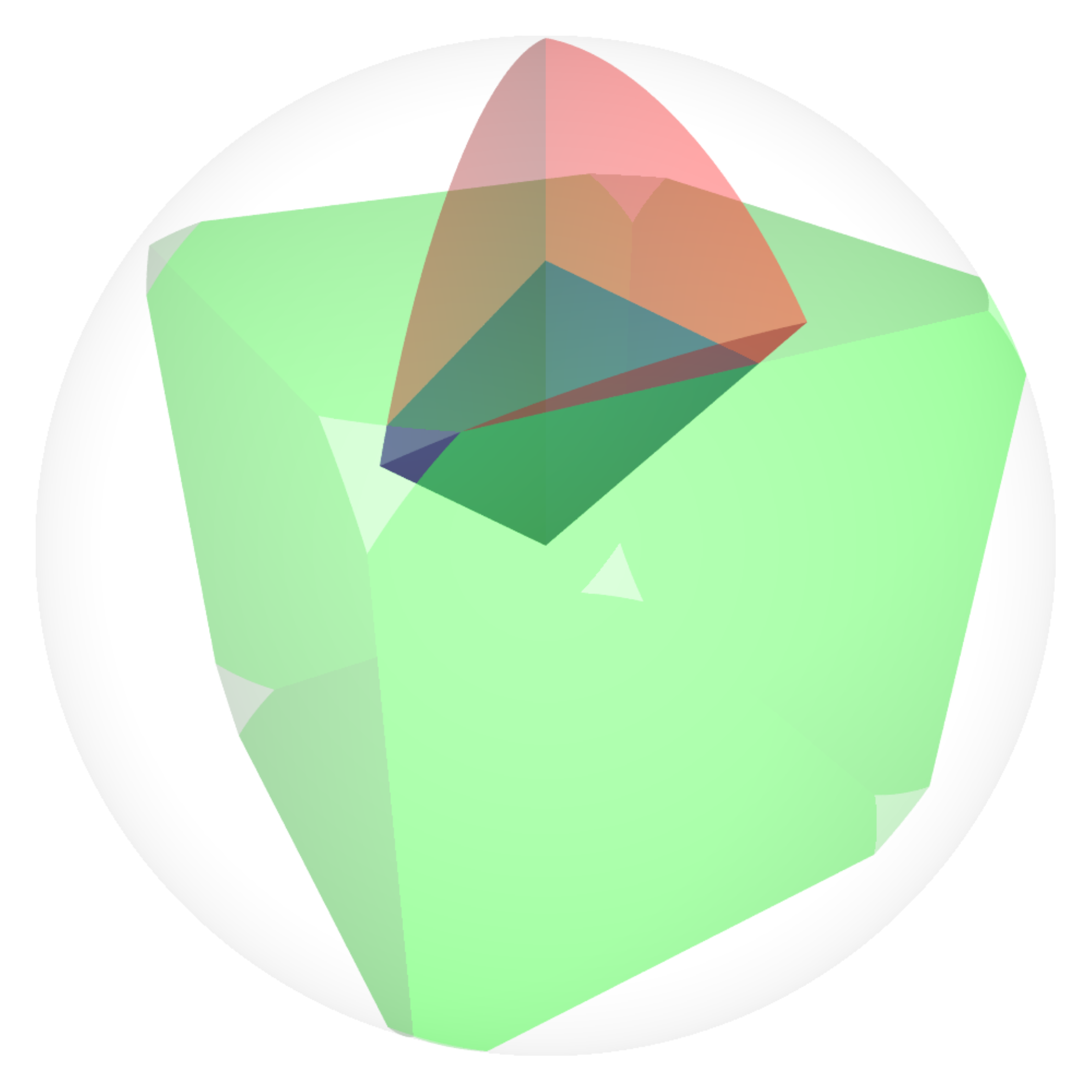}
\caption{The central cell (green), fundamental simplex (blue) and the fundamental wedge (red) for the $\{4,3,7\}$ honeycomb, drawn in the Klein model. The labels mark the four vertices of the simplex, which live at the cell center (0), face (1), edge (2), and vertex (3).}
\label{Fig:wedge}
\end{figure}

As an example to illustrate the algorithm we use to calculate this distance, consider the cell shown in Figure \ref{Fig:wedge}. The cell vertices are hyperideal (point 3 is outside of $\HH^3$), and the cell center is at the origin. In general, the vertices and centers of the cells may be material, ideal, or hyperideal, but the algorithm is the same in all cases. We number the simplex faces according to the numbering of the opposite vertex.  The three faces incident with the cell center (faces 1, 2 and 3) define a radial slice of the honeycomb, which we call the \emph{fundamental wedge}.  

Our algorithm proceeds as follows, starting with the location of a pixel of our image.
\begin{enumerate}
\item Set $d=0$.
\item Apply reflections in faces 1, 2, and 3 until we are in the fundamental wedge. 
\item If we are inside the central cell we are done, return $d$.  If not, reflect across face 0, increment $d$, and go to step (2).
\end{enumerate}

Note that there is no upper bound on the number of reflections needed to move the location of a pixel to the central cell, although in practice it rarely takes many steps. Our implementation stops the calculation if we have not reached the central cell in 4,000 reflections. We found that antialiasing made a huge difference in the quality of small features, so in fact for each pixel we sample 16 nearby points (in a $4\times 4$ grid) and average the colour of all of the points which converged.

We claim that the above algorithm calculates the distance (measured in the cell adjacency graph of the honeycomb) between the cell each pixel is in, and the central cell. This follows from standard results in Coxeter groups, see \cite[chapter 12]{garrett}. We sketch a proof here. 

We work on the level of simplices. We have the fundamental simplex $\sigma_f$, and the simplex $\sigma$ containing the current pixel, and we want to know the distance from one to the other. Each triangular face of a simplex in the honeycomb is part of a geodesic plane. Consider the collection $\mathcal{P}$ of all of these planes. Let $\mathcal{P}_i \subset \mathcal{P}$ be the planes that come from faces conjugate to face $i$ of the fundamental simplex. The simplices $\sigma$ and $\sigma_f$ are separated from each other by some number of geodesic planes in $\mathcal{P}$. The number of planes of $\mathcal{P}_0$ that separate $\sigma$ and $\sigma_f$ is the distance in the cell adjacency graph, and also the output given by the above algorithm. In fact, the number of planes of $\mathcal{P}_i$ that separate $\sigma$ and $\sigma_f$ is the number of reflections in face $i$ needed to get from $\sigma$ to $\sigma_f$. To see this, suppose face $i$ of $\sigma_f$ separates $\sigma$ from $\sigma_f$. If we were to reflect $\sigma_f$ across face $i$ to get $\sigma_f'$, then the number of planes of $\mathcal{P}_i$ that separate $\sigma$ from $\sigma_f'$ goes down by one, and no other number of separating planes has changed. Of course we aren't reflecting $\sigma_f$ across face $i$ -- instead we are reflecting $\sigma$ across that face to get $\sigma'$. But if we then reflect both $\sigma_f$ and $\sigma'$ across face $i$, we get $\sigma_f'$ and $\sigma$, and we are in the previously described situation, so the counts of the number of planes separating the two simplices changes as before.

This shows that the above algorithm counts the number of planes of $\mathcal{P}_0$ that separate $\sigma$ from $\sigma_f$. Also note that reflecting $\sigma_f$ to $\sigma'_f$ is a single step in the \emph{simplex} (as opposed to cell) adjacency graph, so the distance from $\sigma'_f$ to $\sigma$ (which is equal to the distance from $\sigma_f$ to $\sigma'$) is at most one less than the distance from $\sigma_f$ to $\sigma$. In fact, a single reflection either increases or decreases the distance by one -- it cannot leave the distance the same. An induction argument then shows that the number of planes separating $\sigma$ from $\sigma_f$ is the same as the distance in the simplex adjacency graph. The same result also follows only paying attention to reflections in the planes of $\mathcal{P}_0$, corresponding to distance in the cell adjacency graph.

\subsection{Simplex construction}
\label{Sec:simplex_construction}

\begin{figure}[htbp]
\centering 
\subfloat[The boundary pattern. The meta-$\{6,5\}$-tiling, also generated by the red mirrors, is overlaid in black.]{
\includegraphics[width=0.46\textwidth]{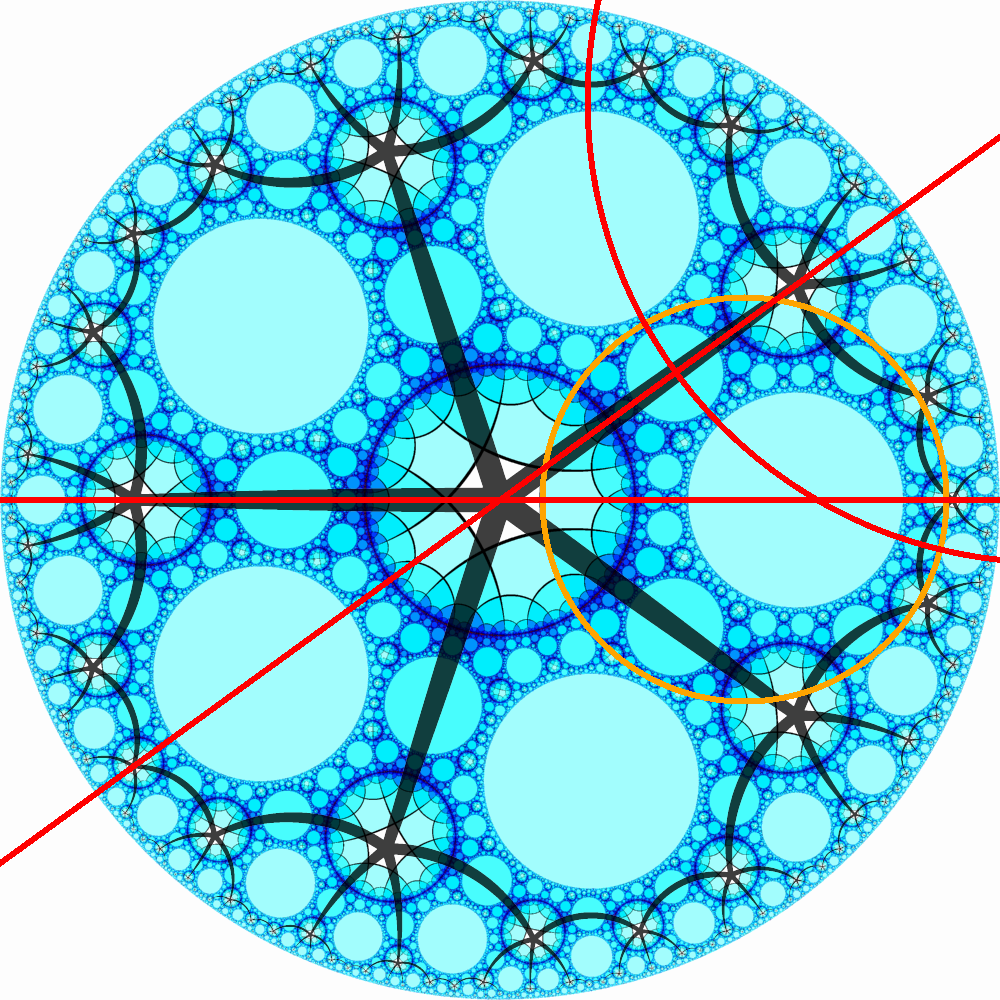}
\label{Fig:654_meta}
}
\quad
\subfloat[Fundamental simplex of the $\{6,5,4\}$ honeycomb in the upper half space model. ]{
\includegraphics[width=0.46\textwidth]{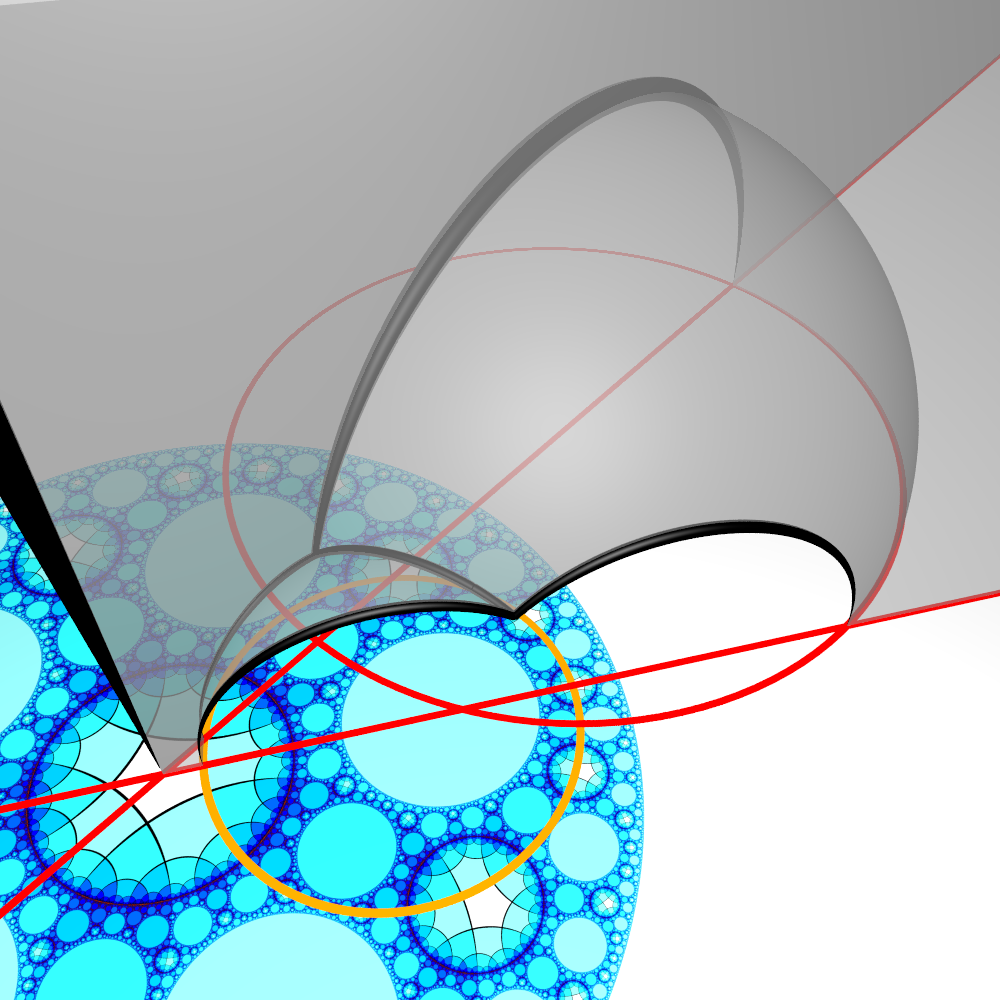}
\label{Fig:654_simplex}
}

\caption{The $\{6,5,4\}$ honeycomb. The intersection of simplex faces 1, 2, and 3 with the boundary plane are drawn red.  The intersection of face 0 is drawn orange.}
\label{Fig:654_construction}
\end{figure}

Calculating the geometry of the fundamental simplex relies on an observation about hyperideal honeycombs, namely that the upper half space pattern of a $\{p,q,r\}$ honeycomb visually looks like a set of solid disks (cell heads) and $\{q,r\}$ disks (hyperideal vertices) arranged in a \emph{meta-$\{p,q\}$-tiling} pattern.  See Figure \ref{Fig:654_meta}.  In the meta-tiling, the solid disks live at meta-tile centers and the $\{q,r\}$ disks live at meta-tile vertices.  The meta-tiling may be spherical (as in the first three columns of Table \ref{Table:3qr}), euclidean (as in Table \ref{Table:ideal_cells}), or hyperbolic (as in the right two columns of Table \ref{Table:3qr}).

Three of the simplex mirrors (faces 1, 2, and 3, drawn in red) reflect the simplex to another simplex within a cell of the honeycomb, and the remaining mirror (face 0, drawn in orange) reflects the simplex to a simplex of a different cell.  On the plane at infinity, the three red mirrors are exactly those of the meta-tiling. The orange mirror is centered on a meta-tile, and its radius can be calculated using trigonometry so that $r$ cells surround each edge. 

Figure \ref{Fig:654_simplex} shows the fundamental simplex of the $\{6,5,4\}$ honeycomb in the upper half space model.  This simplex has two material vertices and two hyperideal vertices.  The six edges of the simplex are drawn in black, with constant hyperbolic thickness.  The simplex edges connected to a hyperideal vertex do not meet within the space, so this simplex itself has two legs, one at the left of the image (inside the cell leg in the center of the boundary pattern), and one at the right (inside the cell head centered at the boundary plane's point-at-infinity).  There is only one edge that connects to both hyperideal vertices, which is the vertical edge at the left.  The other end of this edge is not visible, since it meets the boundary at the boundary plane's point-at-infinity.  There is also only one edge connected to both material vertices, in the central part of the image.  The remaining four edges have one material and one hyperideal endpoint.  The four simplex faces are rendered in semitransparent gray.  Two are portions of planes and two are portions of spheres.  All are geodesic and orthogonal to the boundary.


Infinite terms in the Schl\"afli symbols provide challenges in the same way that they do for two-dimensional tilings.  For example, since an $\infty$-gon does not have a material center, we cannot draw a tile-centered $\{\infty,3\}$ tiling.  It is more natural to construct this tiling as the dual to a $\{3,\infty\}$ tiling.  For three-dimensional honeycombs, we handle various combinations of infinite $p$, $q$, or $r$ similarly, and our orientation choices were made to keep the visual progressions in our tables as consistent as possible.  


\subsection{Edges as Dupin cyclides}
\label{Sec:dupin}

Edges, the one dimensional elements of a honeycomb, are geodesics in the induced geometry.  To realize honeycomb edges in our physical models, we must give them some thickness.  A constant offset from the line of an edge forms a \emph{Clifford surface}.  From within hyperbolic space, thickened edges appear cylindrical, but they take on a banana-like shape when projected to the Poincar\'e ball model (see Figures \ref{Fig:435_535}, \ref{Fig:353}, \ref{Fig:633} and \ref{Fig:444}) or upper half space model (see Figure \ref{Fig:654_simplex}).  These shapes are known as \emph{Dupin cyclides}~\cite{Schrott2006}.  Dupin cyclides are inversions of cylinders, tori, or cones.

A cyclide is the envelope of a smooth one-parameter family of spheres. We construct meshes to represent our cyclides by calculating the centers and radii of regularly spaced spheres in the cyclide, and placing mesh vertices on the tangent circles where each sphere touches the cyclide -- the curve on the sphere furthest from the cyclide's core geodesic. 

Consider a sphere with hyperbolic center $\bm{p_\HH}$ in the Poincar\'e ball model, and hyperbolic radius $r_\HH$.  Because of the distortions of the Poincar\'e ball model, the euclidean center $\bm{p_\EE}$ of this sphere is not generally the same as $\bm{p_\HH}$. The exception is at the origin, where the hyperbolic and euclidean centers coincide, and the euclidean radius is $r_0=\tanh(r_\HH/2)$.  For all other spheres, take $d=\norm{\bm{p_\HH}}$, then a straightforward calculation shows the sphere has euclidean center

\begin{equation}
\bm{p_\EE}=\frac{\bm{p_\HH}(r_0^2-1)}{d^2r_0^2-1}
\end{equation}

and euclidean radius

\begin{equation}
r_\EE=\frac{r_0(d^2-1)}{d^2r_0^2-1}.
\end{equation}

Note that $\bm{p_\EE$} is closer to the origin than $\bm{p_\HH}$. The variable $r_0$ gives us a parameter to control the thickness of edges in our 3D printed model, balancing strength against cost and visibility into the interior of the structure.

\subsection{3D printing issues}
\label{Sec:3d_printing}


In the Poincar\'e ball model, features become very small as we approach the boundary of the ball. The 3D printers we used can print wire-like features with a minimum diameter of around 1mm, a restriction we quickly run into.  We manage this issue in two ways.  

For models with material vertices, such as the $\{6,3,3\}$ honeycomb shown in Figure \ref{Fig:633}, we print only the edges whose thickness remains above this threshold. This first stage of culling leaves ``dangling'' edges, which are connected to the rest of the model only at one end. A second culling removes the dangling edges, a process which can reveal yet more dangling edges. We repeat until every edge is connected at both ends.

For models with ideal vertices, such as the $\{4,4,4\}$  honeycomb shown in Figure \ref{Fig:444}, all edges approach a thickness of zero near the boundary of the ball.  For these, we artificially adjust the edge thickness.  In Figure \ref{Fig:444} the edges have accurate thicknesses until they shrink to 1mm, but remain a constant 1mm thickness beyond that. For the $\{3,3,6\}$ honeycomb in Figure \ref{Fig:336} we used edges of constant thickness in $\RR^3$.

In order to reduce costs, we also hollow out the interior of the cyclides, while keeping their wall thickness within a separate threshold required by the printers.

A particular challenge for these models is mesh unioning. We initially generated models in which each edge was a separate mesh.  For some of the simpler honeycombs, in which a small number of edges meet at a vertex, we were able (with some tinkering) to use tools such as netfabb~\cite{Netfabb} and MeshLab~\cite{Meshlab} to merge the meshes into a single mesh suitable for printing.  For models such as the $\{4,4,4\}$ honeycomb, so many meshes meet at each vertex that those software packages fail.  We wrote a function using the AbFab3D~\cite{AbFab3D} library that converts our meshes into voxels, performs unioning in voxel-space, then converts back to a mesh from the voxel representation.  Computer memory was the limiting factor for the resolution we could achieve using this method.

\subsection{Colouring schemes}
\label{Sec:colouring}

As described in Appendix \ref{Sec:honeycomb_construction}, we colour the cells in Section \ref{Sec:schlafli_cube} according to distance from a central cell in the cell adjacency graph of the honeycomb.  The immediate neighbours of the central cell all have the same colour, as do all of the cells at distance two, and so on.  We cycle through a set of colours, and repeat if necessary.

To generate the set of colours, we employ a simple colouring scheme based on traversing the edges of the RGB colour cube.  With the red, green, and blue colour components of colour normalized between zero and one, the components of the colour are all zero, corresponding to black, at one vertex of this cube.  At the antipode, the components are all one, corresponding to white.  Starting from the white vertex, we move to cyan, then blue, then black, then red, yellow, and back to white. The full colour path traverses a skew hexagon along the edges of the RGB colour cube.

The first three edges of this hexagon are cool colours and the final three edges are warm colours, producing something of a rainbow effect across the entire range.  We chose the rate of movement along the hexagon so that the images most visibly expose the cool range.  This rate is different for each Schl\"afli symbol, in order to best show the features of each image. The rate is chosen so that most of the pixels in the image have colours in the first half of the hexagon path. We calculated the cell-adjacency depth for each pixel, and then a mean depth for the entire picture.  We experimented with a number of different functions of this mean to control the rate of movement, and settled on a rate proportional to the square root of the mean, which we felt gave the best visual result.

\begin{figure}[htbp]
\centering 
\includegraphics[width=0.7\textwidth]{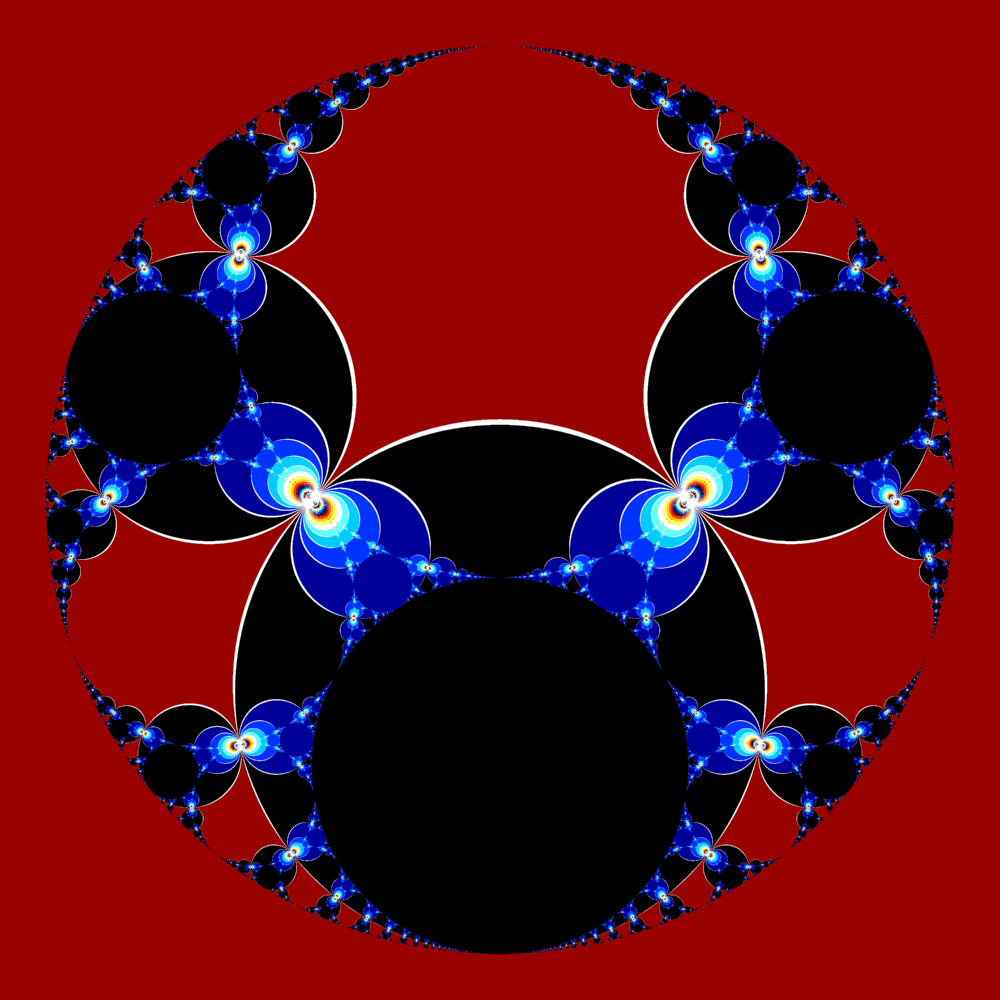}
\caption{The $\{\infty,\infty,\infty\}$ honeycomb.}
\label{Fig:iii}
\end{figure}

There is a choice of six different skew hexagons for a colouring path starting and ending at the white vertex.  Some of the other hexagons produced striking results, and we also experimented by starting at various points along the hexagon.  Figure \ref{Fig:iii} shows an artwork taking advantage of both of these changes.  It starts near the black vertex and moves along the hexagon described above in the opposite direction.

\subsection{Drawing curves on the upper half space boundary}
\label{Sec:drawing_lines}

In Figures \ref{Fig:437_boundary_pattern} through \ref{Fig:43r_progression},  many images in Section \ref{Sec:schlafli_cube}, and Figures \ref{Fig:464} through \ref{Fig:iii}, there is a subtlety in how we draw the cross section of our honeycombs on the boundary of upper half space.  The faces of the cells, even if given a constant thickness within the interior of $\HH^3$, have a cross section on the boundary along a line with zero thickness. In fact, any two distinct points on the boundary are infinitely far apart from each other.

However, for aesthetic reasons it is preferable to draw lines with thickness. Figure \ref{Fig:464} has many disks, each of which appears to contain a $\{6,4\}$ tiling of the Poincar\'e disk model of the hyperbolic plane. Each of these disks is associated to a hyperideal vertex of the honeycomb, where infinitely many cells meet, so we refer to these disks as \emph{vertex disks}. These disks are not two-dimensional tilings in the usual sense.  The apparent edges are intersections of faces with the boundary, and the apparent tiles are intersections of cells with the boundary.  A line of apparent edges travels through an infinite number of the vertex disks, expanding towards the center and shrinking towards the boundary of each vertex disk.

\begin{figure}[htbp]
\centering 
\includegraphics[width=0.6\textwidth]{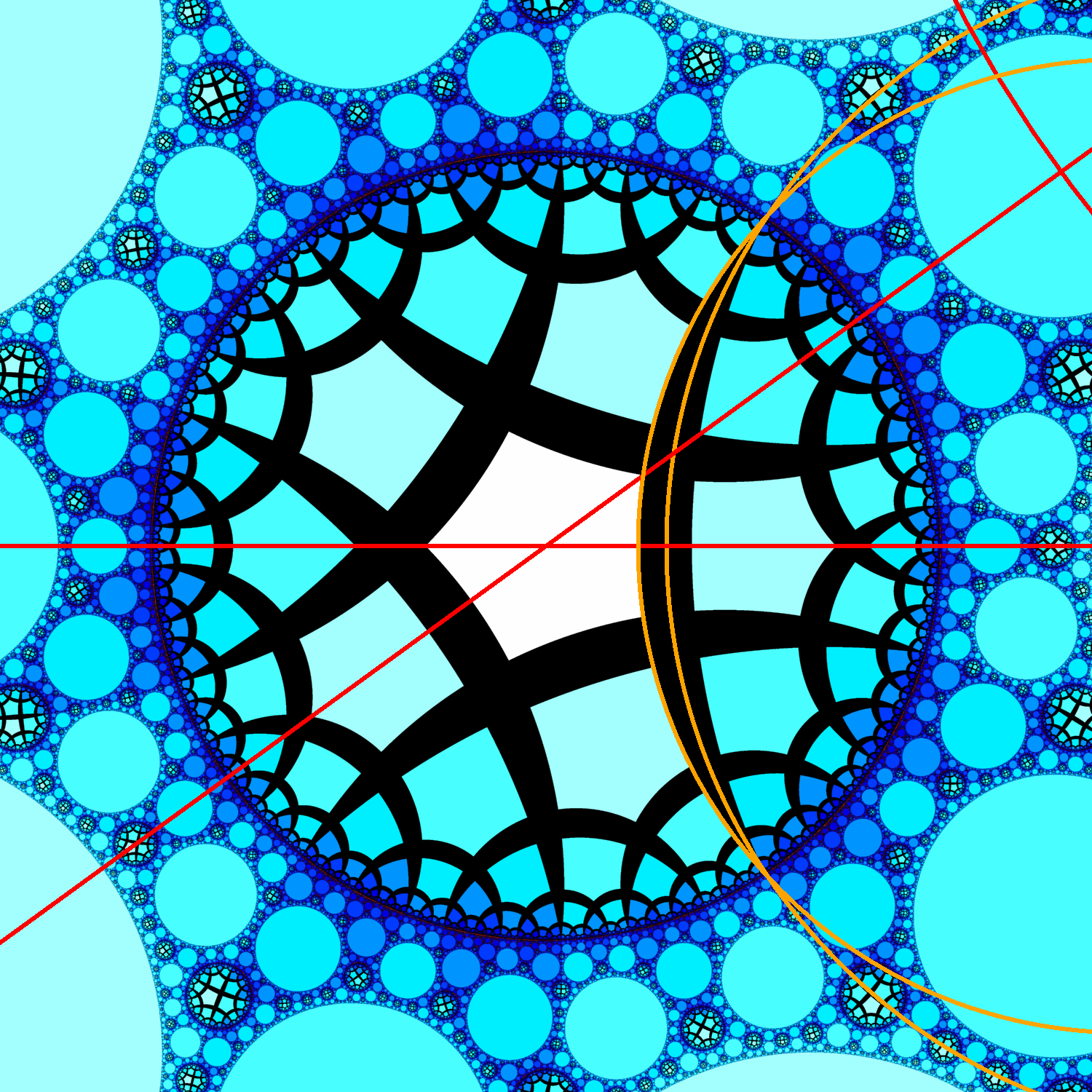}
\caption{A closeup of the central apparent Poincar\'e disk of the $\{6,5,4\}$ honeycomb from Figure \ref{Fig:654_construction}. We have exaggerated the banana thickness and have also rendered the boundary of the offset plane used to shade the bananas black.}
\label{Fig:654_closeup}
\end{figure}

We chose to render the vertex disks as geometrically accurate two-dimensional tilings, with geodesic edges rendered as constant thickness ``bananas'' within the hyperbolic metric in each apparent Poincar\'e disk model of the hyperbolic plane.  The method of achieving this effect involves using a geodesic plane closely related to, but different from the plane containing face 0 of the simplex (see Figure \ref{Fig:wedge}).  The new \emph{offset plane} coincides with the original plane on the boundary circle of a vertex disk, but diverges in the disk interior.  Any points that lie within the simplex but outside of the offset plane are colored black (or in Figure \ref{Fig:iii}, white). In all of our boundary pictures, the radius of each banana is 0.025, measured in the hyperbolic metric of each apparent Poincar\'e disk.  Figure \ref{Fig:654_closeup} shows the boundaries of both the original simplex planes, and of the offset plane used to achieve this effect.  For full details, see the implementation of the function that calculates the offset plane~\cite{hyp_honeycombs_code_bananas}.


\section{Uniqueness of honeycombs for a given Schl\"afli symbol}
\label{Sec:uniqueness}

In Appendix \ref{Sec:honeycomb_construction}, we described how we generate our honeycombs from a Schl\"afli symbol. The question arises whether or not there are other regular honeycombs described by the same Schl\"afli symbol. 

We produce our honeycombs, starting from a fundamental simplex, as shown in Figures \ref{Fig:wedge} and \ref{Fig:654_simplex}. In fact, since we are discussing \emph{regular} honeycombs, the entire honeycomb is determined by this fundamental simplex, under a symmetry group. So, uniqueness of the honeycomb follows from uniqueness of the fundamental simplex.

In \emph{Regular Polytopes}, Coxeter calls such a fundamental simplex (for a regular honeycomb) an \emph{orthoscheme}. The name here references the fact that three of the dihedral angles of the simplex are right angles. 
Coxeter shows~\cite[p. 139]{coxeter1973} (in the spherical case, explicitly), that the other three dihedral angles of the orthoscheme for the Schl\"afli symbol $\{p,q,r\}$ are $\pi/p, \pi/q,$ and $\pi/r$, and so this uniquely determines the simplex. For the hyperbolic case, by following the construction in Figure \ref{Fig:654_construction}, as long as none of the dihedral angles are zero, the four sphere (or plane) faces of the fundamental simplex are similarly determined by the dihedral angles, and so again, we have uniqueness.

\begin{figure}[htbp]
\centering 
\hspace{-.7cm}
\subfloat[]
{
\includegraphics[width=0.43\textwidth]{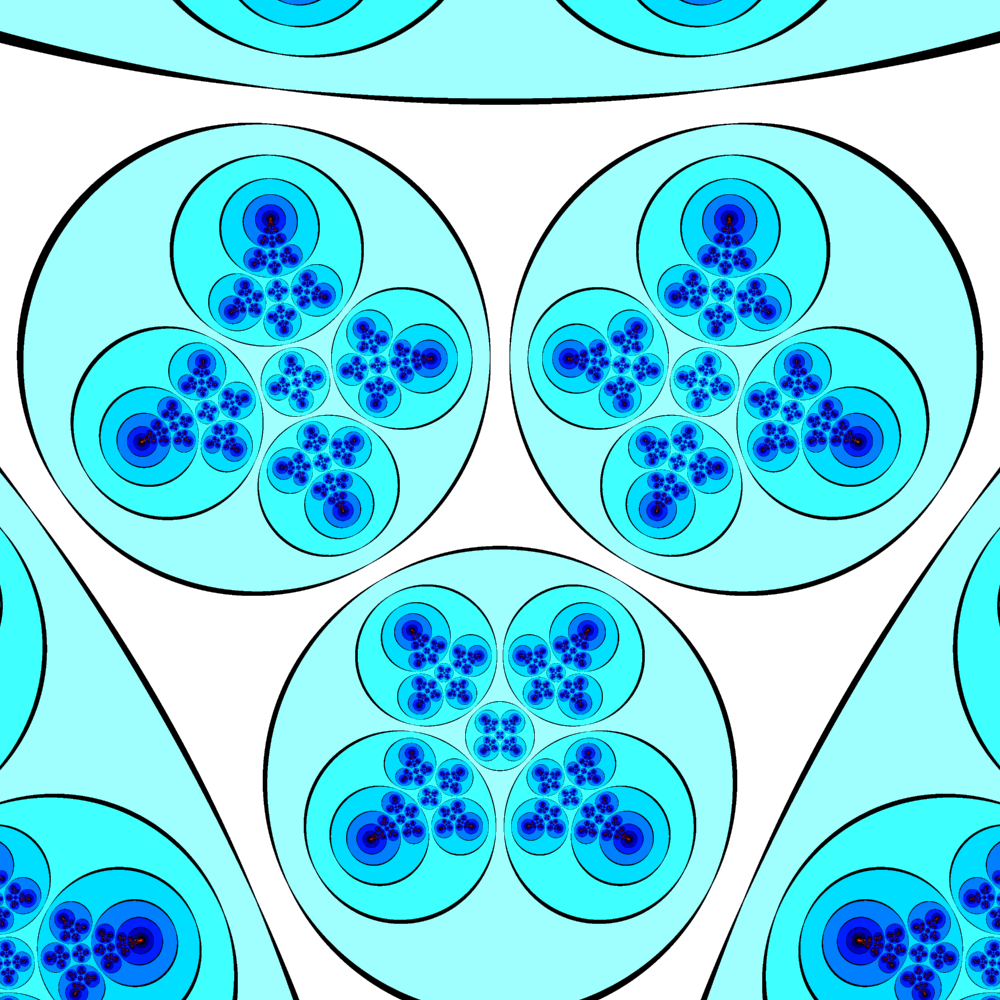}
\label{Fig:43j}
} 
\hspace{.7cm}
\subfloat[]
{
\includegraphics[width=0.43\textwidth]{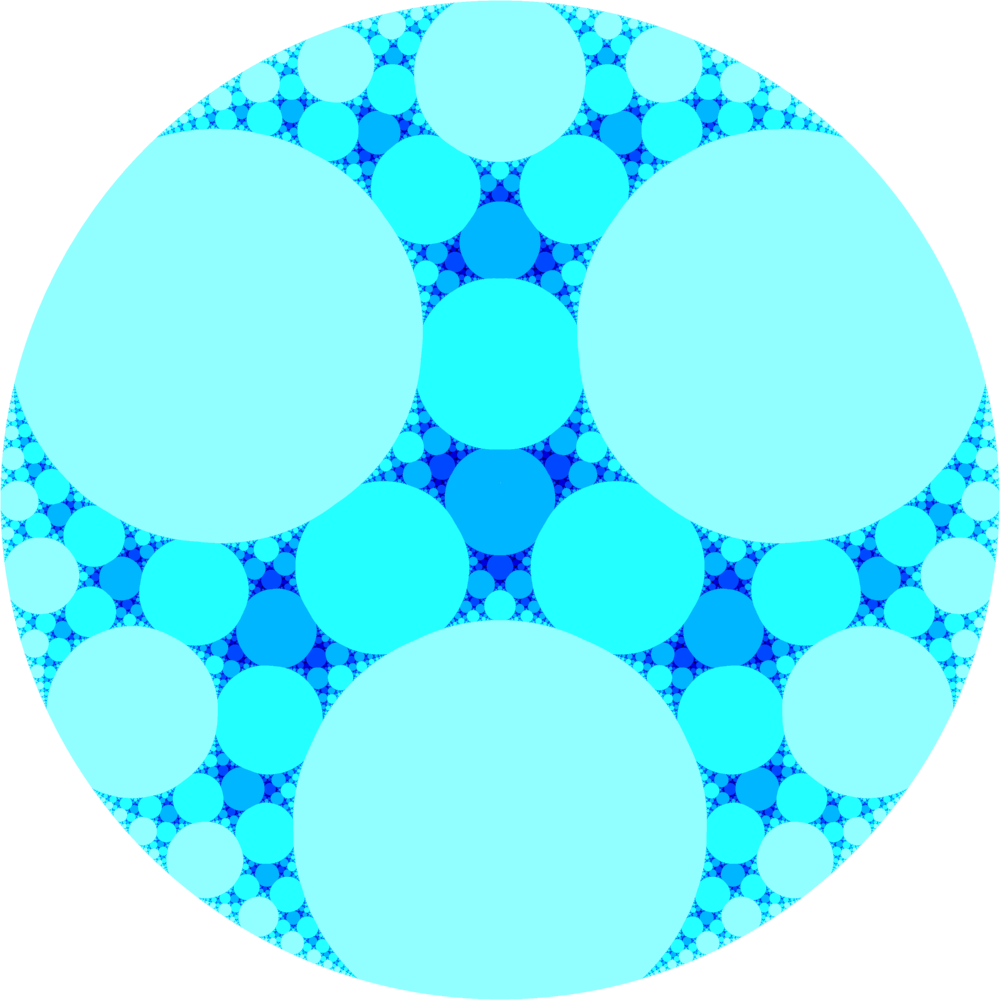}
\label{Fig:j34}
}
\caption{Alternate representations of the dual $\{4,3,\infty\}$ and $\{\infty,3,4\}$ honeycombs.}
\end{figure}

When one of the dihedral angles is zero, corresponding to an infinite term in the Schl\"alfli symbol, we lose uniqueness. For example, see Figure \ref{Fig:43j}, which ``continues'' the sequence shown in Figure \ref{Fig:43r_progression}. In Figure \ref{Fig:43i}, the edges of the cube meet $\HH^3$ only in a single ideal point, while in Figure \ref{Fig:43j} they no longer intersect $\HH^3$ at all, even on the boundary of $\HH^3$. This is still a regular honeycomb, and it is still a somewhat reasonable representation of the Schl\"afli symbol $\{4,3,\infty\}$, having the same symmetry group as the honeycomb shown in Figure \ref{Fig:43i}. One could question however, whether it makes sense to say that infinitely many cells fit around an edge, when the edge in question doesn't exist in the space, even in a limiting sense.

When drawing honeycombs for Schl\"afli symbols with infinite terms (other than in this appendix), we have required that the simplex edges intersect the boundary of $\HH^3$ in at least one point. With this extra condition, the faces of the fundamental simplex are uniquely determined, and we again have a unique honeycomb. 

\begin{figure}[htbp]
\centering 
\hspace{-.7cm}
\subfloat[]
{
\includegraphics[width=0.43\textwidth]{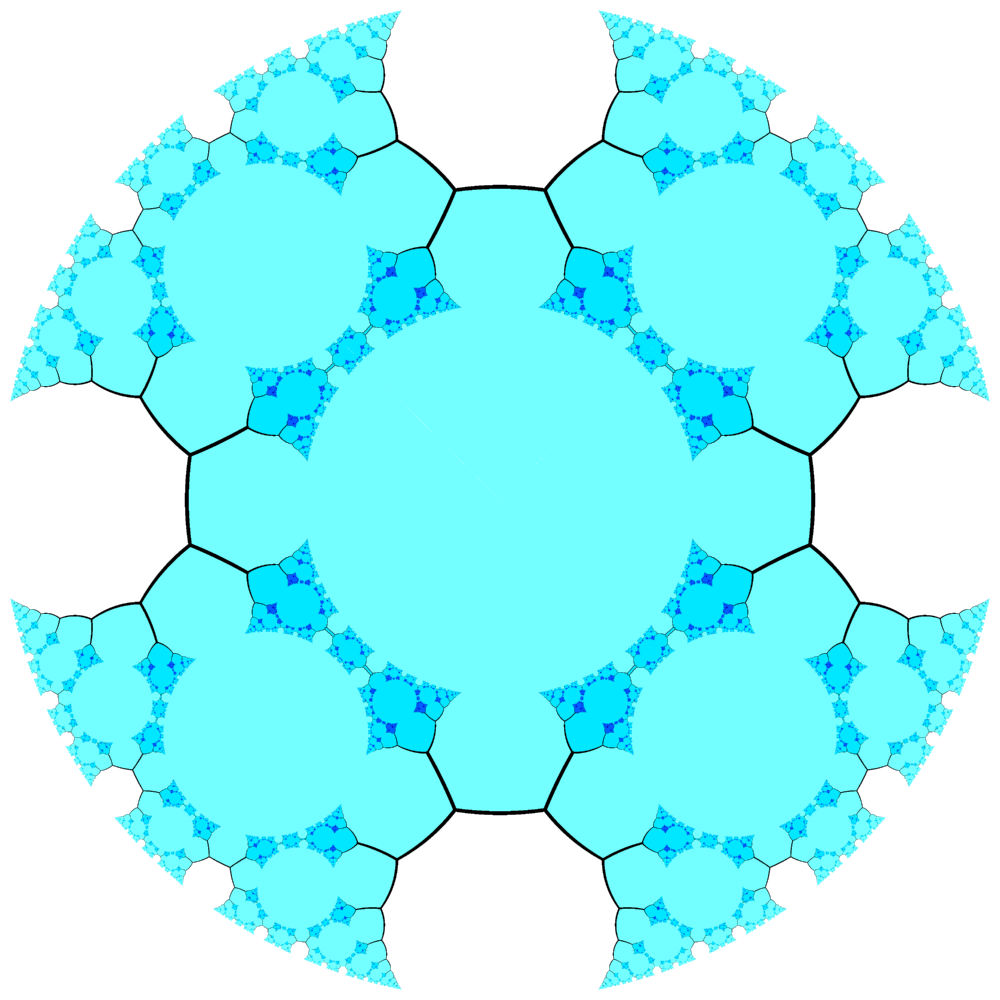}
} 
\hspace{.7cm}
\subfloat[]
{
\includegraphics[width=0.43\textwidth]{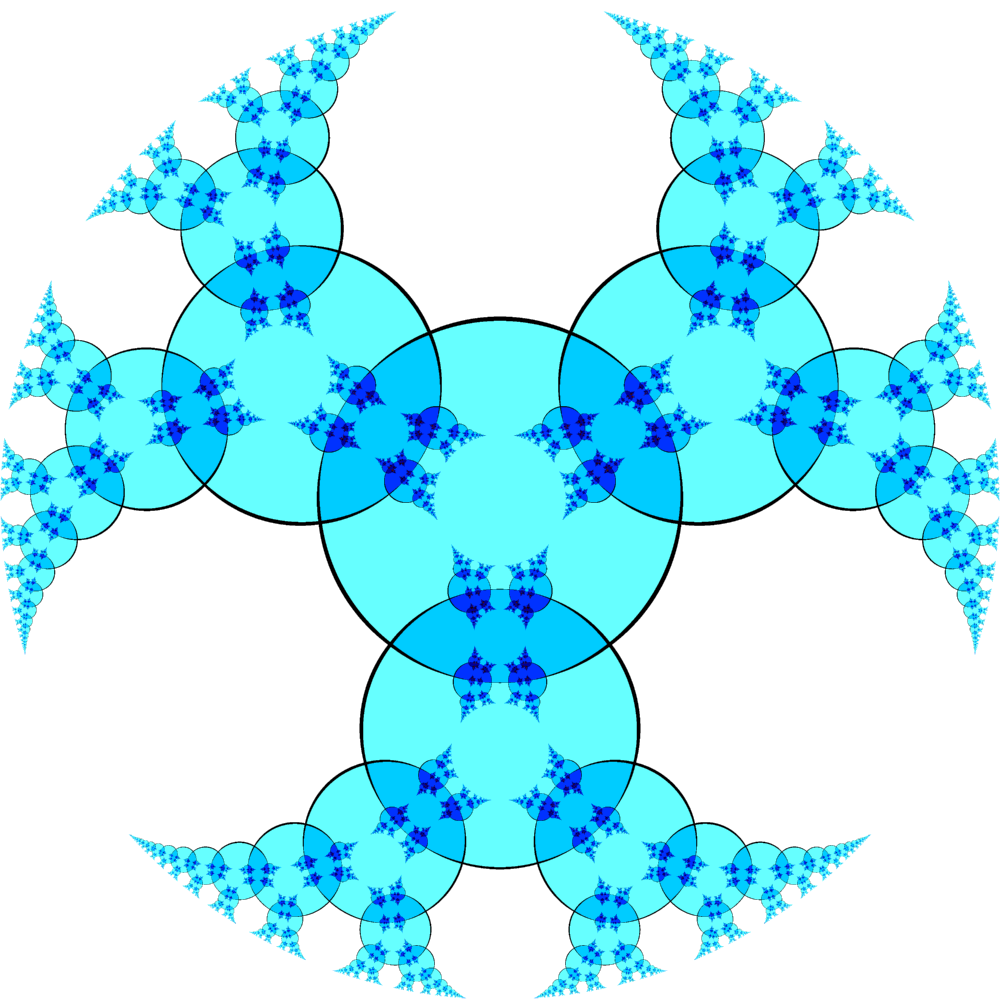}
}
\caption{Alternate representations of the dual $\{4,\infty,3\}$ and $\{3,\infty,4\}$ honeycombs.}
\label{Fig:4j3_3j4}
\end{figure}

The dual $\{\infty,3,4\}$ honeycomb, shown in Figure \ref{Fig:j34}, has the same fundamental simplex as $\{4,3,\infty\}$.  Even though it has hyperideal cells rather than hyperideal vertices, the same condition gives uniqueness.  We require the same edge of the fundamental simplex to have at least one ideal point.  When $q$ is infinite, we can also get a unique honeycomb with this restriction, although a different point of the edge ends up ideal.  Figure \ref{Fig:4j3_3j4} shows alternative honeycombs with $q=\infty$ and the uniqueness restriction lifted.

\begin{figure}[htbp]
\centering 
\hspace{-.7cm}
\subfloat[]
{
\includegraphics[width=0.43\textwidth]{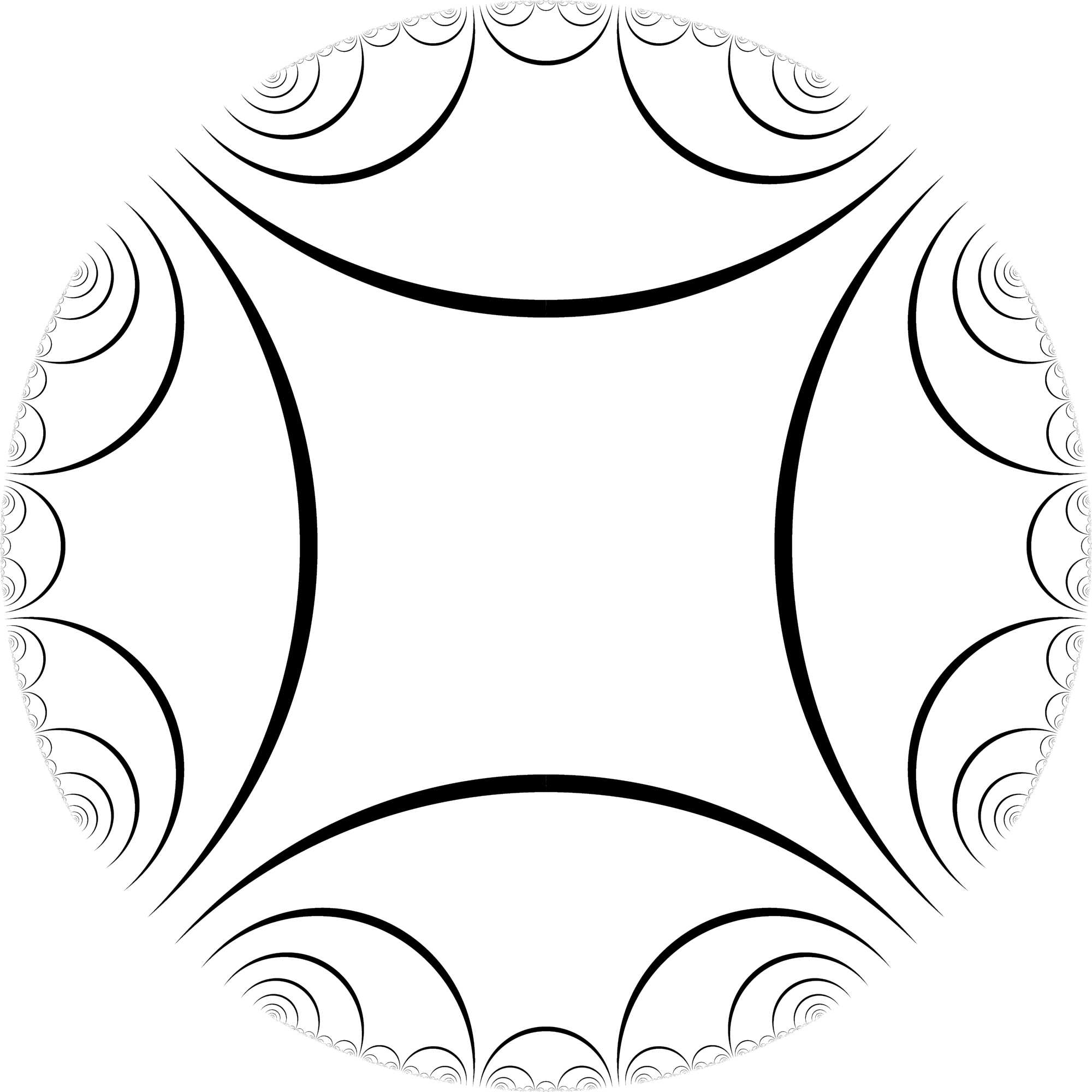}
} 
\hspace{.7cm}
\subfloat[]
{
\includegraphics[width=0.43\textwidth]{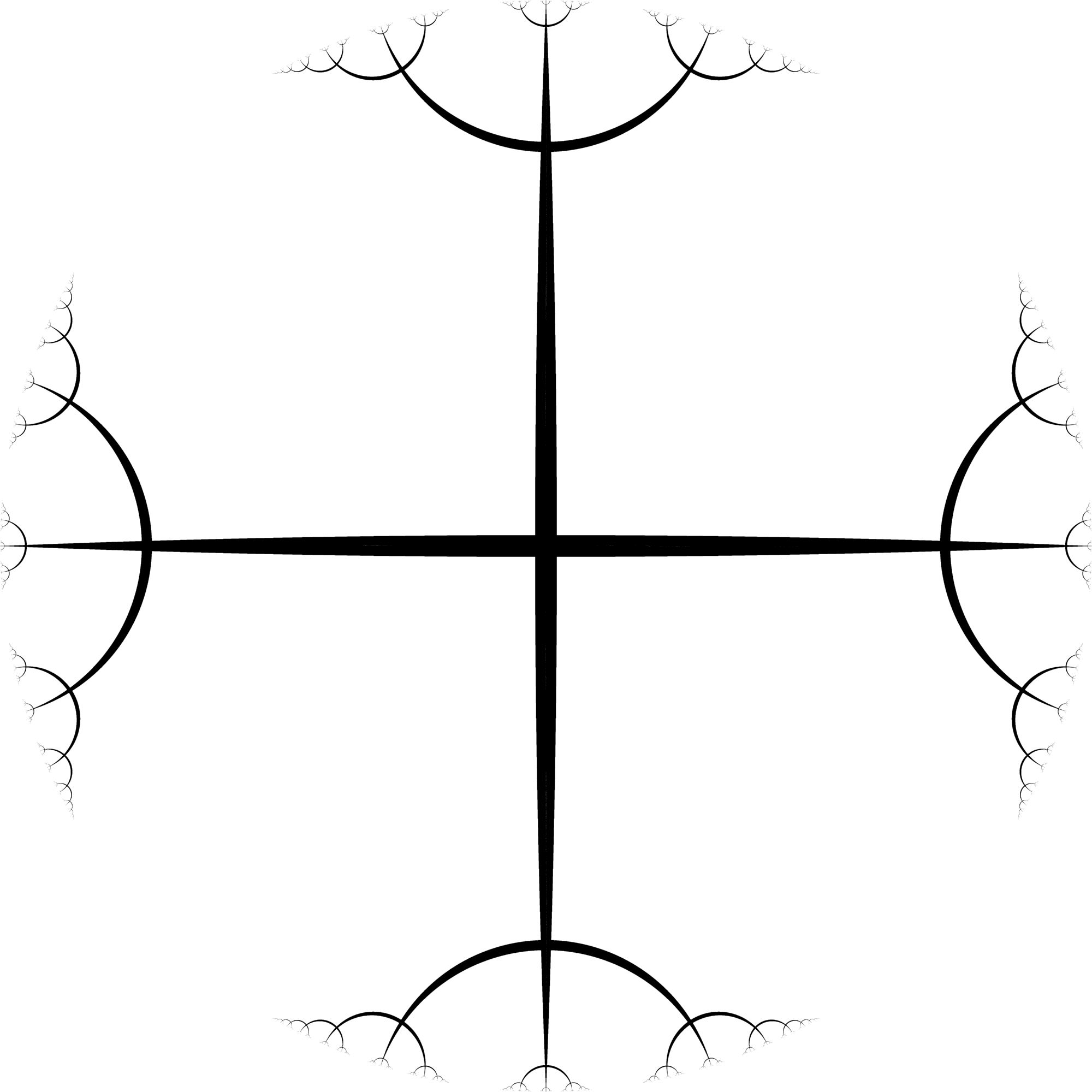}
}
\caption{Alternate representations of the dual $\{4,\infty\}$ and $\{\infty,4\}$ tilings.}
\label{Fig:4j_j4}
\end{figure}

Though we have omitted discussion here, the situation is similar in the two-dimensional case, as depicted in Figure \ref{Fig:4j_j4}.  Here the fundamental triangle angles determine uniqueness when $p$ or $q$ is finite.  When $p$ is infinite, we retain uniqueness by requiring that the vertices remain ideal. Dually, when $q$ is infinite, we retain uniqueness by requiring that the faces remain ideal.  

\RaggedRight 
\bibliographystyle{plain} 
\bibliography{hyperbolic_honeycombs}

\begin{thebibliography}{10}

\bibitem{bulatovbending}
Vladimir Bulatov.
\newblock Bending hyperbolic kaleidoscopes.
\newblock \url{http://bulatov.org/math/1107/}, 2011.

\bibitem{calegari20143}
Danny Calegari and Henry Wilton.
\newblock 3-manifolds everywhere.
\newblock arXiv:1404.7043, 2014.

\bibitem{chen2015lorentzian}
Hao Chen and Jean-Philippe Labb{\'e}.
\newblock Lorentzian {C}oxeter systems and {B}oyd--{M}axwell ball packings.
\newblock {\em Geometriae Dedicata}, 174(1):43--73, 2015.

\bibitem{Coxeter1954}
H.~S.~M. Coxeter.
\newblock Regular honeycombs in hyperbolic space.
\newblock In {\em Proceedings of the International Congress of Mathematicians
  of 1954}, 1954.

\bibitem{coxeter1973}
H.~S.~M. Coxeter.
\newblock {\em Regular Polytopes}.
\newblock Dover Publications, 1973.

\bibitem{Dix2013}
Jacob Dix.
\newblock Unbroken composition.
\newblock \url{http://www.mobilography.net/unbroken-composition/}, 2013.
\newblock Accessed: 2014-01-21.

\bibitem{not_knot}
D.B.A. Epstein and Charlie Gunn.
\newblock Not knot, video {ISBN} 0-86720-240-8 and supplemental book, 1991.

\bibitem{garrett}
Paul~B. Garrett.
\newblock {\em Buildings and Classical Groups}.
\newblock Chapman and Hall/CRC, 1997.

\bibitem{greenberg1993euclidean}
Marvin~J Greenberg.
\newblock {\em Euclidean and non-Euclidean geometries: Development and
  history}.
\newblock Macmillan, 1993.

\bibitem{Cubing2012}
4D~Cubing~Yahoo Group.
\newblock Hyperbolic honeycomb \{7,3,3\}.
\newblock
  \url{https://groups.yahoo.com/neo/groups/4D\_Cubing/conversations/messages/2291},
  2012.
\newblock Accessed: 2015-10-10.

\bibitem{gunn1993discrete}
Charlie Gunn.
\newblock Discrete groups and visualization of three-dimensional manifolds.
\newblock In {\em Proceedings of the 20th annual conference on Computer
  graphics and interactive techniques}, pages 255--262. ACM, 1993.

\bibitem{jos_leys}
Jos Leys.
\newblock Kleinian groups art galleries.
\newblock \url{http://www.josleys.com/galleries.php?catid=7}.
\newblock Accessed: 2016-09-11.

\bibitem{benoit1998multifractals}
Benoit~B. Mandelbrot.
\newblock Multifractals and $1/f$ noise, 1998.

\bibitem{Mumford2002}
David Mumford, Caroline Series, and David Wright.
\newblock {\em Indra's pearls: the vision of {Felix Klein}}.
\newblock Cambridge University Press, 2002.

\bibitem{Needham199902}
Tristan Needham.
\newblock {\em Visual Complex Analysis}.
\newblock Oxford University Press, USA, 1999.

\bibitem{hyp_honeycombs_code_bananas}
Roice Nelson.
\newblock Facet alteration function for apparent 2d tilings.
\newblock
  \url{https://github.com/roice3/Honeycombs/blob/fff18c1b96c6155fbabced2e5cae3ef045cd802a/code/HyperbolicModels/CoxeterImages.cs#L679},
  2015.

\bibitem{hyp_honeycombs_code}
Roice Nelson.
\newblock Software for generating images and models of hyperbolic honeycombs.
\newblock \url{https://github.com/roice3/Honeycombs}, 2015.

\bibitem{Povray}
Persistence of~Vision Pty. Ltd.~(2004).
\newblock Persistence of vision raytracer (version 3.7).
\newblock \url{http://www.povray.org/download/}, 2015.

\bibitem{phillips1992visualizing}
Mark Phillips and Charlie Gunn.
\newblock Visualizing hyperbolic space: Unusual uses of 4x4 matrices.
\newblock In {\em Proceedings of the 1992 symposium on Interactive 3D
  graphics}, pages 209--214. ACM, 1992.

\bibitem{schleimer2012sculptures}
Saul Schleimer and Henry Segerman.
\newblock Sculptures in {$S^3$}.
\newblock In Robert Bosch, Douglas McKenna, and Reza Sarhangi, editors, {\em
  Proceedings of Bridges 2012: Mathematics, Music, Art, Architecture, Culture},
  pages 103--110. Tessellations Publishing, 2012.
\newblock Available online at
  \url{http://archive.bridgesmathart.org/2012/bridges2012-103.html}.

\bibitem{Schrott2006}
Michael Schrott and Boris Odehnal.
\newblock Ortho-circles of {D}upin cyclides.
\newblock {\em J. Geom. Graph}, 10:73--98, 2006.

\bibitem{AbFab3D}
Abfab3d -- toolkit for 3d printing.
\newblock \url{http://abfab3d.com/}, 2014.

\bibitem{Meshlab}
Meshlab -- system for editing 3d triangular meshes.
\newblock \url{http://meshlab.sourceforge.net/}, 2015.

\bibitem{Netfabb}
netfabb -- software for 3d printing.
\newblock \url{https://netfabb.azurewebsites.net}, 2014.

\bibitem{stange2014visualising}
Katherine~E. Stange.
\newblock Visualising the arithmetic of quadratic imaginary fields.
\newblock arXiv:1410.0417, 2014.

\bibitem{Levy199701}
William~P. Thurston.
\newblock {\em Three-Dimensional Geometry and Topology, Vol. 1}.
\newblock Princeton University Press, 1997.

\bibitem{curved_spaces}
Jeff Weeks.
\newblock A flight simulator for multiconnected universes.
\newblock \url{http://geometrygames.org/CurvedSpaces}, 2016.

\end{thebibliography}

\end{document}